\documentclass[11pt]{article}

\usepackage{amsfonts,amsmath,amssymb,amsthm}
\usepackage{graphicx,caption,subcaption}
\usepackage{fullpage}
\usepackage{psfrag}
\usepackage{fancyhdr}
\usepackage{color,enumerate}
\usepackage{amsxtra,amscd}
\usepackage{setspace}
\usepackage{rotating}
\usepackage[none]{hyphenat}
\usepackage{placeins}

\usepackage[hidelinks,breaklinks]{hyperref}
\hypersetup{
    colorlinks=true,
    linkcolor=blue, 
    citecolor=blue, 
}
\usepackage{ragged2e}

\usepackage{multirow,rotating,makecell,longtable,threeparttable}
\usepackage{mathtools}
\usepackage{tikz,pgfplots}
\pgfplotsset{compat=1.18}
\usepackage{apptools}
\AtAppendix{\counterwithin{theorem}{section}}
\AtAppendix{\counterwithin{lemma}{section}}
\AtAppendix{\counterwithin{definition}{section}}

\usepackage{algorithm,algorithmicx,algpseudocode}
\usepackage{booktabs}

\usepackage{esvect}
\usepackage{enumitem}

\newtheorem{theorem}{Theorem}
\newtheorem{proposition}{Proposition}
\newtheorem{lemma}{Lemma}

\newtheorem{corollary}{Corollary}
\newtheorem{assumption}{Assumption}
\newtheorem{definition}{Definition}

\usepackage[numbers,sort&compress]{natbib}

\newcommand{\comment}[1]{}

\newcommand{\field}[1]{\ensuremath{\mathbb{#1}}}
\newcommand{\N}{\ensuremath{\field{N}}} 
\newcommand{\T}{{\top}} 
\newcommand{\R}{\ensuremath{\field{R}}} 
\newcommand{\1}{\vec{1}} 
\newcommand{\I}[1]{\ensuremath{\mathbf{1}{\{#1\}}}} 
\newcommand{\PR}{\ensuremath{\mathbb{P}}} 
\newcommand{\PP}{\ensuremath{\mathsf{P}}} 
\newcommand{\E}{\ensuremath{\mathbb{E}}} 

\newcommand{\sign}{\mathop{\mathrm{sign}}}
\newcommand\independent{\protect\mathpalette{\protect\independenT}{\perp}}
\def\independenT#1#2{\mathrel{\rlap{$#1#2$}\mkern2mu{#1#2}}}
\newcommand{\iidsim}{\stackrel{\mathrm{iid}}{\sim}}

\newcommand{\disteq}{\stackrel{\mathrm{d}}{=}}

\newcommand{\Sup}{\mathrm{supp}}

\newcommand{\rvline}{\hspace*{-\arraycolsep}\vline\hspace*{-\arraycolsep}}

\newcommand{\Escr}{\ensuremath{\mathcal E}}

\newcommand{\Gscr}{\ensuremath{\mathcal G}}

\newcommand{\Lscr}{\ensuremath{\mathcal L}}

\newcommand{\Nscr}{\ensuremath{\mathcal N}}

\newcommand{\dif}{\mathrm{d}}
\DeclareMathOperator{\diag}{diag}

\DeclareMathOperator{\Span}{\mathrm{span}}

\DeclareMathOperator*{\argmin}{\mathrm{argmin}}

\DeclarePairedDelimiter\abs{\lvert}{\rvert}

\DeclarePairedDelimiter\inprod{\langle}{\rangle}

\newcommand{\cov}{\ensuremath{\mathop{\mathrm{cov}}}}
\newcommand{\Var}{\ensuremath{\mathop{\mathrm{var}}}}
\newcommand{\trace}{\ensuremath{\mathop{\mathrm{tr}}}}

\newcommand{\smalop}{o_{\PP}(1)}

\newcommand{\probcov}{\stackrel{\PR}{\to}}
\newcommand{\plim}{\ensuremath{\mathop{\mathrm{plim}}}}

\newcommand{\ro}{\textup{\uppercase\expandafter{\romannumeral1}}%
}
\newcommand{\rt}{\textup{\uppercase\expandafter{\romannumeral2}}%
}
\newcommand{\roi}{\textup{\uppercase\expandafter{\romannumeral1}},\infty}
\newcommand{\rti}{{\textup{\uppercase\expandafter{\romannumeral2}},\infty}%
}

\newcommand{\ix}[1]{{X_{#1}^{\mathrm{iid}}}}
\newcommand{\stac}{\hat{\beta}_{\mathrm{stack}}}
\newcommand{\specop}{o^{*}_{\PR}(\sqrt{p})}

\newcommand{\indiv}[1]{\hat{\beta}_{#1}}
\newcommand{\aver}{\hat{\beta}_{\mathrm{average}}}
\newcommand{\indi}{\mathrm{ind}}

\newcommand{\joi}{\mathrm{joint}}
\newcommand{\ada}{\mathrm{ada}}
\newcommand{\se}{\hat{\beta }_{\rt}}

\DeclarePairedDelimiterX{\infdivx}[2]{(}{)}{%
  #1\;\delimsize\|\;#2%
}

\DeclarePairedDelimiter{\norm}{\lVert}{\rVert}

\definecolor{ysr}{RGB}{0,200,200}
\definecolor{kjr}{RGB}{200,0,200}

\author{Longlin Wang, Yanke Song, Kuanhao Jiang, Pragya Sur\footnote{Corresponding author: pragya@fas.harvard.edu} \\Dept.~of Statistics \\ Harvard University}

\title{Multi-Environment GLAMP: Approximate Message Passing for Transfer Learning with Applications to Lasso-based Estimators}

\begin{document}

\maketitle
\bibliographystyle{unsrtnat}

\begin{abstract}
Approximate Message Passing (AMP) algorithms enable precise characterization of certain classes of random objects in the high-dimensional limit, and have found widespread applications in fields such as signal processing, statistics, and communications. In this work, we introduce Multi-Environment Generalized Long AMP, a novel AMP framework that applies to transfer learning problems with multiple data sources and distribution shifts. We rigorously establish state evolution for multi-environment GLAMP.
We demonstrate the utility of this framework by precisely characterizing the risk of three Lasso-based transfer learning estimators for the first time: the Stacked Lasso, the Model Averaging Estimator, and the Second Step Estimator. We also demonstrate the remarkable finite sample accuracy of our theory via extensive simulations.
  \end{abstract}
  \vspace{-15mm} ~\\
  
  \newpage

\definecolor{ps}{RGB}{200,0,0}
\newcommand{\ps}[1]{\textcolor{ps}{[PS: #1]}}

\vspace{2cm}

\section{Introduction}\label{sec:intro} 

Transfer learning has emerged as a popular paradigm for handling distribution shifts between training and test time. In high-dimensional regression, substantial prior work has derived minimax error bounds demonstrating when auxiliary datasets allow to reduce estimation error and improve convergence rates in transfer learning \citep{bastani2021predicting,li2022transfer,lei2021near,cai2021transfer,zhang2022class,li2022searching, cai2022individual,maity2022meta,duan2023adaptive,li2023targeting,cai2024transfer}. These error bounds are often conservative, derived for worst-case scenarios, and omit constant factors. Consequently, when applied in practice, the risks for transfer learning based estimators may differ vastly, making it difficult to perform realistic comparisons among estimators. This paper addresses this difficulty by introducing a novel approximate message passing framework for tracking the risks of Lasso-inspired estimators in transfer learning. 

Specifically, we consider the 
proportional asymptotics regime, where the number of samples and features both diverge with the ratio converging to a constant. This regime has garnered significant attention in high-dimensional statistics and machine learning due to its ability to capture empirical phenomena in moderate to high-dimensional/overparametrized problems remarkably well \citep{johnstone2009statistical,donoho2009message,bean2013optimal,thrampoulidis2018precise,sur2019likelihood,sur2019modern,salehi2019impact,candes2020phase,barbier2019optimal,malekibridge, mondelli2021approximate,dhifallah2021phase,zhao2022asymptotic,jiang2022new,hastie2022surprises,liang2022precise,celentano2023lasso,song2024hede,luo2024roti,li2025optimal}.
We introduce a new approximate message passing framework, coined Multi-Environment Generalized Long AMP, that allows to characterize the risk of several transfer learning-based estimators. 
Concretely, we show the utility of our framework for the following three estimators: the Stacked Lasso \citep{bastani2021predicting}, the Model Averaging estimator \citep{bastani2021predicting}, and the Second-Step estimator \citep{li2022transfer}. To our knowledge, our paper provides the first specialized framework for characterizing the precise risks of these transfer learning-based estimators in high dimensions.

Approximate message passing (AMP) algorithms were introduced in the compressed sensing literature \citep{donoho2009message} and have seen enormous success in quantifying the risks of  convex regularized estimators in high dimensions \citep{bayati2011lasso,donoho2016high,zdeborova2016statistical,sur2019likelihood,sur2019modern,malekibridge,feng2022unifying,montanari2024friendly}. In its symmetric formulation, AMP takes as input a random design matrix $A \in \R^{N\times N}$, a sequence of deterministic  separable functions $f_t:\R^N\rightarrow \R^N$, and generates a sequence of vectors $x^t \in \R^N$ in the following way:
\begin{equation*}
x^{t+1} = Af_t(x^t) - f_{t-1}(x^{t-1})B_t^\top,
\end{equation*} where $B_t$ has a specific form determined by the choice of the denoising functions $f_t$. Remarkably, AMP iterates admit precise characterizations: for each $t$, as $N \rightarrow \infty$, $x^t$ behaves as a Gaussian variable with a tractable variance, a property known as \textit{state evolution} \citep{bolthausen2014iterative,bayati2011dynamics}. For
a broad range of statistical estimators, one can design suitable AMP iterations that track the estimator of interest (as $t$ additionally diverges to $\infty$), and thereby obtain exact risk characterization of these estimators. This is particularly useful for studying estimators without closed forms as demonstrated in \cite{bayati2011lasso,donoho2016high,sur2019modern}, among others. Since its original introduction, the AMP framework has been extended to Generalized AMP (GAMP) \citep{rangan2011generalized,javanmard2013state}, which allows iteration on matrices $x^t \in \R^{N\times q}$, and Long AMP (LAMP)  \citep{berthier2020state}, which allows $f_t$ to be non-separable, and Generalized Long AMP (GLAMP), which allows both \citep{gerbelot2023graph}.

For transfer learning scenarios, the challenge lies in tackling multiple design matrices from different environments, as well as anisotropic covariances rising from distribution shifts. To address these, we introduce multi-environment GLAMP, an approximate message passing framework that applies under distribution shift settings. Multi-enrivonment GLAMP tracks iterates from different environments simultaneously, while allowing non-separable denoising functions. Remarkably, this combination enables one to track properties of estimators under distribution shifts. 
We establish rigorous state evolution results for multi-environment GLAMP and show how it can be used to derive precise asymptotic risks for transfer learning estimators.

To preview the utility of our framework, we describe our basic setting: we observe $E$ datasets $(y_e,X_e)$ with $y_e \in \mathbb{R}^{n_e}, X_e \in \mathbb{R}^{n_e \times p}$ from potentially different distributions. Each dataset---referred to as an environment---follows a regression model with its own coefficient vector $\beta_e$ (see Assumption~\ref{asp:prelim.transfer.model}). The goal in transfer learning is to leverage information across these environments to improve the estimation of regression coefficients in a specific environment of interest. Multiple approaches have been proposed for this purpose, and we demonstrate the utility of multi-environment GLAMP by analyzing three representative estimators: (i) Stacked Lasso \citep{bastani2021predicting}, which serves as an early fusion method by combining all datasets at the start of learning. It solves an empirical risk minimization problem that takes a weighted aggregation of datasets with a usual lasso penalty; (ii) Model Averaging Estimator \citep{bastani2021predicting}, which represents a late fusion approach by combining information at the end of learning. It computes individual lasso estimators from the environments and then takes their weighted average; (iii) Second-step estimator \citep{li2022transfer}, which follows a pre-training and fine-tuning paradigm. It learns an initial estimator from all but the target environment, and then uses the target data to shrink the final estimator toward this initial one. As we will show in Section \ref{sec:transfer}, multi-environment GLAMP allows characterizing the mean-squared error (MSE) of these estimators by accounting for both model and distribution shifts across environments. It enables a unified analysis of these, and potentially a much wider range of transfer learning methods. We further validate our theoretical results through extensive simulations (see Figures~\ref{fig:stacked_Lasso}, \ref{fig:model_avg_Lasso}, \ref{fig:second_step_Lasso}), demonstrating remarkable accuracy even in finite samples.

The rest of the paper is organized as follows. In Section~\ref{sec:related_work}, we provide a more comprehensive overview of related work in AMP and transfer learning. We formally introduce our multi-environment GLAMP formulation in Section~\ref{sec:multi_glamp}, and demonstrate its applicability for specific estimators in Section~\ref{sec:transfer}. Finally, Section~\ref{sec:discussion} outlines discussions and future directions.

\subsection*{Notations} We summarize below some notations used throughout the manuscript. Let \( \|X\|_2 \) denote the spectral norm of a matrix, and \( \|v\|_2 \) the Euclidean (or \( \ell_2 \)) norm of a vector. More generally, \( \|\cdot\|_q \) denotes the \( \ell_q \) norm of a vector for \( 0 \leq q \leq \infty \). For any set \( A \), let \( |A| \) denote its cardinality. Let \( I_p \) be the \( p \times p \) identity matrix, \( \E[\cdot] \) the expectation operator, and \( \PR[\cdot] \) the probability of an event. For any positive definite matrix \( \Sigma \), define \( \|v\|_{\Sigma} = v^{\top} \Sigma v \) and \( \langle v_1, v_2 \rangle_{\Sigma} = v_1^{\top} \Sigma v_2 \). Moreover, we use \( \lambda_{\max}(\Sigma) \) and \( \lambda_{\min}(\Sigma) \) to denote the largest and smallest eigenvalues of the matrix \( \Sigma \), respectively.
 For \( v \in \mathbb{R}^p \), let \( \diag(v) \in \mathbb{R}^{p \times p} \) denote the diagonal matrix with entries \( [\diag(v)]_{i,j} = v_j \cdot \mathbf{1}\{i = j\} \). For any \( m \in \mathbb{N}_+ \), let \( [m] = \{1, 2, \dots, m\} \). Convergence in probability and in distribution are denoted by \( \xrightarrow{p} \) and \( \xrightarrow{d} \), respectively. We write \( o_{\mathrm{P}}(1) \) for sequences of random variables converging to zero in probability.

\section{Further Related Work}\label{sec:related_work}
\subsection{Approximate Message Passing}
Approximate Message Passing (AMP) was first introduced in \cite{donoho2009message, bayati2011dynamics}. Its origins can be traced back to the TAP equations from statistical physics \citep{thouless1977solution, bolthausen2014iterative}, and to loopy belief propagation from graphical models \citep{montanari2012graphical}. AMP has since achieved notable success in various statistical and applied domains such as linear regression \citep{donoho2009message, bayati2011dynamics}, generalized linear models \citep{rangan2011generalized, sur2019modern,barbier2019optimal}, penalized regression \citep{bayati2011lasso, bayati2013estimating}, low-rank matrix estimation \citep{montanari2015non}, deep learning \citep{pandit2019asymptotics}, and communications \citep{barbier2017approximate}. For a more comprehensive review, see \cite{zdeborova2016statistical,feng2022unifying,montanari2024friendly} and the references therein.

On the technical front, the first rigorous state evolution analysis was provided in \cite{bayati2011dynamics}, relying on a crucial conditioning technique introduced in \cite{bolthausen2014iterative}.
Subsequent generalizations have extended AMP in multiple directions, including extensions to right-rotationally-invariant random matrices \citep{rangan2019vector, ma2017orthogonal,fan2022approximate} and finite sample analyses \citep{rush2018finite,li2022non}. More closely related to our work, Generalized AMP (GAMP) was introduced in \cite{rangan2011generalized, javanmard2013state}, which accommodates iterates of stacked vectors or matrices; Long AMP (LAMP) was proposed in \cite{berthier2020state}, enabling the use of non-separable denoising functions; meanwhile the graph-based AMP proposed in \cite{gerbelot2023graph} combined the strengths of both. This algorithm will be the starting point of our framework, although we do not study graph-based problems in our work. We observe that the algorithm in \cite{gerbelot2023graph} may be viewed as an extension of GAMP and LAMP so as to simultaneously allow matrix-valued iterates and non-separable denoising. In this light, we name it GLAMP.  The focus of this work is to develop a multi-environment version of GLAMP that specifically targets transfer learning.

\subsection{Transfer Learning}\label{subsec:related_transfer}
Substantial prior work exists in transfer learning. For this literature review, we focus specifically on penalized estimators that aggregate data from multiple environments. In Section \ref{sec:transfer}, we demonstrate the utility of our GLAMP framework for three representative estimators: Stacked Lasso, Model Averaging and Second-Step. However, the applicability of GLAMP extends beyond these examples.
For instance, the methods proposed in \cite{li2022transfer,li2023targeting,he2024transfusion}, which are based on similar multi-stage estimation schemes, could be analyzed by adapting our approach for the Second-Step estimator. 
Likewise, the residual-weighting method from \cite{zhao2023residual} can be studied using an analysis similar to that for the Stacked Lasso.
Our framework can also accommodate penalties other than the vanilla Lasso, in the same way that AMP has been used to analyze a wide range of penalized procedures  \cite{bu2020algorithmic,chen2021asymptotic} by selecting appropriate denoising functions $f_t$. This is particularly relevant for the transfer learning methods in \cite{cai2022individual,bellec2021chi,maity2022meta}, which rely on variants of the group Lasso.
Finally, prior work has characterized the precise asymptotic risk of ridge-based transfer learning estimators using random matrix theory \citep{yang2020analysis,tripuraneni2021covariate,song2024generalization,patil2024optimal,mallinar2024minimum}. However, these arguments do not generalize beyond the ridge penalty or to settings where random matrix theory is inapplicable, such as the Lasso-based estimators that we study in this paper.

\section{Multi-environment Generalized Long AMP}\label{sec:multi_glamp}
In this section, we introduce our \emph{Multi-environment GLAMP} framework and relevant convergence results. 

\subsection{Setup}

We first define pseudo-Lipschitz functions, which are used throughout the manuscript.

\begin{definition}[Pseudo-Lipschitz functions]\label{def:prelim.plip}
  For some fixed \( q \in \N_+\) and fixed $k\in\R,\, k\geq 1$:
  \begin{enumerate}[label=\upshape(\alph*),ref=\thedefinition (\alph*)]
    \item\label{def:prelim.plip.vec} We say a class of functions \( \{ \phi _N : (\R^{N})^{q}\to \R\} \) is uniformly pseudo-Lipschitz of order \( k \) iff there exists some \( L \) unrelated to \( N \) such that \( \forall\, x_j,y_j\in\R^{N},j\in[q] \),
    \begin{equation*}
      \abs{\phi _N(x_1,x_2,\cdots ,x_q)-\phi _N(y_1,y_2,\cdots ,y_q)} \leq L \left[ 1 + \sum_{j\in[q]} \left( \frac{\norm{x_j}_2}{\sqrt{N}} \right)^{k-1} + \sum_{j\in[q]} \left( \frac{\norm{y_j}_2}{\sqrt{N}} \right)^{k-1} \right] \frac{\sum_{j\in[q]} \norm{x_j-y_j}_2}{\sqrt{N}}.
    \end{equation*}
    \item\label{def:prelim.plip.mat} We say a class of functions \( \{ \phi _N : \R^{N\times q}\to \R\} \) is uniformly pseudo-Lipschitz of order \( k \) iff there exists some \( L \) unrelated to \( N \) such that \( \forall\, X, Y\in\R^{N\times q} \) with spectral norms \( \norm{X}_2, \norm{Y}_2,\)
    \begin{equation*}
      \abs{\phi _N(X)-\phi _N(Y)} \leq L \left[ 1 + \left( \frac{\norm{X}_2}{\sqrt{N}} \right)^{k-1} + \left( \frac{\norm{Y}_2}{\sqrt{N}} \right)^{k-1} \right] \frac{\norm{X-Y}_2}{\sqrt{N}}.
    \end{equation*}
  \end{enumerate}
\end{definition}

As a remark, Definition~\ref{def:prelim.plip.vec}--with a tuple of vectors as input--and \ref{def:prelim.plip.mat}--with a matrix as input--are compatible. That is, if we write \( X=[x_1|x_2|\cdots |x_q],Y=[y_1|y_2|\cdots |y_q],X,Y\in\R^{N\times q} \), then for any function sequence \( \{\phi _N : (\R^{N})^{q}\to \R\} \) that is pseudo-Lipschitz by Definition~\ref{def:prelim.plip.vec}, it is still uniformly pseudo-Lipschitz by Definition~\ref{def:prelim.plip.mat} when viewed as functions from \( \R^{N\times q}\to \R \), and vice versa. This is because with fixed \( q\in\N_+ \), we have the matrix norm inequalities \( \norm{X}_2^{2} \leq \norm{X}_F^{2} = \sum_{i\in[q]}\norm{x_i}_2^{2} \leq q\cdot \norm{X}_2^{2} \).

Now we are ready to officially define a multi-environment GLAMP framework.

\subsection{Multi-environment GLAMP}\label{sec:multienv}
The multi-environment GLAMP is an AMP framework that allows matrix-valued iterates and non-separable non-linearities. Moreover, it explicitly specifies a multi-environment formulation via a partition structure.

\begin{definition}[Multi-environment GLAMP instance]\label{def:GLAMP.asym.multienv}
Suppose there are \( E \in \mathbb{N}_+ \) environments. Let \( q \) be a fixed and known positive integer with \( q = \sum_{e \in [E]} q_e \). Consider the asymptotic regime where 
\[
\lim_{p \to \infty} \frac{p}{n_e} \to \kappa_e \in (0, +\infty), \quad \forall\, e \in [E].
\]
A \emph{multi-environment GLAMP instance} is a tuple
\[
\left( \bigcup_{e \in [E]} \left\{ X_e, (\Psi_e^t)_{t \in \mathbb{N}}, (\eta_e^t)_{t \in \mathbb{N}}, V_e^0, W_e \right\} \right) \cup \{ B \},
\]
with the following components:
\begin{itemize}
  \item For each \( e \in [E] \), \( X_e \in \mathbb{R}^{n_e \times p} \) is an asymmetric random design matrix;
  \item \( B \in \mathbb{R}^{p \times q} \), and \( W_e \in \mathbb{R}^{n_e \times q} \) encode information that the GLAMP iterates can use;
  \item \( (\Psi_e^t)_{t \in \mathbb{N}} \) is a sequence of mappings, with 
  \[
  \Psi_e^t : \mathbb{R}^{n_e \times q_e} \times \mathbb{R}^{n_e \times q} \to \mathbb{R}^{n_e \times q_e};
  \]
  \item \( (\eta_e^t)_{t \in \mathbb{N}} \) is a sequence of mappings, with
  \[
  \eta_e^t : \left( \times_{e \in [E]} \mathbb{R}^{p \times q_e} \right) \times \mathbb{R}^{p \times q} \to \mathbb{R}^{p \times q_e};
  \]
  \item \( V_e^0 \in \mathbb{R}^{p \times q_e} \) is the initial condition of the asymmetric GLAMP iterates.
\end{itemize}
\end{definition}

This tuple fully determines the multi-environment GLAMP framework, and all subsequent iteration and convergence notations are defined based on this instance. Moreover, when \( V_e^0 \), \( W_e \), and \( B \) are clear in a particular context, we sometimes only use $( \bigcup_{e \in [E]} \{ X_e, (\Psi_e^t)_{t \in \mathbb{N}}, (\eta_e^t)_{t \in \mathbb{N}}\} )$ to denote the instance. Similarly, as we always take $W_e$ and $B$ as the last argument of $\Psi_e^t$ and $\eta_e^0$, respectively, we may ignore them in the formulations.

We list the following assumptions accompanying the multi-environment GLAMP instance:

\begin{assumption}[Multi-environment GLAMP assumptions]\label{asp:GLAMP.asym.multienv.converge}
Consider the GLAMP instance in Definition~\ref{def:GLAMP.asym.multienv}. We impose the following assumptions:

\begin{enumerate}[label=\upshape(\alph*),ref=\theassumption~(\alph*)]

\item For all \( e \in [E] \), the matrix \( X_e \in \mathbb{R}^{n_e \times p} \) has i.i.d.\ entries drawn from \( \mathcal{N}(0, 1/n_e) \).

\item \( \limsup_{p \to \infty} \frac{1}{\sqrt{p}} \|B\|_2 < \infty \). For all \( e \in [E] \), \( \limsup_{n_e \to \infty} \frac{1}{\sqrt{n_e}} \|W_e\|_2 < \infty \).

\item For all \( e \in [E] \), the initialization \( V_e^0 \in \mathbb{R}^{p \times q_e} \) satisfies \( \limsup_{p \to \infty} \frac{1}{p} \|V_e^0\|_2^2 < \infty \).

\item For all \( e \in [E] \) and all \( t \in \mathbb{N} \), the mapping \( \Psi_e^t \) is uniformly Lipschitz in \( n_e \), and \( \eta_e^t \) is uniformly Lipschitz in \( p \). Moreover,
\[
\limsup_{n_e \to \infty} \frac{1}{\sqrt{n_e}} \| \Psi_e^t(0, 0) \|_2 < \infty, \quad
\limsup_{p \to \infty} \frac{1}{\sqrt{p}} \| \eta_e^t(0, \dots, 0) \|_2 < \infty.
\]

\item \label{asp:GLAMP.multi.converge,init} For all \( e \in [E] \),
\[
\frac{1}{p} \left( \eta_e^0(V_1^0, \dots, V_E^0) \right)^\top \eta_e^0(V_1^0, \dots, V_E^0) \to \Sigma^{0}_{(R,e)}
\]
for some constant matrix \( \Sigma^{0}_{(R,e)} \in \mathbb{R}^{q_e \times q_e} \).

\item For all \( e \in [E] \), for any constant covariance matrices \( \Sigma_1, \Sigma_2 \in \mathbb{R}^{q_e \times q_e} \), and for all \( s, r \geq 0 \),
\[
\lim_{n_e \to \infty} \frac{1}{n_e} \mathbb{E} \left[ \Psi_e^s(Z_1, W_e)^\top \Psi_e^r(Z_2, W_e) \right]
\]
exists, where \( Z_1, Z_2 \in \mathbb{R}^{n_e \times q_e} \) have i.i.d.\ rows drawn from \( \mathcal{N}(0, \Sigma_1) \) and \( \mathcal{N}(0, \Sigma_2) \), respectively.

\item For all \( e \in [E] \), for any constant covariance matrices \( \Sigma_1, \Sigma_2 \in \mathbb{R}^{q \times q} \), and for all \( s, r \geq 1 \),
\[
\lim_{p \to \infty} \frac{1}{p} \mathbb{E} \left[ \eta_e^s(Z_1, B)^\top \eta_e^r(Z_2, B) \right], \quad
\lim_{p \to \infty} \frac{1}{p} \mathbb{E} \left[ \eta_e^0(V^0, B)^\top \eta_e^r(Z_2, B) \right]
\]
exist, where \( Z_1, Z_2 \in \mathbb{R}^{p \times q} \) have i.i.d.\ rows drawn from \( \mathcal{N}(0, \Sigma_1) \) and \( \mathcal{N}(0, \Sigma_2) \), respectively.

\end{enumerate}
\end{assumption}

Given the instance in Definition~\ref{def:GLAMP.asym.multienv}, the multi-environment GLAMP evolves according to the iterates defined below:

\begin{definition}[Multi-environment GLAMP iterates]\label{def:GLAMP.asym.multienv.iter}
The multi-environment GLAMP iterates for an instance defined in Definition~\ref{def:GLAMP.asym.multienv}, under Assumption~\ref{asp:GLAMP.asym.multienv.converge}, are the alternating sequence of matrices
\[
\left\{ \left( V_1^t, \dots, V_E^t, R^t \right) : t \in \mathbb{N} \right\},
\]
initialized with \( R_e^0 = X_e \eta_e^0(V_e^0, B) \) for each \( e \in [E] \). For \( t \in \mathbb{N_+} \), the updates are defined as follows:
\begin{equation}\label{eq:asym.multienv}
\begin{alignedat}{2}
& V_e^t = \frac{1}{\kappa_e} X_e^\top \Psi_e^{t-1}(R_e^{t-1}, W_e) 
- \frac{1}{\kappa_e} \eta_e^{t-1}(V_1^{t-1}, \dots, V_E^{t-1}, B) \left[ D_{(R,e)}^{t-1} \right]^\top,\\
 & D_{(R,e)}^{t-1} = \frac{1}{n_e} \mathbb{E} \left[ \sum_{j \in [n_e]} 
\left( \frac{\partial (\Psi_e^{t-1})_{j,\cdot}}{\partial (R_e^{t-1})_{j,\cdot}} \right) (Z_{(R,e)}^{t-1}) \right], \\
& R_e^t = X_e \eta_e^t(V_1^t, \dots, V_E^t, B) 
- \Psi_e^{t-1}(R_e^{t-1}, W_e) \left[ D_{(V,e)}^t \right]^\top,\\
& D_{(V,e)}^t = \frac{1}{p} \mathbb{E} \left[ \sum_{j \in [p]} 
\left( \frac{\partial (\eta_e^t)_{j,\cdot}}{\partial (V_e^t)_{j,\cdot}} \right)(Z_{(V,e)}^t) \right],
\end{alignedat}
\end{equation}
where the Gaussian matrix \( Z_{(R,e)}^{t-1} \in \mathbb{R}^{n_e \times q_e} \) has i.i.d. rows drawn from \( \mathcal{N}\left(0, \kappa_e \Sigma^{t-1}_{(R,e)}\right) \) and the matrix \( Z_{(V,e)}^t \in \mathbb{R}^{p \times q} \) has i.i.d. rows drawn from \(\mathcal{N} \left( 0, \frac{1}{\kappa_e^2} \Sigma^t_{(V,e)} \right)\), with $\Sigma^{t}_{(R,e)}$ and $\Sigma^t_{(V,e)}$ to be defined.
\end{definition}

Accompanying the main GLAMP iterates, we also define iterates on $\Sigma^{t}_{(R,e)}$ and $\Sigma^t_{(V,e)}$ as follows:

\begin{definition}[Multi-environment GLAMP state evolutions]\label{def:GLAMP.asym.multienv.state.evo}
Starting from the initial covariance matrices \( \{ \Sigma^0_{(R,e)} : e \in [E] \} \) as specified in Assumption~\ref{asp:GLAMP.multi.converge,init}, we define a sequence of matrices \( \{ \Sigma^t_{(V,e)}, \Sigma^t_{(R,e)} : t \in \mathbb{N_+},\, e \in [E] \} \) as follows: For each \( e \in [E] \) and \( t \in \mathbb{N_+} \),
\begin{equation*}
\begin{aligned}
\Sigma^t_{(V,e)} &= \lim_{n_e \to \infty} \frac{1}{n_e} \mathbb{E} \left[
\Psi_e^{t-1}(Z_{(R,e)}^{t-1}, W_e)^\top 
\Psi_e^{t-1}(Z_{(R,e)}^{t-1}, W_e)
\right]\\
\Sigma^t_{(R,e)} &= \lim_{p \to \infty} \frac{1}{p} \mathbb{E} \left[
\eta_e^t(Z_{(V,1)}^t, \dots, Z_{(V,E)}^t, B)^\top 
\eta_e^t(Z_{(V,1)}^t, \dots, Z_{(V,E)}^t, B)
\right],
\end{aligned}
\end{equation*}
where $Z_{(R,e)}^{t-1}$ and \( Z^t_{(V,e)}\) are defined above.
\end{definition}

The state evolutions define sequences of covariance matrices that the multi-environment GLAMP iterates converge to, as will be formalized in our main result below.

\subsection{Main Result}

\begin{theorem}[Multi-environment GLAMP convergence]\label{cor:asym.multienv}
For the iterates defined in Definition~\ref{def:GLAMP.asym.multienv.iter} under Assumption~\ref{asp:GLAMP.asym.multienv.converge}, the following convergence results hold:

\begin{itemize}
\item For any \( t \in \mathbb{N_+} \) and any order-\( k \) pseudo-Lipschitz function \( \phi : \left( \prod_{e \in [E]} \mathbb{R}^{p \times q_e} \right)^2 \times \mathbb{R}^{p \times q} \to \mathbb{R} \), with \( k \geq 1 \), we have
\[
\phi(V_1^t, \dots, V_E^t, V_1^0, \dots, V_E^0, B) 
= \mathbb{E} \left[ 
\phi(Z_{(V,1)}^t, \dots, Z_{(V,E)}^t, V_1^0, \dots, V_E^0, B) 
\right] + \smalop,
\]
where for each \( e \in [E] \), the matrix \( Z_{(V,e)}^t \in \mathbb{R}^{p \times q_e} \) has i.i.d.\ rows sampled from 
$\mathcal{N}\left(0, \frac{1}{\kappa_e^2} \Sigma_{(V,e)}^t \right),$
as defined in Definition~\ref{def:GLAMP.asym.multienv.state.evo}, and the \( Z_{(V,e)}^t \) are mutually independent across \( e \in [E] \).

\item For any \( t \in \mathbb{N} \), \( e \in [E] \), and any order-\( k \) pseudo-Lipschitz function \( \phi : (\mathbb{R}^{n_e \times q_e})^2 \to \mathbb{R} \), with \( k \geq 1 \), we have
\[
\phi(R_e^t, W_e) = \mathbb{E} \left[ \phi(Z_{(R,e)}^t, W_e) \right] + \smalop,
\]
where \( Z_{(R,e)}^t \in \mathbb{R}^{n_e \times q_e} \) has i.i.d.\ rows drawn from 
$\mathcal{N}(0, \kappa_e \Sigma_{(R,e)}^t),$
as defined in Definition~\ref{def:GLAMP.asym.multienv.state.evo}.
\end{itemize}
\end{theorem}

Theorem~\ref{cor:asym.multienv} is the main result of the paper, and is proved in Section~\ref{sec:apd.pf.glamp.asym.multienv}.

The multi-environment GLAMP framework is typically linked to statistical estimators in the following manner: $\{X_e\}_{e\in E}$ are related to the design matrices in multiple environments, $\{W_e\}_{e\in E}$ to the response variables, and $(\Psi_e^t,\eta_e^t)_{e\in \E, t \in \mathbb{N}}$ to the specific multi-environment estimator of interest. With proper initializations $\{V_e^0\}_{e\in E}$ and under suitable conditions, taking $t\rightarrow \infty$ additionally above the $N \rightarrow \infty$ asymptotics in Theorem \ref{cor:asym.multienv}, the GLAMP instance converges to some $V^\infty, R^\infty$ that corresponds to our statistical estimator of interest. Below we showcase the linkage with specific examples.

\section{Applications to Transfer Learning}\label{sec:transfer}
We demonstrate the utility of multi-environment GLAMP by investigating a class of Lasso-based estimators introduced in the transfer learning literature. 
We first define our transfer learning setting formally. We consider the problem where at training time data from multiple distributions are available---we call these the environments.\footnote{the environments available during training are often known as sources, and the environment of interest is known as the target, but we follow the naming convention from \cite{duchi2024predictive}.} 

\begin{assumption}[Transfer learning model]\label{asp:prelim.transfer.model} Suppose there are $E \in \N_+$ environments, and consider for each $e \in E$, a sequence of problem instances $\{y_e(n_e),X_e(n_e),\beta_e(n_e),w_e(n_e) \}$ with \( y_e(n_e) \in \R^{n_e} \), \( X_{e}(n_e)\in\R^{n_e\times p} \), \( \beta _e(n_e) \in\R^{p}, w_e\in\R^{n_e} \)  the response vector, the design matrix, the true regression coefficient and the noise, respectively. Further, we assume that $p, n_e \rightarrow \infty$, \( \lim_{p\to \infty}p/n_e = \kappa _e \in (0,+\infty) \). In the sequel we will drop the dependence on $n_e$ for simplicity. 
  \begin{enumerate}[label=\upshape(\alph*),ref=\theassumption (\alph*)]
    \item\label{asp:prelim.transfer.model.nonasym} \( \forall\,e\in[E] \), we assume a linear data generative model: \( y_e = X_e \beta _e + w_e \), where \( y_e \in \R^{n_e} \), \( X_{e}\in\R^{n_e\times p} \), \( \beta _e\in\R^{p}, w_e\in\R^{n_e} \).
    \item\label{asp:prelim.transfer.model.asym.X} \( X_{e} = \ix{e} \Sigma _{e}^{1/2} \), where \((\ix{e})_{i,j}\iidsim \Nscr(0,\frac{1}{n_e})\) are also independent across the environments. $\{\Sigma_e\}$ are deterministic, and \( \exists \, c_1>1 \) s.t. \( c_1^{-1} \leq \inf_{p}\{\lambda _{\min }(\Sigma _{e})\} \leq  \sup_{p}\{\lambda _{\max }(\Sigma _{e}) \}\leq c_1,\, \forall\,e\in[E] \). 

    \item \label{asp:prelim.transfer.model.asym.beta} \( \sup_{e\in[E]}\limsup_{p}\frac{1}{p}\norm{\beta _e}_2^{2} \leq +\infty \), and \( \forall\,e\in[E] \), \( \lim_p \frac{1}{p} \beta _e^{\T} \Sigma _{e} \beta _e \) exists. 
    \item \label{asp:prelim.transfer.model.asym.w}The empirical distribution of \( w_e \) converges in 2-Wasserstein distance to \( W_e \), with \( \E[W_e] = 0 \), \( \E[W_e^{2}] < +\infty \). 
  \end{enumerate}
\end{assumption}

Some comments are in order regarding the assumptions. Although Assumption~\ref{asp:prelim.transfer.model.nonasym} operates under linear models, we expect our multi-environment GLAMP framework to be useful beyond this stylistic setting in the same way as GAMP served useful for non-linear models, c.f.~\cite{sur2019modern,barbier2019optimal}. 
In Assumption \ref{asp:prelim.transfer.model.asym.w}, convergence in 2-Wasserstein distance is equivalent to requiring that the empirical distribution of \( w_e \) converges to \( W_e \) not only in distribution, but also in the first and second moments. 
In our formulation, as detailed in Assumption~\ref{asp:prelim.transfer.model}, we directly work with a deterministic sequence of \( \{ w_e \}  \)  and \( \{ \beta_e \}  \). If they are sampled from a distribution, it must be independent from the design matrices or \( \{ \ix{e}:e\in[E] \}  \), so that conditioning on them does not effect the Gaussianity required for GLAMP. Beyond this independence requirement, GLAMP is sufficiently general to accommodate any correlation or dependence among these sequences, provided that the realized sequences for \( \{ w_e \}  \)  and \( \{ \beta_e \}  \) satisfy Assumption~\ref{asp:prelim.transfer.model}(c), (d). 

With the transfer learning setting defined, we investigate three estimators from the previous literature, and use the multi-environment GLAMP framework to characterize their asymptotic risk under suitable conditions. For the sake of conciseness, we sometimes omit the "multi-environment" prefix and only uses GLAMP to denote our framework.

\subsection{The Stacked Lasso}\label{sec:stack}
We begin with the stacked Lasso estimator--also known as the weighted loss estimator \citep{bastani2021predicting}--which is constructed by running the Lasso on the \( E \) datasets stacked together, with different weights. In the hope that the analysis of the stacked Lasso will serve as a foundation for other estimators, we study a slightly generalized version that allows heterogeneous penalties across variables. 
Suppose we are given a deterministic sequence \( \{ \vec{\lambda }: \vec{\lambda } \in\R^{p},p\in\N \}  \), and the penalty is defined as \( \norm{\diag(\vec{\lambda } )b}_1 = \sum_{j\in[p]} \vec{\lambda }_j \abs{b} _j  \) for any candidate \( b\in\R^{p} \). The estimator \( \stac \) is then defined as:
\begin{equation}\label{eq:def.stack.los}
  \stac = \argmin_{b\in\R^{p}} \Lscr(b),\quad \Lscr(b) = \frac{1}{2}\sum_{e\in[E]}\pi _e \norm*{y_e - \ix{e}\Sigma _e^{1/2} b}_2^{2} +  \norm*{\diag(\vec{\lambda } )b}_1.
\end{equation}
The coefficients \( \{ \pi _e:e\in[E] \}  \) represent the weights assigned to each environment. We naturally assume that \( \pi _e \geq 0,\,\forall\,e\in[E] \) and \( \max\{\pi _e : e\in[E]\} >0 \). 

The stacked Lasso estimator in \eqref{eq:def.stack.los} depends on the random quantities $\{X_e\},\{y_e\}$. To capture its limiting behavior using GLAMP, we define a ``fixed-design" counterpart function \( \eta:(\R^{p})^{E}\to\R^{p} \) as follows:

\begin{equation}\label{eq:def.stack.et}
\begin{aligned}
  \forall\,v_1,\cdots ,v_E\in\R^{p},\, 
  \eta(v_1,\cdots ,v_E) =&~ \argmin_{b\in\R^{p}} \left\{ \frac{1}{2}\sum_{e} \varpi _e \norm*{\Sigma _{e}^{1/2}b - (v_e + \Sigma _{e}^{1/2}\beta _e)}_2^{2} + \theta \cdot  \norm*{\diag(\vec{\lambda } )b}_1 \right\}.
\end{aligned}
\end{equation}

By definition, \( \eta \) depends on the covariances \( \Sigma _{e} \), the true signals \( \{ \beta _e \}  \), as well as $E+1$ deterministic (but yet to be defined)  parameters \( \theta \) and \( \{ \varpi _e \}  \). 

We need the following weak assumptions for the regularity of the stacked Lasso estimator:

\begin{assumption}[Stacked Lasso GLAMP regularity]\label{asp:stack.for.glamp}
Consider the transfer learning setting in Assumption~\ref{asp:prelim.transfer.model}, the stacked Lasso estimator in \eqref{eq:def.stack.los} and its fixed-design counterpart in \eqref{eq:def.stack.et}.
\begin{enumerate}[label=\upshape(\alph*),ref=\theassumption (\alph*)]
  \item\label{asp:stack.for.glamp.lambd} \( \liminf_p \min\{ \vec{\lambda } _{j}:j\in[p]\} > 0 \), and \( \limsup_{p} \frac{1}{p} \norm{\vec{\lambda } }_2^{2} < +\infty \).
  \item\label{asp:stack.for.glamp.limits} Define the following intermediate quantity 
  \begin{equation}\label{eq:def.stacl.et.bar}
    \overline{\eta } = \overline{\eta } (\tau _1,\tau_2,\cdots ,\tau _E) = \eta(\tau_1 Z_1,\tau_2 Z_2,\cdots ,\tau _E Z_E), 
  \end{equation}
  for \( \tau_1,\tau_2,\cdots ,\tau _E >0 \) and \( Z_1,Z_2,\cdots ,Z_E\iidsim \Nscr(0,I_p) \). 
  We assume \( \forall\,\tau_1,\tau_2,\cdots ,\tau _E >0 \), \(\theta>0 \), \( \varpi _e>0 \) and \( \sum_{e}\varpi _e=1 \), the following limits exist: \( \lim_p \frac{1}{p}\E [\norm{\overline{\eta }}_{\Sigma _{e}}^{2}] \), \( \lim_{p} \frac{1}{p}\E[\inprod{\beta _e, \overline{\eta } }_{\Sigma _{e}}] \), \( \forall\,e\in[E] \). 
  \item\label{asp:stack.for.glamp.onsager} Define the following intermediate quantity 
  \begin{equation}\label{eq:def.stacl.E}
    \mathcal{E} ^{(p)}(\tau_1,\cdots ,\tau _E)= \frac{1}{p}\E \left[ \min _{b \in\R^{p}}  \left\{ \frac{1}{2}\sum_{e\in[E]} \varpi _e\norm*{\Sigma _{e}^{1/2}b - (\tau _e Z_e + \Sigma _e^{1/2}\beta _e)}_2^{2} + \theta \norm*{\diag{(\vec{\lambda }) }b}_1 \right\}  \right].
  \end{equation}
  We assume \( \forall\,\tau_1,\tau_2,\cdots ,\tau _E >0\), \(\theta>0 \), \( \varpi _e>0 \) and \( \sum_{e}\varpi _e=1 \), there exists a function \( \Escr(\tau_1,\cdots ,\tau _E) \in C^{1} \), s.t. \( \lim_{p}\Escr^{(p)} = \Escr \), \( \lim_{p}\frac{\partial }{\partial \tau _e}\Escr^{(p)} = \frac{\partial }{\partial \tau _e}\Escr  \). 
\end{enumerate}
\end{assumption}
Assumption~\ref{asp:stack.for.glamp.limits} makes sure the GLAMP iterates corresponding to the stacked Lasso have well-defined asymptotic limit, and Assumption~\ref{asp:stack.for.glamp.onsager} guarantees the well-definedness of the Onsager term, as demonstrated by the following result:

\begin{proposition}\label{prop:stack.onsager}
  Under Assumptions~\ref{asp:prelim.transfer.model}, and \ref{asp:stack.for.glamp},  
  denote \( \overline{\delta } _e^{{(p)}} := \frac{1}{\tau _e p} \E\left[ Z_e^{\T} \Sigma _e^{1/2} \overline{\eta } (\tau_1,\tau_2,\cdots ,\tau _E) \right] \), and let \( S\subset [p] \) be the support of \( \overline{\eta } \in\R^{p} \). Then we have the following equivalent forms of \( \overline{\delta } ^{{(p)}} _e \):
  \begin{equation}
\begin{aligned}\label{eq:prop.stack.forms.of.delta}
  \overline{\delta } _e^{{(p)}} =&~ \frac{1}{p}\sum_{j\in[p]}\E \left[ \frac{\partial (\Sigma _e^{1/2} {\eta } )_j}{\partial (\tau _e Z_e)_j} (\tau_1 Z_1, \tau_2 Z_2,\cdots ,\tau _E Z_E) \right] \\
  =&~ \frac{1}{p}\E \left[ \trace\left(\frac{\partial \eta}{\partial (\tau _e Z_e)}(\tau_1 Z_1, \tau_2 Z_2,\cdots ,\tau _E Z_E) \Sigma _e^{1/2} \right)\right] \\
  =&~  \varpi _e\cdot \frac{1}{p} \E \left[ \trace\bigg(\bigg[ \sum_{e\in[E]} \varpi _e (\Sigma _{e})_{S,S} \bigg]^{-1} (\Sigma _e)_{S,S} \bigg) \right].
\end{aligned}
  \end{equation}
  Moreover, \( \overline{\delta } _e := \lim_{p}\overline{\delta } _e^{{(p)}} =  \frac{1}{\tau _e}\lim_{p} \frac{1}{p}\E\left[ Z_e^{\T} \Sigma _e^{1/2} \overline{\eta } (\tau_1,\tau_2,\cdots ,\tau _E) \right] \) is a well-defined limit. We have \( \frac{1}{p}\E[\norm{\overline{\eta }(\tau_1,\tau_2,\cdots ,\tau _E) }_0]  = \sum_{e\in[E]} \overline{\delta }^{(p)} _e \) and \( \lim_p \frac{1}{p}\E[\norm{\overline{\eta }(\tau_1,\tau_2,\cdots ,\tau _E) }_0] =  \sum_{e\in[E]} \overline{\delta } _e\).
\end{proposition}
Proposition~\ref{prop:stack.onsager} is proved in Appendix Section~\ref{sec:apd.prelim.stack.prop}. 

Under Assumptions~\ref{asp:prelim.transfer.model} and \ref{asp:stack.for.glamp}, we can track the asymptotic convergence of the corresponding GLAMP iterates (see Lemma~\ref{lem:apd.stack.glamp}). We are left with one final piece: ensuring that the corresponding GLAMP indeed tracks the behavior of the stacked Lasso estimator. For that, we require the following assumption:

\begin{assumption}[Stacked Lasso GLAMP fixed point]\label{asp:stack.converge} Consider the transfer learning setting in Assumption~\ref{asp:prelim.transfer.model}, the stacked Lasso estimator in \eqref{eq:def.stack.los} and its fixed-design counterpart in \eqref{eq:def.stack.et}.
  \begin{enumerate}[label=\upshape(\alph*),ref=\theassumption (\alph*)]
    \item\label{asp:stack.converge.fixed.point} There exist a \( (2E+1) \)-tuple \( (\tau _1^{*},\cdots ,\tau _E^{*}, \varpi^* _1,\cdots ,\varpi^* _E,\theta^*) \) with \( \theta^*>0 \), \( \tau ^* _e>0 \), and \( \varpi^* _e\geq 0\ \forall\,e\in[E] \), which solves the following system of $(2E+1)$ equations:
    \begin{equation}
      \begin{cases}\label{eq:stack.fixed.point}
          \tau _e^{2} = \E[W_e^{2}] + \kappa _e \lim_{p} \frac{1}{p}\E\left[\norm*{(\overline{\eta } -\beta _e)}_{\Sigma _{e}}^{2}\right]
          ,\ \forall\,e\in[E],\\
          \pi _e \theta = \varpi _e/(1-\kappa _e \overline{\delta } _e),\ \forall\,e\in[E],\\
          \sum_{e\in[E]} \varpi _e = 1,
      \end{cases}
    \end{equation}
    where \( \overline{\eta }  \) from Equation~\eqref{eq:def.stacl.et.bar} and \( \{ \overline{\delta }  _e \}  \) from Proposition~\ref{prop:stack.onsager} are all functions of the unknown variables. 
    \item\label{asp:stack.converge.cauchy} For the \( (\tau _1^{*},\cdots ,\tau _E^{*}, \varpi^* _1,\cdots ,\varpi^* _E,\theta^*) \) assumed to exist in Assumption~\ref{asp:stack.converge.fixed.point}, we define a function \( H^{{(p)}}:[0,1]^{E}\to[0,1]^{E} \), whose \( e \)-th element \( H_e^{(p)} \) is defined as: 
    \begin{equation}
\begin{aligned}\label{eq:stack.def.He}
  \forall\,\rho \in [0,1]^{E},\, &~H^{(p)}_e(\rho) = H^{(p)}_e(\rho_1,\rho_2,\cdots ,\rho _E) =  \frac{1}{(\tau _e^{*})^{2}}\left\{\E[W_e^{2}] +  \frac{\kappa_e}{p}\E\left[ \inprod{\overline{\eta }_{{(1)}}  - \beta_e,\overline{\eta }_{{(2)}}  - \beta_e}_{\Sigma _e} \right] \right\},
\end{aligned}
    \end{equation}
    where \(  \overline{\eta }_{{(1)}}  = \eta(\tau ^*_1 Z_1,\tau ^*_2 Z_2,\cdots ,\tau ^*_E Z_E) \), \(  \overline{\eta }_{{(2)}}  = \eta(\tau ^*_1 Z_1',\tau ^*_2 Z_2',\cdots ,\tau ^*_E Z_E') \), and the Gaussian variables involved satiesfy: \emph{(i)} \( Z_e,Z'_e \) are jointly normal, \( Z_e,Z'_e\sim \Nscr(0,I_p) \) and \( \cov(Z_e,Z'_e)=\rho_e I_p \); \emph{(ii)} The \( E \) pairs of random variables \( (Z_e,Z'_e), e=1,2,\cdots E  \) are jointly independent. 
    Further, we assume \( H^{(p)} \) converges pointwise to some function \( H:[0,1]^{E}\to[0,1]^{E} \) as \( p\to +\infty \) with \( H\in C^{1} \), and the convergence holds for all the partial derivatives, i.e. \( \lim_p \frac{\partial H_e^{{(p)}}}{\partial \rho _i} = \frac{\partial H_e}{\partial \rho _i} \), \(\forall\, e,i\in[E] \). 
  \end{enumerate}
\end{assumption}
To connect Assumption~\ref{asp:stack.converge} with the stacked Lasso, note that
Assumption~\ref{asp:stack.converge.fixed.point} establishes the existence of a fixed point for the GLAMP iterates that corresponds to the loss function from Equation~\eqref{eq:def.stack.los}, while Assumption~\ref{asp:stack.converge.cauchy} guarantees that the GLAMP iterates  converge to that fixed point. 
Together, these two assumptions guarantee that the GLAMP iterates converge to a fixed point that asymptotically tracks \( \stac \). 

With all assumptions in place, we state our main convergence result for $\stac$.

\begin{theorem}[Stacked Lasso convergence]\label{thm:stack}
  Under Assumptions~\ref{asp:prelim.transfer.model},\ref{asp:stack.for.glamp} and \ref{asp:stack.converge}, suppose we fix \( \tau_e = \tau_e^* \), \( \varpi_e = \varpi_e^* \), and \( \theta = \theta^* \) in the definition of \( \overline{\eta} \).
  Then, for any sequence of order-\( k \) pseudo-Lipschitz functions \( \phi : (\mathbb{R}^p)^{E+1} \to \mathbb{R} \) with \( k \geq 1 \), we have
  \begin{equation*}
    \phi\left(\stac,\beta_1,\cdots ,\beta _E\right) = \E \left[ \phi(\overline{\eta } (\tau_1^{*},\tau _2^{*},\cdots ,\tau _E^{*}),\beta_1,\beta_2,\cdots ,\beta _E) \right] + \smalop.
    \end{equation*}
\end{theorem}
Theorem~\ref{thm:stack} is proved in Appendix Section~\ref{sec:apd.stack}.

\begin{figure}
    \centering
    \includegraphics[width=\linewidth]{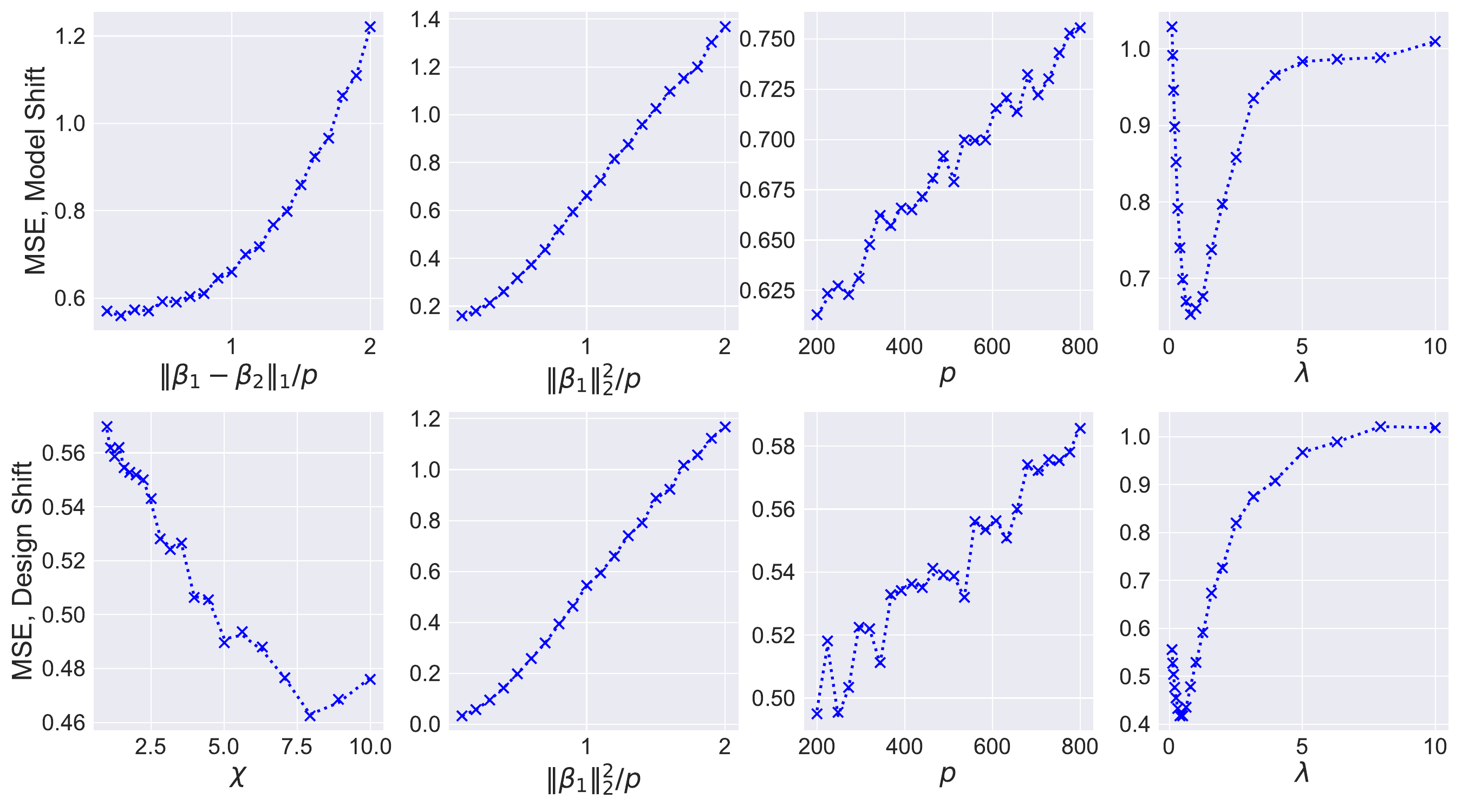}
    \caption{MSE of $\stac$ under different scenarios. The dotted line indicates the theoretical MSE value obtained from our formulae while the cross marks indicate empirical MSE values.  Base setting: $E=2$ (with $1$ as the target environment), $n_1=800$, $n_2=600$, $p=400$, $\Sigma_1=\Sigma_2=I$, $\beta_{1,j}=\beta_{2,j}\iidsim \mathcal{N}(0,\sigma_\beta^2)$ with $\sigma_\beta^2=1$, $w_{e,i}\iidsim \mathcal{N}(0,\sigma_w^2)$, $\sigma_w^2=1$, $\vec{\lambda}_j=1$, $\pi_1=\pi_2=1/2$. For model shift (the first row), $\beta_2=\beta_1 + \tilde\beta$ with $\tilde\beta_j \iidsim \mathcal{N}(0,\tilde\sigma^2\pi/2), \tilde\sigma^2=1$, while $\beta_1$ is taken to be as in the base setting. For design shift (the second row), $\Sigma_2$ differs from the base setting and is taken to be a diagonal matrix with two distinct eigenvalues $\chi,1/\chi$ that have equal multiplicity $p/2$. For each subplot, we change the parameter on the x-axis as shown, fixing others to follow the base setting. Note the remarkable alignment of the theoretical versus empirical values. 
}
    \label{fig:stacked_Lasso}
\end{figure}

We demonstrate the finite-sample accuracy of the asymptotic convergence through numerical examples, with results presented in Figure~\ref{fig:stacked_Lasso}. Specifically, we consider the case with \( E = 2 \) environments. In the base setting, we set the sample sizes to \( n_1 = 800 \), \( n_2 = 600 \), and the feature dimension to \( p = 400 \). The design matrices are isotropic with \( \Sigma_1 = \Sigma_2 = I \).  The signal vectors are generated as \( \beta_{1,j} = \beta_{2,j} \overset{\text{i.i.d.}}{\sim} \mathcal{N}(0, \sigma_\beta^2) \) for all \( j \in [p] \), with \( \sigma_\beta^2 = 1 \), and the noise terms follow \( w_{e,i} \overset{\text{i.i.d.}}{\sim} \mathcal{N}(0, \sigma_w^2) \) for all \( e \in [E] \), \( i \in [n_e] \), with \( \sigma_w^2 = 1 \). We use a homogeneous regularization parameter \( \vec{\lambda} \) with \( \lambda = 1 \), and assign equal weights \( \pi_1 = \pi_2 = 1/2 \). Building on this base case, we consider two types of transfer learning settings:

\begin{itemize}
    \item \textbf{Design/Covariate Shift:} The underlying signal remains the same across environments as in the base setting, but the design matrix covariances differ. Specifically, we fix \( \Sigma_1 = I \) and set \( \Sigma_2 \) to be a diagonal matrix with two distinct eigenvalues \( \chi \) and \( 1/\chi \), each with multiplicity \( p/2 \).
    
    \item \textbf{Model/Signal Shift:} The design matrix covariances remain the same as in the base setting, but the underlying signals differ. Specifically, we set \( \beta_2 = \beta_1 + \tilde{\beta} \), where \( \tilde{\beta}_j \overset{\text{i.i.d.}}{\sim} \mathcal{N}(0, \tilde{\sigma}^2 \pi / 2) \), such that \( \|\tilde{\beta}\|_1 / p \to \tilde{\sigma}^2 \), while $\beta_1$ is taken to be as in the base setting.
\end{itemize}

In each setting, we vary different hyperparameter values and record the mean squared error (MSE) of the stacked Lasso estimator, defined as \( \|\stac - \beta_1\|_2^2 \). Since our formulation assumes fixed signals, we draw a single realization of the signal vector for each experiment.\footnote{which results in mild local non-monotonicity in the plot.} We then perform \( R = 20 \) independent realizations of the random design matrices \( \{X_e\} \) and noise vectors \( \{w_e\} \), along with \( R = 400 \) realizations of \( \{Z_e\} \) from Equation~\eqref{eq:def.stacl.et.bar}, to approximate the expectation.

Figure~\ref{fig:stacked_Lasso} confirms the remarkable finite-sample accuracy of Theorem~\ref{thm:stack}, while also revealing several interesting trends. For example, the U-shaped curve with respect to changes in \( \lambda \) suggests the potential effectiveness of adaptive tuning procedures. Additionally, the observation that the MSE decreases as \( \chi \) increases indicates that heterogeneity across environments can be beneficial—an effect also noted by \cite{song2024generalization}. A more comprehensive study of the MSE and generalization error—which could provide further insights into estimator selection, data pooling strategies, and hyperparameter tuning—is left for future work.

\subsection{The Model Averaging Estimator}\label{sec:aver}

Another straightforward approach to aggregating information from \( E \) environments is to post-combine  \( E \) individual estimators via a weighted average. This method is known as the \emph{model averaging estimator} \citep{bastani2021predicting}:
\begin{align}
  \indiv{e} =&~ \argmin_{b \in \mathbb{R}^{p}} L_e(b), \quad L_e(b) = \frac{1}{2} \left\| y_e - \ix{e} \Sigma_e^{1/2} b \right\|_2^2 + \left\| \diag(\vec{\lambda}_e) b \right\|_1, \label{eq:def.aver.los} \\
  \aver =&~ \sum_{e \in [E]} \pi_e \indiv{e}, \quad \text{where } \sum_{e \in [E]} \pi_e = 1. \label{eq:def.aver}
\end{align}
For the single-environment Lasso estimator \( \indiv{e} \) defined in Equation~\eqref{eq:def.aver.los}, we adopt a heterogeneous \( \ell_1 \) penalty for generality, in the form \( \left\| \diag(\vec{\lambda}_e) b \right\|_1 = \sum_{j \in [p]} (\vec{\lambda}_e)_j |b_j| \). Comparing \( \indiv{e} \) with the stacked estimator \( \stac \) from Equation~\eqref{eq:def.stack.los}, we observe that \( \indiv{e} \) is simply a special case of \( \stac \) when \( E = 1 \), and is therefore covered by the theoretical results for \( \stac \). For the model averaging estimator \( \aver \), which aggregates the individual estimators \( \{ \indiv{e} : e \in [E] \} \) via a weighted sum, the main challenge lies in characterizing the joint convergence of the \( E \) individual estimators.
However, it turns out that the joint convergence of the \( E \) individual estimators \( \indiv{e} \) does not require more than the analogs of Assumptions~\ref{asp:stack.for.glamp} and~\ref{asp:stack.converge} applied to the case \( E = 1 \). For clarity, we restate these requirements as separate assumptions in the context of model averaging. We define the fixed-design counterparts \( \eta_e \) in Equation~\eqref{eq:def.aver.et} in the same way we have defined \( \eta \), as follows:
\begin{equation}\label{eq:def.aver.et}
\begin{aligned}
  \forall\,e \in [E],\,\forall\,v_e \in \mathbb{R}^{p}, \quad
  \eta_e(v_e) =&~ \argmin_{b \in \mathbb{R}^{p}} \left\{ \frac{1}{2} \left\| \Sigma_e^{1/2} b - \left(v_e + \Sigma_e^{1/2} \beta_e \right) \right\|_2^2 
  + \theta_e \cdot \left\| \diag(\vec{\lambda}_e) b \right\|_1 \right\}.
\end{aligned}
\end{equation}

\begin{assumption}[Individual Lasso GLAMP regularity]\label{asp:aver.for.glamp}
Consider the transfer learning setting in Assumption~\ref{asp:prelim.transfer.model}, the individual Lasso estimators in Equation~\eqref{eq:def.aver.los}, and their fixed-design counterparts in Equation~\eqref{eq:def.aver.et}.
\begin{enumerate}[label=\upshape(\alph*),ref=\theassumption~(\alph*)]
  \item\label{asp:aver.for.glamp.lambd} 
  The regularization weights satisfy \( \liminf_p \min\{ (\vec{\lambda}_e)_j : j \in [p] \} > 0 \), and \( \limsup_p \frac{1}{p} \| \vec{\lambda}_e \|_2^2 < \infty \).

  \item\label{asp:aver.for.glamp.limits} 
  For each \( e \in [E] \), define the intermediate quantity \( \overline{\eta}_e \) as
  \begin{equation}\label{eq:def.aver.et.bar}
    \overline{\eta}_e = \overline{\eta}_e(\tau_e) := \eta_e(\tau_e Z_e),
  \end{equation}
  where \( \tau_e > 0 \) and \( Z_e \sim \mathcal{N}(0, I_p) \). Let \( Z_1, Z_2, \dots, Z_E \overset{\text{i.i.d.}}{\sim} \mathcal{N}(0, I_p) \), so that the vectors \( \overline{\eta}_1, \overline{\eta}_2, \dots, \overline{\eta}_E \) are independent. For all \( \tau_e > 0 \) and \( \theta_e > 0 \), we assume the following limits exist:
  \[
  \lim_{p \to \infty} \frac{1}{p} \mathbb{E} \left[ \| \overline{\eta}_e \|_{\Sigma_e}^2 \right], \quad
  \lim_{p \to \infty} \frac{1}{p} \mathbb{E} \left[ \langle \beta_e, \overline{\eta}_e \rangle_{\Sigma_e} \right], \quad \forall e \in [E].
  \]

  \item\label{asp:aver.for.glamp.onsager} 
  For each \( e \in [E] \), define the following intermediate quantity:
  \begin{equation}\label{eq:def.aver.E}
    \mathcal{E}_e^{(p)}(\tau_e) := \frac{1}{p} \mathbb{E} \left[ 
    \min_{b \in \mathbb{R}^p} \left\{ 
    \frac{1}{2} \left\| \Sigma_e^{1/2} b - (\tau_e Z_e + \Sigma_e^{1/2} \beta_e) \right\|_2^2 
    + \theta_e \left\| \diag(\vec{\lambda}_e) b \right\|_1 
    \right\} 
    \right].
  \end{equation}
  We assume that \( \forall\,\theta _e >0 \), there exists a function \( \Escr_e(\tau_e) \in C^{1} \), s.t. \( \lim_{p}\Escr_e^{(p)} = \Escr_e \), \( \lim_{p}\frac{\dif }{\dif \tau _e}\Escr_e^{(p)} = \frac{\partial }{\dif \tau _e}\Escr_e  \). 
\end{enumerate}
\end{assumption}

Assumption~\ref{asp:aver.for.glamp} is a special case of Assumption~\ref{asp:stack.for.glamp} with \( E = 1 \). As a corollary of Proposition~\ref{prop:stack.onsager}, we introduce the following quantities associated with the Onsager terms and use Proposition~\ref{prop:stack.onsager} to establish their well-definedness:
\begin{proposition}\label{prop:aver.onsager}
  Under Assumptions~\ref{asp:prelim.transfer.model}, and \ref{asp:aver.for.glamp}, 
  for each $e\in[E]$, we denote \( \overline{\delta } _{\indi,e}^{{(p)}} := \frac{1}{\tau _e p} \E\left[ Z_e^{\T} \Sigma _e^{1/2} \overline{\eta }_e (\tau_e) \right] \), and let \( S_e\subset [p] \) be the support of \( \overline{\eta }_e \in\R^{p} \). Then we have the following equivalent forms of \( \overline{\delta } ^{{(p)}} _{\indi,e} \):
  \begin{equation}
\begin{aligned}\label{eq:prop.aver.forms.of.delta}
  \overline{\delta } _{\indi,e}^{{(p)}} =&~ \frac{1}{p}\sum_{j\in[p]}\E \left[ \frac{\partial (\Sigma _e^{1/2} {\eta }_e )_j}{\partial (\tau _e Z_e)_j} (\tau_e Z_e) \right] 
  = \frac{1}{p}\E \left[ \trace\left(\frac{\partial \eta _e}{\partial (\tau _e Z_e)}(\tau_e Z_e) \Sigma _e^{1/2} \right)\right] = \frac{1}{p} \E \left[ \abs{S_e}  \right].
\end{aligned}
  \end{equation}
  Moreover, \( \overline{\delta } _{\indi,e} = \lim_{p}\overline{\delta } _{\indi,e}^{{(p)}} =  \frac{1}{\tau _e}\lim_{p} \frac{1}{p}\E\left[ Z_e^{\T} \Sigma _e^{1/2} \overline{\eta }_e (\tau_e) \right] \) is a well-defined limit. 
\end{proposition}
Proposition~\ref{prop:aver.onsager} follows directly as a corollary of Proposition~\ref{prop:stack.onsager}, and its proof is therefore omitted. Finally, we present an assumption that ensures the convergence of the GLAMP iterates to the individual Lasso solutions \( \indiv{e} \), analogous to Assumption~\ref{asp:stack.converge}:

\begin{assumption}[Individual Lasso GLAMP fixed point]\label{asp:aver.converge} 
Consider the transfer learning setting in Assumption~\ref{asp:prelim.transfer.model}, the individual Lasso estimators in Equation~\eqref{eq:def.aver.los}, and their fixed-design counterparts in Equation~\eqref{eq:def.aver.et}.
\begin{enumerate}[label=\upshape(\alph*),ref=\theassumption~(\alph*)]
  
  \item\label{asp:aver.converge.fixed.point} 
  For each \( e \in [E] \), assume the existence of a pair \( (\tau^*_{\indi,e}, \theta^*_{\indi,e}) \) with \( \tau^*_{\indi,e} > 0 \), \( \theta^*_{\indi,e} > 0 \), such that setting \( \tau_e = \tau^*_{\indi,e} \) and \( \theta = \theta^*_{\indi,e} \) solves the following fixed-point system:
  \begin{equation}\label{eq:aver.fixed.point}
  \begin{cases}
    \tau_e^2 = \mathbb{E}[W_e^2] + \kappa_e \lim_{p \to \infty} \frac{1}{p} \mathbb{E} \left[ \| \overline{\eta}_e - \beta_e \|_{\Sigma_e}^2 \right], \\
    \theta = 1 / \left(1 - \kappa_e \overline{\delta}_{\indi,e} \right),
  \end{cases}
  \end{equation}
  where \( \overline{\eta}_e \) is defined in Equation~\eqref{eq:def.aver.et.bar} and \( \overline{\delta}_{\indi,e} \) is as defined in Proposition~\ref{prop:aver.onsager}; both depend implicitly on \( \tau_e \) and \( \theta \).

  \item\label{asp:aver.converge.cauchy} 
  For each \( e \in [E] \), consider the pair \( (\tau^*_{\indi,e}, \theta^*_{\indi,e}) \) from part~\ref{asp:aver.converge.fixed.point}, and define the function \( H^{(p)}_{\indi,e} : [0,1] \to [0,1] \) by
  \begin{equation}\label{eq:aver.def.He}
  \begin{aligned}
    H^{(p)}_{\indi,e}(\rho_e) := \frac{1}{(\tau^*_{\indi,e})^2} \left\{ \mathbb{E}[W_e^2] + \frac{\kappa_e}{p} \mathbb{E} \left[ \langle \overline{\eta}_{e,(1)} - \beta_e, \overline{\eta}_{e,(2)} - \beta_e \rangle_{\Sigma_e} \right] \right\},
  \end{aligned}
  \end{equation}
  where \( \overline{\eta}_{e,(1)} := \eta_e(\tau^*_{\indi,e} Z_e) \), \( \overline{\eta}_{e,(2)} := \eta_e(\tau^*_{\indi,e} Z'_e) \), and \( Z_e, Z'_e \sim \mathcal{N}(0, I_p) \) are jointly Gaussian with \( \operatorname{Cov}(Z_e, Z'_e) = \rho_e I_p \). We assume that \( H^{(p)}_{\indi,e} \) converges pointwise to a function \( H_{\indi,e} \in C^1([0,1]) \) as \( p \to \infty \), and that the convergence also holds for the derivatives:
  \[
  \lim_{p \to \infty} \frac{d H^{(p)}_{\indi,e}}{d \rho_e} = \frac{d H_{\indi,e}}{d \rho_e}.
  \]
\end{enumerate}
\end{assumption}

The main result for the individual Lasso estimators—and consequently for the model average estimator—is stated below:

\begin{theorem}[Individual and model average estimators convergence]\label{thm:aver}
Under Assumptions~\ref{asp:prelim.transfer.model}, \ref{asp:aver.for.glamp}, and \ref{asp:aver.converge}, suppose we fix \( \tau_e = \tau^*_{\indi,e} \) and \( \theta_e = \theta^*_{\indi,e} \) in the definitions of \( \overline{\eta}_e \) from Equation~\eqref{eq:def.aver.et.bar}, for each \( e \in [E] \). Then, for any sequence of order-\( k \) pseudo-Lipschitz functions \( \phi : (\mathbb{R}^p)^{2E} \to \mathbb{R} \), with \( k \geq 1 \), we have:
\begin{equation*}
  \phi\left(\indiv{1}, \dots, \indiv{E}, \beta_1, \dots, \beta_E\right) 
  = \mathbb{E} \left[ 
      \phi\left(
        \overline{\eta}_1(\tau^*_{\indi,1}), \dots, \overline{\eta}_E(\tau^*_{\indi,E}),
        \beta_1, \dots, \beta_E
      \right)
    \right] + o(1).
\end{equation*}

In particular, for the model average estimator \( \aver = \sum_{e \in [E]} \pi_e \indiv{e} \), constructed using weights \( \{ \pi_e : e \in [E] \} \), and for any sequence of uniformly order-\( k \) pseudo-Lipschitz functions \( \phi : (\mathbb{R}^p)^{E+1} \to \mathbb{R} \), with \( k \geq 1 \), we have:
\begin{equation*}
  \phi\left(\aver, \beta_1, \dots, \beta_E\right) 
  = \mathbb{E} \left[ 
      \phi\left( 
        \sum_{e \in [E]} \pi_e \overline{\eta}_e(\tau^*_{\indi,e}), 
        \beta_1, \dots, \beta_E 
      \right)
    \right] + \smalop.
\end{equation*}
\end{theorem}

Theorem~\ref{thm:aver} is proved in Appendix~\ref{sec:apd.aver}.

\begin{figure}
    \centering
    \includegraphics[width=\linewidth]{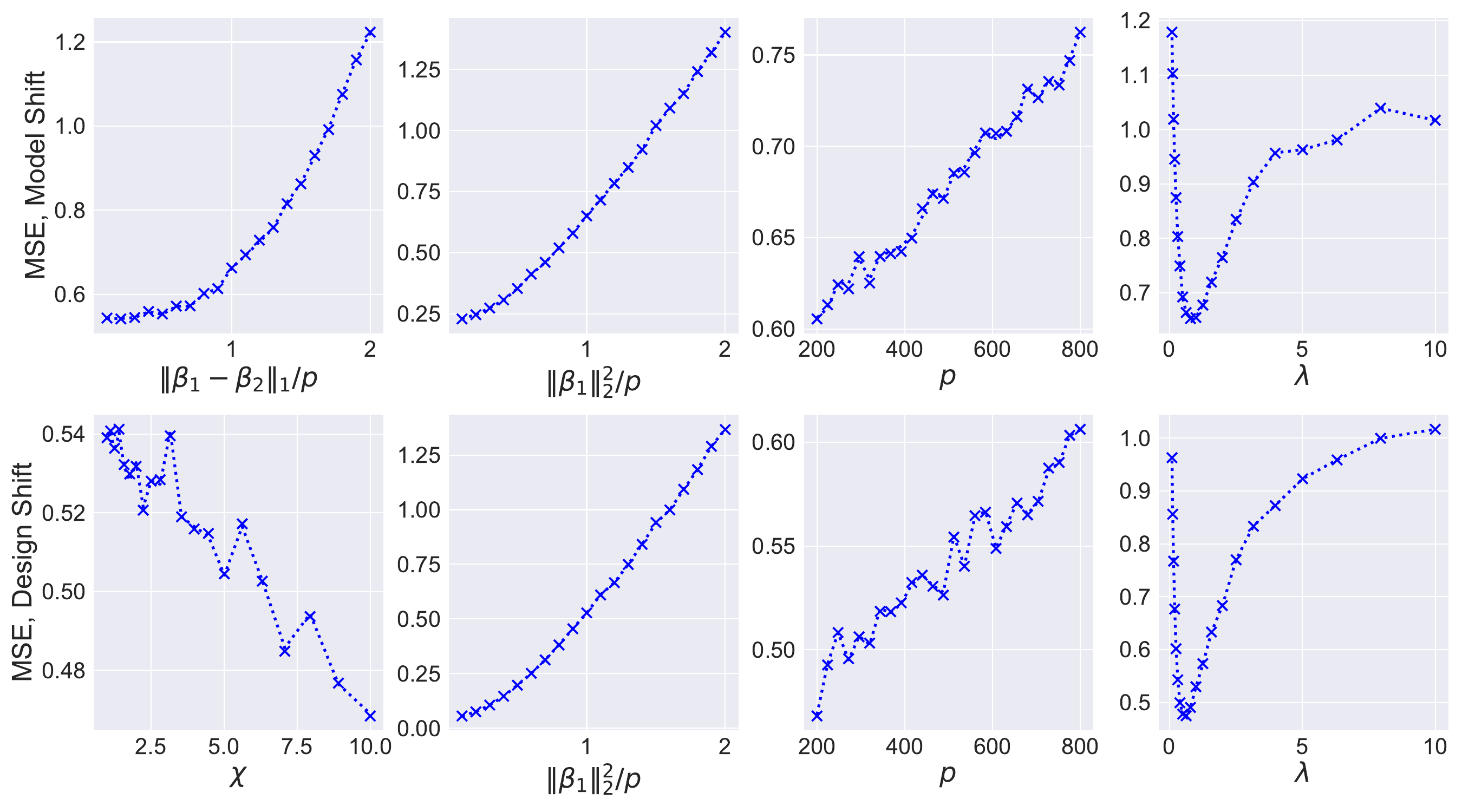}
    \caption{MSE for $\aver$. Same settings as in Figure~\ref{fig:stacked_Lasso}.}
    \label{fig:model_avg_Lasso}
\end{figure}

\subsection{The Second-Step Estimator}

Suppose, without loss of generality, that the first environment is our target environment, and our goal is to construct a good estimator for \( \beta_1 \). While both \( \stac \) and \( \aver \) are designed to aggregate information across environments to aid this task, additional gains can be achieved through target-specific calibration. One example is the \emph{second-step estimator} \citep{bastani2021predicting,li2022transfer}, which refines a first-step estimator by regressing its residual on the target data. We focus on a second-step estimator that uses \( \stac \) from Equation~\eqref{eq:def.stack.los} as the first-step input, while noting that other choices could be analyzed analogously. Given \( \stac \), we define the second-step estimator as:
\begin{gather}\label{eq:def.2nd.los}
  \se = \argmin_{b \in \mathbb{R}^{p}} L_{\rt}(b), \quad
  L_\rt(b) = \frac{1}{2} \left\| y_1 - \ix{1} \Sigma_1^{1/2} b \right\|_2^2 + \mu_{\rt}(b; \stac).
\end{gather}
The penalty function \( \mu : \mathbb{R}^{p} \to \mathbb{R} \) is designed to encourage the second-step estimator to remain close to the first-step estimator \( \stac \). Specifically, we study the following two formulations:

\begin{align}
  \text{(Joint estimator)} \quad 
  & \mu_{\joi}(b ; \stac) = \lambda_{\rt} \left\| b - \stac \right\|_1, \quad \lambda_{\rt} > 0, \label{eq:def.2nd.joi} \\
  \text{(Adaptively weighted estimator)} \quad 
  & \mu_{\ada}(b ; \stac) = \sum_{j \in [p]} \mu\left( \left| (\stac)_j \right| \right) \cdot |b_j| 
  = \left\| D_{\mu} b \right\|_1, \label{eq:def.2nd.ada}
\end{align}
where \( \mu : \mathbb{R}_{\geq 0} \to \mathbb{R}_+ \) is a decreasing function satisfying \( \mu(0) \geq \mu(+\infty) > 0 \), and $$ D_{\mu} = \operatorname{diag}(\mu(|(\stac)_j|) : j \in [p]). $$

The joint estimator with \( \mu_{\rt} = \mu_{\joi} \) was studied in \citet{bastani2021predicting}. The idea is to search within the neighborhood of the first-step estimator \( \stac \), which tends to perform well when \( \beta_1 \) and \( \stac \) share similar supports. The adaptively weighted estimator with \( \mu_{\rt} = \mu_{\ada} \) follows a similar intuition to the adaptive Lasso \citep{zou2006adaptive} and the empirical cross-prior approach proposed in \citet{li2022searching}. By leveraging information from the first step, it assigns smaller penalties to stronger signals, thereby reducing estimation bias. To avoid technical complications, we additionally require the penalty weights to be bounded away from zero and infinity, i.e., \( \mu(0) \geq \mu(+\infty) > 0 \). As a concrete example of such a weight function, one can use:
\[
\mu(x) = 5 + \frac{10}{0.05 + x^2}.
\]

To connect this analysis to the GLAMP framework, we define the fixed-design counterpart of the second-step estimator, which now also depends on a first-step estimator \( \hat{\beta} \). Specifically, define the function \( \xi : \mathbb{R}^p \times \mathbb{R}^p \to \mathbb{R}^p \) such that for any inputs \( v_{\rt}, \hat{\beta} \in \mathbb{R}^p \),
\begin{equation}\label{eq:def.2nd.xi}
\begin{aligned}
  \xi(v_{\rt}, \hat{\beta}) 
  = \argmin_{b \in \mathbb{R}^p} \Bigg\{ 
    \frac{1}{2} \left\| \Sigma_1^{1/2} b - \left[ v_{\rt} + \Sigma_1^{1/2} \beta_1 
    - \kappa_1 \gamma_{\ro} \Sigma_1^{1/2}(\beta_1 - \hat{\beta}) \right] \right\|_2^2 
    + \theta_{\rt} \, \mu_{\rt}(b; \hat{\beta}) 
  \Bigg\}.
\end{aligned}
\end{equation}

We use \( \frac{\partial }{\partial v_{\rt}} \xi \) and \( \frac{\partial }{\partial \hat{\beta}} \xi \) to denote the partial derivatives of \( \xi \) with respect to its first and second arguments, respectively. Note that Assumptions~\ref{asp:stack.for.glamp} and~\ref{asp:stack.converge} are still required here, since \( \stac \) serves as the first-step estimator within this second-step procedure.

\begin{assumption}[Second-step estimator GLAMP regularity]\label{asp:2nd.for.glamp}
Consider the transfer learning setting in Assumption~\ref{asp:prelim.transfer.model}, the stacked Lasso estimator in Equation~\eqref{eq:def.stack.los}, its fixed-design counterpart in Equation~\eqref{eq:def.stack.et}, and the associated Assumptions~\ref{asp:stack.for.glamp} and~\ref{asp:stack.converge}. Fix \( \tau_e = \tau_e^* \), \( \varpi_e = \varpi_e^* \) for all \( e \in [E] \), and \( \theta = \theta^* \) in the definitions of \( \eta \) and \( \overline{\eta} \). Next, consider the second-step estimator in Equation~\eqref{eq:def.2nd.los} and its fixed-design counterpart in Equation~\eqref{eq:def.2nd.xi}.

\begin{enumerate}[label=\upshape(\alph*),ref=\theassumption~(\alph*)]

  \item\label{asp:2nd.for.glamp.limits} 
  Define the intermediate quantity:
  \begin{equation}\label{eq:def.2nd.et.bar}
  \begin{aligned}
    \overline{\xi} = \overline{\xi}(\tau_{\rt}, \zeta) := \xi\left( \tau_{\rt} Z_{\rt},\, \overline{\eta} \right), 
    \quad Z_{\rt} = \zeta Z_1 + \sqrt{1 - \zeta^2} Z_1',
  \end{aligned}
  \end{equation}
  where \( \overline{\eta} = \eta(\tau_1^* Z_1, \dots, \tau_E^* Z_E) \), and \( Z_1, Z_1', Z_2, \dots, Z_E \overset{\text{i.i.d.}}{\sim} \mathcal{N}(0, I_p) \). For any \( \tau_{\rt} > 0 \), \( \zeta \in [-1,1] \), \( \gamma_{\ro} \geq 0 \), and \( \theta_{\rt} > 0 \), we assume the following limits exist:
  \[
  \lim_{p \to \infty} \frac{1}{p} \mathbb{E} \left[ \| \overline{\xi} \|_{\Sigma_1}^2 \right], \quad
  \lim_{p \to \infty} \frac{1}{p} \mathbb{E} \left[ \langle \beta_1, \overline{\xi} \rangle_{\Sigma_1} \right], \quad
  \lim_{p \to \infty} \frac{1}{p} \mathbb{E} \left[ \langle \overline{\eta}, \overline{\xi} \rangle_{\Sigma_1} \right].
  \]

  \item\label{asp:2nd.for.glamp.onsager} 
  Define the following intermediate quantity 
  \begin{equation}\label{eq:def.2nd.E}
  \begin{aligned}
    \mathcal{E}_{\rt}^{(p)} := \frac{1}{p} \mathbb{E} \left[ 
    \min_{b \in \mathbb{R}^p} 
    \left\{ 
      \frac{1}{2} \left\| \Sigma_1^{1/2} b - \left[ \tau_{\rt} Z_{\rt} + \Sigma_1^{1/2} \beta_1 
      - \kappa_1 \gamma_{\ro} \Sigma_1^{1/2} (\beta_1 - \overline{\eta}) \right] \right\|_2^2 
      + \theta_{\rt} \, \mu_{\rt}(b ; \overline{\eta}) 
    \right\} 
    \right],
  \end{aligned}
  \end{equation}
  where \( Z_{\rt} = \zeta Z_1 + \sqrt{1 - \zeta^2} Z_1' \) and \( \overline{\eta} = \eta(\tau_1^* Z_1, \dots, \tau_E^* Z_E) \), with \( Z_1, Z_1', Z_2, \dots, Z_E \overset{\text{i.i.d.}}{\sim} \mathcal{N}(0, I_p) \). We assume that for any fixed \( \gamma_{\ro} \geq 0 \), \( \theta_{\rt} > 0 \), there exists a limiting function \( \mathcal{E}_{\rt}(\tau_{\rt}, \zeta) \in C^1 \) such that:
  \[
  \lim_{p \to \infty} \mathcal{E}_{\rt}^{(p)} = \mathcal{E}_{\rt}, \quad
  \lim_{p \to \infty} \frac{\partial \mathcal{E}_{\rt}^{(p)}}{\partial \tau_{\rt}} = \frac{\partial \mathcal{E}_{\rt}}{\partial \tau_{\rt}}, \quad
  \lim_{p \to \infty} \frac{\partial \mathcal{E}_{\rt}^{(p)}}{\partial \zeta} = \frac{\partial \mathcal{E}_{\rt}}{\partial \zeta}.
  \]
  
\end{enumerate}
\end{assumption}

Assumption~\ref{asp:2nd.for.glamp} guarantees the well-definedness of the Onsager terms:
\begin{proposition}\label{prop:2nd.onsager}
Under Assumptions~\ref{asp:prelim.transfer.model}, \ref{asp:stack.for.glamp}, \ref{asp:stack.converge}, and \ref{asp:2nd.for.glamp}, define
\[
\overline{\gamma}_{\ro}^{(p)} := \frac{1}{\tau_1^* p} \, \mathbb{E} \left[ \left( Z_1 - \frac{\zeta}{\sqrt{1 - \zeta^2}} Z_1' \right)^{\!\top} \Sigma_1^{1/2} \, \overline{\xi} \right], \quad
\overline{\gamma}_{\rt}^{(p)} := \frac{1}{\tau_{\rt} p} \, \mathbb{E} \left[ \left( \frac{1}{\sqrt{1 - \zeta^2}} Z_1' \right)^{\!\top} \Sigma_1^{1/2} \, \overline{\xi} \right].
\]

For the joint estimator \( \mu_{\rt} = \mu_{\joi} \) in Equation~\eqref{eq:def.2nd.xi}, define the support set \( S_{\rt} := \mathrm{supp}\left( \overline{\xi}(\tau_{\rt}, \zeta) - \overline{\eta} \right) \subset [p] \).  
For the adaptively weighted estimator \( \mu_{\rt} = \mu_{\ada} \), define \( S_{\rt} := \mathrm{supp}(\overline{\xi}) \).

Then, in both cases, we have:
\begin{align*}
\overline{\gamma}_{\ro}^{(p)} 
&= \frac{1}{p} \sum_{j \in [p]} \mathbb{E} \left[
\left( \frac{\partial (\Sigma_1^{1/2} \xi)_j}{\partial \hat{\beta}} (\tau_{\rt} Z_{\rt}, \overline{\eta}) \right)^{\!\top}
\frac{\partial \overline{\eta}}{\partial (\tau_1^* Z_1)_j}
\right] = \frac{1}{p} \, \mathbb{E} \left[ \mathrm{Tr} \left( 
\frac{\partial \xi}{\partial \hat{\beta}}(\tau_{\rt} Z_{\rt}, \overline{\eta}) \cdot
\frac{\partial \overline{\eta}}{\partial (\tau_1^* Z_1)} \cdot \Sigma_1^{1/2}
\right) \right], \\
\overline{\gamma}_{\rt}^{(p)} 
&= \frac{1}{p} \sum_{j \in [p]} \mathbb{E} \left[ 
\frac{\partial (\Sigma_1^{1/2} \xi)_j}{\partial (v_{\rt})_j} (\tau_{\rt} Z_{\rt}, \overline{\eta}) 
\right] 
= \frac{1}{p} \, \mathbb{E} \left[ \mathrm{Tr} \left(
\frac{\partial \xi}{\partial v_{\rt}}(\tau_{\rt} Z_{\rt}, \overline{\eta}) \cdot \Sigma_1^{1/2}
\right) \right] 
= \frac{1}{p} \, \mathbb{E} \left[ |S_{\rt}| \right].
\end{align*}

Moreover, the limits \( \lim_{p \to \infty} \overline{\gamma}_{\ro}^{(p)} = \overline{\gamma}_{\ro} \) and \( \lim_{p \to \infty} \overline{\gamma}_{\rt}^{(p)} = \overline{\gamma}_{\rt} \) are well-defined.
\end{proposition}

Proposition~\ref{prop:2nd.onsager} is proved in Appendix~\ref{sec:apd.prelim.2nd.prop}. Readers might wonder what happens if \( \zeta^* = \pm 1 \) in the definitions of \( \overline{\gamma}_{\ro}^{(p)} \) and \( \overline{\gamma}_{\rt}^{(p)} \). It turns out that both terms remain well-defined due to two key facts: \emph{(i)} their equivalent expressions do not involve denominators of the form \( 1 - (\zeta^*)^2 \), and \emph{(ii)} they can be interpreted as limits of the respective expressions as \( \zeta^* \to \pm 1 \). These points are discussed in detail in Appendix~\ref{sec:apd.prelim.2nd.prop}.

Finally, we present the assumption that guarantees convergence of the GLAMP iterates to the second-step estimator \( \se \):

\begin{assumption}[Second-step estimator GLAMP fixed point]\label{asp:2nd.converge}
Consider the transfer learning setting in Assumption~\ref{asp:prelim.transfer.model}, along with the relevant estimators defined in Equations~\eqref{eq:def.stack.los}, \eqref{eq:def.stack.et}, \eqref{eq:def.2nd.los}, and \eqref{eq:def.2nd.xi}. We impose the following conditions in addition to Assumptions~\ref{asp:prelim.transfer.model}, \ref{asp:stack.for.glamp}, \ref{asp:stack.converge}, and \ref{asp:2nd.for.glamp}:

\begin{enumerate}[label=\upshape(\alph*),ref=\theassumption~(\alph*)]

\item\label{asp:2nd.converge.fixed.point} 
Assume that the parameter \( \gamma_{\ro} \) in Equation~\eqref{eq:def.2nd.xi} is chosen such that \( \gamma_{\ro} = \overline{\gamma}_{\ro} \). Further, assume the existence of a triplet \( (\tau_{\rt}^*, \zeta^*, \theta_{\rt}^*) \), where \( \tau_{\rt}^* > 0 \), \( \zeta^* \in [-1, 1] \), and \( \theta_{\rt}^* > 0 \), that solves the following system:
\begin{equation}
\begin{cases}
\tau_{\rt}^2 = (1 - \kappa_1 \overline{\gamma}_{\ro})^2 \mathbb{E}[W_1^2] + \kappa_1 \lim_{p \to \infty} \frac{1}{p} \mathbb{E} \left[ \| \beta_1 - \overline{\xi} - \kappa_1 \overline{\gamma}_{\ro} (\beta_1 - \overline{\eta}) \|_{\Sigma_1}^2 \right], \\
\tau_{\rt} \tau_1^* \zeta = (1 - \kappa_1 \overline{\gamma}_{\ro}) \mathbb{E}[W_1^2] + \kappa_1 \lim_{p \to \infty} \frac{1}{p} \mathbb{E} \left[ \langle \beta_1 - \overline{\eta}, \beta_1 - \overline{\xi} - \kappa_1 \overline{\gamma}_{\ro} (\beta_1 - \overline{\eta}) \rangle_{\Sigma_1} \right], \\
\theta_{\rt} = 1 / (1 - \kappa_1 \overline{\gamma}_{\rt}).
\end{cases}
\end{equation}

\item\label{asp:2nd.converge.cauchy} 
Let \( \epsilon \geq 0 \), and define a family of \( (2E + 2) \)-dimensional Gaussian distributions \( \mathcal{N}^{(\epsilon)} \) as follows. A distribution \( \mathcal{G} \in \mathcal{N}^{(\epsilon)} \) if and only if \( g = [g_1, g_2, \dots, g_{2E+2}]^\top \sim \mathcal{G} \) satisfies:
\begin{enumerate}
  \item[(i)] Marginals: \( g_1, g_2, \dots, g_{E+1} \overset{\text{i.i.d.}}{\sim} \mathcal{N}(0, 1) \), and \( g_{E+2}, \dots, g_{2E+2} \overset{\text{i.i.d.}}{\sim} \mathcal{N}(0, 1) \);
  \item[(ii)] The following \( E \) tuples are mutually independent: \( (g_1, g_2, g_{E+2}, g_{E+3}) \), \( (g_3, g_{E+4}) \), \( (g_4, g_{E+5}) \), ..., \( (g_{E+1}, g_{2E+2}) \);
  \item[(iii)] For all \( i \neq 2 \), with \( 1 \leq i \leq E+1 \), \( \mathrm{Cov}(g_i, g_{i+E+1}) \geq 1 - \epsilon \).
\end{enumerate}

Define the functional \( H^{(p)}_{\rt} : \mathcal{N}^{(\epsilon)} \to \mathbb{R} \) as:
\begin{equation}\label{eq:def.H.2.cal}
\begin{aligned}
H^{(p)}_{\rt}(\mathcal{G}) := 
&\ \frac{(1 - \kappa_1 \overline{\gamma}_{\ro})^2 \mathbb{E}[W_1^2]}{(\tau_{\rt}^*)^2} \\
&+ \frac{\kappa_1}{(\tau_{\rt}^*)^2 p} \mathbb{E} \left[ 
(\beta_1 - \overline{\xi}_{(1)} - \kappa_1 \overline{\gamma}_{\ro} (\beta_1 - \overline{\eta}_{(1)}))^{\!\top} \Sigma_1 
(\beta_1 - \overline{\xi}_{(2)} - \kappa_1 \overline{\gamma}_{\ro} (\beta_1 - \overline{\eta}_{(2)})) 
\right],
\end{aligned}
\end{equation}
where:
\begin{align*}
\overline{\eta}_{(1)} &= \eta(\tau_1^* Z_1, \tau_2^* Z_2, \dots, \tau_E^* Z_E), \;\;\overline{\xi}_{(1)} = \xi\left(\tau_{\rt}^* (\zeta^* Z_1 + \sqrt{1 - (\zeta^*)^2} Z_1'), \overline{\eta}_{(1)} \right), \\
\overline{\eta}_{(2)} &= \eta(\tau_1^* \widetilde{Z}_1, \tau_2^* \widetilde{Z}_2, \dots, \tau_E^* \widetilde{Z}_E), \;\;
\overline{\xi}_{(2)} = \xi\left(\tau_{\rt}^* (\zeta^* \widetilde{Z}_1 + \sqrt{1 - (\zeta^*)^2} \widetilde{Z}_1'), \overline{\eta}_{(2)} \right),
\end{align*}
and \( [Z_1, Z_1', Z_2, \dots, Z_E, \widetilde{Z}_1, \widetilde{Z}_1', \widetilde{Z}_2, \dots, \widetilde{Z}_E] \in \mathbb{R}^{p \times (2E + 2)} \) has i.i.d. rows drawn from \( \mathcal{G} \).

Assume there exists \( \epsilon_H > 0 \) such that for all \( \epsilon \in [0, \epsilon_H] \) and all \( \mathcal{G} \in \mathcal{N}^{(\epsilon)} \), we have:
\[
\lim_{p \to \infty} H^{(p)}_{\rt}(\mathcal{G}) = H_{\rt}(\mathcal{G}),
\]
for some functional \( H_{\rt} : \mathcal{N}^{(\epsilon_H)} \to \mathbb{R} \).

In the special case \( \epsilon = 0 \), each element of \( \mathcal{N}^{(0)} \) is characterized by the correlation \( \rho = \mathrm{Cov}(g_2, g_{E+3}) \), and we slightly overload notation by writing:
\begin{equation}\label{eq:def.H.2}
\forall\, \mathcal{G}_\rho \in \mathcal{N}^{(0)}, \quad
H^{(p)}_{\rt}(\rho) := H^{(p)}_{\rt}(\mathcal{G}_\rho), \quad 
H_{\rt}(\rho) := \lim_{p \to \infty} H^{(p)}_{\rt}(\rho),
\end{equation}
with the additional assumption that:
\[
\lim_{p \to \infty} \frac{d}{d \rho} H^{(p)}_{\rt}(\rho) = \frac{d}{d \rho} H_{\rt}(\rho).
\]

\end{enumerate}
\end{assumption}

Our convergence result for the second-step estimator is stated below:

\begin{theorem}[Second-step estimator convergence]\label{thm:2nd}
Consider the second-step estimator \( \se \) defined in Equation~\eqref{eq:def.2nd.los}. Under Assumptions~\ref{asp:prelim.transfer.model}, \ref{asp:stack.for.glamp}, \ref{asp:stack.converge}, \ref{asp:2nd.for.glamp}, and \ref{asp:2nd.converge}, suppose we fix:
\begin{itemize}
  \item \( \tau_e = \tau_e^* \), \( \varpi_e = \varpi_e^* \) for all \( e \in [E] \), and \( \theta = \theta^* \) as in Assumption~\ref{asp:stack.converge.fixed.point}, used in the definitions of \( \eta \) and \( \overline{\eta} \);
  \item \( \theta_{\rt} = \theta_{\rt}^* \), \( \gamma_{\ro} = \overline{\gamma}_{\ro} \), and \( \zeta = \zeta^* \) as in Assumption~\ref{asp:2nd.converge.fixed.point}, used in the definitions of \( \xi \), \( \overline{\xi} \), and the construction of \( Z_{\rt} \).
\end{itemize}
Then, for any sequence of order-\( k \) pseudo-Lipschitz functions \( \phi : (\mathbb{R}^p)^{E+2} \to \mathbb{R} \), with \( k \geq 1 \), we have:
\[
\phi\left( \se, \stac, \beta_1, \dots, \beta_E \right) 
= \mathbb{E} \left[ \phi\left( \overline{\xi}(\tau_{\rt}^*, \zeta^*), \overline{\eta}(\tau_1^*, \dots, \tau_E^*), \beta_1, \dots, \beta_E \right) \right] +  \smalop.
\]
\end{theorem}

Theorem~\ref{thm:2nd} is proved in Appendix~\ref{sec:apd.2nd}.

Figure~\ref{fig:second_step_Lasso} illustrates the finite-sample accuracy of Theorem~\ref{thm:2nd}, where we evaluate the MSE of the second-step estimator \( \se \) (specificly the joint version as in \eqref{eq:def.2nd.joi}) under the same settings used in Figure~\ref{fig:stacked_Lasso}.

\begin{figure}
    \centering
    \includegraphics[width=\linewidth]{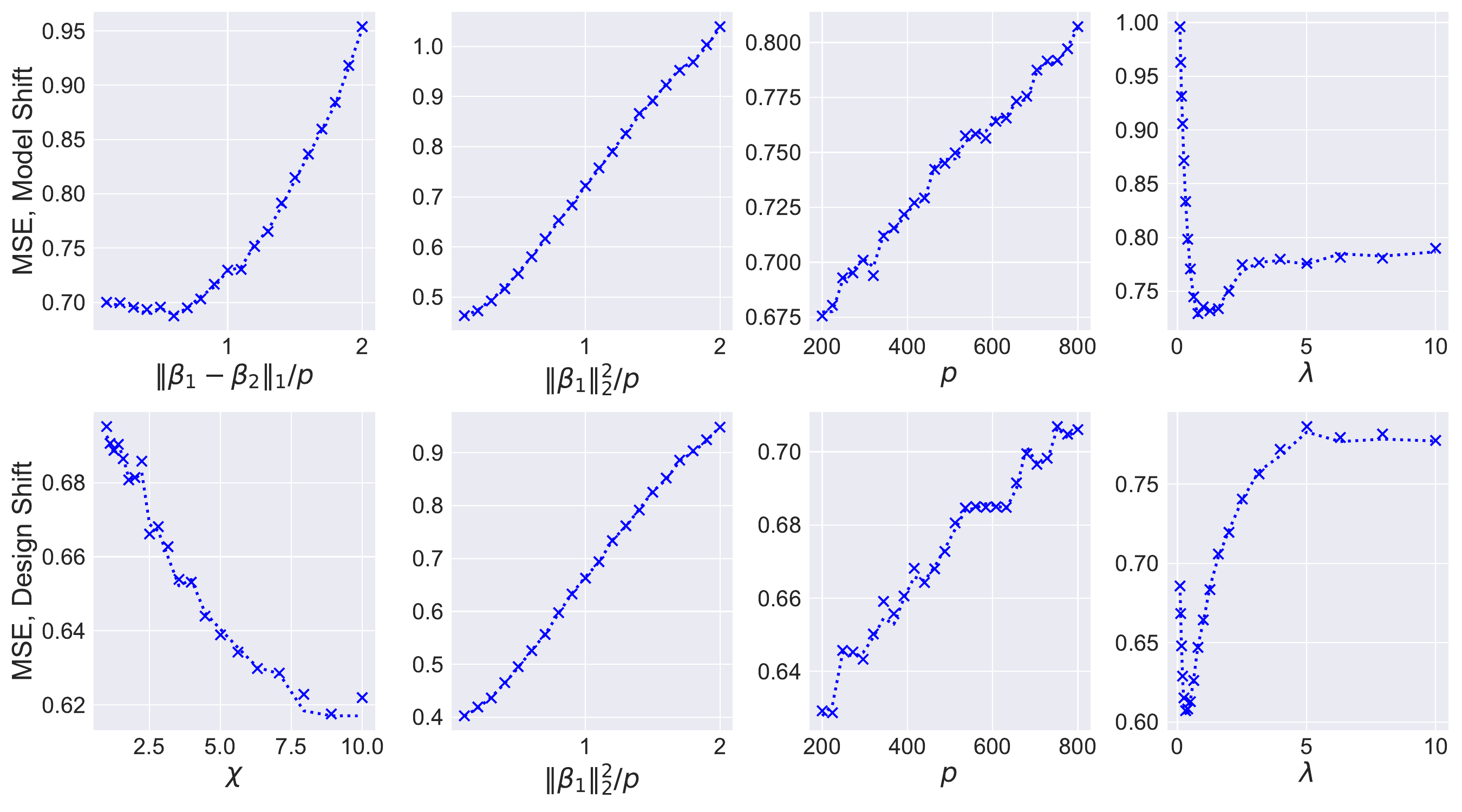}
    \caption{MSE for $\se$. Same settings as in Figure~\ref{fig:stacked_Lasso} with $\lambda_{\rt}=0.5$.}
    \label{fig:second_step_Lasso}
\end{figure}

\section{Discussion and Future Directions}\label{sec:discussion}

In this paper, we introduce \emph{Multi-Environment Generalized Long AMP} (GLAMP), an AMP framework that accommodates both stacked iterations and non-separable denoising functions, while specializing to transfer learning applications. Leveraging the framework, we derive precise asymptotic characterizations for three Lasso-based transfer learning estimators under proportional asymptotics. We conclude by outlining several promising directions for future research:

\begin{enumerate}[label=(\roman*)]
    \item \textbf{Extended Characterizations and Comparisons:} 
    Beyond the three Lasso-based estimators considered in this work, the multi-environment GLAMP framework has the potential to analyze a broader class of estimators, including—but not limited to—those discussed in Section~\ref{subsec:related_transfer}. Additionally, more comprehensive analytical and empirical comparisons could yield deeper insights into the performance of these methods across diverse transfer learning settings and offer more practical guidance for real-world applications.

    \item \textbf{Generalizing Assumptions:} 
    Our analysis assumes Gaussian design matrices; however, in many practical applications, covariates are often discrete or heavy-tailed. Extending existing universality results \citep{bayati2015universality,chen2021universality,hu2022universality,liang2022precise,montanari2022universality,dudeja2023universality,li2023spectrum,wang2024universality,lahiry2024universality,ghane2024universality} or AMP under right rotationally invariant designs \cite{rangan2019vector,fan2022approximate,ma2017orthogonal,li2023spectrum} to the multi-environment GLAMP setting is an exciting avenue for future work.

    \item \textbf{Generalized Linear Models:}
  Finally, we focused on linear models with continuous responses. However, GAMP has been successfully applied to logistic regression in prior work \citep{sur2019modern,barbier2019optimal} and extending our transfer learning applications of multi-environment GLAMP to estimators in GLMs and classification problems \citep{reeve2021adaptive,hanneke2022no,li2024estimation,tian2023transfer,maity2024linear} would be a natural avenue for future research.
\end{enumerate}

\bibliography{refs}
\clearpage

\appendix
\addcontentsline{toc}{section}{Appendices}
\renewcommand{\theequation}{\thesection-\arabic{equation}}
\setcounter{equation}{0}

\renewcommand{\theassumption}{\thesection.\arabic{assumption}}
\renewcommand{\thedefinition}{\thesection.\arabic{definition}}
\renewcommand{\thetheorem}{\thesection.\arabic{theorem}}
\renewcommand{\thelemma}{\thesection.\arabic{lemma}}
\renewcommand{\theproposition}{\thesection.\arabic{proposition}}
\renewcommand{\thecorollary}{\thesection.\arabic{corollary}}

\section{Proof of Theorem~\ref{cor:asym.multienv}}
\label{sec:apd.pf.glamp.asym.multienv}

\subsection{Preliminaries: Symmetric and Asymmetric GLAMP}\label{sec:sy_asy_glamp}

Our proof builds upon the symmetric GLAMP algorithm from \cite{gerbelot2023graph}. We first present it below. Subsequently we present an asymmetric version that is useful to establish validity of our multi-environment GLAMP framework. 

\subsubsection{Symmetric GLAMP}
The Symmetric GLAMP is introduced in \cite[Theorem 2]{gerbelot2023graph} as ``Matrix-valued symmetric AMP iterations with non-separable non-linearities". Here we restate and reformulate it in the same manner as we introduced multi-environment GLAMP in Section~\ref{sec:multi_glamp}.

\begin{definition}[Symmetric GLAMP Instance]\label{def:GLAMP.reg}
  Let $q \in \N_+$ be fixed and known, and consider the asymptotic regime where $N \rightarrow \infty$. A symmetric GLAMP instance is a quadruple \( \{A, (f_{t})_{t \in \mathbb{N}} , x^{0}, Y \}  \), where $A \in \R^{N\times N}$ is a symmetric random design matrix, $Y\in \R^{N\times q}$ encodes information that the GLAMP iterates can use, \( (f_t)_{t \in \mathbb{N}} \) is a collection of mappings with \( f_t:\R^{N\times q}\times \R^{N\times q}\to \R^{N\times q} \), and \( x^{0}  \in \R^{N\times q}\) is the initial condition for the GLAMP iterates. 
\end{definition}

We sometimes omit $x^0$ and $Y$ when they are clear from context. 

\begin{assumption}[Symmetric GLAMP assumptions]\label{asp:glamp.converge}
With the GLAMP instance in Definition~\ref{def:GLAMP.reg}, we assume
\begin{enumerate}[label=\upshape(\alph*),ref=\theassumption (\alph*)]

    \item\label{def:GLAMP.reg.GOE} \( A = G + G^{\T} \) where \( G\in\R^{N\times N} \) has iid Gaussian entries from \( \Nscr(0,\frac{1}{2N}) \). This means 
    \( A \) is drawn from  the Gaussian Orthogonal Ensemble (GOE) of dimension $N$, denoted as \( A\sim \mathrm{GOE}(N) \). 

    \item\label{def:GLAMP.reg.y} \( \limsup_{N}\frac{1}{N}\norm{Y}_2^{2}< +\infty \).
    
    \item\label{def:GLAMP.reg.init} \( \limsup_{N}\frac{1}{N}\norm{x^{0}}_2^{2} < +\infty \). 
    
    \item\label{def:GLAMP.reg.func} 
    For each fixed \( t \in \mathbb{N} \), \( f_t\) is uniformly Lipschitz in \( N \), in the sense that there exists some universal constant \( L_t \) s.t. for any matrices \( M_1,M_2,M_1',M_2'\in\R^{N\times q} \),
      \begin{equation*}
        \norm{f_t(M_1,M_2)-f_t(M_1',M_2')}_2 \leq L_t \norm{M_1-M_1'}_2+L_t\norm{M_2-M_2'}_2, \ \forall \, N.
      \end{equation*}
      Moreover, \( \limsup_{N} \frac{1}{\sqrt{N}}\norm{f_t(0,0)}_2 < +\infty \). 
    
\item \label{asp:glamp.converge.sep} Consider some fixed known \( E\in\N_+ \) and a sequence of \( E \)-partitions (groups) of \( [N] \): \( [N] = \cup_{a\in[E]} C_{a}^{N} \), and \( \lim_{N\to \infty}\abs{C_{a}^{N}} / N = c_a \in (0,1)  \).  For any step \( t \in \mathbb{N} \), \( f_t\) is separable among groups, in the sense that there are functions \( g_{t,a}(\cdot ,\cdot) : \R^{\abs{C_a^{N}}\times q}\times\R^{\abs{C_a^{N}}\times q} \to \R^{\abs{C_a^{N}}\times q}  \) for each \( a\in[E] \), s.t. 

\begin{equation*}
  f_t(x^{t},Y) = \begin{bmatrix}
    g_{t,1}(x^{t}_{C_1^{N}},Y_{C_1^{N}}) \\
    g_{t,2}(x^{t}_{C_2^{N}},Y_{C_2^{N}}) \\
    \vdots \\
    g_{t,E}(x^{t}_{C_E^{N}},Y_{C_p^{N}}) 
  \end{bmatrix}.
\end{equation*}
\item \label{asp:glamp.converge.init} If we write \( q^{0} = f_0(x^{0}) \), then \( \forall\, a\in[E] \),
\begin{equation*}
  \frac{1}{\abs{C_a^{N}} } (q^{0}_{C_a^{N}})^{\T} q^{0}_{C_a^{N}} \to {\Sigma }_a^{0}
\end{equation*}
for some matrices $\{{\Sigma}_a^0\in\R^{q\times q}:a\in[E]\}$. Specifically, 
$$\frac{1}{N} (q^0)^T q^0 \rightarrow \sum_{a\in E} c_a \Sigma_a^0 := \Sigma^0$$
\item \label{asp:glamp.converge.later}  For any constant covariance matrices $\Sigma_1,\Sigma_2\in\R^{q\times q}$, any $s,r > 0$ and any $a\in [E]$, 
\begin{equation*}
  \lim_{N\to \infty}\frac{1}{\abs{C_a^{N}} }\E\left[ g_{s,a}(Z_1,Y_{C_a^{N}})^{\T} g_{r,a}(Z_2,Y_{C_a^{N}}) \right],\quad \lim_{N\to \infty}\frac{1}{\abs{C_a^{N}} }\E\left[ (q^{0}_{C_a^{N}})^{\T} g_{r,a}(Z_2,Y_{C_a^{N}}) \right]
\end{equation*} 
exist, where \( Z_1,Z_2 \in \R^{\abs{C_a^{N}}\times q } \) have iid row vectors from \( \Nscr(0,\Sigma _1) \) and \( \Nscr(0,\Sigma _2) \) respectively. 
\end{enumerate}
\end{assumption}

\begin{definition}[Symmetric GLAMP Iterates]\label{def:GLAMP.iter}
  We define the GLAMP iterates corresponding to an instance \( \{ A, (f_t)_{t \in \mathbb{N}}, x^{0},Y \}  \) as the sequence of matrices \( x^{t}\in\R^{N\times q}:t \in \mathbb{N} \) defined as follows:
  \begin{equation}\label{eq:def.glamp}
    x^{t+1} = A f_t(x^t) -f_{t-1}(x^{t-1}) B_t^{\T}, \ {B}_t = \E \left[ \frac{1}{N} \left( \sum_{j\in[N]} \frac{\partial (f_t)_{j,\cdot }}{\partial x_{j,\cdot }} \right) (Z^{t}) \right],
  \end{equation}
  where \( Z^{t}\in\R^{N\times q} \) has its rows drawn iid from \( \Nscr(0,\Sigma ^t) \) with $\Sigma^t$ to be defined. The term \( {B}_t \) is understood as the average of \( N \) Jacobians. As \( f_t : \R^{N\times q}\to \R^{N\times q} \), \(  \frac{\partial (f_t)_{j,\cdot }}{\partial x_{j,\cdot }} \) takes the \( j \)-th row of \( f_t(\cdot ) \) as a vector function, and consider its Jacobian only w.r.t. the \( j \)-th row of the input \( x^{t} \). 
\end{definition}

\begin{definition}[Symmetric GLAMP State Evolution]\label{def:GLAMP.state.evo}
  Starting from \( \Sigma^0  \) in Assumption~\ref{asp:glamp.converge.init}, we define a sequence of matrices for (\( t \in \mathbb{N} \)):
  \begin{equation*}
    \Sigma ^{t+1} = \lim_{N\to\infty}\frac{1}{N}\E \left[ f_t(Z^t, Y)^{\T}f_t(Z^t,Y)\right],
  \end{equation*}
  where \( Z^{t}\in\R^{N \times q} \) has iid row vectors drawn from \( \Nscr(0,\Sigma ^t) \).   
\end{definition}

\begin{theorem}[Symmetric GLAMP convergence]\label{thm:main}
  Under Assumption~\ref{asp:glamp.converge}, for any \( t \in \mathbb{N_+} \) and any sequence of uniformly order-\( k \) pseudo-Lipschitz functions \( \phi _N : \R^{|C_a^N|\times(3q)}\to \R \) for any \( k\geq 1 \), the GLAMP iterates in Definition~\ref{def:GLAMP.iter} satisfy: 
\begin{equation*}
  \forall\, a\in[E],\ \phi _N([x^{0}_{C_a^{N}}| Y_{C_a^{N}}|x^{t}_{C_a^{N}}]) = \E \left[ \phi _N(x^{0}_{C_a^{N}}| Y_{C_a^{N}}| Z^{t}_a) \right] + \smalop,
\end{equation*}
where \( Z^{t}_a\in\R^{\abs{C_a^{N}} \times q} \) has its row vectors drawn iid from \( \Nscr(0,\Sigma ^t) \). 
\end{theorem}

Theorem \ref{thm:main} is proved in \cite{gerbelot2023graph} as their Theorem 2.

\subsubsection{Asymmetric GLAMP}\label{sec:asymmetric}

Many previous AMP works—including the seminal paper by \citet{donoho2009message}—adopt an alternative but equivalent formulation based on an asymmetric random matrix \( X \), along with suitably defined auxiliary quantities. In this section, we articulate the asymmetric version of GLAMP following that tradition. This asymmetric version is also the prerequisite for the multi-environment version in Section~\ref{sec:multi_glamp}.

\begin{definition}[Asymmetric GLAMP instance]\label{def:GLAMP.asym}
Let \( q \) be a fixed and known positive integer, and consider the asymptotic regime where \( \lim_{N \rightarrow +\infty} p/N \rightarrow \kappa \in (0, +\infty) \). An \emph{asymmetric GLAMP instance} is defined as a tuple
\[
\left\{ X, (\Psi^t)_{t \in \mathbb{N}}, (\eta^t)_{t \in \mathbb{N}}, V^0, B, W \right\},
\]
with the following components:
\begin{itemize}
  \item \( X \in \mathbb{R}^{N \times p} \) is an asymmetric random design matrix;
  \item \( B \in \mathbb{R}^{p \times q} \), \( W \in \mathbb{R}^{N \times q} \) encode information accessible to the GLAMP iterates;
  \item \( (\Psi^t)_{t \in \mathbb{N}} \) and \( (\eta^t)_{t \in \mathbb{N}} \) are sequences of update mappings, where
  \[
  \Psi^t : \mathbb{R}^{N \times q} \times \mathbb{R}^{N \times q} \to \mathbb{R}^{N \times q}, \quad
  \eta^t : \mathbb{R}^{p \times q} \times \mathbb{R}^{p \times q} \to \mathbb{R}^{p \times q};
  \]
  \item \( V^0 \in \mathbb{R}^{p \times q} \) is the initial condition for the GLAMP iterates.
\end{itemize}
\end{definition}

As in prior formulations, the components \( V^0 \), \( B \), and \( W \) are omitted from notation when the context is clear. The following assumptions are associated with the asymmetric GLAMP instance:

\begin{assumption}[Asymmetric GLAMP assumptions]\label{asp:glamp.asym.converge}
Given the GLAMP instance in Definition~\ref{def:GLAMP.asym}, we make the following assumptions:

\begin{enumerate}[label=\upshape(\alph*),ref=\theassumption~(\alph*)]

\item The matrix \( X \in \mathbb{R}^{N \times p} \) has i.i.d.\ entries drawn from \( \mathcal{N}(0, 1/N) \).

\item The matrices \( B \in \mathbb{R}^{p \times q} \) and \( W \in \mathbb{R}^{N \times q} \) satisfy
\[
\limsup_{p \to \infty} \frac{1}{\sqrt{p}} \|B\|_2 < \infty, \quad \limsup_{N \to \infty} \frac{1}{\sqrt{N}} \|W\|_2 < \infty.
\]

\item The initial condition \( V^0 \in \mathbb{R}^{p \times q} \) satisfies
\[
\limsup_{p \to \infty} \frac{1}{p} \|V^0\|_2^2 < \infty.
\]

\item For each fixed \( t \in \mathbb{N} \), the mappings \( \Psi^t \) are uniformly Lipschitz in \( N \), and \( \eta^t \) are uniformly Lipschitz in \( p \). In addition,
\[
\limsup_{N \to \infty} \frac{1}{\sqrt{N}} \| \Psi^t(0, 0) \|_2 < \infty, \quad 
\limsup_{p \to \infty} \frac{1}{\sqrt{p}} \| \eta^t(0, 0) \|_2 < \infty.
\]

\item \label{asp:GLAMP.asym.converge,sep} 
Fix \( E \in \mathbb{N}_+ \), and let \( [N] = \bigcup_{a \in [E]} C_a^N \) be a sequence of partitions (groups) such that \( \lim_{N \to \infty} |C_a^N| / N = c_a \in (0, 1) \). For any \( t \in \mathbb{N} \), the function \( \Psi^t \) is separable across groups, i.e., there exist mappings 
\[
\Psi_a^t : \mathbb{R}^{|C_a^N| \times q} \times \mathbb{R}^{|C_a^N| \times q} \to \mathbb{R}^{|C_a^N| \times q}
\]
such that
\[
\Psi^t(R^t, W) = 
\begin{bmatrix}
\Psi_1^t(R^t_{C_1^N}, W_{C_1^N}) \\
\Psi_2^t(R^t_{C_2^N}, W_{C_2^N}) \\
\vdots \\
\Psi_E^t(R^t_{C_E^N}, W_{C_E^N})
\end{bmatrix}.
\]

\item \label{asp:GLAMP.asym.converge,init}
The initial condition satisfies
\[
\frac{1}{p} \left( \eta^0(V^0) \right)^\top \eta^0(V^0) \to \Sigma_{(R)}^0
\]
for some matrix \( \Sigma_{(R)}^0 \in \mathbb{R}^{q \times q} \).

\item \label{asp:GLAMP.asym.converge.later.psi}
For any fixed covariance matrices \( \Sigma_1, \Sigma_2 \in \mathbb{R}^{q \times q} \), and any \( s, r \geq 0 \), \( a \in [E] \), the following limit exists:
\[
\lim_{N \to \infty} \frac{1}{|C_a^N|} \mathbb{E} \left[ 
\Psi_a^s(Z_1, W_{C_a^N})^\top \Psi_a^r(Z_2, W_{C_a^N})
\right],
\]
where \( Z_1, Z_2 \in \mathbb{R}^{|C_a^N| \times q} \) have i.i.d.\ rows drawn from \( \mathcal{N}(0, \Sigma_1) \) and \( \mathcal{N}(0, \Sigma_2) \), respectively.

\item \label{asp:GLAMP.asym.converge.later.eta}
For any fixed covariance matrices \( \Sigma_1, \Sigma_2 \in \mathbb{R}^{q \times q} \), and any \( s, r \geq 1 \), the following limits exist:
\[
\lim_{p \to \infty} \frac{1}{p} \mathbb{E} \left[ 
\eta^s(Z_1, B)^\top \eta^r(Z_2, B)
\right], \quad
\lim_{p \to \infty} \frac{1}{p} \mathbb{E} \left[ 
\eta^0(V^0, B)^\top \eta^r(Z_2, B)
\right],
\]
where \( Z_1, Z_2 \in \mathbb{R}^{p \times q} \) have i.i.d.\ rows drawn from \( \mathcal{N}(0, \Sigma_1) \) and \( \mathcal{N}(0, \Sigma_2) \), respectively.

\end{enumerate}
\end{assumption}

Given an instance, the asymmetric GLAMP iterates and their corresponding state evolution are defined as follows:

\begin{definition}[Asymmetric GLAMP iterates]\label{def:glamp.asym.iter}
Let \( \{ X, (\Psi^t)_{t \in \mathbb{N}}, (\eta^t)_{t \in \mathbb{N}}, V^0, B, W \} \) be an asymmetric GLAMP instance satisfying Assumption~\ref{asp:glamp.asym.converge}. The asymmetric GLAMP iterates are the alternating sequence of matrices \( \{ (V^t, R^t) : t \in \mathbb{N} \} \), initialized with
\[
R^0 = X \eta^0(V^0, B).
\]
For \( t \in \mathbb{N_+} \), the iterates are defined recursively as:
\begin{equation}\label{eq:def.glamp.asym}
\begin{alignedat}{2}
V^t &= \frac{1}{\kappa} X^\top \Psi^{t-1}(R^{t-1}, W) 
- \frac{1}{\kappa} \eta^{t-1}(V^{t-1}, B) \left[ D_{(R)}^{t-1} \right]^\top, \quad
&D_{(R)}^{t-1} &= \frac{1}{N} \mathbb{E} \left[ \sum_{j \in [N]} \left( \frac{\partial (\Psi^{t-1})_{j,\cdot}}{\partial (R^{t-1})_{j,\cdot}} \right)(Z_{(R)}^{t-1}) \right], \\
R^t &= X \eta^t(V^t, B) 
- \Psi^{t-1}(R^{t-1}, W) \left[ D_{(V)}^t \right]^\top, \quad
&D_{(V)}^t &= \frac{1}{p} \mathbb{E} \left[ \sum_{j \in [p]} \left( \frac{\partial (\eta^t)_{j,\cdot}}{\partial (V^t)_{j,\cdot}} \right)(Z_{(V)}^t) \right],
\end{alignedat}
\end{equation}
where the Gaussian random matrices \( Z_{(R)}^{t-1} \) and \( Z_{(V)}^t \) will be defined in the state evolution equations (Definition~\ref{def:GLAMP.state.evo.asym}).
\end{definition}

\begin{definition}[Asymmetric GLAMP state evolutions]\label{def:GLAMP.state.evo.asym}
Starting from the initialization \( \Sigma^0_{(R)} \) as given in Assumption~\ref{asp:GLAMP.asym.converge,init}, we define a sequence of matrices as follows:

For \( t \in \mathbb{N_+} \),
\begin{equation*}
\Sigma^t_{(V)} = \lim_{N \to \infty} \frac{1}{N} \mathbb{E} \left[ 
\Psi^{t-1}(Z^{t-1}_{(R)}, W)^\top \Psi^{t-1}(Z^{t-1}_{(R)}, W) 
\right],
\end{equation*}
where \( Z^{t-1}_{(R)} \in \mathbb{R}^{N \times q} \) has i.i.d.\ rows drawn from \( \mathcal{N}(0, \kappa\, \Sigma^{t-1}_{(R)}) \).

Next,
\begin{equation*}
\Sigma^t_{(R)} = \lim_{p \to \infty} \frac{1}{p} \mathbb{E} \left[ 
\eta^t(Z^t_{(V)}, B)^\top \eta^t(Z^t_{(V)}, B) 
\right],
\end{equation*}
where \( Z^t_{(V)} \in \mathbb{R}^{p \times q} \) has i.i.d.\ rows drawn from \( \mathcal{N}(0, \frac{1}{\kappa^2} \Sigma^t_{(V)}) \).
\end{definition}

\begin{corollary}[Asymmetric GLAMP convergence]\label{cor:asym}
Let \( t \in \mathbb{N} \). Consider the iterates defined in Definition~\ref{def:GLAMP.asym} under Assumption~\ref{asp:glamp.asym.converge}. Then:

\begin{itemize}
  \item For any order-\( k \) pseudo-Lipschitz function \( \phi : (\mathbb{R}^{p \times q})^3 \to \mathbb{R} \), with \( k \geq 1 \), we have
  \[
  \phi(V^t, V^0, B) = \mathbb{E} \left[ \phi(Z^{t}_{(V)}, V^0, B) \right] + \smalop,
  \]
  where \( Z^t_{(V)} \in \mathbb{R}^{p \times q} \) has i.i.d.\ rows sampled from \( \mathcal{N}(0, \kappa^{-2} \Sigma^t_{(V)}) \), as defined in Definition~\ref{def:GLAMP.state.evo.asym}.
  
  \item For any \( a \in [E] \) and any order-\( k \) pseudo-Lipschitz function \( \phi : (\mathbb{R}^{|C_a^N| \times q})^2 \to \mathbb{R} \), with \( k \geq 1 \), we have
  \[
  \phi(R^t_{C_a^N}, W_{C_a^N}) = \mathbb{E} \left[ \phi(Z^t_{(R)}, W_{C_a^N}) \right] +\smalop,
  \]
  where \( Z^t_{(R)} \in \mathbb{R}^{|C_a^N| \times q} \) has i.i.d.\ rows drawn from \( \mathcal{N}(0, \kappa \Sigma^t_{(R)}) \), also defined in Definition~\ref{def:GLAMP.state.evo.asym}.
\end{itemize}
\end{corollary}

Corollary~\ref{cor:asym} is a direct consequence of Theorem~\ref{thm:main}, and its proof is omitted.

\subsection{Proof of Theorem \ref{cor:asym.multienv}}
\begin{proof}[\textbf{Proof of Theorem~\ref{cor:asym.multienv}}]
  We start from Corollary~\ref{cor:asym}, together with its full set of assumptions, and specify it down to the form of the multi-environment GLAMP. Then we relax some of the assumptions to get Theorem~\ref{cor:asym.multienv}.

  Suppose we have $E$ environments, each with a sample size $n_e$, $N=\sum_{e\in[E]}n_e$. We define  $\kappa_e = \lim_{p}\frac{p}{n_e}$, and \( \delta _e = \sqrt{n_e/N} \). We partition the large Gaussian matrix \( X\in \R^{N\times p} \), $X_{i,j} \iidsim \Nscr(0,1/N)$ into \( X = [\delta_1 X_1| \delta_2 X_2 | \cdots |\delta_E X_E]^{\T} \), where \( X_e \in \R^{n_e \times p} \) and \( X_{i,j}\iidsim \Nscr(0,1/n_e) \). We also define index set \( C_e^{N} = \{ j:\sum_{e'<e}n_{e'} < j \leq \sum_{e'\leq e}n_{e'} \}  \) for \( e\in[E] \). 
  
  We also consider the setting of Corollary~\ref{cor:asym} with \( q=\sum_{e\in[E]}q_e \). We let \( V^{t} = [V_1^{t}|V_2^{t}|\cdots |V_E^{t}] \), where \( V_e^{t}\in\R^{p\times q_e} \) for \( e\in[E] \) and \( t \in \mathbb{N} \).  We also let 
  \begin{equation}\label{eq:pf.multienv.eta}
   \eta ^t (V^{t},B) = 
   \left[\frac{1}{\delta_1}\eta_1^{t}(V_1^{t},\cdots,V_E^{t},B)\,\Big|\,
   \frac{1}{\delta_2}\eta_2^{t}(V_1^{t},\cdots,V_E^{t},B)\,\Big|\,\cdots \,\Big|\,
   \frac{1}{\delta_E}\eta_E^{t}(V_1^{t},\cdots,V_E^{t},B)
   \right].
  \end{equation}
For \( R^{t}\in\R^{N\times q} \), we let 
\begin{equation*}
  R^{t} = \begin{bmatrix}
    \begin{matrix}
      \begin{matrix}
        R^{t}_1 \\ * \\ * \\ \vdots \\ *
      \end{matrix} \quad\ \rvline \quad
      \begin{matrix}
        * \\ R^{t}_2 \\ * \\ \vdots \\ *
      \end{matrix} \quad\ \rvline \quad
      \begin{matrix}
        * \\ * \\ R^{t}_3 \\ \vdots \\  *
      \end{matrix} \quad\ \rvline \quad \cdots  \quad\ \rvline \ \
      \begin{matrix}
        * \\ * \\ * \\ \vdots \\  R^{t}_E
      \end{matrix}  \hspace*{-0.8\arraycolsep}
    \end{matrix} 
    \end{bmatrix}  , 
\end{equation*}
where `\( * \)' means those elements could be nonzero but are suppressed because they will not be used or studied later. In the above block-matrix form, \( R_e^{t}\in\R^{n_e\times q_e} \), \( e\in[E] \) and \( t \in \mathbb{N} \). We also let 
\begin{equation}\label{eq:pf.multienv.psi}
  \Psi ^t(R^{t},W) = \begin{bmatrix}
    \begin{matrix}
      \begin{matrix}
        \delta_1\Psi ^t_1(R^{t}_1,W_1) \\ 0 \\ 0 \\ \vdots \\ 0
      \end{matrix} \quad\ \rvline \quad
      \begin{matrix}
        0 \\ \delta_2\Psi ^t_2(R^{t}_2,W_2) \\ 0 \\ \vdots \\ 0
      \end{matrix} \quad\ \rvline \quad
      \begin{matrix}
        0 \\ 0 \\ \delta_3\Psi ^t_3(R^{t}_3,W_3) \\ \vdots \\ 0
      \end{matrix} \quad\ \rvline \quad \cdots  \quad\ \rvline \ \
      \begin{matrix}
        0 \\ 0 \\ 0 \\ \vdots \\ \delta _E \Psi ^t_E(R^{t}_E,W_E)
      \end{matrix}  \hspace*{-0.8\arraycolsep}
    \end{matrix} 
    \end{bmatrix},
\end{equation}
where \( W_1,W_2,\cdots ,W_E \) are from \( W=[W_1^{\T}|W_2^{\T}|\cdots |W_E^{\T}]^{\T} \); $W_e\in \R^{n_e\times q}$. 

We need to partition the state evolution and the Onsager terms too. For \( Z^{t}_{(V)}\in\R^{p\times q} \) whose rows are iid \( \Nscr(0,\frac{1}{\kappa ^2}\Sigma _{(V)}^{t}) \) corresponding to \( V^{t} \), we partition it into \( Z^{t}_{(V)} = \left[Z^{t}_{(V,1)}|Z^{t}_{(V,2)}|\cdots |Z^{t}_{(V,E)}\right] \), where $Z^{t}_{(V,e)}\in\R^{p\times q_e}$. For \( Z^{t}_{(R)}\in \R^{N\times q} \) whose rows are iid \( \Nscr(0,\kappa \, \Sigma ^t_{(R)}) \) corresponding to \( R^{t} \), we similarly make the partition 
\begin{equation*}
  Z^{t}_{(R)} = \begin{bmatrix}
    \begin{matrix}
      \begin{matrix}
        Z^{t}_{(R,1)} \\ * \\ * \\ \vdots \\ *
      \end{matrix} \quad\ \rvline \quad
      \begin{matrix}
        * \\ Z^{t}_{(R,2)} \\ * \\ \vdots \\ *
      \end{matrix} \quad\ \rvline \quad
      \begin{matrix}
        * \\ * \\ Z^{t}_{(R,3)} \\ \vdots \\  *
      \end{matrix} \quad\ \rvline \quad \cdots  \quad\ \rvline \ \
      \begin{matrix}
        * \\ * \\ * \\ \vdots \\  Z^{t}_{(R,E)}
      \end{matrix}  \hspace*{-0.8\arraycolsep}
    \end{matrix} 
    \end{bmatrix},\ Z^{t}_{(R,e)}\in\R^{n_e\times q_e}.
\end{equation*}

 The Gaussian matrices \( Z^{t}_{(V,1)}, Z^{t}_{(V,2)}, \cdots, Z^{t}_{(V,E)} \) are independent, which results from our definition of the function \( \Psi ^t(R^{t}) \). Specifically each of \( Z^{t}_{(V,e)} \in \R^{p\times q_0}\) has iid rows from \( \Nscr(0,\frac{1}{\kappa _e^{2}}\Sigma _{(V,e)}^{t}) \),
 \begin{equation*}
 \begin{aligned}
   \frac{1}{\kappa ^2}\Sigma _{(V)}^{t} &= \diag\left(\frac{1}{\kappa _1^{2}}\Sigma _{(V,1)}^{t},\, \frac{1}{\kappa _2^{2}}\Sigma _{(V,2)}^{t}, \cdots , \frac{1}{\kappa _E^{2}}\Sigma _{(V,E)}^{t}\right),\\
   \frac{1}{\kappa _e^{2}}\Sigma _{(V,e)}^{t} &=
  \frac{1}{\kappa _e^{2}}\lim \frac{1}{n_e} \E\left[ \Psi^{t-1}_{e}(Z^{t-1}_{(R,e)},W_e)^{\T}\Psi^{t-1}_e(Z^{t-1}_{(R,e)},W_e) \right].
 \end{aligned}
 \end{equation*}
 
 The Gaussian matrices \( Z^{t}_{(R,1)}, Z^{t}_{(R,2)}, \cdots, Z^{t}_{(R,E)} \) are independent, because \( Z^{t}_{(R)} \) has iid rows. In fact, even the correlation within each row of \( Z^{t}_{(R)} \) will not come into play anywhere in the state evolution. Marginally, each of \( Z^{t}_{(R,e)} \in \R^{n_e \times q_e} \) has iid rows from \( \Nscr(0, \kappa _e \Sigma ^{t}_{(R,e)}) \), 
 \begin{equation*}
  \kappa _e \Sigma ^{t}_{(R,e)} = \kappa _e \lim \frac{1}{p} \E\left[ \eta ^{t}_e(Z_{(V,1)}^{t},Z_{(V,2)}^{t},\cdots ,Z_{(V,E)}^{t},B)^{\T} \eta ^{t}_e(Z_{(V,1)}^{t},Z_{(V,2)}^{t},\cdots ,Z_{(V,E)}^{t},B) \right]. 
 \end{equation*}  
For the matrices involved in the Onsager terms, we get after some algebra: (\( Z^{t}_{(V)} = [Z^{t}_{(V,1)}|Z^{t}_{(V,2)}|\cdots |Z^{t}_{(V,E)}] \))
\begin{equation*}
\begin{aligned}
  D^{t-1}_{(R)} =&~ \diag\left( 
    \left\{ \frac{\delta _e}{N}  \E \left[ \sum_{j\in[n_e]} \left(\frac{\partial (\Psi ^{t-1}_e)_{j,\cdot }}{\partial (R^{t-1}_e)_{j,\cdot }}\right)(Z_{(R,e)}^{t-1}) \right] :e\in[E] \right\} 
   \right),\\
  D^{t}_{(V)} =&~ \left[ \frac{1}{p\delta _{e_2}} \E \left[  \sum_{j\in[p]} \frac{\partial (\eta ^t_{e_2})_j}{\partial (V_{e_1}^{t})_j} \left( Z^{t}_{(V)},B \right) \right]  \right]_{e_1\in[E],e_2\in[E]}.
\end{aligned}
\end{equation*}

Now we have all the ingredients to write out the multi-environment GLAMP:
\begin{equation*}
  \begin{alignedat}{2}
    V_e^{t} =&~ \frac{1}{\kappa _e} X_e^{\T} \Psi ^{t-1}_e(R_e^{t-1},W_e) - \frac{1}{\kappa _e} \eta ^{t-1}_e(V_1^{t-1},\cdots ,V_E^{t-1},B) [ D^{t-1}_{(R,e)} ]^{\T},
    &&\quad  
    D_{(R,e)}^{t-1} = \frac{1}{n_e}\E \left[ \sum_{j\in[n_e]} \left(\frac{\partial (\Psi _e ^{t-1})_{j,\cdot }}{\partial (R_e^{t-1})_{j,\cdot }}\right)(Z_{(R,e)}^{t-1}) \right],\\
    R_e^{t} = &~ X_e \eta _e ^{t}(V_1^{t},V_2^{t},\cdots ,V_E^{t},B) -  \Psi ^{t-1}_e(R_e^{t-1},W_e) [D_{(V,e)}^{t}]^{\T},
    &&\quad  
    D_{(V,e)}^{t} = \frac{1}{p}\E \left[ \sum_{j\in[p]} \left(\frac{\partial (\eta_e ^{t})_{j,\cdot }}{\partial (V^{t}_e)_{j,\cdot }}\right)(Z_{(V)}^{t}) \right].
  \end{alignedat}
  \end{equation*}

Under Assumption~\ref{asp:GLAMP.asym.multienv.converge} with all the specification made, we already have convergence results taking the form of those from Theorem~\ref{cor:asym.multienv}. However, Assumption~\ref{asp:glamp.converge} is stronger than Assumption~\ref{asp:GLAMP.asym.multienv.converge}, because the former would have required the following limits to exist 
\begin{equation*}
  \lim_{p\to \infty}\frac{1}{p}\E\left[ \eta _{e_1}^{s}(Z_1,B)^{\T} \eta _{e_2}^{r}(Z_2,B) \right],\quad \lim_{N\to \infty}\frac{1}{p}\E\left[ (\eta_{e_1} ^0({V^{0},B}))^{\T} \eta _{e_2}^r(Z_2,B) \right]
\end{equation*}
for any \( e_1,e_2\in[E] \), any constant covariance matrices $\Sigma_1,\Sigma_2\in\R^{q\times q}$ and any $s,r \geq  1$. In contrast, Assumption~\ref{asp:GLAMP.asym.multienv.converge} only requires the same limits to exist for \( e=e_1=e_2 \). 
In other words, we need to relax the assumptions we are working under. 

Instead of going to the proof details of GLAMP, we use a sub-sequence argument to prove the convergence results under Assumption~\ref{asp:GLAMP.asym.multienv.converge}. 
To show 
\begin{equation*}
\begin{aligned}
  \abs*{\phi(V_1^{t},V_2^{t},\cdots ,V_E^{t},V^{0},B) - \E\left[ \phi(Z_{(V,1)}^{t},Z_{(V,2)}^{t},\cdots ,Z_{(V,E)}^{t},V^{0},B) \right]} \probcov 0
\end{aligned}
\end{equation*}
under Assumption~\ref{asp:GLAMP.asym.multienv.converge}, we take any sub-sequence \( (p',n_1',\cdots ,n_E') \) of \( (p,n_1,n_2,\cdots ,n_E) \) and show it admits a further subsequence that makes the above convergence hold almost surely, with the GLAMP iterates and covariance matrices \( \{\Sigma ^t _{(V,e)} :e\in[E],t \in \mathbb{N_+}\} \) involved 
independent of the choice of the sub-sequences. 

Firstly, for any fixed \( t \), there are only finitely many covariance matrices involved in the state evolution. By Assumption~\ref{asp:GLAMP.asym.multienv.converge}, we know there exists a further sub-sequence \( (p'',n_1'',\cdots ,n_E'') \) such that the limits 
\begin{equation*}
\begin{aligned}
  &~\frac{1}{p} (\eta^{0}(V^{0}))^{\T} \eta^{0}(V^{0}) \to {\Sigma }^{0}_{(R)},\ a.s.,\\
  &~\Sigma ^s_{(V)} = \lim_{N \rightarrow \infty}\frac{1}{N}\E\left[ \Psi ^{s-1}(Z_{(R)}^{s-1},W)^{\T} \Psi ^{s-1}(Z_{(R)}^{s-1},W) \right],\\
  &~\Sigma ^s_{(R)} = \lim_{p}\frac{1}{p}\E\left[ \eta ^{s}(Z_{(V)}^{s},B)^{\T} \eta ^{s}(Z_{(V)}^{s},B) \right]
\end{aligned}
\end{equation*}
exists for \( 1\leq s\leq t \) under the specification in Equation~\eqref{eq:pf.multienv.eta}, \eqref{eq:pf.multienv.psi} along the sub-sequence \( (p'',n_1'',\cdots ,n_E'') \). The $E$ diagonal blocks of the shapes \( q_e\times q_e  \) from \( {\Sigma }^{0}_{(R)}, \Sigma ^s_{(V)}, \Sigma ^s_{(R)}  \) are guaranteed by Assumption~\ref{asp:GLAMP.asym.multienv.converge} to exist independently of the sub-sequence \( (p',n_1',\cdots ,n_e') \), but the off-diagonal blocks can depend on  \( (p'',n_1'',\cdots ,n_e'') \). 

Secondly, along the sub-sequence \( (p'',n_1'',\cdots ,n_e'') \), there exists a further sub-sequence \( (p''',n_1''',\cdots ,n_e''') \) s.t. all the following limits exist along \( (p''',n_1''',\cdots ,n_e''') \)
\begin{equation*}
  \lim_{p\to \infty}\frac{1}{p}\E\left[ \eta _{e_1}^{s}(Z_{(V)}^s,B)^{\T} \eta _{e_2}^{r}(Z_{(V)}^r,B) \right],\quad \lim_{N\to \infty}\frac{1}{p}\E\left[ (\eta_{e_1} ^0({V^{0},B}))^{\T} \eta _{e_2}^r(Z^r_{(V)},B) \right]
\end{equation*}
for all \( e_1\neq e_2 \), \( 1\leq s\leq r\leq t \), and for all  \( Z_{(V)}^{s} \) with iid rows drawn from \( \Nscr(0,\Sigma ^s_{(V)}) \), \( Z_{(V)}^{r} \) with iid rows drawn from \( \Nscr(0,\Sigma ^r_{(V)}) \). 
Along \( (p''',n_1''',\cdots ,n_e''') \), we can invoke Corollary~\ref{cor:asym} to get 
\begin{equation*}
  \begin{aligned}
    \abs*{\phi(V_1^{t},V_2^{t},\cdots ,V_E^{t},V^{0},B) - \E\left[ \phi(Z_{(V,1)}^{t},Z_{(V,2)}^{t},\cdots ,Z_{(V,E)}^{t},V^{0},B) \right]} \to 0
  \end{aligned}
  \end{equation*}
in probability along  \( (p''',n_1''',\cdots ,n_e''') \), and almost surely along a further sub-sequence of \( (p''',n_1''',\cdots ,n_e''') \). 

Lastly, we note that \emph{(i)} the multi-environment GLAMP iterates in Equation~\eqref{eq:asym.multienv} do not change its definition with any special sub-sequence of \( (p,n_1,\cdots ,n_E) \), and \emph{(ii)} the state evolution enough to characterize their behavior in Definition~\ref{def:GLAMP.asym.multienv.state.evo} is fully determined by Assumption~\ref{asp:GLAMP.asym.multienv.converge}, independently of the choice of sub-sequences. In other words, whichever sub-sequence we are working on, the covariance matrices \( \Sigma ^t_{(V,1)},\Sigma ^t_{(V,2)},\cdots ,\Sigma ^t_{(V,E)} \) used to characterize \( Z ^t_{(V,1)}, Z ^t_{(V,2)},\cdots ,Z ^t_{(V,E)} \) stay the same. 

As a result, we have shown under Assumption~\ref{asp:GLAMP.asym.multienv.converge}, 
\begin{equation*}
  \begin{aligned}
    \abs*{\phi(V_1^{t},V_2^{t},\cdots ,V_E^{t},V^{0},B) - \E\left[ \phi(Z_{(V,1)}^{t},Z_{(V,2)}^{t},\cdots ,Z_{(E,2)}^{t},V^{0},B) \right]} \probcov 0.
  \end{aligned}
  \end{equation*}
The same argument applies to \( \abs*{\phi(R^{t}_e,W_e) - \E\left[ \phi(Z^{t}_e,W_e) \right]} \probcov 0 \). This completes the proof of Theorem~\ref{cor:asym.multienv}.  
\end{proof}

\section{Proofs of Auxiliary Results in Section~\ref{sec:transfer}}
\label{sec:apd.prelim}

\subsection{Proof of Proposition~\ref{prop:stack.onsager}}
\label{sec:apd.prelim.stack.prop}

We first show all the alternative forms of \( \overline{\delta }^{(p)} _{e} \) are equivalent. Then we show the limit \( \lim_{p}\overline{\delta }^{(p)} _{e} = \overline{\delta } _{e} \) is well-defined. 

  By the extended Stein's Lemma in our Lemma~\ref{lem:apd.stein}, we immediately have \( \frac{1}{\tau _e p} \E\left[ Z_e^{\T} \Sigma _e^{1/2} \overline{\eta } (\tau_1,\tau_2,\cdots ,\tau _E) \right] = \frac{1}{p}\sum_{j\in[p]}\E \left[ \frac{\partial (\Sigma _e^{1/2} \overline{\eta } )_j}{\partial (\tau _e Z_e)_j} \right] \), which is the first equality of Equation~\eqref{eq:prop.stack.forms.of.delta}. The second equality in Equation~\eqref{eq:prop.stack.forms.of.delta} is due to a rewriting of the summed partial derivatives. Then we prove the third line of Equation~\eqref{eq:prop.stack.forms.of.delta}. 
Recall that \( \overline{\eta } (\tau_{1},\tau _2,\cdots ,\tau_E) = \eta(\tau_1 Z_1,\tau_2 Z_2,\cdots ,\tau_E Z_E) \), 
where the function \( \eta \) is defined in Equation~\eqref{eq:def.stack.et}. Denote the support of \( \overline{\eta }  \) by \( S \). By \( \eta \)'s definition and the KKT condition of the `argmin' part, we know \( \overline{\eta }_S \in\R^{p} \) satisfies
\begin{equation*}
\begin{aligned}
  \sum_{e\in[E]} \varpi _e (\Sigma _e)_{S,S} \overline{\eta } _S - \sum_{e\in[E]} \varpi _e (\Sigma _e^{1/2})_{S,\cdot }(\tau _e Z_e) - \sum_{e\in[E]} \varpi _e (\Sigma _e \beta _e)_{S} + \theta \diag({\vec{\lambda } })\sign(\overline{\eta } _S) = 0,
\end{aligned}
\end{equation*}
while the KKT conditions for \( \overline{\eta } _{S^{c}} \) are slack ones. By the implicit function theorem, 
\begin{equation}\label{eq:apd.prelim.stack.prop.partial.deriv}
  \frac{\partial \overline{\eta } _S}{\partial (\tau _e Z_e)} = \bigg[ \sum_{e\in[E]} \varpi _e (\Sigma _{e})_{S,S} \bigg]^{-1} \varpi _e (\Sigma _e^{1/2})_{S,\cdot }.
\end{equation}
Thus from the second line of Equation~\eqref{eq:prop.stack.forms.of.delta}, we know 
\begin{equation}\label{eq:apd.prelim.stack.prop.partial.deriv.trace}
  \overline{\delta } _e^{{(p)}} = \varpi _e\cdot \frac{1}{p}\E \left[  \trace\bigg(\bigg[ \sum_{e\in[E]} \varpi _e (\Sigma _{e})_{S,S} \bigg]^{-1} (\Sigma _e)_{S,S} \bigg)\right],
\end{equation}
thus proving the third line of Equation~\eqref{eq:prop.stack.forms.of.delta}. 

We are left to prove \( \lim_p \overline{\delta }  _e^{(p)} = \overline{\delta }  _e \). We know from Assumption~\ref{asp:stack.for.glamp.onsager} that
\begin{equation*}
\begin{aligned}
  \frac{\partial \mathcal{E} ^{{(p)}}}{\partial \tau _e} = - \frac{1}{p}\E \left[  \varpi _e Z_e^{\T}\left[\Sigma _e^{1/2}\overline{\eta } - (\tau _e Z_e + \Sigma _e^{1/2}\beta _e) \right] \right] = - \frac{1}{p}  \varpi _e \E \left[ Z_e^{\T} \Sigma _e^{1/2} \overline{\eta }  \right] + \varpi _e \tau _e
\end{aligned}
\end{equation*}
converges to some limit as \( p\to \infty \). Thus \( \lim_{p}\frac{1}{p}\E [ Z_e^{\T} \Sigma _e^{1/2} \overline{\eta }  ] \) has a well-defined limit, proving \( \overline{\delta } _e = \lim_p \overline{\delta } ^{(p)}_e \). 

By the third line of Equation~\eqref{eq:prop.stack.forms.of.delta}, \( \sum_{e} \overline{\delta } _e^{(p)} = \frac{1}{p} \E[\abs{S} ] = \frac{1}{p} \E[\norm{\overline{\eta } }_0] \). As \( p\to \infty \), we know \( \lim_{p}\frac{1}{p}\E[\norm{\overline{\eta } }_0 ] = \sum_{e} \overline{\delta } _e \). This completes the proof of Proposition~\ref{prop:stack.onsager}.

\subsection{Proof of Proposition~\ref{prop:2nd.onsager}}
\label{sec:apd.prelim.2nd.prop}
We first study the equivalent forms of $\overline{\gamma } _{\rt}^{(p)}$ and then \( \overline{\gamma } _{\ro}^{(p)} \). For \( \overline{\gamma } ^{(p)}_{\rt} \), we invoke Stein's lemma (our Lemma~\ref{lem:apd.stein}) to get 
\begin{equation*}
  \overline{\gamma } ^{(p)}_{\rt} = \frac{1}{p}\sum_{j\in[p]}\E \left[ \frac{\partial (\Sigma _{1}^{1/2}\xi)_j}{\partial (v_{\rt})_j}(\tau _{\rt}Z_{\rt}, \overline{\eta }) \right]  
  = \frac{1}{p} \E \left[ \trace \left( \frac{\partial \overline{\xi } }{\partial v_{\rt}} \Sigma _{1}^{1/2}  \right) \right]  
  = \frac{1}{p} \E \left[ \trace \left( \frac{\partial (\overline{\xi } -\overline{\eta } )}{\partial v_{\rt}} \Sigma _{1}^{1/2}  \right) \right].
\end{equation*}
In the case of \( \mu _{\rt}=\mu _{\joi} \), we have the KKT condition of \( \overline{\xi }  \) on the support of \( S_{\rt} = \Sup(\overline{\xi } -\overline{\eta } ) \) in Equation~\eqref{eq:apd.prelim.2nd.kt.joi}:
\begin{equation}\label{eq:apd.prelim.2nd.kt.joi}
  (\Sigma _{1})_{S_\rt,S_\rt}(\overline{\xi } - \overline{\eta } )_{S_\rt} - (\Sigma ^{1/2}_{1})_{S_\rt,\cdot }(\tau _{\rt}Z_{\rt}) - (1 - \kappa _{1}\gamma _{\ro}) (\Sigma _1)_{S_\rt,\cdot } (\beta _{1} - \overline{\eta } ) + \theta _{\rt} \lambda _{\rt} \sign(\overline{\xi } - \overline{\eta } )_{S_\rt}=0.
\end{equation}
Hence by the implicit function theorem,
\begin{equation}\label{eq:apd.prelim.2nd.partial.deriv.joi}
  \frac{\partial (\overline{\xi } -\overline{\eta } )_{S_{\rt}}}{\partial v_{\rt}}=\left[ (\Sigma _{1})_{S_\rt,S_\rt} \right]^{-1} (\Sigma ^{1/2}_{1})_{S_\rt,\cdot } \implies \overline{\gamma } _{\rt}^{(p)} = \frac{1}{p}\E \left[ \trace \left( \left[ (\Sigma _{1})_{S_\rt,S_\rt} \right]^{-1} (\Sigma _{1})_{S_\rt,S_\rt} \right) \right] = \frac{\E[\abs{S_{\rt}} ]}{p}.
\end{equation}
In the case of \( \mu _{\rt} = \mu _{\ada} \), we have the KKT condition of \( \xi \) on the support \( S_{\rt} = \Sup(\overline{\xi } ) \) in Equation~\eqref{eq:apd.prelim.2nd.kt.ada}, where \( D_{\mu} \) has shown up in Equation~\eqref{eq:def.2nd.ada} as a diagonal matrix, \( (D_{\mu})_{j} = \mu(\abs{(\stac)_j} ) \). 
\begin{equation}\label{eq:apd.prelim.2nd.kt.ada}
  (\Sigma _{1})_{S_\rt,S_\rt}\overline{\xi }_{S_\rt} - (\Sigma ^{1/2}_{1})_{S_\rt,\cdot }(\tau _{\rt}Z_{\rt}) - (\Sigma _{1})_{S_{\rt},\cdot }\beta_1 + \kappa _{1}\gamma _{\ro} (\Sigma _1)_{S_\rt,\cdot } (\beta _{1} - \overline{\eta } ) + \theta _{\rt}( D_{\mu})_{S_{\rt},S_{\rt}} \sign(\overline{\xi })_{S_\rt}=0.
\end{equation}
So we similarly have 
\begin{equation}\label{eq:apd.prelim.2nd.partial.deriv.ada}
  \frac{\partial \overline{\xi } _{S_{\rt}}}{\partial v_{\rt}}=\left[ (\Sigma _{1})_{S_\rt,S_\rt} \right]^{-1} (\Sigma ^{1/2}_{1})_{S_\rt,\cdot } \implies \overline{\gamma } _{\rt}^{(p)} = \frac{1}{p}\E \left[ \trace \left( \left[ (\Sigma _{1})_{S_\rt,S_\rt} \right]^{-1} (\Sigma _{1})_{S_\rt,S_\rt} \right) \right] = \frac{\E[\abs{S_{\rt}} ]}{p}.
\end{equation}

After studying \( \overline{\gamma } _{\rt} \), we look at \( \overline{\gamma } ^{(p)}_{\ro} \), 
We first invoke Stein's lemma and the chain rule of derivatives to get 
\begin{equation*}
  \begin{aligned}
    \overline{\gamma } ^{(p)}_{\ro} =&~ \frac{1}{p}\sum_{j\in[p]}\E \left[ \left( \frac{\partial (\Sigma _{1}^{1/2}\xi)_j}{\partial \hat{\beta } }(\tau _{\rt}Z_{\rt}, \overline{\eta }) \right)^{\T} \frac{\partial \overline{\eta } }{\partial (\tau _1^{*}Z_{1})_j}  \right]
    = \frac{1}{p} \E \left[ \trace \left( \frac{\partial \xi}{\partial \hat{\beta } }(\tau _{\rt}Z_{\rt}, \overline{\eta } ) \frac{\partial \overline{\eta } }{\partial (\tau _1^{*}Z_{1})} \Sigma _{1}^{1/2} \right) \right].
  \end{aligned}
    \end{equation*}

    In the case of \( \mu _{\rt} = \mu _{\joi} \), we know 
    \begin{equation*}
      \frac{\partial \xi}{\partial \hat{\beta } }(\tau _{\rt}Z_{\rt}, \overline{\eta } ) \frac{\partial \overline{\eta } }{\partial (\tau _1^{*}Z_{1})}  
      = \frac{\partial (\overline{\xi } - \overline{\eta } )}{\partial \hat{\beta } }(\tau _{\rt}Z_{\rt}, \overline{\eta } ) \frac{\partial \overline{\eta } }{\partial (\tau _1^{*}Z_{1})}  +  \frac{\partial \overline{\eta } }{\partial (\tau _1^{*}Z_{1})},
    \end{equation*}
    and \begin{equation*}
      \frac{\partial (\overline{\xi } -\overline{\eta } )_{S_\rt}}{\partial \overline{\eta } } = - (1 - \kappa _1 \gamma _{\ro}) \left[ (\Sigma _{1})_{S_\rt,S_\rt} \right]^{-1} (\Sigma _{1})_{S_\rt,\cdot }.
    \end{equation*}
Combined with Equation~\eqref{eq:apd.prelim.stack.prop.partial.deriv}, we know 
(with \( S=\Sup(\overline{\eta } ) \), \( \overline{\Sigma } =\sum_{e}\varpi _e^{*}\Sigma _{e} \) from Section~\ref{sec:apd.prelim.stack.prop})
\begin{equation}
\begin{aligned}\label{eq:prop.2step.gamma.I.joi}
  \overline{\gamma } _{\ro}^{(p)} =&~ \overline{\delta } _{1}^{(p)} - \varpi _1^{*}(1 - \kappa _1 \gamma _{\ro})\cdot \frac{1}{p}\E \left[ \trace \left( \left[ (\Sigma _{1})_{S_\rt,S_{\rt}} \right]^{-1} (\Sigma _{1})_{S_\rt,S} \left[ \overline{\Sigma } _{S,S}\right]^{-1} (\Sigma _{1})_{S,S_{\rt}} \right) \right].
\end{aligned}
\end{equation}
In the case of \( \mu _{\rt} = \mu _{\ada} \), we similarly have 
\begin{equation*}
  \frac{\partial \overline{\xi }_{S_\rt}}{\partial \overline{\eta } } = \kappa _1 \gamma _{\ro} \left[ (\Sigma _{1})_{S_\rt,S_\rt} \right]^{-1} (\Sigma _{1})_{S_\rt,\cdot }
\end{equation*} 
and
\begin{equation}
  \begin{aligned}\label{eq:prop.2step.gamma.I.ada}
    \overline{\gamma } _{\ro}^{(p)} =&~  \varpi _1^{*}\kappa _1 \gamma _{\ro} \cdot  \frac{1}{p}\E \left[ \trace \left( \left[ (\Sigma _{1})_{S_\rt,S_{\rt}} \right]^{-1} (\Sigma _{1})_{S_\rt,S} \left[ \overline{\Sigma } _{S,S}\right]^{-1} (\Sigma _{1})_{S,S_{\rt}} \right) \right].
  \end{aligned}
  \end{equation}

 Now we have studied the equivalent forms of \( \overline{\gamma } _{\ro}^{(p)} \), \( \overline{\gamma } _{\rt}^{(p)} \).  We then move on to prove the existence of the limits. First off, we show the limit \( \lim_{p} \overline{\gamma } ^{(p)}_{\ro} = \overline{\gamma } _{\ro} \). We look at \( ({\partial \mathcal{E}^{(p)} _{\rt}})/({\partial \zeta}) \) which admits a limit by Assumption~\ref{asp:2nd.for.glamp.onsager}. After some algebra, 
\begin{equation*}
\begin{aligned}
  \frac{{\partial \mathcal{E}^{(p)} _{\rt}}}{{\partial \zeta}} =&~ -  \frac{\tau _{\rt}}{p}\E \left[ \left(Z_{1}- \zeta Z_{1}'/\sqrt{1-\zeta ^2}\right)^{\T} \Sigma _{1}^{1/2} \overline{\xi }  \right] + \tau _{\rt}\, \kappa _{1}\, \gamma _{\ro}\, \tau _{1}^{*}\, \overline{\delta }^{(p)} _{1} = \tau _{\rt}\tau _{1}^{*} \left( \kappa _{1}\,\overline{\delta } ^{{(p)}}_{1} \gamma _{\ro} - \overline{\gamma } _{\ro}^{(p)} \right).
\end{aligned}
\end{equation*}
By Assumption~\ref{asp:2nd.for.glamp.onsager}, we know there exists some \( \overline{\gamma } _{\ro} \) such that \( \lim_{p} \overline{\gamma } ^{(p)}_{\ro} = \overline{\gamma } _{\ro} \). Next we show there exists some value \( \overline{\gamma } _{\rt} \geq 0 \), s.t. \( \lim_{p} \overline{\gamma } ^{(p)}_{\rt} = \overline{\gamma } _{\rt} \), which is implied by 
\begin{equation*}
\begin{aligned}
  \frac{{\partial \mathcal{E}^{(p)} _{\rt}}}{{\partial \tau _{\rt}}} =&~ \tau _{\rt}(1- \overline{\gamma } _{\rt}^{(p)}) + \zeta \cdot  \tau _{E}^{*} (\kappa _{E}\overline{\delta } _{E}^{(p)}\gamma _{\ro} - \overline{\gamma } _{\ro}^{(p)})
\end{aligned}
\end{equation*}
and Assumption~\ref{asp:2nd.for.glamp.limits}.

\section{Proof of Theorem~\ref{thm:stack}}
\label{sec:apd.stack}

We work under Assumptions~\ref{asp:prelim.transfer.model},\ref{asp:stack.for.glamp},\ref{asp:stack.converge} for the whole Section~\ref{sec:apd.stack}. 
We also fix  \( \tau _e = \tau _e^{*},\ \varpi _e = \varpi _e^{*},\forall\,e\in[E] \)  and \( \theta = \theta ^* \) in the  definitions of the functions \( \eta, \overline{\eta } , \mathcal{E} ^{{(p)}}, \mathcal{E},\overline{\delta } _e^{(p)},\overline{\delta } _e,H^{(p)}, H \) when we work on them. Among the following, Lemmas~\ref{lem:apd.stack.bounded.lasso},~\ref{lem:apd.stack.glamp},~\ref{lem:stack.final.converge} are used to prove Theorem~\ref{thm:stack}, whereas Lemmas ~\ref{lem:apd.stack.cauchy},~\ref{lem:stacked.lasso.subgradient},~\ref{lem:stacked.lasso.bound.on.sc2} are auxiliary lemmas needed to prove Lemma~\ref{lem:stack.final.converge}.

\begin{lemma}[Boundedness of \( \stac \)]\label{lem:apd.stack.bounded.lasso}
  Under Assumption~\ref{asp:prelim.transfer.model}, there exists \( M>0 \) s.t. for \( \stac \) defined in Equation~\eqref{eq:def.stack.los}, \( \limsup_{p}\frac{1}{p}\norm{\stac}_2^{2} \leq M \) a.s. 
\end{lemma}

Lemma~\ref{lem:apd.stack.bounded.lasso} is proved in Section~\ref{sec:apd.pf.stack.bounded.lasso}. 

\begin{definition}[Multi-Environment GLAMP Formulation]\label{def:stack.glamp}
  Under Assumptions~\ref{asp:prelim.transfer.model},\ref{asp:stack.for.glamp},\ref{asp:stack.converge.fixed.point}, we define a multi-environment GLAMP instance as below:
\begin{equation}\label{eq:def.stack.glamp}
\begin{aligned}
  v_e^{t} =&~ \ix{e}^{\T} (y_e - r_e^{t-1}) - \Sigma _e^{1/2} \beta _e + \Sigma _e^{1/2} \eta(v_1^{t-1}, v_2^{t-1}, \cdots ,v_E^{t-1}), \\
  r_e^{t} =&~ \ix{e} \Sigma _e^{1/2} \eta(v_1^{t}, v_2^{t}, \cdots ,v_E^{t}) - \kappa _e \overline{\delta }^{(p)} _e (y_e - r_e^{t-1})
\end{aligned}
\end{equation}
for all \( t \in \mathbb{N_+} \) and \( e\in[E] \). In Equation~\eqref{eq:def.stack.glamp}, \( \overline{\delta } _e^{(p)} = \frac{1}{\tau _e p} \E\left[ Z_e^{\T} \Sigma _e^{1/2} \overline{\eta } (\tau_1^{*},\tau_2^{*},\cdots ,\tau^{*} _E) \right] \) has been in Proposition~\ref{prop:stack.onsager} and specified with \( \tau _e = \tau _e^{*},\varpi _e=\varpi _e^{*}\,\forall\,e\in[E] \) and \( \theta=\theta ^* \). 

As initialization, we define $\eta ^0$ as follows:
Let \( Z_1,Z_2,\cdots ,Z_E\iidsim \Nscr(0,I_p) \) be independent of all the \( \{ \ix{e}: e\in[E] \}  \), and define \( \eta ^{0} = \eta(\tau _1^{*}Z_1,\tau _e^{*}Z_2,\cdots ,\tau _E^{*}Z_E) \). We initialize 
\( r_e^{0} = \ix{e} \Sigma _e^{1/2} \eta^0 \), \( v_e^{1} = \ix{e}^{\T}(y_e - r_e^{0}) - \Sigma _e ^{1/2} \beta _e + \Sigma _e^{1/2} \eta^0 \), \( \forall\,e\in[E] \).

As a shorthand, define \( \eta^t = \eta(v_1^{t},v_2^{t},\cdots ,v_E^{t}) \) for \( t \in \mathbb{N_+} \). 
\end{definition}

\begin{definition}[State Evolution]\label{def:apd.stack.stat.evo} 
  Let \( \Sigma _{(V,e)} \in (\R^{\N_+})\times(\R^{\N_+}) \) be infinite-dimensional matrices defined as below: \( \forall\,e\in[E] \),

  For all \( i\in\N_+ \), \( (\Sigma _{(V,e)})_{i,i}=(\tau _e^{*})^{2} \).

  For all \( i<j \), \( i,j\in\N_+ \), \( (\Sigma _{(V,e)})_{i,j}=(\Sigma _{(V,e)})_{j,i}=(\tau _e^{*})^{2} \rho _e ^{i,j} \), where \(\rho_e^{i,j}\in[0,1]\), \( \rho ^{i,j} = (\rho ^{i,j}_{1},\rho ^{i,j}_{2},\cdots ,\rho ^{i,j}_{E})\in[0,1]^{E} \). For the definition of \( \{\rho ^{i,j}\} \): \emph{(i)} For all \( j\geq 2 \), \( \rho ^{1,j}=H(\vec{0}_E) \), where \( \vec{0}_E\in\R^{E} \)  is a vector of all zeros; \emph{(ii)} For all \( 2\leq i<j \), \( \rho ^{i,j}=H(\rho ^{i-1,j-1}) \). \( H \) has been defined in Assumption~\ref{asp:stack.converge.cauchy}.

  In particular, we define the \( t \)-dimensional top-left block of \( \Sigma _{(V,e)} \) as \( \Sigma ^t _{(V,e)} = ( \Sigma _{(V,e)})_{[t],[t]} \), \( t \in \mathbb{N_+} \). 
\end{definition}

\begin{lemma}[GLAMP Convergence]\label{lem:apd.stack.glamp}
  Under Assumptions~\ref{asp:prelim.transfer.model},\ref{asp:stack.for.glamp},\ref{asp:stack.converge}, and Definitions~\ref{def:stack.glamp},\ref{def:apd.stack.stat.evo},
 \( \forall\,t \in \mathbb{N_+} \), \( \forall\, \) sequence of order-\( k \) pseudo-Lipschitz functions \( \phi:(\R^{p})^{(t+1)E+1}\to\R \), \( k\geq 1 \),
 \begin{equation*}
\begin{alignedat}{2}
  &\phi\left(v_1^{t},\cdots,v_E^{t},v_1^{t-1},\cdots,v_E^{t-1},\cdots ,v_1^{1},\cdots,v_E^{1},\eta^0,\beta _1,\beta_2,\cdots ,\beta _E \right)&& \\
 &= \E \left[ \phi\left( Z^{t}_1,\cdots ,Z^{t}_E,Z^{t-1}_1,\cdots ,Z^{t-1}_E,\cdots ,Z^{1}_1,\cdots ,Z^{1}_E,\eta(Z^{0}),\beta_1,\beta _2,\cdots ,\beta _E \right) \right] + \smalop.
\end{alignedat}
\end{equation*}
 In particular,
 \( \forall\, \) sequence of order-\( k \) pseudo-Lipschitz functions \( \phi:(\R^{p})^{t+E+1}\to\R \), \( k\geq 1 \),
\begin{equation*}
  \phi\left(\eta^{t},\eta^{t-1},\cdots ,\eta^0,\beta _1,\beta_2,\cdots ,\beta _E \right) = \E \left[ \phi\left( \eta(Z^{t}),\eta(Z^{t-1}),\cdots ,\eta(Z^{0}),\beta_1,\beta _2,\cdots ,\beta _E \right) \right] + \smalop.
\end{equation*}
In the expcetations on the RHS, for any \( 0\leq i\leq t \), we have used \( Z^{i}=[Z^{i}_1|Z^{i}_2|\cdots |Z^{i}_E]\in\R^{p\times E}\) being a Gaussian matrix, and \( \eta(Z^{i}) = \eta(Z^{i}_1,Z^{i}_2,\cdots,Z^{i}_E) \). For the joint distribution of these matrices: \emph{(i)} \( Z^{0}\disteq Z\cdot\diag(\tau _1^{*},\tau _2^{*},\cdots ,\tau _E^{*}) \) where \( Z\in\R^{p\times E} \) is a matrix of iid \( \Nscr(0,1) \) elements; \emph{(ii)} \( Z^{(0)} \independent \{ Z^{i}:i\geq 1 \}  \); \emph{(iii)} For \( Z^{i}=[Z^{i}_1|Z^{i}_2|\cdots |Z^{i}_E]\in\R^{p\times E}\), if we re-arrange the columns of $\{Z^i\}_{i=1}^t$ into $\{Z_e\}_{e=1}^E$ where \( Z_{e}=[Z_e^{1}|Z_e^{2}|\cdots |Z_e^{t}]\in\R^{p\times t} \), then \( Z_e \) has iid rows drawn from \( \Sigma ^t_{(V,e)}\in\R^{t\times t} \) defined in Definition~\ref{def:apd.stack.stat.evo}, and \( \{ Z_e:e\in[E] \}  \) are \( E \) independent Gaussian matrices. 
\end{lemma}

Lemma~\ref{lem:apd.stack.glamp} is proved in Section~\ref{sec:apd.pf.stack.glamp}. 

\begin{lemma}[The Cauchy Property]\label{lem:apd.stack.cauchy}
  Under Assumptions~\ref{asp:prelim.transfer.model},\ref{asp:stack.for.glamp},\ref{asp:stack.converge}, and Definitions~\ref{def:stack.glamp},\ref{def:apd.stack.stat.evo}, 
\( \forall\,\epsilon>0 \), \( \exists \,T=T(\epsilon) \). s.t. \( \forall\,j>i\geq T \), $\forall e \in [E]$,
  \begin{equation*}
\begin{alignedat}{2}
  &\lim_{p\to \infty}\PR \left[ \frac{1}{p}\norm*{v_e^{j}-v_e^{i}}_2^{2} \leq \epsilon \right]=1,\quad &&\lim_{p\to \infty}\PR \left[  \frac{1}{n_e}\norm*{r_e^{j}-r_e^{i}}_2^{2} \leq \epsilon\right] =1, \\
  &\lim_{p\to \infty}\PR\bigg[ \frac{1}{p}\norm*{\eta  ^{j}-\eta ^i}_2^{2} \leq \epsilon \bigg] = 1,\quad  &&\max_{e\in[E]}\{1 - \rho _e^{i,j} \} \leq \epsilon.
\end{alignedat}
  \end{equation*}
\end{lemma}

Lemma~\ref{lem:apd.stack.cauchy} is proved in Section~\ref{sec:apd.pf.stacked.lasso.cauchy}. 

\begin{lemma}[Sub-Gradient Going to Zero] \label{lem:stacked.lasso.subgradient}
  Under Assumptions~\ref{asp:prelim.transfer.model},\ref{asp:stack.for.glamp},\ref{asp:stack.converge}, and Definitions~\ref{def:stack.glamp},\ref{def:apd.stack.stat.evo}, 
  define a vector \( s^{t}\in\R^{p} \) by 
  \begin{equation}\label{eq:def.stack.subgradient}
  s^{t} = \frac{1}{\theta ^*}\left[ \diag(\vec{\lambda } )^{-1} \right]\sum_{e\in[E]} \varpi _e^{*} \left( \Sigma _e^{1/2} v_e^{t} + \Sigma _e(\beta _e - \eta^t) \right).
  \end{equation}
  Then \( s^{t}\in \partial \norm{\eta^t}_1 \). Moreover, if we define \( \nabla L(\eta^t) = - \sum_{e\in[E]} \pi _e \Sigma_e^{1/2} \ix{e}^{\T}(y_e - \ix{e} \Sigma _e^{1/2}\eta^t) + \lambda s^{t} \) as a sub-gradient of \( L \) defined in Equation~\eqref{eq:def.stack.los}, then \( \forall\, \epsilon>0 \), \( \exists\, T \) s.t. \( \forall\,t\geq T \), \( \lim_{p}\PR\left[ \frac{1}{p}\norm{\nabla L(\eta^t)}_2^{2} < \epsilon \right] = 1 \). 
\end{lemma}

Lemma~\ref{lem:stacked.lasso.subgradient} is proved in Section~\ref{sec:apd.pf.stacked.lasso.subgradient}.

\begin{lemma}[Bounding the Sub-Gradient's Magnitude]\label{lem:stacked.lasso.bound.on.sc2}
  For any constant \( c\in(0,1) \), define \( S^{t}(c) = \{ j\in[p] : \abs{s^{t}_j} \geq 1 - c  \}  \) where \( s^{t} \) has been defined in Equation~\eqref{eq:def.stack.subgradient}. 
  Under Assumptions~\ref{asp:prelim.transfer.model},\ref{asp:stack.for.glamp},\ref{asp:stack.converge}, and Definitions~\ref{def:stack.glamp},\ref{def:apd.stack.stat.evo},  
  for any \( \epsilon>0 \), there exists \( c = c(\epsilon) \in(0,1) \) invariant to the time step \( t \) s.t. for any time step \( t \in \mathbb{N_+} \), we have \( \lim_{p} \PR\left[\abs{S^{t}(c)}/p \leq \sum_{e\in[E]} \overline{\delta } _e + \epsilon  \right] = 1 \).
\end{lemma}

Lemma~\ref{lem:stacked.lasso.bound.on.sc2} is proved in Section~\ref{sec:apd.pf.stacked.lasso.bound.on.sc2}. 

\begin{lemma}\label{lem:stack.final.converge}
  Under Assumptions~\ref{asp:prelim.transfer.model},\ref{asp:stack.for.glamp},\ref{asp:stack.converge}, and Definitions~\ref{def:stack.glamp},\ref{def:apd.stack.stat.evo}, for any $\epsilon>0$, 
  \( \exists \, T \) s.t. \( \forall\,t\geq T \), \( \lim_{p} \PR \left[ \frac{1}{p}\norm{\eta^t-\stac}_2^{2} < \epsilon\right] = 1  \).  
\end{lemma}

Lemma~\ref{lem:stack.final.converge} is proved in Section~\ref{sec:apd.pf.stacked.lasso.final.converge}. 

\begin{proof}[\textbf{Proof of Theorem~\ref{thm:stack}}]
  We work with small \( \epsilon \) and \( t\geq T \) required by Lemma~\ref{lem:stack.final.converge}. 
  By Lemma~\ref{lem:apd.stack.bounded.lasso} and Lemma~\ref{lem:apd.stack.glamp}, we have 
 \( \lim_{p}\PR\left[\max\{\frac{1}{p}\norm{\stac}_2^{2}, \frac{1}{p}\norm{\eta^{t}}_2^{2}\} \leq M \right]=1 \). By Assumption~\ref{asp:prelim.transfer.model.asym.beta}, we can choose \( M \) large enough s.t. \( \limsup_{p}\frac{1}{p}\norm{\beta _e}_2^{2} < M\ \forall\,e\in[E] \) as well. So we can choose \( \epsilon \) small enough s.t. with probability approaching \( 1 \) as \( p\to \infty \) at a fixed \( t \), for any sequence of order-\( k \) pseudo-Lipschitz functions \( \phi:(\R^{p})^{E+1}\to\R \) \( (k\geq 1) \),
 \begin{equation*}
\begin{aligned}
  &~\abs*{\phi(\eta^t,\beta_1,\cdots ,\beta _E) - \phi(\stac,\beta_1,\cdots ,\beta _E)} \\
  \leq&~ L \left[ 1 + \sum_{e\in[E]} \left( \frac{\norm{\beta _e}_2}{\sqrt{p}} \right)^{k-1} + \left( \frac{\norm{\eta^t}_2}{\sqrt{p}} \right)^{k-1} + \left( \frac{\norm{\stac}_2}{\sqrt{p}} \right)^{k-1} \right] 
  \frac{\norm{\eta^t - \stac}_2}{\sqrt{p}} \leq \frac{\epsilon}{2}.
\end{aligned}
 \end{equation*}
 Also with probability going to 1 as $p\to \infty$, we have \( \abs{\phi(\eta^t,\beta_1,\cdots ,\beta _E) - \E\left[\phi(\overline{\eta } ,\beta_1,\cdots ,\beta _E)\right]} < \epsilon/2 \) by Lemma~\ref{lem:apd.stack.glamp}. Thus for any small enough \( \epsilon \), we have \( \abs{\phi(\stac,\beta_1,\cdots ,\beta _E) - \E\left[\phi(\overline{\eta } ,\beta_1,\cdots ,\beta _E)\right]} < \epsilon \) with probability going to 1 as $p\to \infty$. This proves Theorem~\ref{thm:stack}. 
\end{proof}

\subsection{Proof of Auxiliary Lemmas in Section~\ref{sec:apd.stack}}\label{sec:apd.pf.stack.lems}

\subsubsection{Proof of Lemma~\ref{lem:apd.stack.bounded.lasso}} \label{sec:apd.pf.stack.bounded.lasso}
  
The proof follows the same steps as those of Lemma 3.2, \citet{bayati2011lasso}. First of all, 
\begin{equation*}
\begin{aligned}
\limsup_{p} L(\stac)/p \leq \limsup_{p} L(0)/p \leq \sum_{e\in[E]} \pi _e \cdot \frac{\norm{y_e}_2^{2}}{2p} \leq \sum_{e\in[E]} \pi _e \cdot \left[ \frac{\norm{w_e}_2^{2}}{p} + \frac{\sigma _{\max }(\ix{e})^{2}c_1^{2} \norm{\beta _e}_2^{2}}{p} \right].
\end{aligned}
\end{equation*}
By Lemma~\ref{lem:apd.sing.val.limit.concent}, \( \lim_{p}\sigma _{\max}(\ix{e}) = 1 + \sqrt{\kappa _e},a.s. \). By Assumption~\ref{asp:prelim.transfer.model}, we know \( \limsup_{p} L(\stac)/p < +\infty \). Now if 
\( \kappa _e<1 \) for any \( e:\pi _e>0 \), then we know \( \liminf_{p} \sigma _{\min }(\ix{e} \Sigma _e^{1/2})\geq (1-\sqrt{\kappa _e})/c_1^{1/2} > 0 \) by Lemma~\ref{lem:apd.sing.val.limit.concent} and Assumption~\ref{asp:prelim.transfer.model}, and thus \( \frac{1}{p}\norm{\stac}_2^{2} \leq B_1\cdot \frac{1}{p}\norm{\ix{e}\Sigma _{e}^{1/2}\stac} _2^{2}\ a.s. \) for some constant \( B_1> c_1/(1-\sqrt{\kappa _e})^{2} \) as \( p\to \infty \),
\begin{equation}
\begin{aligned}\label{eq:apd.pf.stack.bounded.mid}
   \frac{1}{p}\norm{\ix{e}\Sigma _{e}^{1/2}\stac} _2^{2} \leq \frac{2}{p}\left[ \norm{y_e - \ix{e}\Sigma _{e}^{1/2}\stac} _2^{2} + \norm{y_e}_2^{2} \right]
  \leq  \frac{2}{\pi _e} \left[ \frac{L(\stac)}{p} + \frac{L(0)}{p} \right],
\end{aligned}
\end{equation}

and we have proved Lemma~\ref{lem:apd.stack.bounded.lasso}. From now on, we focus on the the case of \( \kappa _e \geq 1 \) for all the environments \( e : \pi _e > 0\).
  
  For any \( e : \pi _e > 0 \), let \( P_e^{\|} \) be  the projection matrix onto \( \mathrm{\ker }(\ix{e} \Sigma _e^{1/2}) \) and \( P_{e}^{\bot} \) onto \( \mathrm{\ker }(\ix{e} \Sigma _e^{1/2})^{\bot} \). Then \( \frac{1}{p}\norm{\stac}_2^{2} \leq \frac{1}{p}\norm{P_e^{\|}\stac}_2^{2} + \frac{1}{p} \norm{P_e^{\bot}\stac}_2^{2} \). At any finite \( p \), \( \mathrm{\ker }(\ix{e} \Sigma _e^{1/2}) \) is still a uniformly random space, so Kashin's theorem (c.f. Lemma~\ref{thm:kashin}) applies just like in the proof of Lemma 3.2, \citet{bayati2011lasso}. For some constant $B_1$, the relationship below holds almost surely as \( p\to \infty \), 
  \begin{equation*}
\begin{aligned}
  \frac{1}{p}\norm{\stac}_2^{2} \leq B_1 \left( \frac{1}{p}\norm{P_e^{\|}\stac}_1 \right)^{2} + \frac{1}{p} \norm{P_e^{\bot}\stac}_2^{2} \leq 2B_1 \left( \frac{1}{p}\norm{\stac}_1  \right)^{2} + \frac{2B_1+1}{p}\norm{P_e^{\bot}\stac}_2^{2}.
\end{aligned}
  \end{equation*}
  The first inequality in the above line uses Kashin's theorem, and the second has omitted some details because they are the same as those in \citet{bayati2011lasso}. For the former term, \(\limsup_{p}\frac{1}{p}\norm{\stac}_1 < +\infty, a.s. \) because \( \limsup_p\frac{1}{p}L(\stac) < +\infty \). For the latter term, since
 \( \mathrm{\ker }(\ix{e} \Sigma _e^{1/2})^{\bot} = \mathrm{range}(\Sigma _e^{1/2}\ix{e}^{\T}) \), \( P_e^{\bot} =\Sigma _e^{1/2}\ix{e}^{\T}(\ix{e}\Sigma _e \ix{e}^{\T})^{-1}\ix{e}\Sigma _e^{1/2}  \). By Lemma~\ref{lem:apd.sing.val.limit.concent} and Assumption~\ref{asp:prelim.transfer.model}, there exists \( B_2>0 \) s.t. 
 \begin{equation*}
  \norm{P_e^{\bot}\stac}_2 \leq B_2\norm{\ix{e} \Sigma _e^{1/2} P_e^{\bot}\stac}_2 = B_2\norm{\ix{e} \Sigma _e^{1/2} \stac}_2, a.s.,
 \end{equation*}
and \( \norm{\ix{e} \Sigma _e^{1/2} \stac}_2 \) has been bounded by Equation~\eqref{eq:apd.pf.stack.bounded.mid}. Now we have shown \( \limsup_{p} \frac{1}{p}\norm{\stac}_2^{2} < +\infty \). 

\subsubsection{Proof of Lemma~\ref{lem:apd.stack.glamp}}\label{sec:apd.pf.stack.glamp}
\begin{proof}[\textbf{Proof of Lemma~\ref{lem:apd.stack.glamp}}]
The way to prove Lemma~\ref{lem:apd.stack.glamp} is to show the GLAMP instance defined in Definition~\ref{def:stack.glamp} is a legitimate multi-environment GLAMP as defined in Definition~\ref{def:GLAMP.asym.multienv}. However, we are faced with one more technicality: we have used external randomness to initialize the iterates. As a result, the original GLAMP would give results conditioning on a realization of \( \eta^0 \), which is not exactly Lemma~\ref{lem:apd.stack.glamp}. 
  
We use a sub-sequence argument to address this technicality. We know from  Lemma~\ref{lem:apd.lem23} and Assumption~\ref{asp:stack.for.glamp} that, (\( \overline{\eta } =\overline{\eta } (\tau_1^{*},\tau _2^{*},\cdots ,\tau _E^{*}) \))
\begin{equation}\label{eq:apd.pf.stack.init.limits}
\begin{alignedat}{2}
\frac{1}{p} \beta _e ^{\T}\Sigma _e \eta^0 =&~ \frac{1}{p}\E\left[ \beta _e ^{\T}\Sigma _e \overline{\eta }  \right]+\smalop &&\implies  \plim_{p\to \infty} \frac{1}{p} \beta _e ^{\T}\Sigma _e \eta^0 = \lim_{p\to \infty}\frac{1}{p}\E\left[ \beta _e ^{\T}\Sigma _e \overline{\eta }  \right], \\
\frac{1}{p} (\eta^0) ^{\T}\Sigma _e \eta^0 =&~ \frac{1}{p}\E\left[ \overline{\eta } ^{\T}\Sigma _e \overline{\eta }  \right]+\smalop &&\implies  \plim_{p\to \infty} \frac{1}{p} (\eta^0) ^{\T}\Sigma _e \eta^0 = \lim_{p\to \infty}\frac{1}{p}\E\left[ \overline{\eta } ^{\T}\Sigma _e \overline{\eta }  \right], \\
\frac{1}{p} (\eta^0) ^{\T}\Sigma _e \E[\overline{\eta } ] =&~ \frac{1}{p}\E\left[ \overline{\eta }_{(1)} ^{\T}\Sigma _e \overline{\eta }_{(2)}  \right]+\smalop &&\implies  \plim_{p\to \infty} \frac{1}{p} (\eta^0) ^{\T}\Sigma _e \E[\overline{\eta } ] = \lim_{p\to \infty}\frac{1}{p}\E\left[ \overline{\eta }_{(1)} ^{\T}\Sigma _e \overline{\eta }_{(2)}  \right],
\end{alignedat}
\end{equation}
where in the last line, \( \overline{\eta } _{(1)} \) and \( \overline{\eta } _{(2)} \) are mutually independent copies of \( \overline{\eta }  \). The existence of the limit \( \lim_p\frac{1}{p}\E\left[ \overline{\eta }_{(1)} ^{\T}\Sigma _e \overline{\eta }_{(2)}  \right] \) is implied by \( \lim_{p}H^{{(p)}} = H \) in Assumption~\ref{asp:stack.converge.cauchy} combined with Assumption~\ref{asp:stack.for.glamp.limits}. 

Recall that a sequence of random variables converges in probability iff any of its sub-sequences admit a further sub-sequence that converges almost surely. 
If we prove Lemma~\ref{lem:apd.stack.glamp} along any sub-sequence of \( (p,n_1,n_2,\cdots ,n_E) \) on which the three limits involving \( \eta^0 \) in Equation~\eqref{eq:apd.pf.stack.init.limits} all hold almost surely (instead of merely in probability), we will have proved Lemma~\ref{lem:apd.stack.glamp} itself. 
This is because any sub-sequence of \( (p,n_1,n_2,\cdots ,n_E) \) firstly admits a further sub-sequence along which all the limits in Equation~\eqref{eq:apd.pf.stack.init.limits} hold almost surely. 
If we prove the desired convergence in Lemma~\ref{lem:apd.stack.glamp} along the further sub-sequence, we will have proved Lemma~\ref{lem:apd.stack.glamp} itself, since the limit is free of the choice of specific sub-sequences. 
  
  With a slight abuse of notation, we still use letters \( (p,n_1,n_2,\cdots ,n_E) \) to denote the sub-sequence along with the three limits involving \( \eta^0 \) in Equation~\eqref{eq:apd.pf.stack.init.limits} all hold almost surely. 
  
  We specify a GLAMP instance in the following way: Let \( \{ q_e:e\in[E] \}  \) all be sufficiently large and fixed integers. Let \( B= [\beta _1|\beta_2|\cdots |\beta _{E}|\eta^0|0|0|\cdots ]\in\R^{p\times q} \), \( W_e=[w_e|0|0|\cdots ]\in\R^{n_e\times q} \). At any time step \( t \in \mathbb{N} \), $\forall\,e\in[E]$,
  \begin{equation*}
  \begin{aligned}
    \ V_e^{t} =&~ \left[v_e^{1}|v_e^{2}|\cdots |v_e^{t}|0|0|\cdots \right], \\
    \eta _e^{t}(V_1^{t},\cdots ,V_{E}^{t})=&~ \Sigma _e^{1/2}\left[\beta _e|\eta^0|\eta(v_1^{1},\cdots ,v_{E}^{1})|\eta(v_1^{2},\cdots ,v_{E}^{2})|\cdots |\eta(v_1^{t},\cdots ,v_{E}^{t})|0|0|\cdots \right],\\
    R_e^{t} =&~ \left[\ix{e}\Sigma _e^{1/2}\beta _e|r_e^{0}|r_e^{1}|\cdots |r_e^{t}|0|0|\cdots \right],\\
    \Psi _e^{t}(R_e^{t}) =&~ \kappa _e\left[ y_e - r_e^{0}|y_e - r_e^{1}|\cdots |y_e - r_e^{t}|0|0|\cdots  \right]\quad (y_e  = \ix{e} \Sigma _e^{1/2}\beta _e + w_e).
  \end{aligned}
  \end{equation*}
  The GLAMP instance starts from \( V_e^{0}=0,\ \forall\,e\in[E] \). 
  
  Now we know the GLAMP iterates take the form 
  \begin{equation*}
    \begin{aligned}
      v_e^{t} =&~ \ix{e}^{\T} (y_e - r_e^{t-1}) - \Sigma _e^{1/2} \beta _e + \Sigma _e^{1/2} \eta(v_1^{t-1}, v_2^{t-1}, \cdots ,v_E^{t-1}), \\
      r_e^{t} =&~ \ix{e} \Sigma _e^{1/2} \eta(v_1^{t}, v_2^{t}, \cdots ,v_E^{t}) - \kappa _e  (y_e - r_e^{t-1}) \cdot \frac{1}{p}\E\left[ \trace \left( \Sigma _e^{1/2} \frac{\partial \eta}{\partial v_e^{t}}(Z_1^{t},Z_2^{t},\cdots, Z_E^{t}) \right) \right].
    \end{aligned}
    \end{equation*}
  where the Gaussian variables used are yet to be characterized by the state evolution. 
  
  To fully write out its form, we need the state evolution first. We discard the all-zero columns in \( V_e^{t} \) and \( R_e^{t} \),  start from \( \Sigma _{(R,e)}^{0} \) and compute the first two matrices for the state evolution. From there,  we prove the whole state evolution for \( V^{t}_e \) with induction.  
  \begin{enumerate}
    \item \( \Sigma _{(R,e)}^{0} = \begin{bmatrix}
      \lim_p \frac{1}{p}\beta _e^{\T}\Sigma _e \beta _e & \lim_p \frac{1}{p}\E[\beta _e^{\T} \Sigma _e \overline{\eta } ]  \\ 
      \lim_p \frac{1}{p}\E[\beta _e^{\T} \Sigma _e \overline{\eta }] &  \lim_p \frac{1}{p}\E[\overline{\eta }^{\T} \Sigma _e \overline{\eta }] \\ 
    \end{bmatrix} \).
    \item \( \Sigma _{(V,e)}^{1} = \lim_p \frac{1}{n_e}\E\left[\norm*{w_e+(Z_{(R,e)}^{0})_1 - (Z_{(R,e)}^{0})_2}_2^{2}\right]  \) where \( Z^{0}_{(R,e)}=[(Z_{(R,e)}^{0})_1|(Z_{(R,e)}^{0})_2]\in\R^{p\times 2}\) has iid rows drawn from \( \Nscr(0,\kappa _e \Sigma _{(R,e)}^{0}) \). Simplifying it with Assumption~\ref{asp:stack.for.glamp} and \ref{asp:stack.converge}, 
    
    \( \Sigma _{(V,e)}^{1} = \E[W_e^{2}] + \kappa _e \lim_p \frac{1}{p}\E[\norm{\overline{\eta } -\beta _e}_{\Sigma _{e}}^{2}] = (\tau _e^{*}) ^{2} \).
    \item  \( \Sigma _{(R,e)}^{1} = \begin{bmatrix}
      \lim_p \frac{1}{p}\beta _e^{\T}\Sigma _e \beta _e & \lim_p \frac{1}{p}\E[\beta _e^{\T} \Sigma _e \overline{\eta } ] & \lim_p \frac{1}{p}\E[\beta _e^{\T} \Sigma _e \overline{\eta } ]   \\ 
      \lim_p \frac{1}{p}\E[\beta _e^{\T} \Sigma _e \overline{\eta }] &  \lim_p \frac{1}{p}\E[\overline{\eta }^{\T} \Sigma _e \overline{\eta }] & \lim_p\frac{1}{p}\E\left[ \overline{\eta }_{(1)} ^{\T}\Sigma _e \overline{\eta }_{(2)}  \right] \\ 
      \lim_p \frac{1}{p}\E[\beta _e^{\T} \Sigma _e \overline{\eta } ] & \lim_p \frac{1}{p}\E\left[ \overline{\eta }_{(1)} ^{\T}\Sigma _e \overline{\eta }_{(2)}  \right] & \lim_p \frac{1}{p}\E[\overline{\eta }^{\T} \Sigma _e \overline{\eta } ],
    \end{bmatrix} \) where \( \overline{\eta } _{(1)} \) and \( \overline{\eta } _{(2)} \) are mutually independent copies of \( \overline{\eta }  \).
    \item \( \Sigma _{(V,e)}^{2} = \begin{bmatrix}
      (\tau _e^{*})^{2} & (\tau _e^{*})^{2}H_e(\vec{0}_E) \\ 
      (\tau _e^{*})^{2}H_e(\vec{0}_E) & (\tau _E^{*})^{2} \\ 
    \end{bmatrix} \). For the diagnoal elements of \( \Sigma _{(V,e)}^{2} \) we have simplified it in the same way as \( \Sigma _{(V,e)}^{1} \). For the off-diagonal term, it is 
    \begin{equation*}
  \begin{aligned}
    &\lim_p \frac{1}{n_e}\E \left[ (w_e+ (Z_{(R,e)}^{0})_1 - (Z_{(R,e)}^{0})_2)^{\T}(w_e+ (Z_{(R,e)}^{0})_1 - (Z_{(R,e)}^{0})_3)\right]\\
     &= \E[W_e^{2}] + \kappa _e \lim_p \frac{1}{p}\E \left[ (\overline{\eta }_{(1)}-\beta _e)^{\T}\Sigma _e(\overline{\eta }_{(2)} -\beta _e) \right] 
     =(\tau _e^{*})^{2}H_e(\vec{0}_E),
  \end{aligned}
    \end{equation*}
    where \( \overline{\eta } _{(1)} \) and \( \overline{\eta } _{(2)} \) are mutually independent copies of \( \overline{\eta }  \).
  \end{enumerate}
  We use induction to prove the whole state evolution of \( V_e^{t} \) takes the form of Definition~\ref{def:apd.stack.stat.evo}. Suppose the form in Definition~\ref{def:apd.stack.stat.evo} has been verified up to and including \( t \); let's prove the case of \( (t+1) \). 
  
  For the diagonal term,  \( (\Sigma _{(V,e)}^{t+1})_{t+1,t+1} = \lim_p \frac{1}{n_e}\E\left[\norm*{w_e+(Z_{(R,e)}^{0})_1 - (Z_{(R,e)}^{t})_{t+2}}_{2}^{2}\right] = (\tau _e^{*})^{2}  \) using the same simplification as that of \( \Sigma _{(V,e)}^{1} \).
    
  For the off-diagonal terms: 
  When \( i=1 \), \( j=t+1 \), we show \( \rho ^{1,t+1}=H(\vec{0}_E) \) by 
  \begin{equation*}
  \begin{aligned}
    (\tau _e^{*})^{2} \rho _e^{1,t+1} =&~\lim_p \frac{1}{n_e}\E \left[ (w_e+ (Z_{(R,e)}^{t})_1 - (Z_{(R,e)}^{t})_2)^{\T}(w_e+ (Z_{(R,e)}^{t})_1 - (Z_{(R,e)}^{t})_{t+2})\right]\\
    =&~ \E[W_e^{2}] + \kappa _e \lim_p \frac{1}{p} (\eta^0-\beta _e)^{\T} \Sigma _e \E \left[ \overline{\eta }-\beta _e \right] \\
    =&~ \E[W_e^{2}] + \kappa _e \lim_p \frac{1}{p}\E \left[ (\overline{\eta }_{(1)}-\beta _e)^{\T}\Sigma _e(\overline{\eta }_{(2)} -\beta _e) \right] \\
    =&~(\tau _e^{*})^{2}H_e(\vec{0}_E),
  \end{aligned}
  \end{equation*}
  where \( \overline{\eta } _{(1)} \) and \( \overline{\eta } _{(2)} \) are mutually independent copies of \( \overline{\eta }  \).
  
  When \( 1<i<j=t+1 \), we show \( \rho ^{i,t+1}=H(\rho ^{i-1,t}) \) by
  \begin{equation*}
    \begin{aligned}
      (\tau _e^{*})^{2} \rho _e^{i,t+1} =&~\lim_p \frac{1}{n_e}\E \left[ (w_e+ (Z_{(R,e)}^{t})_1 - (Z_{(R,e)}^{t})_{i+1})^{\T}(w_e+ (Z_{(R,e)}^{t})_1 - (Z_{(R,e)}^{t})_{t+2})\right]\\
      =&~ \E[W_e^{2}] + \kappa _e \lim_p \frac{1}{p}\E \left[ ({\eta }(Z^{i-1})-\beta _e)^{\T}\Sigma _e({\eta }(Z^{t}) -\beta _e) \right],
    \end{aligned}
    \end{equation*}
  where \( Z^{i-1}=[Z^{i-1}_1|Z^{i-1}_2|\cdots |Z^{i-1}_E] \) and \( Z^{t}=[Z^{t}_1|Z^{t}_2|\cdots |Z^{t}_E] \),
   \( \{ (Z_e^{i-1},Z_e^{t}) :e\in[E] \}  \) are \( E \) independent random tuples, and each row of \( [Z^{i-1}_e|Z^{t}_e] \) is sampled iid from 
    \( \Nscr(0,(\tau _e^{*})^{2}[1,\rho _e^{i-1,t} ; \rho _e^{i-1,t}, 1]) \). Thus \( \rho _e^{1,t+1}=H_e(\rho ^{i-1,t}) \). 
  
  Now we have verified the GLAMP iterates take the form of Definition~\ref{def:stack.glamp} and their state evolution is described by Definition~\ref{def:apd.stack.stat.evo}. By GLAMP convergence, we have shown Lemma~\ref{lem:apd.stack.glamp}. 
  \end{proof}
\subsubsection{Proof of Lemma~\ref{lem:apd.stack.cauchy}}
\label{sec:apd.pf.stacked.lasso.cauchy}

We first prove a lemma and then use it to prove Lemma~\ref{lem:apd.stack.cauchy}.

\begin{lemma}\label{lem:apd.stack.study.H}
  For the function \( H \) defined in Assumption~\ref{asp:stack.converge.cauchy} and any \( e\in[E] \), its \( e \)-th element \( H_e=H_e(\rho_1,\cdots ,\rho _E) \) is increasing and strictly convex in \( \rho _i \), \( \forall\, i\in[E] \), while we fix \( \{ \rho_j:j\neq i \}  \). We have 
  \begin{equation*}
    0\leq \frac{\partial H_e}{\partial \rho_i}\leq \frac{\kappa _e (\tau _i^{*})^{2}\varpi _i^{*}}{(\tau _e^{*})^{2}\varpi ^* _e}\overline{\delta }_i\overline{\delta }_e.
  \end{equation*}
\end{lemma}
\begin{proof}[\textbf{Proof of Lemma~\ref{lem:apd.stack.study.H}}]
  The first half of Proposition~\ref{lem:apd.stack.study.H} follows the same argument as that of Lemma~C.1 of \citet{bayati2011lasso}, or Lemma~C.2 of \citet{donoho2016high}. To elaborate, we first assume without loss of generality that \( i=1 \). We represent the function \( H_e \) in an alternative way. Let \( H_e = \lim_p H^{{(p)}}_e \), where \( H^{{(p)}}_e \) is the finite-\( p \) version of the same function, i.e. 
  \begin{equation*}
\begin{aligned}
  H_e^{{(p)}}(\rho_1,\rho_2,\cdots ,\rho _E)=&~\frac{1}{(\tau _e^{*})^{2}}\left\{\E[W_e^{2}] + \kappa_e \frac{1}{p}\E\left[ (\overline{\eta }_{{(1)}} - \beta_e)^{\T}\Sigma _e (\overline{\eta }_{{(2)}} - \beta _e) \right] \right\} \\
  \mathrm{where}\  \overline{\eta }_{{(1)}} = &~ \eta(\tau ^*_1 \sqrt{\rho_1}Z_1 + \tau _1^{*}\sqrt{1-\rho_1} Y_1,\tau ^*_2 \sqrt{\rho_2}Z_2 + \tau _2^{*}\sqrt{1-\rho_2} Y_2,\cdots,\tau ^*_E \sqrt{\rho_E}Z_E + \tau _E^{*}\sqrt{1-\rho_E} Y_E),\\
  \overline{\eta }_{{(2)}} = &~ \eta(\tau ^*_1 \sqrt{\rho_1}Z_1 + \tau _1^{*}\sqrt{1-\rho_1} Y_1',\tau ^*_2 \sqrt{\rho_2}Z_2 + \tau _2^{*}\sqrt{1-\rho_2} Y_2',\cdots,\tau ^*_E \sqrt{\rho_E}Z_E + \tau _E^{*}\sqrt{1-\rho_E} Y_E'),
\end{aligned}
  \end{equation*}
  and \( \{ Z_e \}  \), \( \{ Y_e \}  \), \( \{ Y_e' \}  \) are a total of \( (3E) \) iid \( \Nscr(0,I_p) \) variables. Then \( H_e \) can be further rewritten with a function \( \nu:\R^{p}\to\R^{p} \) as follows: 
  \begin{equation*}
\begin{aligned}
  H_e =&~ \frac{\E[W_e^{2}]}{(\tau _e^{*})^{2}} + \frac{\kappa _e}{(\tau _e^{*})^{2}}\cdot \frac{1}{p}\E \left[ \E \left[  \nu(\sqrt{\rho_1}Z_1+\sqrt{1-\rho_1} Y_1)^{\T} \nu(\sqrt{\rho_1}Z_1+\sqrt{1-\rho_1} Y_1') \middle| Z_2,Z_3,\cdots ,Z_E  \right] \right], \\
  \nu(\cdot ) =&~ \Sigma _e^{1/2}\left\{ \E\left[ \eta(\tau _1^{*}(\,\cdot\,) ,\tau ^*_2 \sqrt{\rho_2}Z_2 + \tau _2^{*}\sqrt{1-\rho_2} Y_2,\cdots,\tau ^*_E \sqrt{\rho_E}Z_E + \tau _E^{*}\sqrt{1-\rho_E} Y_E) \middle|Z_2,Z_3,\cdots ,Z_E \right] \right\}.
\end{aligned}
  \end{equation*} 
  We introduce a total of \( p \) iid stationary Ornstein-Uhlenbeck processes with covariance function \( \exp\left(-t\right) \), stacked as a vector \( X_t\in\R^{p} \), then \( H_e^{(p)} = \frac{\E[W_e^{2}]}{(\tau _e^{*})^{2}} + \frac{\kappa _e}{(\tau _e^{*})^{2}}\cdot \frac{1}{p}\E \left[ \E \left[  \nu(X_0)^{\T} \nu(X_t) \middle| Z_2,Z_3,\cdots ,Z_E  \right] \right] \), \( t = \log (1/\rho_1) \). With the same eigen-decomposition trick as that of  Lemma~C.1 of \citet{bayati2011lasso}, or Lemma~C.2 of \citet{donoho2016high}, we can show \( H_e^{{(p)}} \) is increasing and strictly convex in \( \rho_1\in [0,1] \), and thus so is \( H_e \). We have proved the first half of Proposition~\ref{lem:apd.stack.study.H}. 

  To show \( 0\leq \frac{\partial H_e}{\partial \rho _i} \leq \frac{\kappa _e (\tau _i^{*})^{2}\varpi _i^{*}}{(\tau _e^{*})^{2}\varpi ^* _e}\overline{\delta }_i\overline{\delta }_e  \), we only need to show \( 0\leq \frac{\partial H_e}{\partial \rho _i}(\1_E) \leq \frac{\kappa _e (\tau _i^{*})^{2}\varpi _i^{*}}{(\tau _e^{*})^{2}\varpi ^* _e}\overline{\delta }_i\overline{\delta }_e \), where \( \1_E \) is a vector of all ones in \( R^{E} \). Without loss of generality, still assume \( i=1 \). 

  We first compute the partial derivatives at a general \( \rho = (\rho_1,\cdots ,\rho _E) \), and then set \( \rho=\1_E \).  Suppose the limit can be exchanged with the partial derivative, we have
  \begin{equation*}
\begin{aligned}
  \frac{\partial H_e}{\partial \rho_1}=&~ \frac{\kappa _e \tau_1^{*}}{(\tau _e^{*})^{2}} \lim_p \frac{1}{p} \E \left[ (\overline{\eta }  _{(2)}-\beta _e)^{\T} \Sigma _e J_1(\overline{\eta } _{(1)}) \left( \frac{1}{\sqrt{\rho_1}}Z_1 - \frac{1}{\sqrt{1-\rho_1}} Y_1 \right) \right],
\end{aligned}
  \end{equation*}
  where we have used the exchangeability of \( Y_1 \) and \( Y_1' \) to simplify the form. The notation `$J_1(\overline{\eta } _{(1)})$' means the \( p \)-by-\( p \) Jacobian matrix of \( \overline{\eta } _{(1)} \) w.r.t. \( (\tau_1^{*}\sqrt{\rho_1}Z_1 + \tau_1^{*}\sqrt{1-\rho_1} Y_1) \). In fact,  $J_1(\overline{\eta } _{(1)})$ has been studied at Equation~\eqref{eq:apd.prelim.stack.prop.partial.deriv}, combined with which we can invoke Stein's Lemma to get 
  \begin{equation*}
\begin{aligned}
  \frac{\partial H_e}{\partial \rho_1}=&~ \frac{\kappa _e \tau _1^{*}}{(\tau _e^{*})^{2}} \lim_p \frac{1}{p} \E \left[ (\overline{\eta }  _{(2)}-\beta _e)^{\T} \Sigma _e J_1(\overline{\eta } _{(1)}) \left( \frac{1}{\sqrt{\rho_1}}Z_1 \right) \right] 
  = \frac{\kappa _e (\tau _1^{*})^{2}}{(\tau _e^{*})^{2}} \lim_p \frac{1}{p} \E \left[ \trace \left( J_1(\overline{\eta } _{(2)})^{\T} \Sigma _e J_1(\overline{\eta } _{(1)}) \right) \right].
\end{aligned}
  \end{equation*}
  When \( \rho=\1_{E} \), \( \overline{\eta } _{(1)} = \overline{\eta } _{(2)} = \overline{\eta }  \) and their support also coincides, which we denote as \( S \). From Equation~\eqref{eq:apd.prelim.stack.prop.partial.deriv}, we know 
  \begin{equation*}
\begin{aligned}
  \frac{\partial H_e}{\partial \rho_1}(\1_{E}) =&~ \frac{\kappa _e (\tau _1^{*})^{2}}{(\tau _e^{*})^{2}} \lim_p \frac{1}{p} \E \left[ \trace \left( J_1(\overline{\eta } _{S})^{\T} (\Sigma _e)_{S,S} J_1(\overline{\eta } _{S}) \right) \right] \\
  =&~ \frac{\kappa _e (\tau _1^{*})^{2} (\varpi _1^{*})^{2}}{(\tau _e^{*})^{2} \varpi _e^{*}} \lim_p \frac{1}{p} \E \left[ \trace \left( (\overline{\Sigma }_{S,S} )^{-1} (\varpi _e\Sigma _{e})_{S,S} (\overline{\Sigma }_{S,S} )^{-1} (\varpi _1\Sigma _1)_{S,S} \right) \right]\quad (\overline{\Sigma } =\sum_{e}\varpi _e^{*}\Sigma _e) \\
  =&~ \frac{\kappa _e (\tau _1^{*})^{2} \varpi _1^{*}}{(\tau _e^{*})^{2} \varpi _e^{*}} \lim_p \frac{1}{p} \E \left[ \trace \left( (\overline{\Sigma }_{S,S} )^{-1/2} (\varpi _e^{*}\Sigma _{e})_{S,S} (\overline{\Sigma }_{S,S} )^{-1/2}(\overline{\Sigma }_{S,S} )^{-1/2} (\varpi _1\Sigma _1)_{S,S}(\overline{\Sigma }_{S,S} )^{-1/2} \right)  \right] \\
  \leq &~\frac{\kappa _e (\tau _1^{*})^{2} \varpi _1^{*}}{(\tau _e^{*})^{2} \varpi _e^{*}} \lim_p \frac{1}{p}\E \left[ \trace \left( (\overline{\Sigma }_{S,S} )^{-1} (\varpi _e^{*}\Sigma _{e})_{S,S}  \right)  \right] \cdot \lim_p \frac{1}{p}\E \left[ \trace \left( (\overline{\Sigma }_{S,S} )^{-1} (\varpi _1^{*}\Sigma _{1})_{S,S}  \right)  \right] \\
  =&~ \frac{\kappa _e (\tau _1^{*})^{2} \varpi _1^{*}}{(\tau _e^{*})^{2} \varpi _e^{*}} \overline{\delta } _e \overline{\delta } _1.
\end{aligned}
  \end{equation*}
  In the inequality above, we have used a trace inequality that \( \trace(AB)\leq \lambda _{\max }(A) \trace(B)\leq \trace(A)\trace(B) \) for \( A,B \) both positive semi-definite \citep[see e.g. Theorem 1 of][]{fang1994inequalities}. So far we have proved Proposition~\ref{lem:apd.stack.study.H}.
\end{proof}

\begin{proof}[\textbf{Proof of Lemma~\ref{lem:apd.stack.cauchy}}]

  We start by studying \( \plim_{p\to \infty}\frac{1}{p}\norm*{v_e^{t+1}-v_e^{t}}_2^{2} \).
  Using Lemma~\ref{lem:apd.stack.glamp}, we know \( \plim_{p\to \infty}\frac{1}{p}\norm*{v_e^{t+1}-v_e^{t}}_2^{2} = 2 (\tau _e^{*})^{2} (1 - \rho _e^{t,t+1}) \). We plan to show \( \rho _e^{t,t+1} \) converges exponentially fast to 1 for each \( e\in[E] \) with an argument of fixed point iteration. 
  
  By the multi-variate mean-value theorem, for any \( \rho _{(1)}\in[0,1]^{E} \) and \( \rho _{(2)} = H(\rho _{(1)}) \), there exists some point \( \rho _{\mathrm{tmp}} \), which is a convex combination of \( \rho _{(1)},\rho _{(2)} \), s.t. \( \1_E-\rho _{(2)} = H(\1_{E}) - H(\rho _{(1)}) = J(H(\rho _{\mathrm{tmp}}))(\1_{E} - \rho_{(1)}) \). \( J(H(\rho _{\mathrm{tmp}}))\in\R^{E\times E} \) is the Jacobian matrix. By Lemma~\ref{lem:apd.stack.study.H}, we know 
\begin{gather*}
0\leq (\tau _e^{*})^{2}\varpi _e^{*}(1-\rho _{(2),e}) \leq \kappa _e \overline{\delta } _e \sum_{i\in[E]} \overline{\delta }_i (\tau _i^{*})^{2} \varpi _i^{*} (1-\rho _{(1),i}) \leq \kappa _e \overline{\delta } _e \max_{e\in[E]}\{ \overline{\delta } _e \} \sum_{i\in[E]} (\tau _i^{*})^{2} \varpi _i^{*} (1-\rho _{(1),i}) \\
\implies 0 \leq \sum_{e\in[E]} (\tau _e^{*})^{2}\varpi _e^{*}(1-\rho _{(2),e}) \leq \max_{e\in[E]}\left\{ \kappa _e \overline{\delta } _e \right\} \max _{e\in[E]} \left\{ \overline{\delta } _e \right\} \sum_{e\in[E]} (\tau _e^{*})^{2}\varpi _e^{*}(1-\rho _{(1),e}).
\end{gather*}
By Assumption~\ref{asp:stack.converge.fixed.point}, we know \( \kappa _e \overline{\delta } _e <1 \); by Proposition~\ref{prop:stack.onsager}, we know \( \sum_{e\in[E]} \overline{\delta } _e = \frac{1}{p}\E[\norm{\overline{\eta } }_0] \leq 1 \). As a result, \( \max_{e\in[E]}\left\{ \kappa _e \overline{\delta } _e \right\} \max _{e\in[E]} \left\{ \overline{\delta } _e \right\} <1 \). By a fixed-point iteration argument, we know \( \rho _{e}^{t,t+1} \) converges exponentially fast to 1. As a result, \( \forall\,\epsilon>0 \), there exists \( T \) s.t. \( \forall\,j>i\geq T \), \( \max _{e\in[E]} \left\{ 1-\rho _e^{i,j} \right\} < \epsilon \). 

The rest of Lemma~\ref{lem:apd.stack.cauchy} is easily implied by what we have proved. Firstly, using Lemma~\ref{lem:apd.stack.glamp}, we know \( \plim_{p\to \infty}\frac{1}{p}\norm*{v_e^{j}-v_e^{i}}_2^{2} = 2 (\tau _e^{*})^{2} (1 - \rho _e^{i,j}) \). So the Cauchy property also holds for \( \{ v_e^{t} \}  \). Second, \( \eta^t = \eta(v_1^{t},\cdots ,v_E^{t}) \) is a Lipschitz function in \( (v_1^{t},v_e^{t},\cdots ,v_{E}^{t}) \), so the Cauchy property of \( \{\eta ^t \}\) follows. Third, by the corresponding version of Lemma~\ref{lem:apd.stack.glamp} that tracks the convergence of $r^t$ (which we omit for simplicity), we have \( \plim\frac{1}{n_e}\norm{r_e^{t}-r_e^{t-1}}_2^{2} = 2\kappa _e (\lim_p \frac{1}{p} \E[\overline{\eta }^{\T}\Sigma _e \overline{\eta } ] - \lim_p \frac{1}{p} \E[\eta(Z^{t})^{\T}\Sigma _e \eta(Z^{t-1})]) = 2 (\tau _e^{*})^{2}(1 - H_e(\rho ^{t-1,t})) = 2 (\tau _e^{*})^{2}(1 - \rho_e^{t,t+1}) \). Hence the Cauchy property of \( \{ r^{t}_e \}  \) also follows. 

\end{proof}

\subsubsection{Proof of Lemma~\ref{lem:stacked.lasso.subgradient}}\label{sec:apd.pf.stacked.lasso.subgradient}

  By the definition of \( \eta \) in Equation~\eqref{eq:def.stack.et}, we know there exists a sub-gradient of \( \norm{\eta^t}_1 \), denoted as \( \partial \norm{\eta^t}_1  \), that satisfies \( \eta^t = \eta(v_1^{t},\cdots ,v_E^{t}) \), \( \sum_{e} \varpi _e^{*} \Sigma _e^{1/2}[\Sigma _e^{1/2} \eta^t - (v_e^{t}+ \Sigma _e^{1/2}\beta _e) ] + \theta ^* \diag(\vec{\lambda } ) \partial \norm{\eta^t}_1 = 0 \). Solving for this \( \partial \norm{\eta^t}_1 \) is gives us the expression of \( s^{t} \) in Equation. 

  For brevity of the proof, we new introduce a new symbol \( \specop \) in the spirit of Landau notation. As two general properties of \( \specop \): \emph{(i)} Some sequence \( X^{(p),t} = \specop  \) iff  \( \forall\, \epsilon>0 \), \( \exists\, T \) s.t. \( \forall\,t\geq T \), \( \lim_{p}\PR\left[ \frac{1}{\sqrt{p}}\norm{X^{{(p),t}}}_2 < \epsilon \right] = 1 \). Apparently \( \specop + \specop = \specop \) and \( \mathit{Constant}\cdot \specop = \specop \). \emph{(ii)} If some \( a^{{(p)}}=o(1) \) as \( p\to \infty \), and \( Y^{(p),t} \) satisfies \( \lim_p \PR \left[ \frac{1}{\sqrt{p}}\norm{Y^{{(p),t}}}_2 \leq M \right] = 1\ \forall\,t \in \mathbb{N_+} \)  for some \( M \) free of \( t \), then \( a^{(p)} Y^{{(p),t}} = \specop \). 
  
  We know from Lemma~\ref{lem:apd.stack.cauchy} that \( r_e^{t}-r_e^{t-1}=\specop \), \( \hat{\beta } ^{t+1}-\eta^t=\specop \). By Definition~\ref{def:stack.glamp}, \( (1-\kappa _e \overline{\delta }^{(p)} _e)(y_e - r_e^{t-1}) = y_e - \ix{e} \Sigma _e^{1/2} \eta^t + \specop \). Lemma~\ref{lem:apd.sing.val.limit.concent} implies  \( \lim_{p} \sigma _{\max }(\ix{e}) = 1 + \sqrt{\kappa _e} \), a.s.  Since \( \norm{y_e - r_e^{t-1}}_2 \leq c_1 \sigma _{\max }(\ix{e}) \norm{\beta _e}_2 + \norm{w_e}_2 + \norm{r_e^{t-1}}_2  \), we know there exists \( M>0 \) s.t. \( \lim_{p}\PR \left[ \frac{1}{\sqrt{p}}\norm{y_e - r_e^{t-1}} \leq M \right]=1 \) for all \( t \in \mathbb{N_+} \). Combined with 
\( \overline{\delta } _e^{(p)} = \overline{\delta }  _e + o(1) \), we know \( (1-\kappa _e \overline{\delta } _e)(y_e - r_e^{t-1}) = y_e - \ix{e} \Sigma _e^{1/2} \eta^t + \specop \). As a result,
\begin{equation*}
\begin{aligned}
  \nabla L(\eta^t) =&~ - \sum_{e\in[E]} \pi _e \Sigma_e^{1/2} \ix{e}^{\T}(y_e - \ix{e} \Sigma _e^{1/2}\eta^t) + \diag(\vec{\lambda } ) s^{t} \\
  =&~ - \sum_{e\in[E]} \pi _e(1-\kappa _e \overline{\delta } _e) \Sigma_e^{1/2} \ix{e}^{\T}(y_e - r_e^{t-1}) + \diag(\vec{\lambda } ) s^{t} + \specop.
\end{aligned}
\end{equation*}
Recall the system of equations~\eqref{eq:stack.fixed.point} from Assumption~\ref{asp:stack.converge.fixed.point}, and the first line of GLAMP iterates in Equation~\eqref{eq:def.stack.glamp}, and then we have 
\begin{equation*}
\begin{aligned}
    \theta ^* \nabla L(\eta^t) =&~ - \sum_{e\in[E]} \varpi _e^{*} \Sigma_e^{1/2} (v_e^{t} + \Sigma _e^{1/2}\beta _e - \Sigma _e^{1/2}\eta^{t-1}) + \theta ^* \diag(\vec{\lambda } ) s^{t} + \specop = \specop, \\
\end{aligned}
\end{equation*}
where the second equality has used the definition of \( s^{t} \) in Equation~\eqref{eq:def.stack.subgradient}.  This completes the proof. 

\subsubsection{Proof of Lemma~\ref{lem:stacked.lasso.bound.on.sc2}}\label{sec:apd.pf.stacked.lasso.bound.on.sc2}
We need to prepare some definitions, and Lemma~\ref{lem:stacked.lasso.indicator.mean}, \ref{lem:stacked.lasso.indicator.converge} to prove Lemma~\ref{lem:stacked.lasso.bound.on.sc2}. 

For any constant \( c\in(0,1) \), define \( S^{t}(c) = \{ j\in[p] : \abs{s^{t}_j} \geq 1 - c  \}  \) where \( s^{t} \) has been defined in Equation~\eqref{eq:def.stack.subgradient}. We also define the `idealized' versions of them indicated by overlines: 
\begin{equation}\label{eq:def.stack.ideal.subgradient}
  \overline{s}\in\R^{p},\ \overline{s} = \frac{1}{\theta ^*}\left[ \diag(\vec{\lambda } )^{-1} \right]\sum_{e\in[E]} \varpi _e^{*}( \tau _e^{*}\Sigma _e^{1/2} Z_e  + \Sigma _e(\beta _e - \eta(\tau _1^{*}Z_1,\cdots ,\tau _E^{*}Z_E)) ),
\end{equation}
where \( Z_1,Z_2,\cdots ,Z_E \iidsim \Nscr(0,I_p) \). 
We still use \( \overline{\eta } = \eta(\tau _1^{*}Z_1,\cdots ,\tau _E^{*}Z_E) \) as a shorthand, but note that we will be using the same \( Z_1,Z_2,\cdots ,Z_E \) throughout the proof. We define \( \overline{S}(c) = \{ j\in[p] : \abs{\overline{s} _j} \geq 1 - c \}   \). We also define a set of re-scaled quantities indicated with tildes: Let \( \lambda _{\min ,\mathrm{all}} = \min _{e\in[E]} \{\lambda _{\min }(\Sigma _e)\} \). We define \( \widetilde\Sigma _e = \Sigma _e / \lambda _{\min ,\mathrm{all}}  \), \( \widetilde\tau ^*_e = \tau _e^{*} / \sqrt{\lambda _{\min, \mathrm{all}}} \), \( \widetilde\theta ^* = \theta ^* / {\lambda _{\min,\mathrm{all} }} \). We also define one covariance without subscript, \( \widetilde{\Sigma } = \sum_{e\in[E]} \varpi _e^{*}\widetilde{\Sigma } _e \). 

Define a new quantity \( \tilde{\eta}^t = \eta^t+ \widetilde{\theta } ^{*}\diag(\vec{\lambda})s^{t} \) and its idealized counterpart 
 \( \overline{\widetilde{\eta}} \in\R^{p} \), \( \overline{\widetilde{\eta}} = \overline{\eta } + \widetilde{\theta }  ^* \diag(\vec{\lambda } ) \overline{s}  \). 
About \( \tilde{\eta}^t \) and \( \overline{\widetilde{\eta}}  \), we notice \( S^{t}(c) = \{ j\in[p] : \abs{\tilde{\eta}^t_j} \geq \widetilde{\theta }^* \vec{\lambda } _j  (1-c) \}  \), \( \overline{S} (c) = \{ j\in[p] : \abs{\overline{\widetilde{\eta}} ^{t}_j} \geq \widetilde{\theta }^*  \vec{\lambda }  _j  (1-c) \}  \). This is because \( s^{t} \) and \( \overline{s}  \) are respectively the sub-gradients of \( \norm{\eta^t}_1 \)  and \( \norm{\overline{\eta } }_1 \). 

Define \( \widetilde{Z}  = \sum_{e\in[E]} \varpi _e^{*} \widetilde{\tau }^{*} _e (\widetilde{\Sigma } )^{-1/2} \widetilde{\Sigma } _e^{1/2} Z_e \in \R^{p} \), then \( \widetilde{Z} \sim \Nscr(0, \sum_{e\in[E]}(\varpi _e^{*}\widetilde{\tau } _e^{*})^{2}\widetilde{\Sigma } ^{-1/2}\widetilde{\Sigma } _e \widetilde{\Sigma } ^{-1/2}) \). For the covariance matrix of \( \widetilde{Z} \), we have
\begin{equation}
\begin{aligned}\label{eq:apd.stack.Ztilde.eigen.val}
  \liminf_{p} \lambda _{\min }\bigg(\sum_{e\in[E]}(\varpi _e^{*}\widetilde{\tau } _e^{*})^{2}\widetilde{\Sigma } ^{-1/2}\widetilde{\Sigma } _e \widetilde{\Sigma } ^{-1/2}\bigg) \geq &~ \frac{1}{c_1^{2}}{\sum_{e\in[E]}(\varpi _e^{*}\widetilde{\tau } _e^{*})^{2}} > 0. \\
\end{aligned}
\end{equation}
We note that \( \overline{\eta }  \) depend on \( (Z_1,Z_2,\cdots ,Z_E) \) only through \( \widetilde{Z} \). Moreover, \( \widetilde{\Sigma }^{1/2} \overline{\eta }  \) is 1-Lipschitz w.r.t. \( \widetilde{Z} \). To see why these hold true, we rewrite \( \overline{\eta }  \) as below and invoke the Lipschitz property of proximal operators.
\begin{equation*}
\begin{aligned}
  \overline{\eta } =&~ \argmin_{b} \bigg\{ \sum_{e\in[E]} \frac{1}{2}\varpi _e^{*} \norm*{\widetilde{\Sigma } _e^{1/2} b - \widetilde{\tau } _e^{*}Z_e - \widetilde{\Sigma } _e^{1/2} \beta _e }_2^{2} + \widetilde{\theta } ^{*}\norm*{\diag(\vec{\lambda } )b}_1 \bigg\} \\
  =&~ \argmin_{\beta} \bigg\{ \frac{1}{2} \norm{\widetilde{\Sigma }^{1/2} b}_2^{2} - \frac{1}{2}\cdot 2 ( \widetilde{\Sigma }^{1/2} b )^{\T} \widetilde{Z} - \frac{1}{2}\cdot 2 \sum_{e\in[E]} b ^\T \widetilde{\Sigma }  \beta _e + \widetilde{\theta } ^{*}\norm*{\diag(\vec{\lambda } )b}_1 \bigg\}.
\end{aligned}
\end{equation*}
We also express \( \overline{\widetilde{\eta}}  \) in terms of \( \widetilde{Z} \). we get 
\begin{equation}
\begin{aligned}\label{eq:def.y}
  \overline{\widetilde{\eta}}  = &~ \sum_{e\in[E]} \varpi _e^{*} \widetilde{\Sigma } _{e} \beta _e + \widetilde{\Sigma } ^{1/2} \left[ \widetilde{Z} - (I_p - \widetilde{\Sigma } ^{-1})\widetilde{\Sigma } ^{1/2} \overline{\eta }  \right].
\end{aligned}
\end{equation} 

For each \( j\in[p] \), define \( \widetilde{\sigma } _{j}\in\R^{p} \) as the \( j \)-th row vector of \( \widetilde{\Sigma }^{1/2}  \) (written vertically despite being a row vector).  
We note that \( \norm{\widetilde{\sigma } _{j}}_2 \geq 1 \). This is because \( 1 \leq \lambda _{\min }(\widetilde{\Sigma } ^{1/2}) = \inf_{\norm{u}_2=1 } \norm{\widetilde{\Sigma } ^{1/2} u}_2 \leq \norm{P_{-j}^{\bot} \widetilde{\sigma }  _{j}}_2 \leq \norm{\widetilde{\sigma } _{j}}_2 \) where \( P^{\bot}_{-j}\in\R^{p\times p} \) is the projection matrix onto the orthogonal complement of \( \Span(\widetilde{\sigma }  _{i}:i\neq j) \), and \( u \) could be taken as \( u = P^{\bot}_{-j} \widetilde{\sigma }  _{j} / \norm{P^{\bot}_{-j} \widetilde{\sigma }  _{j}}_2 \). 

Let \( \widetilde{Z}_{j} = \widetilde{\sigma } _{j}^{\T} \widetilde{Z} \in \R \). We know from Equation~\eqref{eq:apd.stack.Ztilde.eigen.val} and \( \norm{\widetilde{\sigma } _j}_2\geq 1 \) that \( \liminf_{p}\Var(\widetilde{Z} _j) \geq \sum_{e}(\varpi _e^{*}\widetilde{\tau } _e^{*})^{2} / c_1^{2} > 0 \).

We try to write \( \overline{\widetilde{\eta}}_j  \) as a function of \( \widetilde{Z}_{j} \) and study that function. Define \( P_j^{\bot} \) as the projection matrix onto the orthogonal complement of \( \Span(\widetilde{\sigma } _{j}) \). Notice that \( \widetilde{Z}  = \frac{\widetilde{\sigma }  _{j}^{\T} \widetilde{Z} }{\norm{\widetilde{\sigma }  _{j}}_2^{2}} \widetilde{\sigma } _{j} + P^{\bot}_{j}\widetilde{Z}  =   \frac{\widetilde{Z}_{j}}{\norm{\widetilde{\sigma }  _{j}}_2^{2}} \widetilde{\sigma } _{j} + P^{\bot}_{j}\widetilde{Z}  \), where \( \widetilde{Z}_{j}\independent P^{\bot}_{j}\widetilde{Z}  \). 
We define \( f_{j}(\widetilde{Z}_{j}) = (I- \widetilde{\Sigma } ^{-1}) \widetilde{\Sigma } ^{1/2} \overline{\eta }  \). Then \( f_{j}(\widetilde{Z}_{j}) \) is a function of \( \widetilde{Z}_{j} \) but also determined by \( P^{\bot}_{j}\widetilde{Z} \) as well as \( \{ Z_{e'} :e'\neq e \}  \), but these Gaussian variables are all independent from \( \widetilde{Z}_{j} \). Conditioning on them, \( \widetilde{\Sigma } ^{1/2}\overline{\eta }  \) is \( \norm{\widetilde{\sigma } _{j}}_2^{-1} \)-Lipschitz in \( \widetilde{Z} _j \). 
Notice by the construction of \( \widetilde{\Sigma} \),
\begin{equation*}
\begin{aligned}
  \lambda _{\max}(\widetilde{\Sigma }^{-1} ) =&~ \lambda _{\min }(\widetilde{\Sigma } )^{-1} \leq \min _{e}\{ \lambda _{\min }(\Sigma _e) \} / [\sum_{e}\varpi _e^{*}\lambda _{\min }(\Sigma _e)] \leq 1,\\
  \lambda _{\min}(\widetilde{\Sigma }^{-1} ) =&~ \lambda _{\max }(\widetilde{\Sigma } )^{-1} \geq \min _{e}\{ \lambda _{\min }(\Sigma _e) \} / \max _{e} \left\{ \lambda _{\max }(\Sigma _{e}) \right\}.
\end{aligned}
\end{equation*}
We define the ratio \( \kappa _{cond}^{-1} = \min _{e}\{ \lambda _{\min }(\Sigma _e) \} / \max _{e} \{ \lambda _{\max }(\Sigma _e) \} \). Then  \( f_j(\cdot ) \)  is \( (1-\kappa _{cond}^{-1})\norm{\widetilde{\sigma } _{j}}_2^{-1} \)-Lipschitz. 

We then look at how \( \overline{\widetilde{\eta}} _j \) is related to \( \widetilde{Z} _j \), and write \( \overline{\widetilde{\eta}} _j = h_{j}(\widetilde{Z} _{j}) \). By the form of \( \overline{\widetilde{\eta}}  \) in Equation~\eqref{eq:def.y}, we know \( h_{j} \) is also determined jointly by \( \widetilde{Z} _j \) and \( P^{\bot}_{j}\widetilde{Z}  \), \( \{ Z_{e'}:e'\neq e \}  \). Given any \( P^{\bot}_{j}\widetilde{Z} \) and \( \{ Z_{e'}:e'\neq e \}  \), the function \( h_j(\cdot ) \) is Lipschitz continuous, and satisfies  
\begin{equation}
\begin{aligned}\label{eq:apd.stack.property.of.hj}
  \abs{ h_j(x_1)-h_j(x_2) } \geq &~  \abs{x_1-x_2} -  \abs{\widetilde{\sigma } _j^{\T}\left(f_{j}(x_1)-f_{j}(x_2)\right)} 
  \geq \kappa ^{-1}_{cond} \abs{x_1-x_2}.
\end{aligned}
\end{equation}
With Equation~\eqref{eq:apd.stack.property.of.hj} prepared, we need a few more lemmas:
\begin{lemma}\label{lem:stacked.lasso.indicator.mean}
  Let \( \{ I_a\subset \R:a>0 \}  \) be a collection of intervals indexed by \( a>0 \), satisfying \( \lim_{a\downarrow 0}m(I_a)=0 \) where \( m(\cdot ) \) is the Lebesgue measure. Then \( \lim_{a\downarrow 0} \sup_{p}\E\left[\frac{1}{p}\sum_{j\in[p]}\I{\overline{\widetilde{\eta}} _j\in I_a}\right] = 0 \). 
\end{lemma}
\begin{proof}[\textbf{Proof of Lemma~\ref{lem:stacked.lasso.indicator.mean}}] 
  \( \PR [\overline{\widetilde{\eta}} _j \in I_a ] = \PR [h(\widetilde{Z} _j)\in I_a ] = \E\left[ \PR[h(\widetilde{Z} _j)\in I_a|P_j^{\bot}\widetilde{Z},\{Z_{e'}:e'\neq e\} ] \right]   \). For any fixed \( P_j^{\bot}\widetilde{Z} \) and \( \{Z_{e'}:e'\neq e\} \), Equation~\eqref{eq:apd.stack.property.of.hj} implies there exists an interval \( {I}  \), s.t. \( m({I}) \leq \kappa _{cond}\, m(I_a) \) and \( \{h(\widetilde{Z} _j)\in I_a\}|_{P_j^{\bot}\widetilde{Z},\{Z_{e'}:e'\neq e\}} = \{\widetilde{Z} _j\in I\} \). This is because otherwise there must be \( x_1\neq x_2 \) s.t. \( h_j(x_1)=h_j(x_2) \) but \( \abs{x_1-x_2}>0  \), contradicting Equation~\eqref{eq:apd.stack.property.of.hj}. As a result, 
  \begin{equation*}
\begin{aligned}
  \E\left[ \PR[h(\widetilde{Z} _j)\in I_a|P_j^{\bot}\widetilde{Z},\{Z_{e'}:e'\neq e\} ] \right] \leq \sup_{I:m(I)\leq \kappa _{cond}m(I_a)} \PR \left[ \widetilde{Z} _j \in I \right] \leq \frac{\kappa _{cond}\,m(I_a)}{\sqrt{2\pi ( \sum_{e}(\varpi _e^{*}\widetilde{\tau } _e^{*})^{2} / c_1^{2})}},
\end{aligned}
  \end{equation*}
where the second inequality has simply used the upper bound on the density of a normal distribution.
Lemma~\ref{lem:stacked.lasso.indicator.mean} follows right from the above line. 
\end{proof}

\begin{lemma}\label{lem:stacked.lasso.indicator.converge}
  For any fixed \( t \in \mathbb{N_+} \) and \( c\in(0,1) \), \( \abs{S^{t}(c)} / p = \E[\abs{\overline{S} (c)} / p] + \smalop   \). 
\end{lemma}
\begin{proof}[\textbf{Proof of Lemma~\ref{lem:stacked.lasso.indicator.converge}}]
We fix \( c \) for the rest of the proof. We define the notation \( (x)_+ = \max \{ x,0 \},\forall\, x\in\R \). For $0<a<\widetilde{\theta } ^{*}(1-c)/3$, define a Lipschitz function \( \phi ^{(1)}_1(x) : \R\to\R \): \( \phi ^{(1)}_1(x) = [1 - (\widetilde{\theta } ^{*}(1-c) - \abs{x} )_+ / a]_+ \). Then \( \abs{\phi ^{(1)}_1(x) - \I{\abs{x} \geq \widetilde{\theta } ^{*}(1-c)}}  \leq \I{\abs{x} \in (\widetilde{\theta } ^{*}(1-c)-a,\widetilde{\theta } ^{*}(1-c)) } \). 
  
For any $x\in\R$, define \( \mathrm{dist}(x)  \) as the distance of \( \abs{x} \) to the interval \( (\widetilde{\theta } ^{*}(1-c)-a,\widetilde{\theta } ^{*}(1-c)) \), i.e. \( \mathrm{dist}(x) = (\abs{x}  - \widetilde{\theta } ^{*}(1-c))_+ + (\widetilde{\theta } ^{*}(1-c)-a-\abs{x}  )_+ \). Define another Lipschitz function \( \phi ^{(1)}_2(x) : \R\to\R \): \( \phi ^{(1)}_2(x) = (1 - \mathrm{dist}(x)/a)_+ \). Then \(  \I{\abs{x} \in (\widetilde{\theta } ^{*}(1-c)-a,\widetilde{\theta } ^{*}(1-c)) } \leq \phi ^{(1)}_2(x) \leq  \I{\abs{x} \in (\widetilde{\theta } ^{*}(1-c)-2a,\widetilde{\theta } ^{*}(1-c)+a) } \). For $i=1,2$, define \( \phi _i^{{(p)}}:\R^{p}\to \R \): \( \forall\,x\in\R^{p} \), \( \phi _i^{{(p)}}(x) = \sum_{j\in[p]} \phi ^{(1)}_i (x_j) \). By Lemma~\ref{lem:apd.pseudo}, both \( \{ \phi _1^{{(p)}}:p\geq 1 \}  \) and \( \{ \phi _{2}^{{(p)}}:p\geq 1 \}  \) are uniformly Lipschitz functions. 
  
It is already known that  \( S^{t}(c) = \{ j\in[p] : \abs{\tilde{\eta}^t_j} \geq \widetilde{\theta }^* \vec{\lambda } _j(1-c) \}  \), \( \widetilde{S} (c) = \{ j\in[p] : \abs{\overline{\widetilde{\eta}} ^{t}_j} \geq \widetilde{\theta }^* \vec{\lambda } _j (1-c) \}  \). 
By Lemma~\ref{lem:apd.stack.glamp}, for any already chosen and fixed \( a \), as \( p\to \infty \), 
\begin{equation*}
\begin{aligned}
  \abs*{{\abs{S^{t}(c)}}/{p} -  \E[{\abs{\overline{S} (c)}}/{p}]} \leq&~ \abs{\phi _1^{(p)}(\tilde{\eta}^t) - \E[\phi _1^{(p)}(\overline{\widetilde{\eta}} )]} + \phi _2^{(p)}(\tilde{\eta}^t) + \E[\phi _2^{(p)}(\overline{\widetilde{\eta}} )] \leq 2 \E[\phi _2^{(p)}(\overline{\widetilde{\eta}} )] + \smalop \\
  \leq &~ \frac{2}{p}\sum_{j\in[p]}\E[ \I{\abs{\overline{\widetilde{\eta}} _j} \in (\widetilde{\theta } ^{*}(1-c)-2a,\widetilde{\theta } ^{*}(1-c)+a) }] + \smalop. 
\end{aligned}
  \end{equation*}
  By Lemma~\ref{lem:stacked.lasso.indicator.mean}, for any \( \epsilon>0 \), we can fix a small \( a \) s.t. \( \sup_{p}\frac{1}{p}\sum_{j\in[p]}\E[ \I{\abs{\overline{\widetilde{\eta}} _j} \in (\widetilde{\theta } ^{*}(1-c)-2a,\widetilde{\theta } ^{*}(1-c)+a) }] < \epsilon/4 \) before taking \( p\to \infty \).  Thus Lemma~\ref{lem:stacked.lasso.indicator.converge} is proved. 
\end{proof}  
  
  \begin{proof}[\textbf{Proof of Lemma~\ref{lem:stacked.lasso.bound.on.sc2}}]
For any chosen and fixed $c\in(0,1)$, we know from Lemma~\ref{lem:stacked.lasso.indicator.converge}, \( \abs{S^{t}(c)} / p = \E[\abs{\overline{S} (c)} / p] + \smalop \) as \( p\to \infty \).
Furthermore, 
\begin{equation*}
\begin{aligned}
  \E[\abs{\overline{S} (c)} / p] =&~ \frac{1}{p}\sum_{j\in[p]} \E \left[ \I{\abs{\overline{\widetilde{\eta}} _j} > \widetilde{\theta } ^{*}\vec{\lambda } _j } \right] + \frac{1}{p}\sum_{j\in[p]} \E \left[ \I{\abs{\overline{\widetilde{\eta}} _j} \in (\widetilde{\theta } ^{*}\vec{\lambda } _j(1-c), \widetilde{\theta } ^{*}\vec{\lambda } _j] } \right]. 
\end{aligned}
\end{equation*}
Lemma~\ref{lem:stacked.lasso.indicator.mean} implies $\forall\,\epsilon>0$, \( \exists\,c\in(0,1) \) small enough s.t. \( \sup_p\frac{1}{p}\sum_{j\in[p]} \E [ \I{\abs{\overline{\widetilde{\eta}} _j} \in (\widetilde{\theta } ^{*}\vec{\lambda } _j(1-c), \widetilde{\theta } ^{*}\vec{\lambda } _j) } ] < \epsilon \). By the definition of \( \overline{\widetilde{\eta}}  \), 
\begin{equation*}
\begin{aligned}
  \frac{1}{p}\sum_{j\in[p]} \E \left[ \I{\abs{\overline{\widetilde{\eta}} _j} > \widetilde{\theta } ^{*} \vec{\lambda } _j} \right] =&~ \frac{1}{p} \E \left[ \norm{\overline{\eta }}_0  \right] = \sum_{e\in[E]} \overline{\delta }  _e^{(p)} \to \sum_{e\in[E]} \overline{\delta }  _e,
\end{aligned}
\end{equation*}
where the second equality is by Equation~\eqref{eq:apd.prelim.stack.prop.partial.deriv.trace} and the convergence is by Proposition~\ref{prop:stack.onsager}. 

By the choice of \( c \) so far, we know it only depends on the `idealized' quantities such as \( \overline{\eta }  \) and \( \overline{\widetilde{\eta}}  \), not any specific time step \( t \). As a result,  \(\forall\, \epsilon>0 \), there exists \( c = c(\epsilon) \in(0,1) \) s.t. fixing any \( c'\in(0,c] \) and any \( t \in \mathbb{N_+} \), we have \( \lim_{p} \PR\left[\abs{S^{t}(c')}/p \leq \sum_{e\in[E]} \overline{\delta } _e + \epsilon  \right] = 1 \). This completes the proof. 
  \end{proof}

\subsubsection{Proof of Lemma~\ref{lem:stack.final.converge}}\label{sec:apd.pf.stacked.lasso.final.converge}

We need to introduce Lemma~\ref{lem:stacked.lasso.sparse.singval} as preparation for the proof of Lemma~\ref{lem:stack.final.converge}.

First off, we assume without loss of generality that \( \pi _e >0 \) for all \( e\in[E] \); that is, all the environments \( (y_e, \ix{e}\Sigma _{e}^{1/2}) \) that have contributed to \( \stac \) have been assigned positive weights in the loss function in Equation~\eqref{eq:def.stack.los}. This does not lose generality because by Assumption~\ref{asp:prelim.transfer.model} and \ref{asp:stack.converge.fixed.point}, \( \pi _e = 0 \) iff \( \varpi _e^{*} = 0 \); in other words, when one environment is assigned no weight in the loss function in Equation~\eqref{eq:def.stack.los}, it does not enter the GLAMP formulation either, so we can just ignore it.  

Secondly, let \( N = \sum_{e\in[E]} n_e \), \( \kappa = \lim_{p} p/N = (\sum_{e\in[E]}\kappa _E^{-1})^{-1} \), and define  
\begin{equation}\label{eq:apd.def.stack.allstacked}
  y_{\mathrm{stack}} = \begin{bmatrix}
    \sqrt{\pi_1}y_1 \\
    \sqrt{\pi_2}y_2 \\
    \vdots \\
    \sqrt{\pi_E}y_E 
  \end{bmatrix},
  X_{\mathrm{stack}} = \begin{bmatrix}
    \sqrt{\pi_1} \ix{1} \Sigma _1^{1/2} \\
    \sqrt{\pi_2}\ix{2} \Sigma _2^{1/2}\\
    \vdots \\
    \sqrt{\pi_E}\ix{E} \Sigma _E^{1/2} 
  \end{bmatrix} = \begin{bmatrix}
    \sqrt{\frac{n_1}{N}} \ix{1} \cdot \sqrt{\frac{\pi_1 N}{n_1}} \Sigma _1^{1/2} \\
    \sqrt{\frac{n_2}{N}} \ix{2} \cdot \sqrt{\frac{\pi_2 N}{n_2}}\Sigma _2^{1/2} \\
    \vdots \\
    \sqrt{\frac{n_E}{N}} \ix{E} \cdot \sqrt{\frac{\pi_E N}{n_E}} \Sigma _E^{1/2}.
  \end{bmatrix}
\end{equation}
Now we can rewrite the loss function in Equation~\eqref{eq:def.stack.los} as \( L(b) = \frac{1}{2}\norm{y_{\mathrm{stack}} - X_{\mathrm{stack}}b}_2^{2} +  \norm{\diag(\vec{\lambda } )b}_1 \).

Lastly, define the sparse singular eigenvalues of a matrix. For any matrix \( X \) having \( p \) columns,  any subset \( S \subset [p] \), define \( \kappa _{-}(X,S) = \inf \left\{ \norm{Xu}_2:\Sup(u)\subset S,\norm{u}_2=1 \right\} \)
and the \( k \)-sparse singular value 
\begin{equation}\label{eq:apd.def.sparse.sigval}
  \kappa _{-}(X,k) = \min _{S\subset[p],\abs{S}\leq k} \kappa_{-}(X,S). 
\end{equation}
Of course \( \kappa _{-}(X,k)  \) is decreasing in \( k \) for fixed \( X \); when \( k\in\N_+ \), the definition in Equation~\eqref{eq:apd.def.sparse.sigval} can be equivalently written as \( \kappa _{-}(X,k) = \min _{S\subset[p],\abs{S} = k} \kappa _{-}(X,S) \). 

\begin{lemma}[Lower Bound on the Sparse Singular Values]\label{lem:stacked.lasso.sparse.singval}
  Under Assumption~\ref{asp:prelim.transfer.model} and the definition in Equation~\eqref{eq:apd.def.stack.allstacked}, there exists small enough constants \( \epsilon _0>0 \) and $c_0>0$, s.t. \( \kappa \cdot  (\sum_{e\in[E]}\overline{\delta } _e + \epsilon _0) < 1 \), and \( \forall\,\epsilon \in(0,\epsilon _0) \), 
\begin{equation*}
  \lim_{p\to \infty} \PR \left[ \kappa _{-}\left(X_{\mathrm{stack}}, p(\textstyle\sum_{e\in[E]}\overline{\delta } _e+\epsilon)\right) > c_0 \right]  = 1.
\end{equation*}
\end{lemma}
\begin{proof}[\textbf{Proof of Lemma~\ref{lem:stacked.lasso.sparse.singval}}]
  We use \( X_{\mathrm{temp}} \in \R^{N\times p} \) to denote a matrix of iid \( \Nscr(0,1/N) \) elements, and divide it into blocks as \( (X_{\mathrm{temp}})^{\T} = \left[(X_{\mathrm{temp},1})^{\T}|(X_{\mathrm{temp},2})^{\T}|\cdots |(X_{\mathrm{temp},E})^{\T}\right] \), where \( X_{\mathrm{temp},e}\in\R^{n_e\times p} \). For \( X_{\mathrm{stack}} \) in Equation~\eqref{eq:apd.def.stack.allstacked}, and any index set \( S\subset[p] \), the following equality in distribution holds:
  \begin{equation*}
    (X_{\mathrm{stack}})_{\cdot ,S} \disteq  \begin{bmatrix}
(X_{\mathrm{temp},1})_{\cdot ,S} \left( {\frac{\pi_1 N}{n_1}} (\Sigma _1)_{S,S} \right)^{1/2} \\
(X_{\mathrm{temp},2})_{\cdot ,S} \left( {\frac{\pi_2 N}{n_2}}(\Sigma _2)_{S,S}  \right)^{1/2}\\
      \vdots \\
(X_{\mathrm{temp},E})_{\cdot ,S} \left( {\frac{\pi_E N}{n_E}} (\Sigma _E)_{S,S} \right)^{1/2}.
    \end{bmatrix}
  \end{equation*}
To see why the equality in distribution holds, note that each block of \( (X_{\mathrm{stack}})_{\cdot ,S} \) satisfies \( [\sqrt{\frac{n_e}{N}} \ix{e} \cdot \sqrt{\frac{\pi_e N}{n_e}}\Sigma _e^{1/2}]_{\cdot ,S} \disteq [X_{\mathrm{temp},e}  (\sqrt{\frac{\pi_e N}{n_e}}\Sigma _e)^{1/2}]_{\cdot ,S} = X_{\mathrm{temp},e} [(\sqrt{\frac{\pi_e N}{n_e}}\Sigma _e)^{1/2}]_{\cdot ,S} \in\R^{N\times \abs{ S} } \). 
Each row of \( X_{\mathrm{temp},e} [(\sqrt{\frac{\pi_e N}{n_e}}\Sigma _e)^{1/2}]_{\cdot ,S} \) is an iid draw from \( \Nscr(0,\frac{\pi_e N}{n_e}(\Sigma _e)_{S,S}) \), so the block is equal in distribution to \( (X_{\mathrm{temp},e})_{\cdot ,S} ( {\frac{\pi_e N}{n_e}}(\Sigma _e)_{S,S}  )^{1/2} \).

Since \( \liminf_{p\to \infty} \min _{e\in[E]}\{ \lambda _{\min }(\frac{\pi _e N}{n_e}\Sigma _e)\} \geq \min _{e}\{ \pi _e \kappa _e \} / (c_1 \kappa)  \), there exists some non-random large enough \( p_0 \) s.t. \( \forall\,p\geq p_0 \), \( \min _{S \subset [p],e\in[E]}\{ \lambda _{\min }(\frac{\pi _e N}{n_e}(\Sigma _e)_{S,S})\} \geq \min _{e}\{ \pi _e \kappa _e\}/(2c_1 \kappa) \). For \( p\geq p_0 \), we can invoke Lemma~\ref{lem:apd.sing.val.bounds} to get \( \PR[\sigma _{\min }\left((X_{\mathrm{stack}})_{\cdot ,S}\right) \leq c] \leq \PR[(\min _{e}\{ \pi _e \kappa _e\}/(2c_1 \kappa))^{1/2} \sigma _{\min }((X_{\mathrm{temp}})_{\cdot ,S}) \leq c],\ \forall\,c>0 \). 
When we think about all the possible index sets \( S\subset[p], \abs{S} = k\leq p  \) for some \( k\in\N_+ \), 
\begin{equation*}
  \PR\left[ \kappa _{-}(X_{\mathrm{stack}},k) \leq c \right] \leq \binom{p}{k} \PR \left[ \sigma _{\min }((X_{\mathrm{temp}})_{\cdot ,S}) \leq \frac{c}{(\min _{e}\{ \pi _e \kappa _e\}/(2c_1 \kappa))^{1/2}} \right],\ \forall\,c>0,
\end{equation*}
where we have used a union bound over all possible \( S\subset[p],\abs{S}=k  \) and the fact that \( (X_{\mathrm{temp}})_{\cdot ,S} \) has the same marginal distribution for any \( S\subset[p], \abs{S}=k  \). 

Everything from here follows exactly the same steps as those of Lemma 5, \citet{huang2022lasso} (specifically, after its Equation~(63)). With the same proof, we arrive at the following fact: For \( k\in\N_+ \), if \( \lim_{N} \frac{k}{N} = \omega \in(0,+\infty) \), and additionally \( \omega \in (0,1) \) only when \( \kappa \geq 1 \), then there exists some small constant \( c' \) and large \( p_0 \) s.t. fixing any small \( c\in(0,c') \), the non-asymptotic bound holds \( \forall\,p\geq p_0 \):
  \begin{equation*}
    \PR\left[ \kappa _{-}(X_{\mathrm{stack}},k) \leq c \right] \leq C_1 \exp\left(-C_2\cdot p\right)\quad \text{for some constants }C_1,C_2>0.
  \end{equation*}

  We argue $\lfloor p(\sum_{e\in[E]}\overline{\delta } _e + \epsilon) \rfloor$ qualifies as the choice of `$k$' in the above fact for small enough \( \epsilon>0 \). 
  \( \lim_p p(\sum_{e\in[E]}\overline{\delta } _e + \epsilon) / N = \kappa \, (\sum_{e}\overline{\delta } _e + \epsilon)  \).  Assumption~\ref{asp:stack.converge.fixed.point} implies \( \overline{\delta } _e < \kappa _e^{-1} \), so \( \kappa \sum_{e} \overline{\delta } _e < 1 \). We can make \( \epsilon \) small enough s.t. \( \kappa (\sum_{e} \overline{\delta } _e + \epsilon) \) is strictly smaller than one when it is needed for \( \kappa\geq 1 \). The proof is then complete after taking \( p\to \infty \).
\end{proof}

We are finally ready to prove Lemma~\ref{lem:stack.final.converge}.

\begin{proof}[\textbf{Proof of Lemma~\ref{lem:stack.final.converge}}]
  \newcommand{\yst}{y_{\mathrm{stack}}}
  \newcommand{\xst}{X_{\mathrm{stack}}}
  \newcommand{\Xst}{X_{\mathrm{stack}}}
  We still assume without loss of generality that \( \pi _e >0 \) for all \( e\in[E] \).
Recall the loss function defined in Equation~\eqref{eq:def.stack.los} and re-expressed as \( L(b) = \frac{1}{2} \norm{\yst - \xst \beta}_2^{2} +  \norm{\diag(\vec{\lambda } )b}_1 \). Define \( \Delta \eta^t = \stac - \eta^t \). Since \( \stac \) minimizes the loss, \( L(\eta^t + \Delta \eta^t) \leq L(\eta^t) \). 

Applying Lemma~\ref{lem:apd.sing.val.bounds}{(b)} and then Lemma~\ref{lem:apd.sing.val.limit.concent}{(a)} to \( \xst \) defined in Equation~\eqref{eq:apd.def.stack.allstacked}, we know  \( \sigma _{\max}(\xst) \leq 2(\max_{e}\{ \pi _e \kappa _e \} c_1 / \kappa)^{1/2}  (1 + \sqrt{\kappa}) \) holds almost surely as \( p\to \infty \).

We now define the high-probability event we will be working on throughout the proof. We choose \( \epsilon _0\in(0,1 ) \) small enough and \( M > 1 \) large enough s.t. \( \forall \, \epsilon \in(0,\epsilon _0) \), \( \exists\,T \) s.t. \( \forall\,t\geq T \), at the time step \( t \), the following conditions are satisfied simultaneously with probability going to 1 as \( p\to \infty \): 
\begin{enumerate}
  \item[\it(i)]  \( \max\{\frac{1}{\sqrt{p}}\norm{\stac}_2, \frac{1}{\sqrt{p}}\norm{\eta^t}_2, \frac{1}{\sqrt{p}}\norm{\vec{\lambda } }_2\} \leq M \) (by Lemma~\ref{lem:apd.stack.bounded.lasso}, Lemma~\ref{lem:apd.stack.glamp} and Assumption~\ref{asp:stack.for.glamp.lambd});
  \item[\it(ii)] \( \abs{ S^{t}(c_2)} \leq p(\sum_{e}\overline{\delta } _e + {\epsilon}) \) for some \( c_2=c_2(\epsilon) \in (0,1) \) (by Lemma~\ref{lem:stacked.lasso.bound.on.sc2} with \( c_2 \) depending on the specific \( \epsilon \in(0,\epsilon _0) \) but invariant to the choice of the time step \( t \in \mathbb{N_+} \));
  \item[\it(iii)] \( \kappa _{-}(X_{\mathrm{stack}}, p(\sum_{e\in[E]}\overline{\delta } _e+2{\epsilon})) > c_3 \) for some $c_3 = c_3(\epsilon _0)\in(0,1)$ (by Lemma~\ref{lem:stacked.lasso.sparse.singval} with \( c_3 \) depending on the upper bound \( \epsilon _0 \));
  \item[\it(iv)] \( \frac{1}{\sqrt{p}}\norm{\nabla L(\eta^t)}_2 < \epsilon ^2 \cdot  c_2^{2} c_3^{4} / (2M^{2}) \) (by Lemma~\ref{lem:stacked.lasso.subgradient}; the only condition that determines the choice of \( T \));
  \item[\it(v)] \( \sigma _{\max}(\xst) \leq c_4 = 2(\max_{e}\{ \pi _e \kappa _e \} c_1 / \kappa)^{1/2}  (1 + \sqrt{\kappa}) \). 
\end{enumerate}
In particular, the use of Lemma~\ref{lem:stacked.lasso.sparse.singval} in condition \emph{(iii)} has also required \( \kappa \,(\sum_{e\in[E]}\overline{\delta } _e+2{\epsilon _0}) < 1 \) which will be needed later. We also let \( c_5>0 \) be a value that satisfies \( \liminf_p \min _{j\in[p]} \{ \vec{\lambda } _j \} > c_5  \) by Assumption~\ref{asp:stack.for.glamp.lambd}. 

  Let $S = \Sup(\eta^t) \subset [p]$. Recall \( s^{t} \) defined in Equation~\eqref{eq:def.stack.subgradient} from Lemma~\ref{lem:stacked.lasso.subgradient} is a sub-gradient of \( \norm{\eta^t}_1 \), and \( \nabla L(\eta^t) \) in Lemma~\ref{lem:stacked.lasso.subgradient} is defined using \( s^{t} \). 
  We have the following inequalities:
  \begin{equation*}
  \begin{aligned}
   0\geq &~ \frac{1}{p} \left[ L(\eta^t+\Delta \eta^t) - L(\eta^t) \right] \\
   \geq &~ \frac{1}{p} \left( \norm{\diag(\vec{\lambda } _{S})(\eta^{t} _S+\Delta \eta^t_S)}_1 - \norm{\diag(\vec{\lambda }_S )\eta^t_S}_1 \right) + \frac{1}{p} \norm{\diag(\vec{\lambda } _{S^{c}})\Delta \eta^t _{S^{c}}}_1 + \frac{1}{2p} \norm{\xst \Delta \eta^t}_2^{2}  \\
   &~ -  \frac{1}{p} (\yst - \xst \eta^t)^{\T} \Xst \Delta \eta^t \\
   \geq &~ \underbrace{\frac{1}{p} \left( \norm{\diag(\vec{\lambda } _{S})(\eta^{t} _S+\Delta \eta^t_S)}_1 - \norm{\diag(\vec{\lambda }_S )\eta^t_S}_1 - \diag(\vec{\lambda } _S)\sign(\eta^{t}_S)^{\T} \Delta \eta^t _{S}\right) }_{(1)} \\
   &~ + \underbrace{\frac{1}{p}\left( \norm{\diag(\vec{\lambda } _{S^{c}})\Delta \eta^t _{S^{c}}}_1 - [\diag(\vec{\lambda } _{S^{c}})s^{t}_{S^{c}}]^{^{\T}} \Delta \eta^t_{S^{c}} \right)}_{(2)} 
   + \underbrace{\frac{1}{2p} \norm{\Xst\Delta \eta^t}_2^{2}}_{(3)} \\
   &~ + \frac{1}{p}(\nabla L(\eta^t))^{\T} \Delta \eta^t.
  \end{aligned}
  \end{equation*}
  The terms \( (1) \), \( (2) \) and \( (3) \) are all non-negative. In particular for   \( (1) \), we can use the following inequality for any two real numbers \( a,b \): \( \abs{a+b}\geq \abs{a} + b\cdot \sign(a)  \). Since the terms \( (1) \), \( (2) \) and \( (3) \) are all non-negative, we must have \( (\nabla L(\eta^t))^{\T} \Delta \eta^t \leq 0 \). By conditions \emph{(i)}, \emph{(iv)} of the high-probability event and Cauchy's inequality, we know    
  \begin{equation*}
    \lim_{p\to \infty}\PR\left[(1) + (2) + (3) \leq \frac{1}{p}\abs{(\nabla L(\eta^t))^{\T} \Delta \eta^t } \leq \epsilon ^2 \cdot c_2^{2} c_3^{4} / M \right] = 1.
  \end{equation*}
  Respectively for terms (2) and (3), the above tells us on the high probability event, 
  \begin{equation}
  \begin{aligned}\label{eq:apd.stack.convexity.1}
    0\leq&~ \frac{1}{p}\norm{\diag(\vec{\lambda } _{S^{c}})\Delta \eta^t _{S^{c}}}_1 - \frac{1}{p}[\diag(\vec{\lambda } _{S^{c}})s^{t}_{S^{c}}]^{\T} \Delta \eta^t _{S^{c}} \leq  \epsilon ^2 \cdot c_2^{2} c_3^{4} / M, \\
    0\leq&~ \frac{1}{p} \norm{\Xst \Delta \eta^t}_2^{2} \leq 2\epsilon ^2 \cdot c_2^{2} c_3^{4} / M.
  \end{aligned}
  \end{equation}

  Take the SVD of \( \Xst :\, \Xst = \sum_{i=1}^{N \land p} \sigma _{i} u_{i} v_{i}^{\T} \), assuming \( \sigma _{1} \geq \sigma _{2} \geq \cdots \geq \sigma _{N \land p} \). Let \( V_{(1)} = \Span(v_{i} : \sigma _{i}\geq c_{3}/2) \) and \( V_{(2)} \) be the orthogonal complement of \( V_{(1)} \) in \( \R^{p} \). Let \( P_{(1)},P_{(2)} \) be the projection matrices of \( V_{(1)} \) and \( V_{(2)} \). Denote \( X_{(1)} = \xst P_{(1)}  \) and \( X_{(1)} = \xst P_{(2)}  \). 
   Then \(  X_{(1)}^{\T} X_{(2)} = 0 \). Define \( \Delta \eta_{(1)} = P_{(1)} \Delta \eta^t\), \( \Delta \eta_{(2)} = P_{(2)} \Delta \eta^t\). The second line of Equation~\eqref{eq:apd.stack.convexity.1} implies
  \begin{equation*}
    \frac{1}{p} \norm{X_{(1)}\Delta \eta_{(1)}}_2^{2} \leq  2\epsilon ^2 \cdot c_2^{2} c_3^{4}/M,\quad 
    \frac{1}{p} \norm{X_{(2)}\Delta \eta_{(2)}}_2^{2} \leq  2\epsilon ^2 \cdot c_2^{2} c_3^{4}/M.
  \end{equation*}
  For $\Delta \eta_{(1)}$ corresponding to the larger singular values, we consequently have 
  \begin{equation}\label{eq:apd.stack.convexity.2}
    \frac{1}{p} \norm{\Delta \eta_{(1)}}_2^{2} \leq 8 \epsilon ^2 \cdot c_2^{2} c_3^{4} / (c_3^{2}) \leq  8 \epsilon ^2 \cdot c_2^{2}c_3^{2} / M,
  \end{equation} 
  and we are left to bound \( \frac{1}{p}\norm{\Delta \eta_{(2)}}_2^{2} \). 
  
  We use the first line of Equation~\eqref{eq:apd.stack.convexity.1} to bound \( \frac{1}{p}\norm{\Delta \eta_{(2)}}_2^{2} \). For \( \Delta \eta_{(1)} \), we have 
  \begin{equation*}
    \frac{1}{p} [\diag(\vec{\lambda}_{S^{c}})s^{t}_{S^{c}}]^{\T} (\Delta \eta_{(1)})_{S^{c}} \leq \frac{1}{p}\norm{\diag(\vec{\lambda } _{S^{c}})(\Delta \eta _{(1)})_{S^{c}}}_1 \leq \frac{1}{\sqrt{p}} \norm{\vec{\lambda } _{S^{c}}} \cdot \frac{1}{\sqrt{p}} \norm{\Delta \hat{\beta } _{(1)}}_2  \leq 2 \sqrt{2} \epsilon  c_2 c_3,
  \end{equation*}
  where we have used Cauchy's inequality and the bound on \( \frac{1}{\sqrt{p}}\norm{\vec{\lambda } }_2 \) from \emph{(i)} of the high-probability event.
Since \( \Delta \eta_{(2)} = \Delta \eta ^t - \Delta \eta _{(1)} \), we know 
\begin{equation*}
  \frac{1}{p}\norm{\diag(\vec{\lambda }  _{S^{c}})(\Delta \eta_{(2)})_{S^{c}}}_1 - \frac{1}{p} [\diag(\vec{\lambda } _{S^{c}})s^{t}_{S^{c}}]^{\T} (\Delta \eta_{(2)})_{S^{c}} \leq \epsilon ^2\cdot c_2^{2} c_3^{4} + 4 \sqrt{2} \epsilon c_2 c_3,
\end{equation*}
  where we have used the fact that $M>1$ from (i). Recall that for \( S^{t}(c_2) \) defined in Lemma~\ref{lem:stacked.lasso.bound.on.sc2}, \( S^{t}(c_2)^{c} \subset S^{c} \). Hence,
  \begin{equation*}
\begin{aligned}
  \frac{1}{p}\norm{\diag(\vec{\lambda }  _{S^{t}(c_2)^{c}})(\Delta \eta_{(2)})_{S^{t}(c_2)^{c}}}_1 - \frac{1}{p} [\diag(\vec{\lambda } _{S^{t}(c_2)^{c}})s^{t}_{S^{c}}]^{\T} (\Delta \eta_{(2)})_{S^{t}(c_2)^{c}} \leq&~ \epsilon ^2\cdot c_2^{2} c_3^{4} + 4 \sqrt{2} \epsilon c_2 c_3 \\
  \iff \frac{1}{p}\sum_{j\in S^{t}(c_2)^{c}} \vec{\lambda } _{j} \left[ \abs*{(\Delta \eta _{(2)})_j}  - s^{t}_{j} (\Delta \eta _{(2)})_j \right] \leq &~ \epsilon ^2\cdot c_2^{2} c_3^{4} + 4 \sqrt{2} \epsilon c_2 c_3.
\end{aligned}
  \end{equation*}
  
  By the definition of \( S^{t}(c_2) \), \( \forall\,j\in S^{t}(c_2)^{c}\), \( \abs{s^{t}_j} < 1-c_2  \). Also recall we have defined \( c_5 \) to be a constant s.t. \( \liminf_p \min _{j\in[p]} \{ \vec{\lambda } _j \} > c_5  \). As a result,
  \begin{gather}
\frac{1}{p}\norm{(\Delta \eta_{(2)})_{S^{t}(c_2)^{c}}}_1 \leq \frac{1}{c_2 c_5}\left( \epsilon ^2 \cdot c_2^{2} c_3^{4}/M + 4 \sqrt{2} \epsilon  c_2 c_3 \right) \leq ( 1 + 4 \sqrt{2} )\epsilon c_3 / c_5. \label{eq:apd.stack.convexity.3}
  \end{gather}

Recall we have chosen \( \kappa \,(\sum_{e}\overline{\delta } _e + 2{\epsilon}) < 1 \) as required by Lemma~\ref{lem:stacked.lasso.sparse.singval} used for condition \emph{(iii)} of the high-probability event. Condition \emph{(ii)} guarantees \( \abs{S^{t}(c_2)} \leq p(\sum_{e}\overline{\delta } _e + {\epsilon})  \). Then we immediately know \( \abs{S^{t}(c_2)^{c}} > p {\epsilon}  \), because otherwise \( p = \abs{S^{t}(c_2)^{c}} + \abs{S^{t}(c_2)} \leq p(\sum_{e}\overline{\delta } _e+2{\epsilon}) < p \) would give a contradiction. We partition \( S^{t}(c_2)^{c} \) disjointly into \( S^{t}(c_2)^{c} = \cup_{l=1}^{K} S_l \), with each \( S_l \) satisfying \( \abs{S_l} \in [p {\epsilon}/2,p {\epsilon}] \), and they are constructed in a descending order in the sense that for any \( i\in S_{l} \) and \( j\in S_{l+1} \), \( \abs{({\Delta } \eta_{(2)})_i} \geq \abs{(\Delta \eta_{(2)})_j}  \). By the disjoint partition, 
\begin{equation*}
  \frac{1}{p}\norm{\Delta \eta_{(2)}}_2^{2} =   \frac{1}{p}\norm*{\left(\Delta \eta_{(2)}\right)_{\cup_{l=2}^{K}S_l}}_2^{2} + \frac{1}{p}\norm*{\left(\Delta \eta_{(2)}\right)_{S^{t}(c_2)\cup S_1}}_2^{2}.
\end{equation*}
We bound the two terms on the RHS respectively. For the first term, 
  \begin{equation}
  \begin{aligned}\label{eq:apd.laso.convexity.4}
    \frac{1}{p}\norm*{\left(\Delta \eta_{(2)}\right)_{\cup_{l=2}^{K}S_l}}_2^{2} = &~ 
    \frac{1}{p} \sum_{l=2}^{K} \norm*{\left( \Delta \eta_{(2)} \right)_{S_l}}_2^{2} 
    \leq  \frac{1}{p} \sum_{l=2}^{K} \abs{S_l} \cdot \left( \frac{\norm{(\Delta \eta_{(2)})_{S_{l-1}}}_1}{\abs{S_{l-1}} } \right)^{2} \\
    \leq &~ \frac{1}{p}\cdot \frac{4}{p \epsilon} \sum_{l=2}^{K} \norm{(\Delta \eta_{(2)})_{S_{l-1}}}_1^{2} 
    \leq  \frac{1}{p}\cdot \frac{4}{p {\epsilon}} \left[\sum_{l=2}^{K} \norm{(\Delta \eta_{(2)})_{S_{l-1}}}_1\right]^{2} \\
    \leq &~  \frac{4}{p^{2} {\epsilon}} \norm{(\Delta \eta_{(2)})_{S^{t}(c_2)^{c}}}_1^{2} \leq 4 \epsilon  c_3^{2} (1 + 4 \sqrt{2})^{2} / c_5^{2},
  \end{aligned}
  \end{equation}
where the first inequality is because \( \forall\,i\in S_l \), \( \abs{(\Delta \eta_{(2)})_i}^{2} \leq ({\norm{(\Delta \eta_{(2)})_{S_{l-1}}}_1}/{\abs{S_{l-1}} })^{2} \); the second inequality is because \( \abs{S_l} \in [p {\epsilon}/2,p {\epsilon}] \) for all \( l\in[K] \); the last inequality is by Equation~\eqref{eq:apd.stack.convexity.3}.  

  For the second term, we apply conditions \emph{(iii)}, \emph{(v)}  of the high-probability event:
  \begin{equation*}
  \begin{aligned}
    \frac{1}{p} \norm{(\Delta \eta_{(2)})_{S^{t}(c_2)\cup S_1}}_2^{2} \leq&~ \frac{1}{p c_3^{2}}  \norm*{(\Xst)_{\cdot ,S^{t}(c_2)\cup S_1}(\Delta \eta_{(2)})_{S^{t}(c_2)\cup S_1} }_2^{2} \\
    \leq &~ \frac{2}{p  c_3^{2}}  \left[ \norm*{ \Xst \Delta \eta_{(2)} }_2^{2}  
    + \norm*{ (\Xst)_{\cdot ,\cup_{l=2}^{K}S_l} (\Delta \eta_{(2)})_{_{\cup_{l=2}^{K}S_l}} }_2^{2} \right] \\
    \leq &~ \frac{2}{p c_3^{2}} \left[  \frac{c_3^{2}}{4} \left( \norm*{(\Delta \eta_{(2)})_{S^{t}(c_2)\cup S_1}}_2^{2} + \norm*{(\Delta \eta_{(2)})_{\cup_{l=2}^{K}S_l}}_2^{2} \right)  
    + c_4^{2}\, \norm*{\left(\Delta \eta_{(2)}\right)_{\cup_{l=2}^{K}S_l}}_2^{2} \right] \\
    =&~ \frac{1}{2p} \norm{(\Delta \eta_{(2)})_{S^{t}(c_2)\cup S_1}}_2^{2} + \left(\frac{1}{2} + \frac{2 c_4^{2}}{c_3^{2}}\right) \frac{1}{p} \norm{(\Delta \eta_{(2)})_{\cup_{l=2}^{K}S_l}}_2^{2},
  \end{aligned}
\end{equation*}
where the first inequality has used condition \emph{(iii)}; the second inequality is by the triangle inequality; the terms involving `$c_3$' after the third inequality is because \( \Delta \eta_{(2)} \) is \( \Delta \eta^t \) projected onto \( V_{(2)} \) corresponding to the smaller or null singular values of \( \xst \). Further, we plug in the bound from Equation~\eqref{eq:apd.laso.convexity.4},   
\begin{equation*}
  \begin{aligned}
  \frac{1}{p} \norm{(\Delta \eta_{(2)})_{S^{t}(c_2)\cup S_1}}_2^{2} \leq&~ \left(1 + \frac{4 c_4^{2}}{c_3^{2}}\right) \frac{1}{p} \norm{(\Delta \eta_{(2)})_{\cup_{l=2}^{K}S_l}}_2^{2}
  \leq 4 \epsilon \left(c_3^{2} + 4 c_4^{2}\right) (1 + 4 \sqrt{2})^{2} / c_5^{2},
\end{aligned}
\end{equation*}
which, put together with Equation~\eqref{eq:apd.laso.convexity.4} and Equation~\eqref{eq:apd.stack.convexity.2}, yields
\begin{equation*}
\begin{aligned}
   &~\frac{1}{p}\norm{\Delta \eta_{(2)}}_2^{2} \leq 
   8 \epsilon \left(c_3^{2} + 2 c_4^{2}\right) (1 + 4 \sqrt{2})^{2}/ c_5^{2}  \\
   \implies&~ \frac{1}{p} \norm{\Delta\eta^t}_2^{2} \leq \left[ {8}\epsilon c_2^{2}{c_3^{2}} + 8 \left(c_3^{2} + 2 c_4^{2}\right) (1 + 4 \sqrt{2})^{2}/c_5^{2} \right] \cdot {\epsilon} \leq \left[ {8}+ 8 \left(1 + 2 c_4^{2}\right) (1 + 4 \sqrt{2})^{2}/c_5^{2} \right] \cdot {\epsilon}.
  \end{aligned}
  \end{equation*}
Recall that \( c_4= 2(\max_{e}\{ \pi _e \kappa _e \} c_1 / \kappa)^{1/2}  (1 + \sqrt{\kappa}) \), \( c_5 < \liminf_{p} \min_{j\in[p]} \{ \vec{\lambda } _j \} \). Both are  fixed constants free of \( \epsilon \) or \( p \), \( t \). This completes the proof of Lemma~\ref{lem:stack.final.converge}. To re-iterate what we have proved: There exists some small \( \epsilon _0\in(0,1) \), s.t. \( \forall\,\epsilon \in(0,\epsilon _0) \), \( \exists \, T \) s.t. \( \forall\,t\geq T \), \( \lim_{p} \PR \left[ \frac{1}{p}\norm{\eta^t-\stac}_2^{2} < \epsilon\right] = 1  \).
\end{proof}

\section{Proof of Theorem~\ref{thm:aver}}
\label{sec:apd.aver}

For each environment \( e\in[E] \), \( \hat{\beta } _{e} \) defined in Equation~\eqref{eq:def.aver.los} is a special case of the stacked Lasso \( \stac \), because if we set `\( E=1 \)' in the stacked Lasso, we naturally get the Lasso estimator fit on a single environment. This allows us to use the tools already prepared for the stacked Lasso to study every individual estimator \( \hat{\beta } _{e} \), \( e\in[E] \). 

Specially for the GLAMP part, we let `\( E=1 \)' and get the special cases of Definitions~\ref{def:stack.glamp} and \ref{def:apd.stack.stat.evo} as the following Definitions~\ref{def:aver.glamp} and \ref{def:apd.aver.stat.evo}:
\begin{definition}[Multi-Environment GLAMP Formulation]\label{def:aver.glamp}
  Under Assumptions~\ref{asp:prelim.transfer.model},\ref{asp:aver.for.glamp},\ref{asp:aver.converge.fixed.point}, we define AMP iterates for each environment as below:
\begin{equation}\label{eq:def.aver.glamp}
\begin{aligned}
  v_{\indi,e}^{t} =&~ \ix{e}^{\T} (y_e - r_{\indi,e}^{t-1}) - \Sigma _e^{1/2} \beta _e + \Sigma _e^{1/2} \eta_e(v_{\indi,e}^{t-1}), \\
  r_{\indi,e}^{t} =&~ \ix{e} \Sigma _e^{1/2} \eta_e(v_{\indi,e}^{t}) - \kappa _e \overline{\delta }^{(p)} _{\indi,e} (y_e - r_{\indi,e}^{t-1})
\end{aligned}
\end{equation}
for all \( t \in \mathbb{N_+} \) and \( e\in[E] \). 
In Equation~\eqref{eq:def.aver.glamp}, \( \overline{\delta } _{\indi,e}^{(p)} = \frac{1}{\tau _{\indi,e} p} \E\left[ Z_e^{\T} \Sigma _e^{1/2} \overline{\eta }_e (\tau_e) \right] \) has been in Proposition~\ref{prop:aver.onsager} and specified with \( \tau _{\indi,e} = \tau _{\indi,e}^{*}\,\forall\,e\in[E] \) and \( \theta_e=\theta_{\indi,e} ^* \). 

As initialization, we define $\eta_{e} ^0$ as follows:
Let \( Z_e\sim \Nscr(0,(\tau _{\indi,e}^{*})^{2}I_p),\,\forall\,e\in[E] \) be independent Gaussian variables, and they are also independent of \( \{ \ix{e}: e\in[E] \}  \); define \( \eta_{\indi,e} ^{0} = \eta_e(Z_e),\, \forall\,e\in[E] \). We initialize 
\( r_{\indi,e}^{0} = \ix{e} \Sigma _e^{1/2} \eta_{\indi,e}^0 \), \( v_{\indi,e}^{1} = \ix{e}^{\T}(y_e - r_{\indi,e}^{0}) - \Sigma _e ^{1/2} \beta _e + \Sigma _e^{1/2} \eta_{\indi,e}^0 \), \( \forall\,e\in[E] \).

As a shorthand, define \( \eta_{\indi,e}^t = \eta_e(v_e^{t}) \) for \( t \in \mathbb{N_+} \). 
\end{definition}

\begin{definition}[State Evolution]\label{def:apd.aver.stat.evo} 
Let $\{\Sigma _{(\indi,V,e)}:e\in[E]\}$ be \( E \) infinite dimensional matrices. They are defined separately without using each other. \( \forall\,e\in[E] \),

For all \( i\in\N_+ \), \( (\Sigma _{(\indi,V,e)})_{i,i}=(\tau _{\indi,e}^{*})^{2} \).

For all \( i<j \), \( i,j\in\N_+ \), \( (\Sigma _{(\indi,V,e)})_{i,j}=(\Sigma _{(\indi,V,e)})_{j,i}=(\tau _{\indi,e}^{*})^{2} \rho _{\indi,e} ^{i,j} \), where \(\rho_{\indi,e}^{i,j}\in[0,1]\). For the definition of \( \{\rho ^{i,j}_{\indi,e}\} \): \emph{(i)} For all \( j\geq 2 \), \( \rho ^{1,j}_{\indi,e}=H_{\indi,e}(0) \); \emph{(ii)} For all \( 2\leq i<j \), \( \rho _{\indi,e} ^{i,j}=H_{\indi,e}(\rho _{\indi,e}^{i-1,j-1}) \). \( H_{\indi,e} \) has been defined in Assumption~\ref{asp:aver.converge.cauchy}.

In particular, we define the \( t \)-dimensional top-left block of \( \Sigma _{(\indi,V,e)} \) as \( \Sigma ^t _{(\indi,V,e)} = ( \Sigma _{(\indi,V,e)})_{[t],[t]} \), \( t \in \mathbb{N_+} \). 
\end{definition}

Definitions~\ref{def:aver.glamp} and \ref{def:apd.aver.stat.evo} are just the GLAMP formulation for the stacked Lasso applied to each individual environment and put together. As a corollary of Lemma~\ref{lem:stack.final.converge}, we already know \( \eta _{e}(v^{t}_{\indi,e}) \) converges to \( \hat{\beta } _e \) in the sense that 
\begin{equation}\label{eq:apd.aver.cor.for.past}
  \forall\,\epsilon>0,\ \exists\,T, s.t. \forall\,t\geq T,\ \lim_{p\to \infty}\PR \left[ \frac{1}{p}\norm*{\eta _{e}(v_{\indi,e}^{t}) - \hat{\beta } _{e}}_2^{2} < \epsilon \right] = 1.
\end{equation}

However, results of the stacked Lasso applied to each \( \hat{\beta } _e \) only gives us the asymptotic behavior within one environment. In other words, Theorem~\ref{thm:stack} only implies for a single environment that \( \abs*{\phi(\overline{\eta } _{e},\beta _{e}) - \phi(\hat{\beta } _{e},\beta _e)} = \smalop  \) for any pseudo-Lipschitz function \( \phi \). But this is not enough for the proof of Theorem~\ref{thm:aver}, which is the joint convergence of all the environments. It turns out that the GLAMP iterates in Definition~\ref{def:aver.glamp} can also be cast in a multi-environment GLAMP instance, which gives us the desired  joint convergence of all environments. This is shown by Lemma~\ref{lem:apd.aver.glamp}:

\begin{lemma}[GLAMP Convergence]\label{lem:apd.aver.glamp}
  Under Assumptions~\ref{asp:prelim.transfer.model},\ref{asp:aver.for.glamp},\ref{asp:aver.converge}, and Definitions~\ref{def:aver.glamp},\ref{def:apd.aver.stat.evo},
 \( \forall\,t \in \mathbb{N_+} \), \( \forall\, \) sequence of order-\( k \) pseudo-Lipschitz functions \( \phi:(\R^{p})^{(t+2)E}\to\R \), \( k\geq 1 \),
 \begin{equation*}
\begin{alignedat}{2}
  &\phi\left(v_{\indi,1}^{t},\cdots,v_{\indi,E}^{t},v_{\indi,1}^{t-1},\cdots,v_{\indi,E}^{t-1},\cdots ,v_{\indi,1}^{1},\cdots,v_{\indi,E}^{1},\eta^0_{1},\cdots ,\eta^0_{E},\beta _1,\beta_2,\cdots ,\beta _E \right)&& \\
 &= \E \left[ \phi\left( Z^{t}_1,\cdots ,Z^{t}_E,Z^{t-1}_1,\cdots ,Z^{t-1}_E,\cdots ,Z^{1}_1,\cdots ,Z^{1}_E,\eta _{1}(Z_1^{0}),\cdots ,\eta _{E}(Z_E^{0}),\beta_1,\beta _2,\cdots ,\beta _E \right) \right] + \smalop.
\end{alignedat}
\end{equation*}
In the expectation on the RHS, we have \emph{(i)} \( \{ Z_{e}^{i}:e\in[E],0\leq i\leq t \}  \) are jointly Gaussian; \emph{(ii)} \( \{ Z_{e}^{0}:e\in[E] \} \independent \{ Z_{e}^{i}:e\in[E],i\geq 1 \}  \) and \( \{ Z_{e}^{0}:e\in[E] \}\) follows the distribution described in Definition~\ref{def:apd.aver.stat.evo}; \emph{(iii)} For the other Gaussian vectors, the matrices \( [Z_e^{1}|Z_e^{2}|\cdots |Z_e^{t}]\in\R^{p\times t} \) has iid rows drawn from \( \Nscr(0,\Sigma _{(\indi,V,e)}^{t}) \), and the \( E \) matrices are mutually independent. 
\end{lemma}

Lemma~\ref{lem:apd.aver.glamp} is proved in Section~\ref{sec:apd.pf.aver.glamp}. Now we are ready to prove Theorem~\ref{thm:aver}. 

\begin{proof}[\textbf{Proof of Theorem~\ref{thm:aver}}]
Applying Lemma~\ref{lem:apd.stack.bounded.lasso} to each \( \hat{\beta } _{e} \), we know \( \exists\, M \) s.t. \( \sup_{e}\limsup_{p} \frac{1}{p}\norm{\hat{\beta } _{e}}_2^{2} \leq M \) a.s..
Combined with Equation~\eqref{eq:apd.aver.cor.for.past}, we can choose $M$ large enough, s.t. for any order-\( k \) pseudo-Lipschitz function \( \phi : \R^{2E}\to\R \), \( k\geq 1 \), \( \forall\,\epsilon > 0 \), \( \exists\, T  \) s.t. for \( t\geq T \), 
\begin{equation}
  \begin{aligned}\label{eq:apd.aver.intermed}
    &~\abs*{\phi \left(\eta_{1}(v^{t}_{\indi,1}),\cdots ,\eta _{E}(v^{t}_{\indi,E}),\beta_1,\cdots ,\beta _E\right) - \phi(\hat{\beta } _{1},\cdots ,\hat{\beta } _{E},\beta_1,\cdots ,\beta _E)} \\
    \leq&~ L \sum_{e\in[E]}\left[ 1 + \left( \frac{\norm{\beta _e}_2}{\sqrt{p}} \right)^{k-1} + \left( \frac{\norm{\eta _{e}^t(v^{t}_{\indi,e})}_2}{\sqrt{p}} \right)^{k-1} + \left( \frac{\norm{\hat{\beta } _{e}}_2}{\sqrt{p}} \right)^{k-1} \right] 
    \sum_{e\in[E]}\frac{\norm{\eta _{e}(v^{t}_{\indi,e}) - \hat{\beta } _e}_2}{\sqrt{p}} \\
    \leq&~ 3L\cdot E \cdot M^{(k-1)/2} \sum_{e\in[E]}\frac{\norm{\eta _{e}(v^{t}_{\indi,e}) - \hat{\beta } _e}_2}{\sqrt{p}} \leq \epsilon
  \end{aligned}
   \end{equation}
holds with probability approaching one as \( p \to \infty \). 

For the same \( \epsilon \) and \( t \), Lemma~\ref{lem:apd.aver.glamp} implies 
\begin{equation*}
  \abs*{\phi \left(\eta_{1}(v^{t}_{\indi,1}),\cdots ,\eta _{e}(v^{t}_{\indi,E}),\beta_1,\cdots ,\beta _E\right) - \E \left[ \phi(\overline{\eta } _{1},\cdots ,\overline{\eta } _{E},\beta_1,\cdots ,\beta _E) \right]} \leq \epsilon
\end{equation*}
with probability approaching one as \( p \to \infty \). Combining the above line with Equation~\eqref{eq:apd.aver.intermed} proves Theorem~\ref{thm:aver}. 
\end{proof}

\subsection{Proof of Lemma~\ref{lem:apd.aver.glamp}}\label{sec:apd.pf.aver.glamp}

\begin{proof}[\textbf{Proof of Lemma~\ref{lem:apd.aver.glamp}}]
Just like the proof of Lemma~\ref{lem:apd.stack.glamp}, we specify a multi-environment GLAMP instance that has the form and state evolution as described in Definitions~\ref{def:aver.glamp} and \ref{def:apd.aver.stat.evo}. Also analogously to the proof of Lemma~\ref{lem:apd.stack.glamp}, we can use a sub-sequence based argument to address the technicality caused by using random initialization in the GLAMP instance. We then move on to describe the instance itself. 

Same as the GLAMP instance for the stack Lasso estimator, let \( \{ q_e:e\in[E] \}  \) all be sufficiently large and fixed integers, let \( B= [\beta _1|\beta_2|\cdots |\beta _{E}|\eta^0|0|0|\cdots ]\in\R^{p\times q} \), and \( W_e=[w_e|0|0|\cdots ]\in\R^{n_e\times q} \). We add `ind' to the subscripts to differentiate from the iterates for the stacked Lasso; starting from \( V_{\indi,e}^{0}=0 \) and \( R_{\indi,e}^{0} = [\ix{e}\Sigma _{e}^{1/2}\beta _{e}|\ix{e}\Sigma _{e}^{1/2}\eta _{\indi,e}^{0}|0|\cdots ] \), we specify the GLAMP instance as follows: (\( t \in \mathbb{N_+} \))
  \begin{equation*}
  \begin{aligned}
   \forall\,e\in[E],\ V_{\indi,e}^{t} =&~ \left[v_{\indi,e}^{1}|v_{\indi,e}^{2}|\cdots |v_{\indi,e}^{t}|0|0|\cdots \right], \\
   \eta _{\indi,e}^{t}(V_{\indi,e}^{t}) =&~ \Sigma _e^{1/2}\left[\beta _e|\eta _{\indi,e}^0|\eta _{e}(v_{\indi,e}^{1})|\eta _{e}(v_{\indi,e}^{2})|\cdots |\eta _{e}(v_{\indi,e}^{t})|0|0|\cdots \right],\\
    R_{\indi,e}^{t} =&~ [\ix{e}\Sigma _e^{1/2}\beta _e|r_{\indi,e}^{0}|r_{\indi,e}^{1}|\cdots |r_{\indi,e}^{t}|0|0|\cdots ],\\
    \Psi _{\indi,e}^{t}(R_e^{t}) =&~ \kappa _e[ y_e - r_{\indi,e}^{0}|y_e - r_{\indi,e}^{1}|\cdots |y_e - r_{\indi,e}^{t}|0|0|\cdots  ]\quad (y_e  = \ix{e} \Sigma _e^{1/2}\beta _e + w_e).
  \end{aligned}
  \end{equation*}
   
  Now we know the GLAMP iterates take the form 
  \begin{equation*}
    \begin{aligned}
      v_{\indi,e}^{t} =&~ \ix{e}^{\T} (y_e - r_{\indi,e}^{t-1}) - \Sigma _e^{1/2} \beta _e + \Sigma _e^{1/2} \eta_e(v_{\indi,e}^{t-1}), \\
      r_{\indi,e}^{t} =&~ \ix{e} \Sigma _e^{1/2} \eta_e(v_{\indi,e}^{t}) - \kappa _e (y_e - r_{\indi,e}^{t-1})\cdot \frac{1}{p} \E \left[ \trace \left( \Sigma _{e}^{1/2} \frac{\partial \eta _E}{\partial v_{\indi,e}^{t}}(Z_{e}^{t}) \right) \right], 
    \end{aligned}
    \end{equation*}
  where the Gaussian variables used are characterized by the state evolution which we are yet to study. 
  
  To fully write out the form of the GLAMP iterates, we need to find the state evolution first. We discard the all-zero columns in \( V_e^{t} \) and \( R_e^{t} \),  start from \( \Sigma _{(\indi,R,e)}^{0} \) and compute the first two covariance matrices for the state evolution. From there,  we prove the whole state evolution for \( V^{t}_e \) with induction.  
  \begin{enumerate}
    \item \( \Sigma _{(\indi,R,e)}^{0} = \begin{bmatrix}
      \lim_p \frac{1}{p}\beta _e^{\T}\Sigma _e \beta _e & \lim_p \frac{1}{p}\E[\beta _e^{\T} \Sigma _e \overline{\eta }_{e} ]  \\ 
      \lim_p \frac{1}{p}\E[\beta _e^{\T} \Sigma _e \overline{\eta }_{e}] &  \lim_p \frac{1}{p}\E[\overline{\eta }_e^{\T} \Sigma _e \overline{\eta }_e] \\ 
    \end{bmatrix} \).
    \item \( \Sigma _{(\indi,V,e)}^{1} = \lim_p \frac{1}{n_e}\E\left[\norm*{w_e+(Z_{(\indi,R,e)}^{0})_1 - (Z_{(\indi,R,e)}^{0})_2}_2^{2}\right]  \) where \( Z^{0}_{(\indi,R,e)}=[(Z_{(\indi,R,e)}^{0})_1|(Z_{(\indi,R,e)}^{0})_2]\in\R^{p\times 2}\) has iid rows drawn from \( \Nscr(0,\kappa _e \Sigma _{(\indi,R,e)}^{0}) \). Simplifying it with Assumption~\ref{asp:aver.for.glamp} and \ref{asp:aver.converge}, 
    
    \( \Sigma _{(\indi,V,e)}^{1} = \E[W_e^{2}] + \kappa _e \lim_p \frac{1}{p}\E[\norm{\overline{\eta }_{e} -\beta _e}_{\Sigma _{(\indi,e)}}^{2}] = (\tau _{\indi,e}^{*}) ^{2} \).
    \item  \( \Sigma _{(\indi,R,e)}^{1} = \begin{bmatrix}
      \lim_p \frac{1}{p}\beta _e^{\T}\Sigma _e \beta _e & 
      \lim_p \frac{1}{p}\E[\beta _e^{\T} \Sigma _e \overline{\eta }_e ] & 
      \lim_p \frac{1}{p}\E[\beta _e^{\T} \Sigma _e \overline{\eta }_e ]   \\ 
      \lim_p \frac{1}{p}\E[\beta _e^{\T} \Sigma _e \overline{\eta }_e] &  
      \lim_p \frac{1}{p}\E[\overline{\eta }_{e}^{\T} \Sigma _e \overline{\eta }_{e}] & 
      \lim_p\frac{1}{p}\E\left[ \overline{\eta }_{e,(1)} ^{\T}\Sigma _e \overline{\eta }_{e,(2)}  \right] \\ 
      \lim_p \frac{1}{p}\E[\beta _e^{\T} \Sigma _e \overline{\eta } ] & \lim_p \frac{1}{p}\E\left[ \overline{\eta }_{e,(1)} ^{\T}\Sigma _e \overline{\eta }_{e,(2)}  \right] & 
      \lim_p \frac{1}{p}\E[\overline{\eta }_e^{\T} \Sigma _e \overline{\eta }_e ],
    \end{bmatrix} \) 
    where \( \overline{\eta } _{e,(1)} \) and \( \overline{\eta } _{e,(2)} \) are mutually independent copies of \( \overline{\eta }_e  \).
    \item \( \Sigma _{(\indi,V,e)}^{2} = \begin{bmatrix}
      (\tau _{\indi,e}^{*})^{2} & (\tau _{\indi,e}^{*})^{2}H_{\indi,e}({0}) \\ 
      (\tau _{\indi,e}^{*})^{2}H_{\indi,e}({0}) & (\tau _{\indi,e}^{*})^{2} \\ 
    \end{bmatrix} \). 
    For the diagnoal elements of \( \Sigma _{(\indi,V,e)}^{2} \) we have simplified it in the same way as \( \Sigma _{(\indi,V,e)}^{1} \). For the off-diagonal term, it is 
    \begin{equation*}
  \begin{aligned}
    &\lim_p \frac{1}{n_e}\E \left[ (w_e+ (Z_{(\indi,R,e)}^{0})_1 - (Z_{(\indi,R,e)}^{0})_2)^{\T}(w_e+ (Z_{(\indi,R,e)}^{0})_1 - (Z_{(\indi,R,e)}^{0})_3)\right]\\
     &= \E[W_e^{2}] + \kappa _e \lim_p \frac{1}{p}\E \left[ (\overline{\eta }_{(1)}-\beta _e)^{\T}\Sigma _e(\overline{\eta }_{(2)} -\beta _e) \right] 
     =(\tau _{\indi,e}^{*})^{2}H_{\indi,e}({0}),
  \end{aligned}
    \end{equation*}
    where \( \overline{\eta } _{e,(1)} \) and \( \overline{\eta } _{e,(2)} \) are mutually independent copies of \( \overline{\eta }_e  \).
  \end{enumerate}
  We use induction to prove the whole state evolution of \( V_{\indi,e}^{t} \) takes the form of Definition~\ref{def:apd.aver.stat.evo}. Suppose the form in Definition~\ref{def:apd.aver.stat.evo} has been verified up to and including \( t \); let's prove the case of \( (t+1) \). 
  
  For the diagonal term,  \( (\Sigma _{(\indi,V,e)}^{t+1})_{t+1,t+1} = \lim_p \frac{1}{n_e}\E\left[\norm*{w_e+(Z_{(\indi,R,e)}^{0})_1 - (Z_{(\indi,R,e)}^{t})_{t+2}}_{2}^{2}\right] = (\tau _{\indi,e}^{*})^{2}  \) using the same simplification as that of \( \Sigma _{(\indi,V,e)}^{1} \).
    
  For the off-diagonal terms: 
  When \( i=1 \), \( j=t+1 \), we show \( \rho _{\indi,e} ^{1,t+1}=H_{\indi,e}(0) \) by 
  \begin{equation*}
  \begin{aligned}
    (\tau _{\indi,e}^{*})^{2} \rho _{\indi,e}^{1,t+1} =&~\lim_p \frac{1}{n_e}\E \left[ (w_e+ (Z_{(\indi,R,e)}^{t})_1 - (Z_{(\indi,R,e)}^{t})_2)^{\T}(w_e+ (Z_{(\indi,R,e)}^{t})_1 - (Z_{(\indi,R,e)}^{t})_{t+2})\right]\\
    =&~ \E[W_e^{2}] + \kappa _e \lim_p \frac{1}{p} (\eta _{\indi,e}^0-\beta _e)^{\T} \Sigma _e \E \left[ \overline{\eta }_e-\beta _e \right] \\
    =&~ \E[W_e^{2}] + \kappa _e \lim_p \frac{1}{p}\E \left[ (\overline{\eta }_{e,(1)}-\beta _e)^{\T}\Sigma _e(\overline{\eta }_{e,(2)} -\beta _e) \right] \\
    =&~(\tau _{\indi,e}^{*})^{2}H_{\indi,e}(0),
  \end{aligned}
  \end{equation*}
  where \( \overline{\eta } _{e,(1)} \) and \( \overline{\eta } _{e,(2)} \) are mutually independent copies of \( \overline{\eta }_{e}  \).
  
  When \( 1<i<j=t+1 \), we show \( \rho _{\indi,e} ^{i,t+1}=H(\rho_{\indi,e} ^{i-1,t}) \) by
  \begin{equation*}
    \begin{aligned}
      (\tau _{\indi,e}^{*})^{2} \rho _{\indi,e}^{i,t+1} =&~\lim_p \frac{1}{n_e}\E \left[ (w_e+ (Z_{(\indi,R,e)}^{t})_1 - (Z_{(\indi,R,e)}^{t})_{i+1})^{\T}(w_e+ (Z_{(\indi,R,e)}^{t})_1 - (Z_{(\indi,R,e)}^{t})_{t+2})\right]\\
      =&~ \E[W_e^{2}] + \kappa _e \lim_p \frac{1}{p}\E \left[ ({\eta }_e(Z^{i-1}_e)-\beta _e)^{\T}\Sigma _e({\eta }_{e}(Z^{t}_{e}) -\beta _e) \right],
    \end{aligned}
    \end{equation*}
  where \( Z^{i-1}_{e}\in\R^{p} \) and \( Z^{t}_e \in\R^{p} \) are jointly Gaussian, with each row of \( [Z^{i-1}_e|Z^{t}_e] \) being sampled iid from 
    \( \Nscr(0,(\tau _{\indi,e}^{*})^{2}[1,\rho _{\indi,e}^{i-1,t} ; \rho _{\indi,e}^{i-1,t}, 1]) \). Thus \( \rho _{\indi,e}^{1,t+1}=H_{\indi,e}(\rho _{\indi,e} ^{i-1,t}) \). 
  
  Now we have verified the GLAMP iterates take the form of Definition~\ref{def:aver.glamp} and their state evolution is described by Definition~\ref{def:apd.aver.stat.evo}. By GLAMP convergence, we have shown Lemma~\ref{lem:apd.aver.glamp}. 
  \end{proof}

\section{Proof of Theorem~\ref{thm:2nd}}
\label{sec:apd.2nd}

\begin{lemma}[Boundedness of \( \hat{\beta } _{\rt} \)]\label{lem:2nd.bounded}
  Under Assumption~\ref{asp:prelim.transfer.model}, there exists \( M > 0 \) s.t. for \( \hat{\beta } _{\rt} \) defined in Equation~\eqref{eq:def.2nd.los}, \( \limsup_{p} \frac{1}{p}\norm{\se}_2^{2} \leq M \), a.s. 
  \end{lemma}
  
  Lemma~\ref{lem:2nd.bounded} is proved in Section~\ref{sec:apd.pf.2nd.bounded}.

\begin{definition}[Augmented GLAMP Formulation]\label{def:2nd.amp}
    Under Assumptions~\ref{asp:prelim.transfer.model},\ref{asp:stack.for.glamp},\ref{asp:stack.converge},\ref{asp:2nd.for.glamp},\ref{asp:2nd.converge}, on top of Definition~\ref{def:stack.glamp} and \ref{def:apd.stack.stat.evo}, we define an augmented multi-environment GLAMP instance. With parameter \( t_{0}\in\N_+ \), 
    \begin{equation}\label{eq:def.2nd.glamp}
      \begin{aligned}
  &\begin{aligned}
    \mathrm{Part}\ \ro:\ 
    \forall\,t \in \mathbb{N_+},e\in[E],\ v_e^{t} =&~ \ix{e}^{\T} (y_e - r_e^{t-1}) - \Sigma _e^{1/2} \beta _e + \Sigma _e^{1/2} \eta(v_1^{t-1}, v_2^{t-1}, \cdots ,v_1^{t-1}), \\
    r_e^{t} =&~ \ix{e} \Sigma _e^{1/2} \eta(v_1^{t}, v_2^{t}, \cdots ,v_1^{t}) - \kappa _e \overline{\delta }^{(p)} _e (y_e - r_e^{t-1}) .
  \end{aligned}\\
  &\begin{aligned}
    \mathrm{Part}\ \rt:\ 
        v_{\rt}^{0} =&~ \zeta ^* \cdot \tau _{\rt}^{*} v_{1}^{t_0}/\tau _{1}^{*}+\sqrt{1-(\zeta ^*)^{2}}\cdot \tau _{\rt}^{*} Z_{1}^{\mathrm{init}},
        \\
        r_{\rt}^{0} =&~  \ix{1}\Sigma_{1}^{1/2} \xi\left(v_{\rt}^{0},\,\eta(v_1^{t_0}, \cdots ,v_1^{t_0}) \right) - (\overline{\gamma } _{\rt} \tau _{\rt}^{*}\zeta ^* / \tau ^*_{1} + \overline{\gamma } _{\ro}) \cdot  \kappa_1(y_1 - r_{1}^{t_0-1});\\
        \forall\,t \in \mathbb{N_+},\; v_{\rt}^{t} =&~ \ix{1}^{\T}(y_1-r_{\rt}^{t-1}-\kappa_1 \overline{\gamma } _{\ro}(y_1 - r_1^{t+t_0-1})) - (1-\kappa_1\overline{\gamma } ^{(p)}_{\ro})\Sigma_{1}^{1/2}\beta _1 - \kappa_1 \overline{\gamma } _{\ro}^{(p)} \Sigma_{1}^{1/2} \hat{\beta }^{t+t_0-1} \\
        &~ + \Sigma_{1}^{1/2} \xi(v_{\rt}^{t-1},\eta(v_1^{t+t_0-1},\cdots ,v_{1}^{t+t_0-1})) , \\
      r_{\rt}^{t}=&~ \ix{1} \Sigma_{1}^{1/2}\xi(v_{\rt}^{t},\eta(v_1^{t+t_0},\cdots ,v_{1}^{t+t_0})) - \kappa_1 \overline{\gamma } _{\ro}^{(p)} (y_1 - r_1^{t+t_0-1}) \\&~- \kappa_1 \overline{\gamma } _{\rt}^{(p)} (y_1 - r_{\rt}^{t-1}-\kappa_1 \overline{\gamma } _{\ro}(y_1 - r_1^{t+t_0-1})).
      \end{aligned}
  \end{aligned}
      \end{equation}
    In Equation~\eqref{eq:def.2nd.glamp}, \( \overline{\gamma } ^{(p)}_{\ro} \), \( \overline{\gamma } ^{(p)}_{\rt} \) are as defined in Proposition~\ref{prop:2nd.onsager} and further specified at \( \tau _{\rt} = \tau _{\rt}^{*},\ \zeta = \zeta ^* \). The initialization of Part I in Equation~\eqref{eq:def.2nd.glamp} is unchanged from that of Equation~\eqref{eq:def.stack.glamp} in Definition~\ref{def:stack.glamp}. The initialization of Part II has used an external $Z_{1}^{\mathrm{init}}\in\R^{p}$, \( Z_{1}^{\mathrm{init}}\sim N(0,I_p) \) and \( Z_{1}^{\mathrm{init}}\independent \{ \ix{1},\cdots ,\ix{E}, \hat{\beta } ^{0} \}  \).  
  \end{definition}
   
  \begin{lemma}[GLAMP Convergence]\label{lem:2nd.glamp}
    Under Assumptions~\ref{asp:prelim.transfer.model},\ref{asp:stack.for.glamp},\ref{asp:stack.converge},\ref{asp:2nd.for.glamp},\ref{asp:2nd.converge}, \( \exists \,t_0 \) sufficiently large, s.t.  
   \( \forall\,t \in \mathbb{N_+} \), \( \forall\, \) sequence of order-\( k \) pseudo-Lipschitz functions \( \phi:(\R^{p})^{(t+1)\times (E+1)+E}\to\R \), \( k\geq 1 \),
   \begin{equation*}
  \begin{alignedat}{2}
    &\phi\left(v_{\rt}^{t},v_{\rt}^{t-1},\cdots ,v_{\rt}^{1},v_{\rt}^{0},v_1^{t+t_0},\cdots,v_E^{t+t_0},v_1^{t+t_0-1},\cdots,v_E^{t+t_0-1},\cdots ,v_1^{t_0},\cdots,v_E^{t_0},\beta _1,\beta_2,\cdots ,\beta_E \right) + \smalop&& \\
   &= \E \left[ \phi\left( Z_{\rt}^{t},Z_{\rt}^{t-1},\cdots ,Z_{\rt}^{1},Z_{\rt}^{0}, Z^{t+t_0}_1,\cdots ,Z^{t+t_0}_E,Z^{t+t_0-1}_1,\cdots ,Z^{t+t_0-1}_E,\cdots ,Z^{t_0}_1,\cdots ,Z^{t_0}_E,\beta_1,\beta _2,\cdots ,\beta _E \right) \right].
  \end{alignedat}
  \end{equation*}
  For the Gaussian variables used in the limit: \emph{(i)} All of them are jointly Gaussian with a fully determined distribution, and \( \{ Z_{e}^{t_0+i}:e\in[E],0\leq i\leq t \}  \) are as they are defined in Lemma~\ref{lem:apd.stack.glamp}; \emph{(ii)} \( \{ Z_{\rt}^{i}:0\leq i\leq t \} \cup \{ Z_{1}^{t_0+i}:0\leq i\leq t \}  \) are independent from ther rest of the Gaussian variables; \emph{(iii)} \( \{ Z_{\rt}^{i}:0\leq i\leq t \} \cup \{ Z_{1}^{t_0+i}:0\leq i\leq t \}  \) marginally satisfies: \( [Z_{\rt}^{i}|Z_{1}^{i+t_0}] \) has iid rows drawn from \( N(0,[(\tau ^*_{\rt})^{2},\,\tau ^*_{\rt}\tau ^*_{1}\zeta ^*;\, \tau ^*_{\rt}\tau ^*_{1}\zeta ^*,\, (\tau ^*_{1})^{2}]) \) \( \forall\,0\leq i\leq t \).
  \end{lemma}

  Lemma~\ref{lem:2nd.glamp} is proved in Section~\ref{sec:apd.pf.lem:2nd.glamp}.

  \begin{definition}[Induced AMP]\label{def:2nd.induced.glamp}
    Under Assumptions~\ref{asp:prelim.transfer.model},\ref{asp:stack.for.glamp},\ref{asp:stack.converge},\ref{asp:2nd.for.glamp},\ref{asp:2nd.converge}, on top of Definition~\ref{def:stack.glamp} and \ref{def:2nd.amp},  we define an AMP instance induced by Part II of Definition~\ref{def:2nd.amp}.
  We firstly define two objects using \( \stac \) already studied: 
  \begin{equation*}
    v_1^{\infty} = \frac{\ix{1}^{\T}(y_1 - \ix{1}\Sigma_1^{1/2}\stac)}{1-\kappa_1 \overline{\delta }_1 } + \Sigma _1^{1/2}(\stac - \beta _{E}),\quad r_1^{\infty} = y_1 - \frac{y_1- \ix{1}\Sigma_1^{1/2}\stac}{1-\kappa_1 \overline{\delta }_ {1}}.
  \end{equation*}
  Then we define the following AMP iterates as a limiting version of Part II of Definition~\ref{def:2nd.amp}. 
    \begin{equation}\label{eq:def.2nd.induced.amp}
      \begin{aligned}
        v_{\rti}^{0} =&~ \zeta ^* \cdot \tau _{\rt}^{*} v_1^{\infty}/\tau _{E}^{*}+\sqrt{1-(\zeta ^*)^{2}}\cdot \tau _{\rt}^{*} Z_{1}^{\mathrm{init}},
        \\
        r_{\rti}^{0} =&~  \ix{1}\Sigma_1^{1/2} \xi\left(v_{\rti}^{0},\stac \right) - (\overline{\gamma } _{\rt} \tau _{\rt}^{*}\zeta ^* / \tau ^*_{E} + \overline{\gamma } _{\ro}) \kappa_1(y_1 - r_1^{\infty});\\
        \forall\,t \in \mathbb{N_+},\; v_{\rti}^{t} =&~ \ix{1}^{\T}(y_1-r_{\rti}^{t-1}-\kappa_1 \overline{\gamma } _{\ro}(y_1 - r_1^{\infty})) - (1-\kappa_1\overline{\gamma } ^{(p)}_{\ro})\Sigma_1^{1/2}\beta _{1} - \kappa_1 \overline{\gamma } _{\ro}^{(p)} \Sigma_1^{1/2} \stac \\
        &~ + \Sigma_1^{1/2} \xi(v_{\rti}^{t-1},\stac) , \\
      r_{\rti}^{t}=&~ \ix{1} \Sigma_1^{1/2}\xi(v_{\rti}^{t},\stac) - \kappa_1 \overline{\gamma } _{\ro}^{(p)} (y_1 - r_1^{\infty}) - \kappa_1 \overline{\gamma } _{\rt}^{(p)} (y_1 - r_{\rti}^{t-1}-\kappa_1 \overline{\gamma } _{\ro}(y_1 - r_1^{\infty})).
      \end{aligned}
      \end{equation}
   The initial \( v_{\rti}^{0} \) has used  $Z_{1}^{\mathrm{init}}\in\R^{p}$, \( Z_{1}^{\mathrm{init}}\sim N(0,I_p) \) and \(Z_{1}^{\mathrm{init}}\independent \{ \ix{1},\cdots ,\ix{E}, \eta ^{0} \}  \). As a shorthand, define \( \xi^{t} = \xi(v_{\rti}^{t},\stac),\; \forall\,t \in \mathbb{N} \). 
  \end{definition}

  \begin{definition}[State Evolution]\label{def:2nd.induced.evo} 
    Under Assumptions~\ref{asp:prelim.transfer.model},\ref{asp:stack.for.glamp},\ref{asp:stack.converge},\ref{asp:2nd.for.glamp},\ref{asp:2nd.converge}, we define \( \Sigma _{(V,\rt)} \in (\R^{\N_{\geq 0}})\times(\R^{\N_{\geq 0}}) \) as an infinite-dimensional matrix whose indices start from zero: 
    
      For all \( i\in\N_{\geq 0} \), \( (\Sigma _{(V,\rt)})_{i,i}=1 \). 
    
      For all \( i<j \), \( (\Sigma _{(V,\rt)})_{i,j}=(\Sigma _{(V,\rt)})_{j,i}= \rho _\rt ^{i,j}\in\R \). For the definition of \( \{\rho _{\rt} ^{i,j}\} \): \emph{(i)} For all \( j\geq 1 \), \( \rho _{\rt} ^{0,j}=0 \); \emph{(ii)} For all \( 1\leq i<j \), \( \rho_{\rt} ^{i,j}=(H_{\rt}(\rho_{\rt} ^{i-1,j-1}) - (\zeta ^*)^{2})/(1 - (\zeta ^*)^{2}) \), or \( \rho_{\rt} ^{i,j}=1 \) in the case of \( \abs{\zeta ^*}=1 \). \( H_{\rt} \) has been defined in Assumption~\ref{asp:2nd.converge.cauchy}.
    
      In particular, we define \( \Sigma ^t _{(V,\rt)}  \) to be the \( (t+1) \)-dimensional top-left block of \( \Sigma _{(V,\rt)} \), \( \Sigma ^t _{(V,\rt)} = ( \Sigma _{(V,\rt)})_{\{ 0 \} \cup [t],\{ 0 \} \cup[t]} \), \( t \in \mathbb{N_+} \). 
    \end{definition}

\begin{lemma}[Induced AMP Convergence]\label{lem:2nd.amp.plus}
  Under Assumptions~\ref{asp:prelim.transfer.model},\ref{asp:stack.for.glamp},\ref{asp:stack.converge},\ref{asp:2nd.for.glamp},\ref{asp:2nd.converge}, \( \forall\,t \in \mathbb{N_+} \), \( \forall\, \) sequence of order-\( k \) pseudo-Lipschitz functions \( \phi:(\R^p)^{(t+E+3)}\to \R \),
  \begin{equation*}
\begin{aligned}
  &\begin{alignedat}{2}
    &\phi\left(v_{\rti}^{t},v_{\rti}^{t-1},\cdots ,v_{\rti}^{1},v_{\rti}^{0},v_1^{\infty},\stac,\beta _1,\beta_2,\cdots ,\beta _E \right) + \smalop&& \\
   &= \E \left[ \phi\left( Z_{\rti}^{t} , Z_{\rti}^{t-1}, \cdots , Z_{\rti}^{1}, Z_{\rti}^{0},  Z_{1},\eta(Z_1,\cdots ,Z_{E}) ,\beta_1,\beta _2,\cdots ,\beta _E \right) \right], 
  \end{alignedat} \\
  \text{where }& Z_{\rti}^{i} = \tau _{\rt}^{*} \zeta ^* Z_1 + \tau _\rt^{*} \sqrt{1-(\zeta ^*)^{2}} Z_{1,\infty}^{i}, \; 0\leq i\leq t ;\  \overline{\eta } = \eta(\tau _1^{*}Z_{1,\infty},\cdots ,\tau _{E}^{*}Z_{E,\infty}),
\end{aligned}
    \end{equation*} 
  where 
  \emph{(i)} all of \( \{ Z_1,Z_2,\cdots ,Z_{E},Z_{1,\infty}^{0},\cdots ,Z_{1,\infty}^{t} \}  \) all jointly Gaussian vectors in \( R^{p} \); 
  \emph{(ii)} 
  \( Z_{e}\sim N(0,(\tau ^*_e)^{2}I_p) \) independently for \( e\in[E] \); 
  \emph{(iii)} The Gaussian matrix \( [Z_{1,\infty}^{0}|Z_{1,\infty}^{1}|\cdots |Z_{1,\infty}^{t}] \in\R^{p\times (t+1)} \) is independent of \( \{Z_{1},\cdots ,Z_{E}\} \) and has iid rows drawn from \( N(0,\Sigma ^t_{(V,\rt)}) \), \( \Sigma ^t_{(V,\rt)}\in\R^{(t+1)\times(t+1)} \) from Definition~\ref{def:2nd.induced.evo}.
  
 In particular,  \( \forall\, \) sequence of order-\( k \) pseudo-Lipschitz functions \( \phi:(\R^p)^{(t+E+2)}\to \R \),
  \begin{equation*}
\begin{alignedat}{2}
        &\phi\left(\xi^{t},\xi^{t-1},\cdots ,\xi^{1},\xi^{0},\stac,\beta _1,\beta_2,\cdots ,\beta _E \right) + \smalop&& \\
       &= \E \left[ \phi\left( \xi(Z_{\rti}^{t},\overline{\eta } ) , \xi(Z_{\rti}^{t-1},\overline{\eta } ), \cdots , \xi(Z_{\rti}^{1},\overline{\eta } ),  \xi(Z_{\rti}^{0},\overline{\eta } ), \overline{\eta } ,\beta_1,\beta _2,\cdots ,\beta _E \right) \right].
      \end{alignedat} 
        \end{equation*} 
\end{lemma}

Lemma~\ref{lem:2nd.amp.plus} is the convergence result of the induced AMP iterates from Definition~\ref{def:2nd.induced.glamp}. Since the induced AMP is not a direct specification from any existing AMP theory, we will prove it by showing both sides are well-approximated by Part II of the GLAMP from Defintion~\ref{def:2nd.amp} as \( t_0 \to \infty \). The proof is deferred to Section~\ref{sec:apd.pf.2nd.amp.plus}.

After Lemma~\ref{lem:2nd.amp.plus}, we introduce Lemma~\ref{lem:2nd.cauchy} which establishes a Cauchy property for the induced AMP iterates in Definition~\ref{def:2nd.induced.glamp}. 
\begin{lemma}\label{lem:2nd.cauchy}
  Under Assumptions~\ref{asp:prelim.transfer.model},\ref{asp:stack.for.glamp},\ref{asp:stack.converge},\ref{asp:2nd.for.glamp},\ref{asp:2nd.converge}, and Definition~\ref{def:2nd.induced.glamp},\ref{def:2nd.induced.evo},
  \begin{equation*}
      \limsup_{t\to \infty}\sup_{t\leq i<j}\plim_{p\to \infty}\frac{1}{p}\norm*{v_{\rti}^{i}-v_\rti^{j}}_2^{2} =0,   \quad \limsup_{t\to \infty} \sup_{t\leq i<j} \abs{1 - \rho _\rt^{i,j} }  =0,
      \end{equation*}
  and for any $\epsilon>0$, there exists some \( T = T(\epsilon) \), \( \forall\,j>i\geq T \), \( \lim_{p\to \infty}\PR\left[\frac{1}{p}\norm*{\xi^{i}-\xi ^j}_2^{2} \leq \epsilon \right] = 1 \). 
\end{lemma}

Lemma~\ref{lem:2nd.cauchy} is proved in Section~\ref{sec:apd.pf.2nd.cauchy}. 

With the Cauchy property established, we can prove that the induced AMP iterates asymptoticsally satisfy the KKT condition for \( \se \) in its loss function. 
\begin{lemma}[Sub-Gradient Going to Zero] \label{lem:2nd.subgradient}
  Under Assumptions~\ref{asp:prelim.transfer.model},\ref{asp:stack.for.glamp},\ref{asp:stack.converge},\ref{asp:2nd.for.glamp},\ref{asp:2nd.converge}, and Definitions~\ref{def:2nd.induced.glamp},\ref{def:2nd.induced.evo},  define 
  \begin{equation}\label{eq:def.2nd.subgradient}
    s^{t}_{\rt} = \frac{1}{\theta _{\rt} ^*}  \left[ \Sigma_1^{1/2} v_\rti^{t} + (1-\kappa_1 \overline{\gamma } _{\ro}^{(p)}) \Sigma_1(\beta _1 - \stac) - \Sigma_1(\xi^{t} - \stac)\right].
  \end{equation}
  Then \( s_{\rt}^{t} \in \nabla \mu _{\rt}(\xi ^t ; \stac) \), that is, \( s_{\rt}^{t} \) is a sub-gradient of the penalty \( \mu _{\rt} \) in the loss function \( L_{\rt} \) in Equation~\eqref{eq:def.2nd.los}, either its the joint estimator (\( \mu _{\rt} = \mu _{\joi} \)) or the adaptively weighted estimator (\( \mu _{\rt} = \mu _{\ada} \)). 

   Further, we use \( s^{t}_{\rt} \) to define a member of \( L_{\rt} \)'s subgradient  as \( \nabla L_{\rt}(\xi^{t}) = - \Sigma_1^{1/2}\ix{1}^{\T}(y_1 - \ix{1} \Sigma_1^{1/2}\xi^{t}) +  s_{\rt}^{t} \). We have: \( \forall\,\epsilon>0 \), \( \exists \,T \), s.t. \( \forall \, t\geq T \), \( \lim_{p}\PR[ \frac{1}{p} \norm{\nabla L_{\rt}^{t}(\xi^{t})}_2^{2} \leq \epsilon  ] = 1   \).
\end{lemma}

Lemma~\ref{lem:2nd.subgradient} is proved in Section~\ref{sec:apd.pf.2nd.subgradient}.

\begin{lemma}[Bounding the Sub-Gradient's Magnitude]\label{lem:2nd.bound.on.sc2}
  For any constant \( c\in(0,1) \), define a set \( S^{t}_{\rt}(c) \) in the following way:
  \begin{equation}\label{eq:def.2nd.sc}
    S^{t}_{\rt}(c) = \begin{cases}
        \{ j\in[p] : \abs{(s^{t}_{\rt})_j} \geq \lambda _{\rt}(1-c)  \} 
        \text{ in the case of the joint estimator, }\mu _{\rt}=\mu _{\joi}, \\
        \{ j\in[p] : \abs{(s^{t}_{\rt})_j} \geq \mu(\abs{(\stac)_j})\cdot (1-c)  \} \\ 
        \hspace{2cm}\text{ in the case of the adaptively weighted estimator, }\mu _{\rt}=\mu _{\ada}. 
    \end{cases}
  \end{equation}
  Under Assumptions~\ref{asp:prelim.transfer.model},\ref{asp:stack.for.glamp},\ref{asp:stack.converge},\ref{asp:2nd.for.glamp},\ref{asp:2nd.converge}, and Definitions~\ref{def:2nd.induced.glamp},\ref{def:2nd.induced.evo}, 
  for any \( \epsilon>0 \), there exists \( c = c(\epsilon) \in(0,1) \) invariant to the time step \(  t \in \mathbb{N_+} \), s.t. \(\forall\, t \in \mathbb{N_+} \), \( \lim_{p} \PR\left[\abs{S_{\rt}^{t}(c)}/p \leq \overline{\gamma }  _{\rt} + \epsilon  \right] = 1 \).
\end{lemma}

Lemma~\ref{lem:2nd.bound.on.sc2} is proved in Section~\ref{sec:apd.pf.2nd.bound.on.sc2}. 

\begin{lemma}\label{lem:2nd.final.converge}
  Under Assumptions~\ref{asp:prelim.transfer.model},\ref{asp:stack.for.glamp},\ref{asp:stack.converge},\ref{asp:2nd.for.glamp},\ref{asp:2nd.converge}, and Definitions~\ref{def:2nd.induced.glamp},\ref{def:2nd.induced.evo}, 
  for any $\epsilon>0$,  
  \( \exists \, T \) s.t. \( \forall\,t\geq T \), \( \lim_{p} \PR \left[ \frac{1}{p}\norm{\xi^{t}-\se}_2^{2} < \epsilon\right] = 1  \).  
\end{lemma}

Lemma~\ref{lem:2nd.final.converge} is proved in Section~\ref{sec:apd.pf.2nd.final.converge}. 

With all the preparation so far, Theorem~\ref{thm:2nd} easily follows from Lemma~\ref{lem:2nd.amp.plus} and Lemma~\ref{lem:2nd.final.converge}, and we omit the details here.

\subsection{Proof of Auxiliary Lemmas in Section~\ref{sec:apd.2nd}}
\label{sec:apd.pf.2nd.lems}

\subsubsection{Proof of Lemma~\ref{lem:2nd.bounded}}\label{sec:apd.pf.2nd.bounded}

The proof follows the exact same steps as those of our Lemma~\ref{lem:apd.stack.bounded.lasso}, or Lemma 3.2 of \citet{bayati2011lasso}, so we only point out some key steps this time. 
By Assumption~\ref{asp:prelim.transfer.model} and Lemma~\ref{lem:apd.stack.bounded.lasso}, we know there exists \( M_1>0 \) s.t.
\( \limsup_{p} L_{\rt}(\se)/p \leq \limsup_{p} L_{\rt}(0)/p \leq  M_1  \). 

If  \( \kappa _1<1 \), then we  have \( \liminf_{p} \sigma _{\min }(\ix{1} \Sigma _1^{1/2})\geq (1-\sqrt{\kappa _1})/c_1^{1/2} > 0 \) using Lemma~\ref{lem:apd.sing.val.limit.concent} and Assumption~\ref{asp:prelim.transfer.model}. In this case, we would have proved Lemma~\ref{lem:2nd.bounded}. For the rest of the proof, we focus on the case of \( \kappa _1 \geq 1 \).
  
Let \( P_1^{\|} \) be  the projection matrix onto \( \mathrm{\ker }(\ix{1} \Sigma _1^{1/2}) \) and \( P_{1}^{\bot} \) onto \( \mathrm{\ker }(\ix{1} \Sigma _1^{1/2})^{\bot} \). 
Following the same argument as that of Lemma~\ref{lem:apd.stack.bounded.lasso} and Lemma 3.2 of \citet{bayati2011lasso}, we know from Kashin's theorem that
there exists some constant $B_1$, s.t. almost surely as \( p\to \infty \), 
 \begin{equation*}
\begin{aligned}
 \frac{1}{p}\norm{\se}_2^{2} \leq B_1 \left( \frac{1}{p}\norm{P_1^{\|}\se}_1 \right)^{2} + \frac{1}{p} \norm{P_1^{\bot}\se}_2^{2} \leq 2B_1 \left( \frac{1}{p}\norm{\se}_1  \right)^{2} + \frac{2B_1+1}{p}\norm{P_1^{\bot}\se}_2^{2},
\end{aligned}
 \end{equation*}
 where \( \frac{1}{p}\norm{\se}_1  \) is  bounded by some non-random constant $M_2$ almost surely since \( \frac{1}{p}L_{\rt}(\se) \) is bounded. 
 
 Because \( \mathrm{\ker }(\ix{1} \Sigma _1^{1/2})^{\bot} = \mathrm{range}(\Sigma _1^{1/2}\ix{1}^{\T}) \), \( P_1^{\bot} =\Sigma _1^{1/2}\ix{1}^{\T}(\ix{1}\Sigma _1 \ix{1}^{\T})^{-1}\ix{1}\Sigma _1^{1/2}  \). By Lemma~\ref{lem:apd.sing.val.limit.concent} and Assumption~\ref{asp:prelim.transfer.model}, there exists \( B_2>0 \) s.t. \( \norm{P_1^{\bot}\se}_2 \leq B_2\norm{\ix{1} \Sigma _1^{1/2} P_1^{\bot}\se}_2 = B_2\norm{\ix{1} \Sigma _1^{1/2} \se}_2 \). From here, we can follow the exact same steps as those of Lemma 3.2 of \citet{bayati2011lasso} to prove our Lemma~\ref{lem:2nd.bounded}. 

 \subsubsection{Proof of Lemma~\ref{lem:2nd.glamp}}\label{sec:apd.pf.lem:2nd.glamp}
 We first verify the form of the GLAMP instance defined in Equation~\eqref{eq:def.2nd.glamp} of Lemma~\ref{lem:2nd.glamp} can be cast as a multi-environment GLAMP instance. To do that, we simply augment the GLAMP instance already specified in the proof of Lemma~\ref{lem:apd.stack.glamp}. We still choose \( \{ q_e \}  \) to be integers large enough. We partition the iterates \( V_1^{t} \) and \( R_1^{t} \) into two blocks, whose separation we will emphasize with `$\|$' later.  For the steps \( 0\leq t < t_0 \), everything remains the same as the GLAMP in the proof of Lemma~\ref{lem:apd.stack.glamp}. At the particular step \( t_0 \), we change \( \eta _1^{t_0} \) and \( \Psi _1^{t_0} \) to initialize a new phase as follows:
\begin{equation*}
  \begin{aligned}
    &~\forall\,e\in[E],\ V_e^{t_0} = \left[v_e^{1}|v_e^{2}|\cdots |v_e^{t_0}|0|\cdots \right], \\
    &~\eta _1^{t_0}(V_1^{t_0},\cdots ,V_{E}^{t_0})= \Sigma_1^{1/2}\Big[\beta_1|\eta ^{0}|\eta(v_1^{1},\cdots ,v_{E}^{1})|\eta(v_1^{2},\cdots ,v_{E}^{2})|\cdots |\eta(v_1^{t},\cdots ,v_{E}^{t_0})|0|\cdots\Big\|\\
    &~\hspace{4.35cm} \xi(v_{\rt}^{0}, \eta(v_1^{t_0},\cdots ,v_{E}^{t_0}))|0|\cdots  \Big],\text{ where } v_{\rt}^{0}=\zeta ^* \tau _{\rt}^{*} v_{1}^{t_0}/\tau _{1}^{*}+\sqrt{1-(\zeta ^*)^{2}} \tau _{\rt}^{*} Z_{1}^{\mathrm{init}},\\
    &~\forall\,e > 1,\ \eta _e^{t_0}(V_1^{t_0},\cdots ,V_{E}^{t_0})= \Sigma _e^{1/2}\left[\beta _e|\eta ^{0}|\eta(v_1^{1},\cdots ,v_{E}^{1})|\eta(v_1^{2},\cdots ,v_{E}^{2})|\cdots |\eta(v_1^{t_0},\cdots ,v_{E}^{t_0})|0|\cdots \right].
  \end{aligned}
  \end{equation*}
After the change in \( \eta _1^{t_0} \), the GLAMP iterates for the first environment will have an augmented section starting from \( R_1^{t_0} \):
  \begin{equation*}
    \begin{aligned}
      &~R_1^{t_0} = \left[\ix{1}\Sigma_1^{1/2}\beta_1|r_1^{0}|r_1^{1}|\cdots |r_1^{t}|0|\cdots \|r_{\rt}^{0}|0|\cdots  \right],\\
      &~\forall\,e>1,\ R_e^{t_0} = \left[\ix{e}\Sigma _e^{1/2}\beta _e|r_e^{0}|r_e^{1}|\cdots |r_e^{t}|0|\cdots \right],\\
      &~\Psi_1^{t_0}(R_1^{t_0}) = \kappa_1\left[ y_1 - r_1^{0}|y_1 - r_1^{1}|\cdots |y_1 - r_1^{t}|0|\cdots\middle\|y_1 - r_{\rt}^{0}-\kappa_1 \overline{\gamma } _{\ro}(y_1 - r_1^{t_0})|0|\cdots   \right],\\
      &~\forall\,e>1,\ \Psi _e^{t_0}(R_e^{t_0}) = \kappa _e\left[ y_e - r_e^{0}|y_e - r_e^{1}|\cdots |y_e - r_e^{t}|0|\cdots  \right].
    \end{aligned}
    \end{equation*}

After such initialization of Part II, for any \( t\geq t_0+1 \), the GLAMP proceeds as 
\begin{equation*}
  \begin{aligned}
    &~V_1^{t} = \left[v_1^{1}|v_1^{2}|\cdots |v_1^{t}|0|\cdots\middle\|v_{\rt}^{1}|\cdots |v_{\rt}^{t-t_0}|0|\cdots  \right], \\
    &~\forall\,e>1,\ V_e^{t} = \left[v_e^{1}|v_e^{2}|\cdots |v_e^{t}|0|\cdots \right], \\
    &~\eta _1^{t}(V_1^{t},\cdots ,V_{E}^{t})= \Sigma_1^{1/2}\Big[\beta_1|\eta ^0|\eta(v_1^{1},\cdots ,v_{E}^{1})|\eta(v_1^{2},\cdots ,v_{E}^{2})|\cdots |\eta(v_1^{t},\cdots ,v_{E}^{t})|0|\cdots \|,\\
    &~\hspace{4cm} \xi\left(v_{\rt}^{0}, \eta(v_1^{t_0},\cdots ,v_{E}^{t_0})\right)|\cdots |\xi\left(v_{\rt}^{t-t_0}, \eta(v_1^{t},\cdots ,v_{E}^{t})\right)|0|\cdots  \Big],\\
    &~ \forall\,e>1,\ \eta _e^{t}(V_1^{t},\cdots ,V_{E}^{t})= \Sigma _e^{1/2}\left[\beta _e|\eta ^0|\eta(v_1^{1},\cdots ,v_{E}^{1})|\eta(v_1^{2},\cdots ,v_{E}^{2})|\cdots |\eta(v_1^{t},\cdots ,v_{E}^{t})|0|\cdots \right],\\
    &~R_1^{t} = \left[\ix{1}\Sigma_1^{1/2}\beta_1|r_1^{0}|r_1^{1}|\cdots |r_1^{t}|0|\cdots \|r_{\rt}^{0}|\cdots |r_{\rt}^{t-t_0}|0|\cdots  \right],\\
    &~\forall\,e>1,\ R_e^{t} = \left[\ix{e}\Sigma _e^{1/2}\beta _e|r_e^{0}|r_e^{1}|\cdots |r_e^{t}|0|\cdots \right],\\
    &~ \Psi_1^{t}(R_1^{t}) = \kappa_1\big[ y_1 - r_1^{0}|y_1 - r_1^{1}|\cdots |y_1 - r_1^{t}|0|\cdots \|,\\
    &~\hspace{2.5cm} y_1 - r_{\rt}^{0}-\kappa_1 \overline{\gamma } _{\ro}(y_1 - r_1^{t_0})|\cdots |y_1 - r_{\rt}^{t-t_0}-\kappa_1 \overline{\gamma } _{\ro}(y_1 - r_1^{t})|0|\cdots   \big],\\
    &~ \forall\,e>1,\ \Psi _e^{t}(R_e^{t}) = \kappa _e\left[ y_e - r_e^{0}|y_e - r_e^{1}|\cdots |y_e - r_e^{t}|0|\cdots  \right].
  \end{aligned}
  \end{equation*}
We are left to verify that the state evolution is well-defined, and satisfies the marginal properties in Lemma~\ref{lem:2nd.glamp}. For the asymptotic variance of \( v_{\rt}^{1} \) and its covariance with \( v_{1}^{t_0+1} \), they are 
\begin{equation*}
\begin{aligned}
  \plim_{p}&\frac{1}{p}\norm{v_{\rt}^{1}}_2^{2}\\ 
  =&~ (1-\kappa_1 \overline{\gamma } _{\ro})^{2}\E[W_1^{2}] + \kappa_1 \lim_{p}\frac{1}{p}\E \left[ \norm*{\Sigma_1^{1/2}(\beta_1-\xi(\tau _{\rt}^{*}[\zeta ^* Z^{t_0}_1 + \sqrt{1-(\zeta ^*)^{2}}Z_{1}^{\mathrm{init}}],\eta(Z^{t_0})) - \kappa_1 \overline{\gamma } _{\ro} (\beta_1-\eta(Z^{t_0})) ) }_2^{2}\right] \\
  =&~ (\tau _{\rt}^{*})^{2},\\
  \plim_{p}&\frac{1}{p}(v_{\rt}^{1})^{\T} (v_1^{t_0+1})\\
   =&~
  (1-\kappa_1 \overline{\gamma } _{\ro})\E[W_1^{2}] \\ 
  &~ + \kappa_1 \lim_{p}\frac{1}{p}\E \bigg[(\beta_1-\eta(Z^{t_0}))\Sigma_1(\beta_1-\xi(\tau _{\rt}^{*}[\zeta ^* Z^{t_0}_1 + \sqrt{1-(\zeta ^*)^{2}}Z_{1}^{\mathrm{init}}],\eta(Z^{t_0})) - \kappa_1 \overline{\gamma } _{\ro} (\beta_1-\eta(Z^{t_0})) ) \bigg] \\
  =&~\tau _1^{*}\tau _{\rt}^{*} \zeta ^*,
\end{aligned}
\end{equation*}
both by Assumption~\ref{asp:2nd.converge.fixed.point}. As a result, the coefficients of the Onsager terms are just \( \overline{\gamma } ^{(p)}_{\ro},\overline{\gamma } ^{(p)}_{\ro} \) defined in Proposition~\ref{prop:2nd.onsager} with parameters specified as \( \tau _\rt = \tau _{\rt}^{*}, \zeta = \zeta ^* \).  We can write out the GLAMP iterates as they are in Lemma~\ref{lem:2nd.glamp}. 

We still need to verify the state evolution characterizing the correlation between different time steps. The quantities involved in this are uniquely determined by the multi-environment GLAMP instance, but we choose not to fully write them out for brevity, We only verify those quantities are well defined since they all come as limits as \( p\to \infty \). 

By Assumption~\ref{asp:2nd.converge.cauchy}, this needs us to choose a large enough \( t_0 \), so that the GLAMP instance defined in Equation~\eqref{eq:def.2nd.glamp} has well-defined state evolution. From Lemmas~\ref{lem:apd.stack.glamp} and \ref{lem:apd.stack.cauchy}, we know for any \( \epsilon >0 \), there exists \( T=T(\epsilon) \) s.t. 
\( \sup_{j>i\geq T,\, e\in[E]}\{ 1-\rho ^{i,j}_e \} < \epsilon  \). Hence we can choose \( t_0 > T(\epsilon _H) \) for the \( \epsilon _H \) required by Assumption~\ref{asp:2nd.converge.cauchy}, and then the limits involved in the state evolution are thereby well-defined.  To be more specific, for any \( t_0\leq i<j \), 
\begin{equation*}
\begin{aligned}
  \plim_{p}\frac{1}{p}(v_{\rt}^{i})^\T v_{\rt}^{j} =&~  (1-\kappa_1 \overline{\gamma } _{\ro})^{2}\E[W_1^{2}] + \\
  &~\hspace{-1.5cm}\kappa_1 \lim_{p}\frac{1}{p}\E \Big[ (\beta_1-\xi(\tau _{\rt}^{*} Z_{\rt}^{i-1} ,\eta(Z^{t_0+i-1}_1,\cdots , Z^{t_0+i-1}_E)) - \kappa_1 \overline{\gamma } _{\ro} (\beta_1-\eta(Z^{t_0+i-1}_1,\cdots , Z^{t_0+i-1}_E)) )\\ 
  &~\hspace{-0.25cm}\Sigma_1 (\beta_1-\xi(\tau _{\rt}^{*} Z_{\rt}^{j-1},\eta(Z^{t_0+i-1}_1,\cdots , Z^{t_0+i-1}_E)) - \kappa_1 \overline{\gamma } _{\ro} (\beta_1-\eta(Z^{t_0+i-1}_1,\cdots , Z^{t_0+i-1}_E)) ) \Big],\\
  \text{where } Z_{\rt}^{i-1} =&~ \zeta ^*Z_1^{t_0+i-1} + \sqrt{1-(\zeta ^*)^{2}}(Z_{1}')^{i-1}\\
  Z_{\rt}^{j-1} =&~ \zeta ^*Z_1^{t_0+j-1} + \sqrt{1-(\zeta ^*)^{2}}(Z_{1}')^{j-1}.
\end{aligned}
\end{equation*}
The Gaussian matrix \( [Z_1^{t_0+i-1}|(Z_{1}')^{i-1}|Z_2^{t_0+i-1}|\cdots |Z_{E}^{t_0+i-1}|Z_1^{t_0+j-1}|(Z_{1}')^{j-1}|Z_2^{t_0+j-1}|\cdots |Z_{E}^{t_0+j-1}]\in\R^{p\times 2(E+1)} \) has iid rows from some \( \mathcal{G} \in \mathcal{N} ^{(\epsilon _{H})} \). Hence Assumption~\ref{asp:2nd.converge.cauchy} guarantees the above limit is well-defined. 

\subsubsection{Proof of Lemma~\ref{lem:2nd.amp.plus}}
\label{sec:apd.pf.2nd.amp.plus}

We need preparation before showing Lemma~\ref{lem:2nd.amp.plus}. Firstly, we present Lemma~\ref{lem:2nd.study.H}, which implies the matrix \( \Sigma _{(V,\rt)} \) from Definition~\ref{def:2nd.induced.evo} is well-defined and its elements converge to 1 in the lower-right direction. 

\begin{lemma}\label{lem:2nd.study.H}
  The function \( H_{\rt}:[0,1]\to\R \) defined in Assumption~\ref{asp:2nd.converge.cauchy} is increasing and strictly convex in \( \rho \in[0,1] \). Moreover, \( 0\leq \frac{\mathrm{d}H_{\rt}}{\mathrm{d}\rho}\leq \kappa _{1}\overline{\gamma } _{\rt} [1 - (\zeta ^*)^{2}] \), and \( (\zeta ^*)^{2} \leq  H_{\rt}(0) \leq H_{\rt}(1) =1 \). 
\end{lemma}
\begin{proof}[\textbf{Proof of Lemma~\ref{lem:2nd.study.H}}]
  The first half of Lemma~\ref{lem:2nd.study.H} follows the same argument as that of our Lemma~\ref{lem:apd.stack.study.H}, or the Lemma~C.1 of \citet{bayati2011lasso}, or the Lemma~C.2 of \citet{donoho2016high}. In our case, for any \( p \), we can re-express \( H^{(p)}_{\rt} \) with a function \( \nu:\R^{p}\to \R^{p} \):
  \begin{align*}
    H_{\rt}^{(p)}(\rho) =&~  \frac{(1 - \kappa_1\,\overline{\gamma } _{\ro})^{2} \E[W_1^{2}]}{(\tau _{\rt}^{*})^{2}} + \frac{\kappa_1 }{(\tau _{\rt}^{*})^{2}p}\E
    \left[ \E \left[ \nu\left(\sqrt{\rho} Z_{1}' + \sqrt{1-\rho } Z_{1}''\right)^{\T}\nu\left(\sqrt{\rho} Z_{1}' + \sqrt{1-\rho } Z_{1}'''\right) \middle| Z_1,Z_2,\cdots ,Z_{E}\right] \right],
    \\
\nu(\cdot ) =&~ \Sigma_1^{1/2} 
    \left\{ \left[\beta _1 - \xi\left(\tau _{\rt}^{*}\cdot \left[\zeta ^* Z_{1}/\tau ^*_1 + \sqrt{1-(\zeta ^*) ^2} (\,\cdot\, )\right],\, \overline{\eta } \right) - \kappa_1 \overline{\gamma } _{\ro}(\beta _1 - \overline{\eta } )\right] \right\},
  \end{align*}
where \( \overline{\eta } = \eta(Z_1,\cdots ,Z_E) \), \( Z_1,\cdots ,Z_{E},Z_1',Z_1'',Z_{1}''' \) are independent Gaussian vectors in \( R^{p} \), \( Z_{e}\sim N(0,(\tau _e^{*})^{2}I_p) \), \( Z_{1}',Z_{1}'',Z_{1}'''\iidsim N(0,I_p) \). 

  We introduce a total of \( p \) iid stationary Ornstein-Uhlenbeck processes with covariance function \( \exp\left(-t\right) \), stacked as a vector \( X_t\in\R^{p} \), then 
  \begin{equation*}
    H_\rt^{(p)} = \frac{(1-\kappa_1 \overline{\gamma } _{\ro})^{2}\E[W_1^{2}]}{(\tau _\rt^{*})^{2}} + \frac{\kappa_1}{(\tau _\rt^{*})^{2} p} \E \left[ \E \left[  \nu(X_0)^{\T} \nu(X_t) \middle| Z_1,\cdots ,Z_E  \right] \right]
  \end{equation*}
 with \( t = \log (1/\rho) \). With the same eigen-decomposition trick as that of  Lemma~C.1 of \citet{bayati2011lasso}, or Lemma~C.2 of \citet{donoho2016high}, we can show \( H_\rt^{{(p)}} \) is increasing and strictly convex in \( \rho\in [0,1] \), and thus so is \( H_\rt \). We have proved the first half of Lemma~\ref{lem:2nd.study.H}. 

  \( H_\rt(1) = 1 \) follows from the first equation in Assumption~\ref{asp:2nd.converge.fixed.point}. As a sufficient condition of \( 0\leq \frac{\mathrm{d}H_{\rt}}{\mathrm{d}\rho}\leq \kappa_1\overline{\gamma } _{\rt} [1 - (\zeta ^*)^{2}] \), we prove \( \frac{\mathrm{d}H_\rt}{\mathrm{d} \rho}(1) = \kappa_1 \overline{\gamma } _{\rt}[1-(\zeta ^*)^{2}] \). 

  We first compute the partial derivatives at general \( \rho \in[0,1] \), and then set \( \rho=1 \).  Supposing the limit can be exchanged with the partial derivative, we have
  \begin{equation*}
\begin{aligned}
  \frac{\mathrm{d} H_\rt}{\mathrm{d} \rho}=&~ \frac{-\kappa_1 \sqrt{1 - (\zeta ^*)^{2}}}{\tau _\rt^{*}} \lim_p \frac{1}{p} \E \left[ (\beta _1 - \overline{\xi } _{(2)} - \kappa_1 \overline{\gamma } _1 (\beta _1 - \overline{\eta } ))^{\T} \Sigma_1 J_1(\overline{\xi } _{(1)}) \left( \frac{1}{\sqrt{\rho}}Z_1' - \frac{1}{\sqrt{1-\rho}} Z_1'' \right) \right],
\end{aligned}
  \end{equation*}
  where we have used the exchangeability of \( \overline{\xi } _{(1)} \) and \( \overline{\xi } _{(2)} \) to simplify the form. The notation `$J_1(\overline{\xi } _{(1)})$' means the \( p \)-by-\( p \) Jacobian matrix of \( \overline{\xi } _{(1)} \) w.r.t. its very first argument from Equation~\eqref{eq:def.H.2.cal}. We then invoke Stein's Lemma:
  \begin{equation*}
\begin{aligned}
  \frac{\mathrm{d} H_\rt}{\mathrm{d} \rho}=&~ \frac{-\kappa_1 \sqrt{1 - (\zeta ^*)^{2}}}{\tau _\rt^{*} \sqrt{\rho}} \lim_p \frac{1}{p} \E \left[ (\beta _1 - \overline{\xi } _{(2)} - \kappa_1 \overline{\gamma } _1 (\beta _1 - \overline{\eta } ))^{\T} \Sigma_1 J_1(\overline{\xi } _{(1)}) Z_1' \right]\\ 
  =&~ \kappa_1 [1 - (\zeta ^*)^{2}] \lim_p \frac{1}{p} \E \left[ \trace \left( J_1(\overline{\xi } _{(2)})^{\T} \Sigma_1 J_1(\overline{\xi } _{(1)})  \right)\right].
\end{aligned}
  \end{equation*}
  When \( \rho=1 \), \( \overline{\xi } _{(1)} = \overline{\xi } _{(2)} = \overline{\xi }  \). 
  For either \( S_{\rt} = \Sup(\overline{\xi } -\overline{\eta } ) \) in the case of \( \mu _{\rt}=\mu _{\joi} \), or \( S_{\rt} = \Sup(\overline{\xi }) \) in the case of \( \mu _{\rt}=\mu _{\ada} \), we have \( [J_1(\overline{\xi })]_{S^{c}_{\rt}} = 0 \), so 
  \begin{equation*}
\begin{aligned}
  \frac{\mathrm{d} H_\rt}{\mathrm{d} \rho}=&~ \kappa_1 [1 - (\zeta ^*)^{2}] \lim_p \frac{1}{p} \E \left[ \trace \left( (J_1(\overline{\xi } _{(2)})_{S_{\rt}})^{\T} (\Sigma_1)_{S_\rt,S_\rt} (J_1(\overline{\xi } _{(1)}))_{S_\rt} \right) \right]  \\
  =&~ \kappa_1 [1 - (\zeta ^*)^{2}] \lim_p \frac{1}{p} \E \left[ \trace \left( (\Sigma_1^{1/2})_{\cdot ,S_\rt} [(\Sigma_1)_{S_\rt,S_\rt}]^{-1}  (\Sigma_1^{1/2})_{S_\rt,\cdot } \right) \right] \\
  =&~ \kappa_1 \overline{\gamma } _{\rt} [1 - (\zeta ^*)^{2}], 
\end{aligned}
  \end{equation*}
  where the second and third equalities have both used Equation~\eqref{eq:apd.prelim.2nd.partial.deriv.joi} for the joint estimator, or Equation~\eqref{eq:apd.prelim.2nd.partial.deriv.ada} for the adaptively weighted estimator. 

Finally we prove \( (\zeta ^*)^{2} \leq H_\rt(0)\leq H_{\rt}(1) \). 
The side of \( H_\rt(0)\leq H_{\rt}(1) \) has been argued and follows from  the definition of \( (\tau ^*_{\rt})^{2} \) in Assumption~\ref{asp:2nd.converge.fixed.point}. 
For the rest of the claim,  
recall that we have written \( \overline{\eta } = \eta(Z_1,\cdots ,Z_{E}) \), and \( \overline{\xi } = \xi(\tau _{\rt}^{*}(\zeta ^* Z_{1}/\tau _1^{*} + \sqrt{1-(\zeta ^*)^{2}}Z_{1}'),\, \overline{\eta } ) \), where \( Z_1,\cdots ,Z_{E},Z_{1}' \) are independent, \( Z_1'\sim N(0,I_p) \). 
Define \( \breve{\xi} = \E[\overline{\xi } |Z_1,Z_2,\cdots ,Z_{E}] \). In other words, the integral is taken only w.r.t. \( Z_{1}' \). After some algebra, we have
\begin{equation*}
\begin{aligned}
  (\tau _{1}^{*}\tau _{\rt}^{*})^{2}H_{\rt}(0) = &~ 
  (1-\kappa_1\overline{\gamma } _{\ro})^{2}\E[W_1^{2}]^{2} + \kappa_1 \E[W_1^{2}] \lim_{p}\frac{1}{p} \E\left[\norm*{\Sigma_1^{1/2} \left(\beta _1 - \breve{\xi}- \kappa_1 \overline{\gamma } _{\ro}(\beta _{1} - \overline{\eta } )\right) }_2^{2}\right] \\
  &~ + \kappa_1^{2} \lim_{p}\frac{1}{p} \E\left[\norm*{\Sigma_1^{1/2} \left(\beta_1 - \breve{\xi}- \kappa_1 \overline{\gamma } _{\ro}(\beta _{1} - \overline{\eta } )\right) }_2^{2}\right] \lim_{p}\frac{1}{p} \E\left[\norm*{\Sigma_1^{1/2} \left(\beta _{1} - \overline{\eta }\right) }_2^{2}\right] \\
  &~ + \kappa_1 (1-\kappa_1\overline{\gamma } _{\ro})^{2} \E[W_1^{2}] \lim_{p}\frac{1}{p} \E\left[\norm*{\Sigma_1^{1/2} \left(\beta _{1} - \overline{\eta }\right) }_2^{2}\right],  
\end{aligned}
\end{equation*}
\begin{equation*}
  \begin{aligned}
    (\tau _{1}^{*}\tau _{\rt}^{*})^{2}H_{\rt}(1) = &~ 
    (1-\kappa_1\overline{\gamma } _{\ro})^{2}\E[W_1^{2}]^{2} + \kappa_1 \E[W_1^{2}] \lim_{p}\frac{1}{p} \E\left[\norm*{\Sigma_1^{1/2} \left(\beta_1 - \overline{\xi}- \kappa_1 \overline{\gamma } _{\ro}(\beta _{E} - \overline{\eta } )\right) }_2^{2}\right] \\
    &~ + \kappa_1^{2} \lim_{p}\frac{1}{p} \E\left[\norm*{\Sigma_1^{1/2} \left(\beta_1- \overline{\xi}- \kappa_1 \overline{\gamma } _{\ro}(\beta _{1} - \overline{\eta } )\right) }_2^{2}\right] \lim_{p}\frac{1}{p} \E\left[\norm*{\Sigma_1^{1/2} \left(\beta _{1} - \overline{\eta }\right) }_2^{2}\right] \\
    &~ + \kappa_1 (1-\kappa_1\overline{\gamma } _{\ro})^{2} \E[W_1^{2}] \lim_{p}\frac{1}{p} \E\left[\norm*{\Sigma_1^{1/2} \left(\beta _{1} - \overline{\eta }\right) }_2^{2}\right],
  \end{aligned}
  \end{equation*}
  \begin{equation*}
\begin{aligned}
  (\tau _1^{*}\tau_\rt^* \zeta ^*)^{2} =&~ (1-\kappa_1\overline{\gamma } _{\ro})^{2}\E[W_1^{2}]^{2} + \kappa_1^{2} \lim_{p}  \frac{1}{p^{2}}\E \left[ (\beta _{1} - \overline{\eta } ) \Sigma_1 \left[\beta _{1} - \breve{\xi }  - \kappa_1 \overline{\gamma } _{\ro}(\beta _{1}-\overline{\eta } )\right] \right]^{2} \\
  &~ + 2 \kappa_1 (1-\kappa_1 \overline{\gamma } _{\ro}) \E[W_1^{2}] \lim_{p}  \frac{1}{p}\E \left[ (\beta _{1} - \overline{\eta } ) \Sigma_1 \left[\beta _{1} - \breve{\xi } - \kappa_1 \overline{\gamma } _{\ro}(\beta _{1}-\overline{\eta } )\right] \right]^{2}.
\end{aligned}
\end{equation*}
By the law of iterated expectations, \( H_{\rt}(1) \geq H_{\rt}(0) \). By Cauchy's inequality and the AM-GM inequality, we also see the following relationships:
\begin{equation*}
\begin{aligned}
  H_{\rt}(0) \geq (\zeta ^*)^{2};& \quad H_{\rt}(0) = (\zeta ^*)^{2} \iff \lim_{p} \frac{1}{p}\E\left[\norm{\breve\xi - \overline{\eta } }_2^{2}\right] = 0, \\
  H_{\rt}(1) \geq (\zeta ^*)^{2};& \quad H_{\rt}(1) = (\zeta ^*)^{2} \iff \lim_{p} \frac{1}{p}\E\left[\norm{\overline\xi - \overline{\eta } }_2^{2}\right] = 0.
\end{aligned}
\end{equation*}
This completes the proof of Lemma~\ref{lem:2nd.study.H}. 
\end{proof}

We need yet another lemma before proving Lemma~\ref{lem:2nd.amp.plus}. Lemma~\ref{lem:2nd.amp.prep.1st.stage} establishes the convergence of \( (v_{1}^{t},r_{1}^{t},\eta(v_{1}^{t},\cdots ,v_{E}^{t})) \) of Part I of the GLAMP from Definition~\ref{def:2nd.amp} towards \( (v_{1}^{\infty},r_{1}^{\infty},\stac) \). 

\begin{lemma}\label{lem:2nd.amp.prep.1st.stage}
In addition to \( v_{1}^{\infty} \) and \( r_{1}^{\infty} \), define for any \( e\in[E] \)
\begin{equation*}
  v_{e}^{\infty} = \frac{\ix{e}^{\T}(y_{e} - \ix{e}\Sigma _{e}^{1/2}\stac)}{1-\kappa _e \overline{\delta }_e } + \Sigma _e^{1/2}(\stac - \beta _{e}),\quad r_{e}^{\infty} = y_{e} - \frac{y_{e}- \ix{e}\Sigma _{e}^{1/2}\stac}{1-\kappa _e \overline{\delta }_ {e}}.
\end{equation*}
\( \forall\,\epsilon > 0 \), \( \exists \, T_0 \), s.t. \( \forall\,t_0\geq T_0,\,t \in \mathbb{N},\,e\in[E] \), \( \lim_{p\to \infty}\PR[ \frac{1}{p}\norm{\eta(v_1^{t+t_0},\cdots ,v_{E}^{t+t_0})-\stac}_2^{2} < \epsilon ] = 1 \), and
\begin{equation*}
  \begin{alignedat}{2}
    &\lim_{p\to \infty}\PR\bigg[ \frac{1}{p}\norm*{v_e^{t+t_0}-v_e^{\infty}}_2^{2} < \epsilon \bigg] = 1, \  \lim_{p\to \infty}\PR\bigg[ \frac{1}{n_e}\norm*{r_e^{t+t_0}-r_e^{\infty}}_2^{2} < \epsilon \bigg] = 1.
  \end{alignedat}
    \end{equation*}
\end{lemma}
\begin{proof}[\textbf{Proof of Lemma~\ref{lem:2nd.amp.prep.1st.stage}}]
  \newcommand{\brevop}{\breve{o}_{\PR}(\sqrt(p))}
The claim about \( \eta(v_1^{t+t_0},\cdots ,v_{E}^{t+t_0}) \) has been proved in Lemma~\ref{lem:stack.final.converge} and is just revisted here. We focus on the rest of Lemma~\ref{lem:2nd.amp.prep.1st.stage}.

Recall from Lemma~\ref{lem:apd.stack.cauchy} that 
for any \( \epsilon>0 \), there exists \( T_0 \) s.t. \( \forall\,j>i\geq T_0,\,\forall\,e\in[E] \), \( \lim_{p}\PR[\frac{1}{p}\norm{v_{e}^{i}-v_{e}^{j}}_2^{2} < \epsilon] = 1 \) and \( \lim_{p}\PR[\frac{1}{n_e}\norm{r_{e}^{i}-r_{e}^{j}}_2^{2} < \epsilon] = 1 \). 
By re-writing the GLAMP iterates from Part I of Definition~\ref{def:2nd.amp}, we know for \( t_0\geq T_0 \) and \( t \in \mathbb{N_+} \),
\begin{equation*}
\begin{aligned}
  y_e - r_e^{t+t_0} =&~ \frac{y_e - \ix{e} \Sigma ^{1/2}_{e} \hat{\beta } ^{t+t_0}}{1 - \kappa_e \overline{\delta } ^{(p)}_{e}} + \frac{\kappa _{e} \overline{\delta } ^{(p)}_{e}}{1 - \kappa_e \overline{\delta } ^{(p)}_{e}}(r_{e}^{t+t_0}-r_{e}^{t+t_0-1}),\\
  v_{e}^{t+t_0+1} =&~ \frac{\ix{e}^{\T}(y_e - \ix{e} \Sigma ^{1/2}_{e} \eta^{t+t_0})}{1 - \kappa_e \overline{\delta } ^{(p)}_{e}} + \frac{\kappa _{e} \overline{\delta } ^{(p)}_{e}}{1 - \kappa_e \overline{\delta } ^{(p)}_{e}}\ix{e}^{\T}(r_{e}^{t+t_0}-r_{e}^{t+t_0-1}) + \Sigma _{e}^{1/2}(\eta^{t+t_0}-\beta _e), \\
  \text{where } \eta ^{t+t_0}=&~\eta(v_{1}^{t+t_0},\cdots ,v_{E}^{t+t_0})
\end{aligned}
\end{equation*}

Recall that \( \overline{\delta } ^{(p)}_{e} \to \overline{\delta } _{e} \), \( \lim_{p} \sigma _{\max }(\ix{e}) \leq 1 + \sqrt{\kappa _e} \) a.s. and \( \lambda _{\max }(\Sigma _{e}) \leq c_1 \). As a result, for any \( \epsilon>0 \), there exists \( T_0 \) large enough s.t. \( \forall\,t_0\geq T_0 \) and \( t \in \mathbb{N} \), 
\begin{equation*}
  \lim_{p}\PR \left[ \frac{1}{n_e}\norm{r_{e}^{t+t_0} - r_{e}^{\infty}}_2^{2} < \epsilon \right] = 1,\quad \lim_{p}\PR \left[ \frac{1}{p}\norm{v_{e}^{t+t_0} - v_{e}^{\infty}}_2^{2} < \epsilon \right] = 1.
\end{equation*}
Thus we have proved Lemma~\ref{lem:2nd.amp.prep.1st.stage}. 
\end{proof}

\begin{lemma}\label{lem:2nd.stronger.wasser}
  Fixing \( t \in \mathbb{N_+} \), for any sequence of order-$k$ pseudo-Lipschitz functions \( \phi:(\R^{p})^{(t+1)\times(E+1)+E}\to \R \), \( k\geq 1 \), define 
  \begin{equation*}
    A = \phi\left( Z_{\rt}^{t},Z_{\rt}^{t-1},\cdots ,Z_{\rt}^{1},Z_{\rt}^{0}, Z^{t+t_0}_1,\cdots ,Z^{t+t_0}_E,Z^{t+t_0-1}_1,\cdots ,Z^{t+t_0-1}_E,\cdots ,Z^{t_0}_1,\cdots ,Z^{t_0}_E,\beta_1,\beta _2,\cdots ,\beta _E \right),
  \end{equation*}
  where the Gaussian variables on the right hand side are exactly from the limit in Lemma~\ref{lem:2nd.glamp}. We define another quantity \( B \):
  \begin{equation*}
    B = \phi\left( Z_{\rti}^{t},Z_{\rti}^{t-1},\cdots ,Z_{\rti}^{1},Z_{\rti}^{0}, Z_{1},\cdots ,Z_{E},Z_{1},\cdots ,Z_{E},\cdots ,Z_{1},\cdots ,Z_{E},\beta_1,\beta _2,\cdots ,\beta _E \right),
  \end{equation*}
  where 
  \emph{(i)} \( \forall\,0\leq i\leq t \), \( Z_{\rti}^{i} = \tau_{\rt}^{*} \zeta ^* (Z_{1}/\tau _{E}^{*}) + \tau _{\rt}^{*} \sqrt{1-(\zeta ^*)} Z_{1,\infty}^{i} \); 
  \emph{(ii)} \( Z_{e}\sim N(0,(\tau _e^{*})^{2}I_p) \) for \( e\in[E] \), \(  \) 
  \emph{(iii)}  The Gaussian matrix  \( [Z_{1,\infty}^{0}|Z_{E,\infty}^{1}|\cdots |Z_{E,\infty}^{t}] \in\R^{p\times (t+1)}\) is independent of \( \{ Z_1,Z_2,\cdots ,Z_{E} \}  \), and has iid rows drawn from \( N(0, \Sigma _{V,\rt}^{t})  \). 
  Then we have 
  \begin{equation*}
    \lim_{t_0\to \infty} \sup_{p\geq 1} \abs*{E[A]-E[B]} = 0. 
  \end{equation*}
\end{lemma}
\begin{proof}[\textbf{Proof of Lemma~\ref{lem:2nd.stronger.wasser}}]
  If we write the Gaussian variables used by $A$ and $B$ into matrices \( Z_A,Z_B\in\R^{(p\times(t+1)(E+1))} \):
\begin{equation*}
\begin{aligned}
  Z_A =&~ [Z_{\rt}^{t}|Z_{\rt}^{t-1}|\cdots |Z_{\rt}^{1}|Z_{\rt}^{0}| Z^{t+t_0}_1|\cdots |Z^{t+t_0}_E|Z^{t+t_0-1}_1|\cdots |Z^{t+t_0-1}_E|\cdots |Z^{t_0}_1|\cdots |Z^{t_0}_E],\\
  Z_{B} =&~ [Z_{\rti}^{t}|Z_{\rti}^{t-1}|\cdots |Z_{\rti}^{1}|Z_{\rti}^{0}| Z_{1}|\cdots |Z_{E}|Z_{1}|\cdots |Z_{E}|\cdots |Z_{1}|\cdots |Z_{E}].
\end{aligned}
\end{equation*}
Then by the definitions of the Gaussian variables from Lemma~\ref{lem:2nd.glamp} and Lemma~\ref{lem:2nd.stronger.wasser}, 
we know both \( Z_A \) and \( Z_{B} \) have iid rows drawn from some \( (t+1)(E+1) \)-dimensional Gaussian distribution, with mean zero and covariance matrices \( \Sigma _{A} \) and \( \Sigma _{B} \). Both covariance matrices are \( (t+1)(E+1) \)-dimensional and do not depend on \( p \). 

We first show \( \lim_{t_0\to \infty}\norm{\Sigma _{A}-\Sigma _B}_2 = 0 \). 

By definition, \( \Sigma _{B}\in\R^{(t+1)(E+1)\times(t+1)(E+1)} \) do not depend on \( t_0 \). On the other hand for \( \Sigma _{A} \), by Lemma~\ref{lem:apd.stack.cauchy}, \( \lim_{t_0\to \infty}\sup_{j>i\geq 0, e\in[E]}\rho ^{t_0+i,t_0+j}_{e} = 1  \). 
Hence for any fixed \( e\in[E] \), if we look at the marginal distribution of each row of \( [Z_{e}^{t+t_0}|Z_{e}^{t+t_0-1}|\cdots |Z_{e}^{t_0}] \in\R^{p\times(t+1)} \), all the \( (t+1) \) elements become perfectly correlated as \( t_0\to \infty \); in other words, the elements in \( \Sigma _{A} \) corresponding to their covariances converge to \( (\tau _e^{*})^{2} \).  
Among different environments, we know for different \( e\in[E] \),  \( \{ [Z_{e}^{t+t_0}|Z_{e}^{t+t_0-1}|\cdots |Z_{e}^{t_0}] :e\in[E]\} \) are independent Gaussian matrices. As a result, we have proved that for the lower-right sub-blocks: 
\begin{equation*}
  \lim_{t_0 \to \infty} (\Sigma _{A})_{(t+2) : (t+1)(E+1),\,(t+2) : (t+1)(E+1)} = (\Sigma _{B})_{(t+2) : (t+1)(E+1),\,(t+2) : (t+1)(E+1)}.
\end{equation*}

For the first \( (t+1) \) rows of \( Z_A \), we know from Lemma~\ref{lem:2nd.glamp} that \( Z_{\rt}^{i} = \tau _{\rt}^{*} \zeta ^* Z_{1}^{t_0+i} / \tau _{1}^{*} + \tau _{\rt}^{*} \sqrt{1-(\zeta ^*)^{2}} (Z_{1}')^{i} \) for some \( (Z_{1}')^{i}\independent Z_{1}^{t_0+i} \),  \( \forall\,0\leq i\leq t \). For convenience, we apply a linear transformation to \( Z_{A} \), take an excerpt at two time steps \( 0\leq i<j\leq t \), and consider the Gaussian matrix \( [\frac{1}{\tau _1^{*}}Z_{1}^{i+t_0}|(Z_{1}')^{i}|\frac{1}{\tau _2^{*}}Z_{2}^{i+t_0}|\cdots |\frac{1}{\tau _E^{*}}Z_{E}^{i+t_0}|\frac{1}{\tau _1^{*}}Z_{1}^{j+t_0}|(Z_{1}')^{j}|\frac{1}{\tau _2^{*}}Z_{2}^{j+t_0}|\cdots |\frac{1}{\tau _E^{*}}Z_{E}^{j+t_0}]\in\R^{p\times 2(E+1)} \). Each row of it follows a $2(E+1)$-dimensional Gaussian distribution, which we denote as \( \mathcal{G} ^{i,j} \). 

We focus on the correlation between the two elements in each iid row of \( [(Z_{1}')^{i}|(Z_{1}')^{j}]\in\R^{p\times 2} \), which we denote as \( \varrho ^{i,j} \).  We want to show \( \lim_{t_0 to \infty} \varrho ^{i,j} = (\Sigma _{(V,\rt)})_{i,j}\), and then the upper-left block of \( \Sigma _{A} \) will also converge to that of \( \Sigma _{B} \). This is because for columns of \( B \), \( Z_{\rt}^{i} = \tau _{\rt}^{*}\zeta ^* Z_{1}/\tau _1^{*} + \tau _{\rt}^{*} \sqrt{1-(\zeta ^*)^{2}} Z^{i}_{1,\infty} \) and \( [Z^{0}_{1,\infty}|Z^{1}_{1,\infty}|\cdots |Z^{t}_{1,\infty}] \) has iid rows from \( N(0,\Sigma _{(V,\rt)}^{t}) \) by construction. 

When \( i=0 \), \( (Z_{1}')^{0} = Z^{\mathrm{init}}_{1} \) is independent from all other variables, so \( \varrho ^{i,j} = 0 = (\Sigma _{(V,\rt)})_{0,j} \). 

When \( i\geq 1 \), \( \varrho ^{i,j} = H_{\rt}^{\epsilon _H}(\Gscr^{i-1.j-1}) \). This is because of the GLAMP convergence in Lemma~\ref{lem:2nd.glamp}. Then for any \( j \), we can show by induction that \( \lim_{t_0}\varrho ^{i,j} = \lim_{t_0\to \infty} H_{\rt}^{\epsilon _{H}}(\Gscr^{i-1.j-1}) = H_{\rt}^{\epsilon _H}(\lim_{t_0}\varrho ^{i-1,j-1}) \), and \( \lim_{t_0}\varrho ^{i,j} = (\Sigma _{(V,\rt)})_{i,j} \). 

As a result, we have proved that for the upper-left sub-blocks:
\begin{equation*}
  \lim_{t_0 \to \infty} (\Sigma _{A})_{[t+1],[t+1]} = (\Sigma _{B})_{[t+1],[t+1]}.
\end{equation*}

For the rest of the elements in \( \Sigma _{A} \), note that the columns of \( [Z_{e}^{t+t_0}|Z_{e}^{t+t_0-1}|\cdots |Z_{e}^{t_0}] \in\R^{p\times(t+1)} \) becomes perfectly correlated as \( t_0 \to \infty \), so those elements also converge to the corresponding elements in \( \Sigma _{B} \)

Thus we have shown \( \lim_{t_0}\norm{\Sigma _{A}-\Sigma _B}_2 = 0 \). From here, following the same argument as Equation (156)-(166) of the proof of Lemma 10, \citet{berthier2020state}, we arrive at 
\begin{equation*}
  \abs{ \E[A]-\E[B]}^{2} \leq \mathit{Constant}\cdot \trace(\Sigma _A + \Sigma _B - 2(\Sigma _B^{1/2}\Sigma _A \Sigma _B^{1/2})^{1/2}) = o(\norm{\Sigma _{A}-\Sigma _B}_2),
\end{equation*}
where the \( \mathit{Constant} \) does not depend on \( p,t_0, \Sigma_A,\Sigma _B \). 
This completes the proof of Lemma~\ref{lem:2nd.stronger.wasser}. 
\end{proof}

Now we are ready to show Lemma~\ref{lem:2nd.amp.plus}.
\begin{proof}[\textbf{Proof of Lemma~\ref{lem:2nd.amp.plus}}]

We first prove by induction on \( t \) that \( \forall\,\epsilon>0 \), \( \exists\,T_0  \) s.t. \( \forall\,t_0\geq T_0 \),
\begin{equation}\label{eq:2nd.2amp.equiv}
  \lim_{p}\PR\left[\frac{1}{p}\norm{v_{\rt}^{t} - v_{\rti}^{t}}_2^{2}\leq \epsilon\right]=1. 
\end{equation}
When \( t=0 \), \( v_{\rt}^{0}= \zeta ^* \cdot \tau _{\rt}^{*} v_1^{t_0}/\tau _{1}^{*}+\sqrt{1-(\zeta ^*)^{2}}\cdot \tau _{\rt}^{*} Z_{\rt}^{\mathrm{init}}\), \( v_{\rti}^{0} =  \zeta ^* \cdot \tau _{\rt}^{*} v_1^{\infty}/\tau _{1}^{*}+\sqrt{1-(\zeta ^*)^{2}}\cdot \tau _{\rt}^{*} Z_{\rt}^{\mathrm{init}} \). Then Equation~\eqref{eq:2nd.2amp.equiv} directly follows from Lemma~\ref{lem:2nd.amp.prep.1st.stage}. Now suppose we have proved Equation~\eqref{eq:2nd.2amp.equiv} for all \( 0\leq i\leq t-1 \). For \( t \),
\begin{equation*}
\begin{aligned}
  v_{\rt}^{t} =&~ \ix{1}^{\T}(y_{1}-r_{\rt}^{t-1}-\kappa _1 \overline{\gamma } _{\ro}(y_{1} - r_{1}^{t+t_0-1})) - (1-\kappa _{1}\overline{\gamma } ^{(p)}_{\ro})\Sigma _{1}^{1/2}\beta _{1} - \kappa _{1} \overline{\gamma } _{\ro}^{(p)} \Sigma _{1}^{1/2} \eta^{t+t_0-1}\\
  &~ + \Sigma _{1}^{1/2} \xi(v_{\rt}^{t-1},\eta^{t+t_0-1}) , \\
  v_{\rti}^{t} =&~ \ix{1}^{\T}(y_{1}-r_{\rti}^{t-1}-\kappa _1 \overline{\gamma } _{\ro}(y_{1} - r_{1}^{\infty})) - (1-\kappa _{1}\overline{\gamma } ^{(p)}_{\ro})\Sigma _{1}^{1/2}\beta _{1} - \kappa _{1} \overline{\gamma } _{\ro}^{(p)} \Sigma _{1}^{1/2} \stac \\
  &~ + \Sigma _{1}^{1/2} \xi(v_{\rti}^{t-1},\stac), 
\end{aligned}
\end{equation*}
and Equation~\eqref{eq:2nd.2amp.equiv} follows from Lemma~\ref{lem:2nd.amp.prep.1st.stage} and the induction hypothesis. 

We know from Lemma~\ref{lem:2nd.glamp} that 
  \begin{equation}\label{eq:apd.2nd.amp.plus}
    \begin{alignedat}{2}
      &\phi\left(v_{\rt}^{t},v_{\rt}^{t-1},\cdots ,v_{\rt}^{1},v_{\rt}^{0},v_1^{t+t_0},\eta ^{t+t_0},\beta _1,\beta_2,\cdots ,\beta _E \right) + \smalop&& \\
     &= \E \left[ \phi\left( Z_{\rt}^{t},Z_{\rt}^{t-1},\cdots ,Z_{\rt}^{1},Z_{\rt}^{0}, Z^{t+t_0}_1,\eta(Z_{1}^{t+t_0},\cdots ,Z_{E}^{t+t_0}),\beta_1,\beta _2,\cdots ,\beta _E \right) \right].
    \end{alignedat}
    \end{equation}
We wish to take \( t_0\to \infty \) on both sides of Equation~\eqref{eq:apd.2nd.amp.plus} to prove Lemma~\ref{lem:2nd.amp.plus}. 

    For the LHS of Equation~\eqref{eq:apd.2nd.amp.plus}, by Equation~\eqref{eq:2nd.2amp.equiv} and Lemma~\ref{lem:2nd.amp.prep.1st.stage}, \( \forall\,\epsilon>0 \), \( \exists \,T_0 \) s.t. \( \forall\, t_0\geq T_0 \) and \( t \in \mathbb{N_+} \), 
\begin{equation*}
  \lim_{p}\PR\left[\abs*{\phi\left(v_{\rt}^{t},\cdots ,v_{\rt}^{0},v_1^{t+t_0},\eta ^{t+t_0},\beta _1,\beta_2,\cdots ,\beta _E \right) - \phi\left(v_{\rti}^{t},\cdots ,v_{\rti}^{0},v_1^{\infty},\stac,\beta _1,\beta_2,\cdots ,\beta _E \right)} \leq \epsilon \right]=1,
\end{equation*}
where the Gaussian variables are as they are defined in Lemma~\ref{lem:2nd.glamp}. 

For the RHS of Equation~\eqref{eq:apd.2nd.amp.plus}, by Lemma~\ref{lem:2nd.stronger.wasser}, \( \forall\,\epsilon>0 \), \( \exists \,T_0 \) s.t. \( \forall\, t_0\geq T_0 \) and \( t \in \mathbb{N_+} \), 
\begin{equation*}
\begin{aligned}
  &~\sup_{p} \Big|\E \left[ \phi\left( Z_{\rt}^{t},Z_{\rt}^{t-1},\cdots ,Z_{\rt}^{1},Z_{\rt}^{0}, Z^{t+t_0}_1,\eta(Z_{1}^{t+t_0},\cdots ,Z_{1}^{t+t_0}),\beta_1,\beta _2,\cdots ,\beta _E \right) \right] \\
  &~ \hspace{3cm} - \E \left[ \phi\left( Z_{\rti}^{t} , Z_{\rti}^{t-1}, \cdots , Z_{\rti}^{1}, Z_{\rti}^{0},  Z_{1},\overline{\eta } ,\beta_1,\beta _2,\cdots ,\beta _E \right) \right] \Big|< \epsilon, 
\end{aligned} 
\end{equation*}
where the Gaussian variables are as they are defined in the term `$B$' in Lemma~\ref{lem:2nd.stronger.wasser}. 

Putting the analysis of both sides together, we know 
\begin{equation*}
  \begin{alignedat}{2}
    &\phi\left(v_{\rti}^{t},v_{\rti}^{t-1},\cdots ,v_{\rti}^{1},v_{\rti}^{0},v_1^{\infty},\stac,\beta _1,\beta_2,\cdots ,\beta _E \right) + \smalop&& \\
   &= \E \left[\phi\left( Z_{\rti}^{t} , Z_{\rti}^{t-1}, \cdots , Z_{\rti}^{1}, Z_{\rti}^{0}, Z_{1},\overline{\eta } ,\beta_1,\beta _2,\cdots ,\beta _E \right) \right], 
  \end{alignedat}
  \end{equation*}
  which completes the proof of Lemma~\ref{lem:2nd.amp.plus}. 
\end{proof}

\subsubsection{Proof of Lemma~\ref{lem:2nd.cauchy}}\label{sec:apd.pf.2nd.cauchy}

Let \( \widetilde{H} _{\rt}(\cdot ) = (H_{\rt}(\cdot )-(\zeta ^*)^{2}) / (1-(\zeta ^*)^{2}) \). By Lemma~\ref{lem:2nd.study.H}, \( \widetilde{H} _{\rt}(\cdot ) \) is increasing, strictly convex, \( \widetilde{H} _{\rt}(1)=1 \), \( \widetilde{H} _{\rt}(0)\geq 0 \). Also, \( 0\leq \frac{\mathrm{d}\widetilde{H} _{\rt}(\rho)}{\mathrm{d}\rho} \leq \kappa _1 \overline{\gamma } _{\rt} < 1 \) by the last equation in Assumption~\ref{asp:2nd.converge.fixed.point}. By Definition~\ref{def:2nd.induced.evo}, 
\begin{equation*}
  \abs{1 - \rho ^{i,j}_{\rt}} \leq (\kappa_1 \overline{\gamma } _{\rt}) \abs{1 -  \rho ^{i-1,j-1}_{\rt}} \leq (\kappa_1 \overline{\gamma } _{\rt})^{i}.  
\end{equation*}
Thus \( \limsup_{T\to \infty}\sup_{T\leq i<j}\abs{1 - \rho _{\rt}^{i,j}} = 0 \). Also, \( \plim_{p}\frac{1}{p}\norm{v_{\rti}^{i}-v_{\rti}^{j}}_2^{2} = 2(\tau _\rt^{*})^{2} (1-\rho ^{i,j}_{\rt}) \), and \( \xi^{t} = \eta(v^{t}_{\rti},\stac) \) is Lipschitz in \( v^{t}_{\rti} \). Thus we have shown Lemma~\ref{lem:2nd.cauchy}.

\subsubsection{Proof of Lemma~\ref{lem:2nd.subgradient}}
\label{sec:apd.pf.2nd.subgradient}
Either \( \mu _{\rt}=\mu _{\joi} \) or \( \mu _{\ada} \), 
we can use the definition of 
\( \xi \) in Eqn~\eqref{eq:def.2nd.xi} to write out the KKT condition of \( \xi^{t} = \xi(v_{\rt}^{t},\stac) \), 
\begin{equation*}
  \Sigma_1(\xi^{t} - \stac)  - \Sigma_1^{1/2} v_{\rti}^{t} - (1 - \kappa_1 \overline{\gamma } _{\ro}) \Sigma_1 (\beta _{1}-\stac) + \theta _{\rt}^{*} \nabla \mu_{\rt}(\xi ^t ; \stac) \ni 0.
\end{equation*}
Set the LHS to exactly zero, and let \( s^{t}_{\rt} \) be the member of sub-gradient that solves the equation. Then we get 
\begin{equation*}
  s_{\rt}^{t} = \frac{1}{\theta _{\rt}^{*}} \left[ \Sigma_1^{1/2} v_{\rti}^{t} + (1 - \kappa_1 \overline{\gamma } _{\ro}) \Sigma_1 (\beta _{1}-\stac) - \Sigma_1(\xi^{t} - \stac) \right] \in \nabla\mu _{\rt}(\xi ^t ; \stac).
\end{equation*}

To show the claim about \( \nabla L_{\rt}(\xi^{t} ) \), we again use the \( \specop \) notation defined in the proof of Lemma~\ref{lem:stacked.lasso.subgradient} in Section~\ref{sec:apd.pf.stacked.lasso.subgradient}. To recall its definition, some sequence \( X^{(p),t} = \specop  \) iff  \( \forall\, \epsilon>0 \), \( \exists\, T \) s.t. \( \forall\,t\geq T \), \( \lim_{p}\PR\left[ \frac{1}{\sqrt{p}}\norm{X^{{(p),t}}}_2 < \epsilon \right] = 1 \). 

By Lemma~\ref{lem:2nd.cauchy}, \( v_{\rti}^{t+1}-v_{\rti}^{t} = \specop \), \( \xi^{t+1} - \xi^{t} = \specop  \). From how the induced AMP iterates are defined in Definition~\ref{def:2nd.induced.glamp}, we know
\begin{equation*}
\begin{aligned}
  &y_1 - r_{\rti}^{t}-\kappa_1 \overline{\gamma } ^{(p)}_{\ro}(y_1 - r^{\infty}_{1}) = y_1 - \ix{1} \Sigma_1 ^{1/2} \xi^t + \kappa_1 \overline{\gamma } _{\rt}^{(p)} (y_1 - r_{\rti}^{t-1}-\kappa_1 \overline{\gamma } _{\ro}(y_1 - r_1^{\infty}))\\
  \implies~&\ix{1}^{\T}\left[ y_1 - r_{\rti}^{t}-\kappa_1 \overline{\gamma } ^{(p)}_{\ro}(y_1 - r^{\infty}_{1}) \right] \\
  &~=  \ix{1}^{\T}(y_1 - \ix{1} \Sigma_1 ^{1/2} \xi^t ) + \kappa_1 \overline{\gamma } ^{(p)}_{\rt} \ix{1}^{\T} \left[ y_1 - r_{\rti}^{t-1}-\kappa_1 \overline{\gamma } _{\ro}(y_1 - r_1^{\infty}) \right] \\
  \implies~&
  v_{\rti}^{t+1} + (1-\kappa_1\overline{\gamma } ^{(p)}_{\ro})\Sigma_1^{1/2}\beta _1 + \kappa_1 \overline{\gamma } _{\ro}^{(p)} \Sigma_1^{1/2} \stac - \Sigma_1^{1/2} \xi^t + \specop \\
  &~ = \ix{1}^{\T}(y_1 - \ix{1} \Sigma_1 ^{1/2} \xi^t ) + \kappa_1 \overline{\gamma } ^{(p)}_{\rt}\left[ v_{\rti}^{t} + (1-\kappa_1\overline{\gamma } ^{(p)}_{\ro})\Sigma_1^{1/2}\beta _1 + \kappa_1 \overline{\gamma } _{\ro}^{(p)} \Sigma_1^{1/2} \stac - \Sigma_1^{1/2} \xi^{t-1} \right] \\
  \implies ~& v_{\rti}^{t} + (1-\kappa_1\overline{\gamma }_{\ro})\Sigma_1^{1/2}\beta _1 + \kappa_1 \overline{\gamma } _{\ro} \Sigma_1^{1/2} \stac - \Sigma_1^{1/2} \xi^t = \frac{\ix{1}^{\T}(y_1 - \ix{1} \Sigma_1 ^{1/2} \xi^t )}{1-\kappa_1 \overline{\gamma } _{\rt}} + \specop,
\end{aligned}
\end{equation*}
where the first equation is true just by plugging in the definition of  \( r^{t}_{\rti} \), the second multiplies both sides with \( \ix{1}^{\T} \), and the third equation has suppressed \( (\overline{\gamma } _{\ro}^{(p)}-\overline{\gamma } _{\ro})\cdot \kappa _1 \ix{1} (y_1-r^{\infty}_1) \) into \( \specop \). 

Compare the LHS of the last equation above to the definition of \( s^{t}_{\rt} \), and recall that \( \theta ^*_{\rt} (1- \kappa_1 \overline{\gamma } _{\rt}) = 1 \) from Assumption~\ref{asp:2nd.converge.fixed.point}. We notice that 
\begin{equation*}
  s^{t}_{\rt}
   = \frac{1}{\theta ^*_{\rt}} \cdot \frac{\Sigma _{1}^{1/2}\ix{1}^{\T}(y_1 - \ix{1} \Sigma_1 ^{1/2} \xi^t )}{1-\kappa_1 \overline{\gamma } _{\rt}} + \specop 
   = \Sigma _{1}\ix{1}^{\T}(y_1 - \ix{1} \Sigma_1 ^{1/2} \xi^t ) + \specop.
\end{equation*}
Hence \( \nabla L_{\rt}(\xi^t) = - \Sigma _{1}\ix{1}^{\T}(y_1 - \ix{1} \Sigma_1 ^{1/2} \xi^t ) + s^{t}_{\rt} + \specop \). This concludes the proof of Lemma~\ref{lem:2nd.subgradient}. 

\subsubsection{Proof of Lemma~\ref{lem:2nd.bound.on.sc2}}
\label{sec:apd.pf.2nd.bound.on.sc2}

The proof of Lemma~\ref{lem:2nd.bound.on.sc2} only slightly differs for \( \mu _{\rt} = \mu _{\joi} \) and \( \mu _{\rt} = \mu _{\ada} \). We will use the former as an example and the latter can be proved in an analogous way. 

We define some re-scaled quantities for the subsequent proof. Let \( \widetilde{\Sigma }_1 = \Sigma_1 / \lambda _{\min }(\Sigma_1) \), \( \widetilde{\theta } ^{*}_{\rt} = \theta ^*_{\rt} / \lambda _{\min }(\Sigma_1) \), \( \widetilde{\tau } ^{*}_{\rt} = \tau ^*_{\rt} / \sqrt{\lambda _{\min }(\Sigma_1)} \). 

We consider \( \overline{\xi }  \) defined in Eqn~\eqref{eq:def.2nd.xi} from Assumption~\ref{asp:2nd.converge.fixed.point}: \( \overline{\xi } = \xi(\tau _{\rt}^{*}(\zeta ^* Z_1/\tau _1^{*} + \sqrt{1 - (\zeta ^*)^{2}} Z_1'), \overline{\eta } ) \), \( \overline{\eta } = \eta(Z_{1},\cdots ,Z_{E}) \)  where \( Z_1',(Z_1/\tau _1^{*}),(Z_2/\tau _2^{*}),\cdots ,(Z_{E}/\tau _E^{*})\iidsim N(0,I_p) \). We will be using the same set of Gaussian variables throughout the proof. We define an idealized version of \( s_{\rt}^{t} \) called \( \overline{s} _{\rt}   \), 
\begin{equation*}
  \overline{s} _{\rt} = \frac{1}{\widetilde{\theta }  ^*_{\rt}}\left[ \widetilde{\tau } _{\rt}^{*}\widetilde{\Sigma }_1^{1/2}(\zeta ^* Z_{1}/\tau _1^{*}+\sqrt{1-(\zeta ^*)^{2}}Z_1') + (1-\kappa_1 \overline{\gamma } _{\ro} )\widetilde{\Sigma }_1(\beta_1 - \overline{\eta }  ) - \widetilde{\Sigma }_{1}(\overline{\xi }  - \overline{\eta } ) \right].
\end{equation*}
Then \( \overline{s} _{\rt} \in \nabla \mu _{\rt}(\overline{\xi } ; \overline{\eta } ) \), or in the case of the joint estimator, \(\overline{s} _{\rt} / \lambda _{\rt} \in \partial \norm{\overline{\xi } -\overline{\eta } }_1 \).  

Define \( \widetilde{\xi }^{t} = \xi^{t} - \stac + \widetilde{\theta } ^{*}_{\rt} s^{t}_{\rt} / \lambda _{\rt} \).
Recall we have defined \( S^{t}_{\rt}(c) = \{ j\in[p] : \abs{(s^{t}_{\rt})_j} \geq \lambda _{\rt}(1-c)  \}  \). Since \( s^{t}_{\rt} / \lambda _{\rt} \in \partial \norm{\xi ^t - \stac}_1 \), \( S^{t}_{\rt}(c) \) can be equivalently written as \( S^{t}_{\rt}(c) = \{ j\in[p] : \abs{(\widetilde{\xi } ^{t})_j} \geq \widetilde{\theta } ^{*}_{\rt}(1-c) \}  \). 

Define the idealized counterpart of \( \widetilde{\xi } ^{t} \) as \( \bar{\tilde{\xi } }  = \overline{\xi }  - \overline{\eta } + \widetilde{\theta } ^{*}_{\rt} \overline{s}_{\rt} / \lambda _{\rt} \). If we define the idealized counterpart of \( S^{t}_{\rt}(c) \) as \( \overline{S}_{\rt}(c) = \{ j\in[p] : \abs{ (\overline{s}_{\rt})_j } \geq \lambda _{\rt}(1-c) \}  \), then the set \( \overline{S} _{\rt}(c) \) also admits an equivalent form using \( \bar{\tilde{\xi }}  \): 
\( \overline{S} _{\rt}(c) = \{ j\in[p] : \abs{(\bar{\tilde{\xi } } )_j} \geq \widetilde{\theta } ^{*}_{\rt}(1-c)\}  \). We rearrange the terms of \( \bar{\tilde{\xi }}  \) to get an equivalent form:
\begin{equation*}
    \bar{\tilde{\xi}} = (1- \kappa_1 \overline{\gamma } _{\ro}) \widetilde{\Sigma }_1 (\beta _1 - \overline{\eta } ) 
    + \widetilde{\tau } ^{*}_{\rt} \zeta ^* \widetilde{\Sigma }_1^{1/2} Z_1 
    + \widetilde{\Sigma } ^{1/2}_1\left[ \widetilde{\tau } ^{*}_{\rt} \sqrt{1 - (\zeta ^*)^{2}}Z_1' - (I_p - \widetilde{\Sigma } ^{-1}_1) \widetilde{\Sigma } ^{1/2}_1(\overline{\xi } -\overline{\eta } ) \right].
\end{equation*}

The rest of the proof follows the same steps as those of Lemma~\ref{lem:stacked.lasso.bound.on.sc2}'s proof. For any $j\in[p]$, use \( \widetilde{\sigma } _{1,j}\in\R^{p} \) to denote the \( j \)-th row of \( \widetilde{\Sigma } ^{1/2}_{1} \) but written as a vertical vector.  By the normlization, we know \( \norm{\widetilde{\sigma } _{1,j}}_2 \geq 1 \). Define \( Z'_{1,j} = \widetilde{\tau } ^{*}_{\rt} \sqrt{1 - (\zeta ^*)^{2}} (\widetilde{\sigma } _{1,j})^{\T} Z'_{1} \). We write \( (\bar{\tilde{\xi }})_{j} = h_{\rt,j}(Z'_{1,j}) \), where \( h_{\rt,j}(\cdot ) \) also implicitly depends on 
\( Z_1,\cdots ,Z_{E} \) and \( P_{1,j}^{\bot}Z_{1}' \), \( P_{1,j}^{\bot} \) being the projection matrix onto the orthogonal complement of \( \Span(\widetilde{\sigma } _{1,j}) \). Note that these other Gaussian variables \( h_{\rt,j}(\cdot ) \) implicitly depends on are all independent of \( Z_{1,j}' \); conditioning on any values of those Gaussian variables, we have 
\begin{equation}\label{eq:apd.2nd.property.of.hj}
  \abs{h_{\rt,j}(x_1) - h_{\rt,j}(x_2)} \geq \frac{\lambda _{\min }(\Sigma _1)}{\lambda _{\max}(\Sigma_1)} \abs{x_1-x_2} 
\end{equation}
in the same way we have got Equation~\eqref{eq:apd.stack.property.of.hj}. We have omitted some details because they are analogous to the corresponding steps in the proof of Lemma~\ref{lem:stacked.lasso.bound.on.sc2}. 

Note that \( h_{\rt,j}(\cdot ) \) is also Lipschitz in its argument. Similarly to Lemma~\ref{lem:stacked.lasso.indicator.mean} used in the proof of Lemma~\ref{lem:stacked.lasso.bound.on.sc2}, we have the following fact written as a lemma:

\begin{lemma}\label{lem:2step.indicator.mean}
  Let \( \{ I_a:a>0 \}  \) be a collection of intervals indexed by \( a>0 \), satisfying \( \lim_{a\downarrow 0} m(I_a) = 0\) where \( m(\cdot ) \) is the Lebesgue measure. Then \( \lim_{a\downarrow 0} \sup_{p} \E \left[ \frac{1}{p}\sum_{j\in[p]} \I{ (\bar{\tilde{\xi } } )_j \in I_a }  \right] = 0\). 
\end{lemma}
We need yet Lemma~\ref{lem:2step.indicator.converge} before the proof of Lemma~\ref{lem:2nd.bound.on.sc2} itself. It is analogous to Lemma~\ref{lem:stacked.lasso.indicator.converge} for the proof of Lemma~\ref{lem:stacked.lasso.bound.on.sc2}.

\begin{lemma}\label{lem:2step.indicator.converge}
  For  any \( t \in \mathbb{N_+} \), any \( c\in(0,1) \), \( \abs{S^{t}_{\rt}(c)} / p = \E[\abs{\overline{S} _{\rt}(c)} / p] + \smalop  \). 
\end{lemma}

\begin{proof}[\textbf{Proof of Lemma~\ref{lem:2step.indicator.converge}}]
  We fix \( c \) for the rest of the proof. We define the notation \( (x)_+ = \max \{ x,0 \},\forall\, x\in\R \). For $0<a<\widetilde{\theta } ^{*}_{\rt}(1-c)/3$, define a Lipschitz function \( \phi ^{(1)}_{\rt,1}(x) : \R\to\R \): \( \phi ^{(1)}_{\rt,1}(x) = [1 - (\widetilde{\theta }_{\rt} ^{*}(1-c) - \abs{x} )_+ / a]_+ \). Then \( \abs{\phi ^{(1)}_{\rt,1}(x) - \1\{\abs{x} \geq \widetilde{\theta } _{\rt}^{*}(1-c)\}}  \leq \1\{\abs{x} \in (\widetilde{\theta }_{\rt} ^{*}(1-c)-a,\widetilde{\theta }_{\rt} ^{*}(1-c)) \} \). 
  
  For any $x\in\R$, define \( \mathrm{dist}(x)  \) as the distance of \( \abs{x} \) to the interval \( (\widetilde{\theta }_{\rt} ^{*}(1-c)-a,\widetilde{\theta }_{\rt} ^{*}(1-c)) \), i.e. \( \mathrm{dist}(x) = (\abs{x}  - \widetilde{\theta }_{\rt} ^{*}(1-c))_+ + (\widetilde{\theta }_{\rt} ^{*}(1-c)-a-\abs{x}  )_+ \). 
  Define another Lipschitz function \( \phi ^{(1)}_{\rt,2}(x) : \R\to\R \): \( \phi ^{(1)}_{\rt,2}(x) = (1 - \mathrm{dist}(x)/a)_+ \). Then \(  \I{\abs{x} \in (\widetilde{\theta }_{\rt} ^{*}(1-c)-a,\widetilde{\theta }_\rt ^{*}(1-c)) } \leq \phi ^{(1)}_{\rt,2}(x) \leq  \1\{\abs{x} \in (\widetilde{\theta }_{\rt} ^{*}(1-c)-2a,\widetilde{\theta }_{\rt} ^{*}(1-c)+a) \} \). For $i=1,2$, define \( \phi _{\rt,i}^{{(p)}}:\R^{p}\to \R \): \( \forall\,x\in\R^{p} \), \( \phi _{\rt,i}^{{(p)}}(x) = \sum_{j\in[p]} \phi ^{(1)}_{\rt,i} (x_j) \). By Lemma~\ref{lem:apd.pseudo}, both \( \{ \phi _{\rt,1}^{{(p)}}:p\geq 1 \}  \) and \( \{ \phi _{\rt,2}^{{(p)}}:p\geq 1 \}  \) are uniformly Lipschitz functions. 

  It is already known that  \( S^{t}_{\rt}(c) = \{ j\in[p] : \abs{(\widetilde{\xi }^{t})_j} \geq \widetilde{\theta } ^{*}_{\rt}(1-c) \}  \), \( \overline{S} _{\rt}(c) = \{ j\in[p] : \abs{(\bar{\tilde{\xi }})_j} \geq \widetilde{\theta } ^{*}_{\rt}(1-c) \}  \). 
  By Lemma~\ref{lem:2nd.amp.plus}, for any already chosen and fixed \( a \), as \( p\to \infty \), 
\begin{equation*}
\begin{aligned}
  \abs*{{\abs{S_\rt^{t}(c)}}/{p} -  \E[{\abs{\overline{S}_\rt (c)}}/{p}]} \leq&~ \abs{\phi _{\rt,1}^{(p)}(\widetilde{\xi }^{t}) - \E[\phi _{\rt,1}^{(p)}(\bar{\tilde{\xi }})]} + \phi _{\rt,2}^{(p)}(\widetilde{\xi }^{t}) + \E[\phi _{\rt,2}^{(p)}(\bar{\tilde{\xi } })] \leq 2 \E[\phi _{\rt,2}^{(p)}(\bar{\tilde{\xi } })] + \smalop \\
  \leq &~ \frac{2}{p}\sum_{j\in[p]}\E[ \I{\abs{(\bar{\tilde{\xi }}) _j} \in (\widetilde{\theta }_{\rt} ^{*}(1-c)-2a,\widetilde{\theta }_{\rt} ^{*}(1-c)+a) }] + \smalop. 
\end{aligned}
  \end{equation*}
  By Lemma~\ref{lem:2step.indicator.mean}, for any \( \epsilon>0 \), we can fix a small \( a \) s.t. $$ \sup_{p}\sum_{j\in[p]}\E[ \I{\abs{(\bar{\tilde{\xi } }) _j} \in (\widetilde{\theta }_{\rt} ^{*}(1-c)-2a,\widetilde{\theta }_{\rt} ^{*}(1-c)+a) }] / p < \epsilon/4 $$ before taking \( p\to \infty \).  Thus Lemma~\ref{lem:2step.indicator.converge} is proved. 
\end{proof}

Finally we are ready to prove Lemma~\ref{lem:2nd.bound.on.sc2}. 
\begin{proof}[\textbf{Proof of Lemma~\ref{lem:2nd.bound.on.sc2}}]
  For any chosen and fixed $c\in(0,1)$, by Lemma~\ref{lem:2step.indicator.converge}, \( \abs{S_\rt^{t}(c)} / p = \E[\abs{\overline{S}_{\rt} (c)} / p] + \smalop \) as \( p\to \infty \).
Furthermore, 
\begin{equation*}
\begin{aligned}
  \E[\abs{\overline{S}_{\rt} (c)} / p] =&~ \frac{1}{p}\sum_{j\in[p]} \E \left[ \I{\abs{(\bar{\tilde{\xi }}) _j} > \widetilde{\theta }_{\rt} ^{*} } \right] + \frac{1}{p}\sum_{j\in[p]} \E \left[ \I{\abs{(\bar{\tilde{\xi }}) _j} \in (\widetilde{\theta }_{\rt} ^{*}(1-c), \widetilde{\theta }_{\rt} ^{*}] } \right].
\end{aligned}
\end{equation*}
By Lemma~\ref{lem:2step.indicator.mean}, $\forall\,\epsilon>0$,  \(\exists \, c\in(0,1) \) close enough to 0, s.t. \( \sup_p\sum_{j\in[p]} \E [ \I{\abs{(\bar{\tilde{\xi } }) _j} \in (\widetilde{\theta } ^{*}(1-c), \widetilde{\theta } ^{*}] } ] / p < \epsilon \). By the definition of \( \bar{\tilde{\xi }}  \), 
\begin{equation*}
\begin{aligned}
  \frac{1}{p}\sum_{j\in[p]} \E \left[ \I{\abs{(\bar{\tilde{\xi }} ) _j} > \widetilde{\theta }_{\rt} ^{*} } \right] =&~ \frac{1}{p} \E \left[ \norm{\overline{\xi } - \overline{\eta } }_0  \right] = \overline{\gamma } _{\rt}^{(p)} \to \overline{\gamma } _{\rt},
\end{aligned}
\end{equation*}
where the second equality is by Eqn~\eqref{eq:apd.prelim.2nd.partial.deriv.joi} and the convergence is by Proposition~\ref{prop:2nd.onsager}.

By the choice of \( c \) so far, we know it only depends on the `idealized' quantities such as \( \overline{\eta }  \) and \( \bar{\tilde{\xi } }  \), not any specific time step \( t \).  As a result,  \(\forall\, \epsilon>0 \), there exists \( c = c(\epsilon) \in(0,1) \) s.t. fixing any \( c'\in(0,c(\epsilon)] \) and any \( t \in \mathbb{N_+} \), we have \( \lim_{p} \PR\left[\abs{S_\rt^{t}(c')}/p \leq \overline{\gamma } _{\rt} + \epsilon  \right] = 1 \). This completes the proof. 
\end{proof}

\subsubsection{Proof of Lemma~\ref{lem:2nd.final.converge}}
\label{sec:apd.pf.2nd.final.converge}

The proof of Lemma~\ref{lem:2nd.final.converge} is analogous to that of Lemma~\ref{lem:stack.final.converge} in Section~\ref{sec:apd.pf.stacked.lasso.final.converge}, so we may skip certain details. Just like Lemma~\ref{lem:2nd.bound.on.sc2}, the proof for the joint estimator (\( \mu _{\rt} = \mu _{\joi} \)) is only slightly different to the case of the adaptively weighted estimator (\( \mu _{\rt} = \mu _{\ada} \)). We do the former proof  and the latter can be derived in an analogous way.

Recall how we have defined the sparse singular eigenvalues of a matrix in Section~\ref{sec:apd.pf.stacked.lasso.final.converge}: For any matrix \( X \) having \( p \) columns,  any subset \( S \subset [p] \), define \( \kappa _{-}(X,S) = \inf \left\{ \norm{Xu}_2:\Sup(u)\subset S,\norm{u}_2=1 \right\} \)
and the \( k \)-sparse singular value 
\begin{equation*}
  \kappa _{-}(X,k) = \min _{S\subset[p],\abs{S}\leq k} \kappa _{-}(X,S) = \min _{S\subset[p],\abs{S} = k} \kappa _{-}(X,S) .
\end{equation*}

\begin{lemma}[Lower Bound on the Sparse Singular Values]\label{lem:2nd.sparse.singval}
  Under Assumption~\ref{asp:prelim.transfer.model}, there exists small enough constants \( \epsilon _0>0 \) and $c_0>0$, s.t. \( \kappa_1\cdot  (\overline{\gamma } _{\rt} + \epsilon _0) < 1 \), and \( \forall\,\epsilon \in(0,\epsilon _0) \), 
\begin{equation*}
  \lim_{p\to \infty} \PR \left[ \kappa _{-}\left(\ix{1} \Sigma_1^{1/2}, p(\overline{\gamma } _{\rt}+\epsilon)\right) > c_0 \right]  = 1.
\end{equation*}
\end{lemma}

The proof of Lemma~\ref{lem:2nd.sparse.singval} is analogous to Lemma~\ref{lem:stacked.lasso.sparse.singval} so it is omitted here. 

\begin{proof}[\textbf{Proof of Lemma~\ref{lem:2nd.final.converge}}]
  \newcommand{\yst}{y_1}
  \newcommand{\xst}{\ix{1} \Sigma_1^{1/2}}
  \newcommand{\Xst}{\ix{1} \Sigma_1^{1/2}}

Define \( \Delta \xi ^{t} = \se  - \xi^{t} \), where \( \se  \) is the minimizer of \( L_{\rt}(\cdot ) \) defined in Eqn~\eqref{eq:def.2nd.los}.  \( L_{\rt}(\xi^t + \Delta \xi^t) \leq L_{\rt}(\xi ^t) \). 
Applying Lemma~\ref{lem:apd.sing.val.bounds}\emph{(b)} and then Lemma~\ref{lem:apd.sing.val.limit.concent}\emph{(a)} to \( \ix{1} \Sigma_1^{1/2} \), we know  \( \sigma _{\max}(\ix{1} \Sigma_1^{1/2}) \leq 2(c_1 / \kappa_1)^{1/2}  (1 + \sqrt{\kappa_1}) \)  almost surely as \( p\to \infty \), in which \( c_1 \) is the constant from Assumption~\ref{asp:prelim.transfer.model}. 

We now define the high-probability event we will be working on throughout the proof. We choose $M>1$ large enough and \( \epsilon _0\in(0,1 ) \) small enough such that  \( \kappa_1 \,(\overline{\gamma } _{\rt}+2{\epsilon _0}) < 1 \) to satisfy the requirement of Lemma~\ref{lem:2nd.sparse.singval} s.t. \( \forall \, \epsilon \in(0,\epsilon _0) \), \( \exists\,T \) s.t. \( \forall\,t\geq T \), at the time step \( t \), the following conditions are satisfied simultaneously with probability going to 1 as \( p\to \infty \): 
\begin{enumerate}
  \item[\it(i)]  \( \max\{\frac{1}{\sqrt{p}}\norm{\se}_2, \frac{1}{\sqrt{p}}\norm{\xi^t}_2\} \leq M \) (by Lemma~\ref{lem:2nd.bounded} and Lemma~\ref{lem:2nd.amp.plus});
  \item[\it(ii)] \( \abs{ S_{\rt}^{t}(c_2)} \leq p(\overline{\gamma } _{\rt} + {\epsilon}) \) for some \( c_2=c_2(\epsilon) \in (0,1) \) (by Lemma~\ref{lem:2nd.bound.on.sc2} with \( c_2 \) depending only on the specific \( \epsilon \in(0,\epsilon _0) \) but invariant to the choice of the time step \( t \in \mathbb{N_+} \));
  \item[\it(iii)] \( \kappa _{-}(\ix{1} \Sigma_1^{1/2}, p( \overline{\gamma } _{\rt}+2{\epsilon})) > c_3 \) for some $c_3 = c_3(\epsilon _0)\in(0,1)$ (by Lemma~\ref{lem:2nd.sparse.singval} with \( c_3 \) depending on the upper bound \( \epsilon _0 \));
  \item[\it(iv)] \( \frac{1}{\sqrt{p}}\norm{\nabla L_{\rt}(\xi^t)}_2 < \epsilon ^2 \cdot  c_2^{2} c_3^{4} / (4M) \) (by Lemma~\ref{lem:2nd.subgradient}; the only condition that determines the choice of \( T \));
  \item[\it(v)] \( \sigma _{\max}(\ix{1} \Sigma_1^{1/2}) \leq c_4 = 2( c_1 / \kappa_1)^{1/2}  (1 + \sqrt{\kappa_1}) \). 
\end{enumerate}

  We are working under \( \mu _{\rt} = \mu _{\joi} \), so \( (\se - \stac) \) is sparse. Accordingly, let $S_{\rt} = \Sup(\xi^t - \stac) \subset [p]$. 
  We have the following inequalities:
  \begin{equation*}
  \begin{aligned}
   0\geq &~ \frac{1}{p} \left[ L_{\rt}(\xi^t+\Delta \xi^t) - L_{\rt}(\xi^{t}) \right] \\
   \geq &~ \underbrace{\frac{\lambda _\rt}{p} \left( \norm{(\xi^t+\Delta \xi^t)_{S_\rt}}_1 - \norm{(\xi^t) _{S_\rt}}_1 - \sign((\xi^t)_{S_{\rt}})^{\T} (\Delta \xi^{t}) _{S_\rt}\right) }_{(1)} \\
   &~
   + \underbrace{\frac{\lambda _{\rt}}{p}\left( \norm{(\Delta \xi ^{t}) _{S_{\rt}^{c}}}_1 - [(s_{\rt}^{t})_{S_{\rt}^{c}}]^{^{\T}} (\Delta \xi^{t})_{S_{\rt}^{c}} \right)}_{(2)}  
   + \underbrace{\frac{1}{2p} \norm{\Xst\Delta \xi^t}_2^{2}}_{(3)} 
    + \frac{1}{p}(\nabla L_{\rt}(\xi^{t}))^{\T} \Delta \xi^{t}.
  \end{aligned}
  \end{equation*}
The terms \( (1) \), \( (2) \) and \( (3) \)  are all non-negative in the same way that we have asserted so in the proof of Lemma~\ref{lem:stack.final.converge}. Then we must have \( (\nabla L_{\rt}(\xi^{t}))^{\T} \Delta \xi^{t} \leq 0 \). By conditions \emph{(i)}, \emph{(iv)} of the high-probability event and Cauchy's inequality, we know   
  \begin{equation*}
    \lim_{p\to \infty}\PR\left[(1) + (2) + (3) \leq \frac{1}{p}\abs{(\nabla L_{\rt}(\xi^t))^{\T} \Delta \xi^t } \leq \epsilon ^2 \cdot c_2^{2} c_3^{4} \right] = 1.
  \end{equation*}
  Respectively for terms (2) and (3), the above tells us on the high probability event, 
  \begin{equation}
  \begin{aligned}\label{eq:apd.2step.laso.convexity.1}
    0\leq&~ \frac{1}{p}\norm{(\Delta \xi^{t}) _{S_{\rt}^{c}}}_1 - \frac{1}{p}[(s_{\rt}^{t})_{S_{\rt}^{c}}]^{^{\T}} (\Delta \xi^{t})_{S_{\rt}^{c}}  \leq  \epsilon ^2 \cdot c_2^{2} c_3^{4}/ \lambda _{\rt} \\
    0\leq&~ \frac{1}{p} \norm{\Xst \Delta \xi^t}_2^{2} \leq 2\epsilon ^2 \cdot c_2^{2} c_3^{4}.
  \end{aligned}
  \end{equation}

  Take the SVD of \( \Xst :\, \Xst = \sum_{i=1}^{n_{1} \land p} \sigma _{1,i} u_{1,i} v_{1,i}^{\T} \), with \( \sigma _{1,1} \geq \sigma _{1,2} \geq \cdots \geq \sigma _{1, n_{1}\land p} \). Let \( V_{(1)} = \Span(v_{1,i} : \sigma _{1,i}\geq c_{3}/2) \) and \( V_{(2)} \) be the orthogonal complement of \( V_{(1)} \) in \( \R^{p} \). Let \( P_{(1)},P_{(2)} \) be the projection matrices onto \( V_{(1)} \) and \( V_{(2)} \). Denote \( X_{(1)} = \xst P_{(1)}  \) and \( X_{(1)} = \xst P_{(2)}  \). 
   Then \(  X_{(1)}^{\T} X_{(2)} = 0 \). Define \( \Delta \xi_{(1)} = P_{(1)} \Delta \xi^t\), \( \Delta \xi_{(2)} = P_{(2)} \Delta \xi^t\). The second line of Eqn~\eqref{eq:apd.2step.laso.convexity.1} implies
  \begin{equation*}
    \frac{1}{p} \norm{X_{(1)}\Delta \xi_{(1)}}_2^{2} \leq  2\epsilon ^2 \cdot c_2^{2} c_3^{4},\quad 
    \frac{1}{p} \norm{X_{(2)}\Delta \xi_{(2)}}_2^{2} \leq  2\epsilon ^2 \cdot c_2^{2} c_3^{4}.
  \end{equation*}
  For $\Delta \xi_{(1)}$ corresponding to the larger singular values, we have 
  \begin{equation}\label{eq:apd.2step.laso.convexity.2}
    \frac{1}{p} \norm{\Delta \xi_{(1)}}_2^{2} \leq 8 \epsilon ^2 \cdot c_2^{2} c_3^{4} / (c_3^{2}) = 8 \epsilon ^2 \cdot c_2^{2}c_3^{2}.
  \end{equation} 
  
  We are left to bound \( \frac{1}{p}\norm{\Delta \xi_{(2)}}_2^{2} \). 
  We utilize the first line of Eqn~\eqref{eq:apd.2step.laso.convexity.1}. For \( \Delta \xi_{(1)} \), \(   \frac{1}{p} [(s_{\rt}^{t})_{S_{\rt}^{c}}]^{\T} (\Delta \xi_{(1)})_{S_\rt^{c}} \leq \frac{1}{p}\norm{(\Delta \xi_{(1)})_{S_{\rt}^{c}}}_1 \leq \frac{1}{\sqrt{p}} \norm{\Delta \xi_{(1)}}_2 \leq 2 \sqrt{2} \epsilon  c_2 c_3 \). Since \( \Delta \xi_{(2)} = \Delta \xi ^t- \Delta\xi _{(1)} \), we know \( \frac{1}{p}\norm{(\Delta \xi_{(2)})_{S_\rt^{c}}}_1 - \frac{1}{p} [(s_\rt^{t})_{S^{c}}]^{\T} (\Delta \xi_{(2)})_{S_\rt^{c}} \leq \epsilon ^2\cdot (c_2 {c_3})^{2}/ \lambda _{\rt} + 4 \sqrt{2} \epsilon c_2 c_3 \). 
  Recall that \( S^{t}_{\rt}(c_2) \) defined in Lemma~\ref{lem:2nd.bound.on.sc2} is a subset of \( S_{\rt}^{c} \). By the definition of \( S^{t}_{\rt}(c_2) \), 
  \begin{gather*}
    \frac{1}{p}\norm{(\Delta \xi_{(2)})_{S_{\rt}^{t}(c_2)^{c}}}_1 - \frac{1}{p} (s^{t}_{S_\rt^{t}(c_2)^{c}})^{\T} (\Delta \xi_{(2)})_{S_\rt^{t}(c_2)^{c}} \leq \epsilon ^2 \cdot c_2^{2} c_3^{4}/ \lambda _{\rt} + 4 \sqrt{2} \epsilon c_2 c_3 \notag \\
    \implies \frac{1}{p} \norm{(\Delta \xi_{(2)})_{S_{\rt}^{t}(c_2)^{c}}}_1 \leq \frac{1}{c_2}\left( \epsilon ^2 \cdot c_2^{2} c_3^{2}/ \lambda _\rt + 4 \sqrt{2} \epsilon  c_2 c_3 \right) \leq ( \lambda _\rt ^{-1} + 4 \sqrt{2} )\epsilon c_3. 
  \end{gather*}

Recall we have chosen \( \kappa_1 \,( \overline{\gamma } _1 + 2{\epsilon}) \leq  \kappa_1 \,( \overline{\gamma } _1 + 2\epsilon_0) < 1 \) as required by Lemma~\ref{lem:2nd.sparse.singval} used for condition \emph{(iii)} of the high-probability event. Condition \emph{(ii)} guarantees \( \abs{S_{\rt}^{t}(c_2)} \leq p( \overline{\gamma } _{\rt} + {\epsilon})  \). Then we immediately know \( \abs{S^{t}_{\rt}(c_2)^{c}} > p {\epsilon}  \), because otherwise \( p = \abs{S^{t}_{\rt}(c_2)^{c}} + \abs{S_{\rt}^{t}(c_2)} \leq p( \overline{\gamma } _{\rt}+2{\epsilon}) < p \) would give a contradiction. We partition \( S^{t}_{\rt}(c_2)^{c} \) disjointly into \( S_\rt^{t}(c_2)^{c} = \cup_{l=1}^{K} S_l \), with each \( S_l \) satisfying \( \abs{S_l} \in [p {\epsilon}/2,p {\epsilon}] \), and they are constructed in a descending order in the sense that for any \( i\in S_{l} \) and \( j\in S_{l+1} \), \( \abs{({\Delta } \xi_{(2)})_i} \geq \abs{(\Delta \hat\beta _{(2)})_j}  \). By the disjoint partition, 
\begin{equation*}
  \frac{1}{p}\norm*{\Delta \xi_{(2)}}_2^{2} =   \frac{1}{p}\norm*{\left(\Delta \xi_{(2)}\right)_{\cup_{l=2}^{K}S_l}}_2^{2} + \frac{1}{p}\norm*{\left(\Delta \xi_{(2)}\right)_{S^{t}(c_2)\cup S_1}}_2^{2}.
\end{equation*}
We bound the two terms on the RHS respectively. For the first term, same as Eqn~\eqref{eq:apd.laso.convexity.4}, we eventually get
\begin{equation}
\begin{aligned}\label{eq:apd.2step.laso.convexity.4}
\frac{1}{p}\norm*{\left(\Delta \xi_{(2)}\right)_{\cup_{l=2}^{K}S_l}}_2^{2}  \leq 4 \epsilon  c_3^{2} (\lambda ^{-1}_{\rt} + 4 \sqrt{2})^{2}. 
\end{aligned}
\end{equation}
  
For the second term, we apply conditions \emph{(iii)}, \emph{(v)}  of the high-probability event same as the steps in the proof of Lemma~\ref{lem:stack.final.converge} to get:
\begin{equation*}
\begin{aligned}
\frac{1}{p} \norm{(\Delta \xi_{(2)})_{S_\rt^{t}(c_2)\cup S_1}}_2^{2} \leq&~  \frac{1}{2p} \norm{(\Delta \xi_{(2)})_{S_\rt^{t}(c_2)\cup S_1}}_2^{2} + \left(\frac{1}{2} + \frac{2 c_4^{2}}{c_3^{2}}\right) \frac{1}{p} \norm{(\Delta \xi_{(2)})_{\cup_{l=2}^{K}S_l}}_2^{2} \\  
\leq&~ \left(1 + \frac{4 c_4^{2}}{c_3^{2}}\right) \frac{1}{p} \norm{(\Delta \xi_{(2)})_{\cup_{l=2}^{K}S_l}}_2^{2}.
\end{aligned}
\end{equation*}
Further, we plug in the bound from Eqn~\eqref{eq:apd.2step.laso.convexity.4},   
\begin{equation*}
  \begin{aligned}
  \frac{1}{p} \norm{(\Delta \xi_{(2)})_{S_\rt^{t}(c_2)\cup S_1}}_2^{2} \leq&~ \left(1 + \frac{4 c_4^{2}}{c_3^{2}}\right) \frac{1}{p} \norm{(\Delta \xi_{(2)})_{\cup_{l=2}^{K}S_l}}_2^{2}
  \leq 4 \epsilon \left(c_3^{2} + 4 c_4^{2}\right) (\lambda _{\rt} ^{-1} + 4 \sqrt{2})^{2}, 
\end{aligned}
\end{equation*}
which, put together with Eqn~\eqref{eq:apd.2step.laso.convexity.4} and Eqn~\eqref{eq:apd.2step.laso.convexity.2}, yields
\begin{equation*}
\begin{aligned}
   &~\frac{1}{p}\norm{\Delta \xi_{(2)}}_2^{2} \leq 
   8 \epsilon \left(c_3^{2} + 2 c_4^{2}\right) (\lambda _\rt ^{-1} + 4 \sqrt{2})^{2}   \\
   \implies&~ \frac{1}{p} \norm{\Delta\xi^t}_2^{2} \leq \left[ {8}\epsilon c_2^{2}{c_3^{2}} + 8 \left(c_3^{2} + 2 c_4^{2}\right) (\lambda _\rt ^{-1} + 4 \sqrt{2})^{2} \right] \cdot {\epsilon} \leq \left[ {8}+ 8 \left(1 + 2 c_4^{2}\right) (\lambda_\rt ^{-1} + 4 \sqrt{2})^{2} \right] \cdot {\epsilon}.
  \end{aligned}
  \end{equation*}
Recall that \( c_4= 2( c_1 / \kappa_1)^{1/2}  (1 + \sqrt{\kappa_1}) \) is a fixed constant free of \( \epsilon \). 

To re-iterate what we have proved: There exists some small \( \epsilon _0\in(0,1) \), s.t. \( \forall\,\epsilon \in(0,\epsilon _0) \), \( \exists \, T \) s.t. \( \forall\,t\geq T \), \( \lim_{p} \PR \left[ \frac{1}{p}\norm{\xi^t-\se}_2^{2} < \epsilon\right] = 1  \). This completes the proof of Lemma~\ref{lem:2nd.final.converge}.  
\end{proof}

\section{Supporting Technical Lemmas}
\label{sec:apd.support}

\begin{lemma}[Relationships Between Pseudo-Lipschitz Functions] \label{lem:apd.pseudo}
  Let $q$ be a fixed known constant. Suppose we have  \( \psi:\R^{q}\to\R \) that is pseudo-Lipschitz of order \( k \) for some \( k\geq 1 \), according to Definition~\ref{def:prelim.plip.vec}. Define a family of functions \( \{ \psi _N :N=1,2,\cdots  \}  \) by
  \begin{equation*}
    \psi _N:(\R^{N})^{q}\to \R;\quad \psi _N(x_1,x_2,\cdots ,x_q) = \frac{1}{N} \sum_{j\in[N]} \psi\left((x_1)_j,(x_2)_j,\cdots ,(x_q)_j\right).
  \end{equation*}
  Then, \( \{ \psi_N \}  \)  is uniformly pseudo-Lipschitz iff \( 1\leq k\leq 2 \); when it is, it is uniformly pseudo-Lipschitz of order \( k \). 
\end{lemma}

\begin{proof}[\textbf{Proof of Lemma~\ref{lem:apd.pseudo}}]
  When \( k\leq 2 \), we prove \( \{\psi _N\} \)  constructed in such a way must be pseudo-Lipschitz. Since we have argued Defintions~\ref{def:prelim.plip.vec} and \ref{def:prelim.plip.mat} are compatible, we can write \( X = [x_1|x_2|\cdots |x_q]\in\R^{N\times q} \) and \(  Y = [y_1|y_2|\cdots |y_q]\in\R^{N\times q} \), slightly abuse notation by letting \( \psi _N \) be \( \R^{N\times q}\to \R \) and verify Definition~\ref{def:prelim.plip.mat}. Then 
  \begin{equation*}
\begin{aligned}
    \abs{\psi _N(X) - \psi_N(Y)} \leq &~ \frac{L}{N}\sum_{j\in[N]} \left({1 + \norm{X_{j,\cdot }}_2^{k-1} + \norm{Y_{j,\cdot }}_2^{k-1}} \right)\cdot \norm*{X_{j,\cdot }-Y_{j,\cdot }}_2 \\
    \leq &~ L \left[\frac{1}{N}\sum_{j\in[N]} \left({1 + \norm{X_{j,\cdot }}_2^{k-1} + \norm{Y_{j,\cdot }}_2^{k-1}} \right)^{2}\right]^{1/2}\cdot  \frac{\norm*{X-Y}_F}{\sqrt{N}}\quad \text{(Cauchy's inequality)}.
\end{aligned}
  \end{equation*}
By the relationship between the spectral norm and the Frobenius norm, \( \norm*{X-Y}_F \leq \sqrt{q} \cdot \norm*{X-Y}_2 \). For the sum of squares in the above inequality, we know 
\begin{equation*}
\begin{aligned}
  \frac{1}{N}\sum_{j\in[N]} \left({1 + \norm{X_{j,\cdot }}_2^{k-1} + \norm{Y_{j,\cdot }}_2^{k-1}} \right)^{2} \leq&~  \frac{3}{N}\sum_{j\in[N]} \left[ 1 + (\norm{X_{j,\cdot }}_2^{2})^{k-1} + (\norm{Y_{j,\cdot }}_2^{2})^{k-1} \right] \\
  = &~ 3 + \frac{3\sum_{j\in[N]}(\norm{X_{j,\cdot }}_2^{2})^{k-1}}{N} + \frac{3\sum_{j\in[N]}(\norm{Y_{j,\cdot }}_2^{2})^{k-1}}{N}\\
  \leq &~ 3 + 3\left( \frac{\norm{X}_F^{2}}{N} \right)^{k-1} + 3 \left( \frac{\norm{Y}_F^{2}}{N} \right)^{k-1} \quad \text{(Jensen's inequality)} \\
  \leq &~ 3 q^{k-1} \left[ 1 + \left( \frac{\norm{X}_2^{2}}{N} \right)^{k-1} + \left( \frac{\norm{Y}_2^{2}}{N} \right)^{k-1} \right] \\
  \leq &~ 3 q^{k-1} \left[ 1 + \left( \frac{\norm{X}_2}{\sqrt{N}} \right)^{k-1} + \left( \frac{\norm{Y}_2}{\sqrt{N}} \right)^{k-1} \right]^{2}.
\end{aligned}
\end{equation*}
Plug this into the inequality we have had for \( \abs{\psi _N(X) - \psi_N(Y)} \), and then we have proved 
\begin{equation*}
  \abs{\psi _N(X) - \psi_N(Y)} \leq L \sqrt{3 q^{k-1}}\left[ 1 + \left( \frac{\norm{X}_2}{\sqrt{N}} \right)^{k-1} + \left( \frac{\norm{Y}_2}{\sqrt{N}} \right)^{k-1} \right] \frac{\norm{X-Y}_2}{\sqrt{N}},
\end{equation*}
meaning \( \{\phi _N\} \)  is uniformly pseudo-Lipschitz. 

When \( k>2 \), as long as there is some \( \psi \) exactly attaining the rate of the order-$k$ pseudo-Lipschitz bound, we can give a counter-example against \( \{ \psi _N \}   \) being uniformly pseudo-Lipschitz. For example, let \( q=1 \) and \( \psi(x) = [(x)_+]^{k} \). We can evaluate the functions at \( x = [N,0,0,\cdots ,0]^{\T}\in\R^{N} \) and \( y=\vec{0}\in\R^{N} \). We will notice that as \( N\to\infty \), no finite \( L \) can meet the requirement for uniform pseudo-Lipschitz-ness.   
\end{proof}

\begin{lemma}\label{lem:apd.lem23}
  Let \( \Sigma\succcurlyeq 0 \) be a covariance matrix, \( Z \in \R^{N\times q} \) and \( Z^{\T} \) has column vectors drawn iid from \( \Nscr(0,\Sigma) \). For any order-$k$ uniformly pseudo-Lipschitz functions \( \{ \phi _N :\R^{N\times q}\to \R \}  \) according to Definition~\ref{def:prelim.plip.mat}, \( \phi _N(Z) = \E[\phi _N(Z)] + \smalop \) as \( N\to \infty \).   
\end{lemma}
\begin{proof}[\textbf{Proof of Lemma~\ref{lem:apd.lem23}}]
Without loss of generality, assume \( \Sigma \) is invertible. Otherwise, some dimension from \( \Nscr(0,\Sigma) \) can be written as a linear combination of the others, and we can just cut  the dimension down from $q$ by only using the other dimensions to eliminate the degeneracy. 

We use induction on \( q \). For \( q= 1 \), it is proved in Lemma 23, \citet{berthier2020state}. Now suppose we have proved the results for \( q=q_0 \) and are considering \( (q_0+1) \). Equivalently, let's write \( Z = [Z_1|Z_2] \) where \( Z_1\in\R^{N\times q_0} \), \( Z_2\in\R^{N} \), and each row of \( Z = [Z_1|Z_2]\) is an iid draw from \( \Nscr(0,\overline{\Sigma}) \), where
\begin{equation*}
  \overline{\Sigma} = \begin{bmatrix}
    \Sigma & \alpha \\
    \alpha ^\T & \beta.
  \end{bmatrix}
\end{equation*}
Then 
\begin{equation*}
\begin{aligned}
  \phi _N(Z_1,Z_2) |_{Z_1} \disteq&~ \phi _N(Z_1,\sqrt{\beta-\alpha ^\T \Sigma ^{-1}\alpha} \tilde{Z}_2 +  Z_1 \Sigma ^{-1} \alpha  )\\
  =&~ \E_{\tilde{Z}_2\sim \Nscr(0,I_N)}\left[ \phi _N(Z_1,\sqrt{\beta-\alpha ^\T \Sigma ^{-1}\alpha} \tilde{Z}_2 +  Z_1 \Sigma ^{-1} \alpha  ) \right] + \smalop.
\end{aligned}
\end{equation*}
Then, using the induction hypothesis,
\begin{equation*}
  \begin{aligned}
    \phi _N(Z_1,Z_2)
    =&~ \E_{{Z_1\independent\tilde{Z}_2}}\left[ \phi _N(Z_1,\sqrt{\beta-\alpha ^\T \Sigma ^{-1}\alpha} \tilde{Z}_2 +  Z_1 \Sigma ^{-1} \alpha  ) \right] + \smalop.
  \end{aligned}
  \end{equation*}
In the argument, each row of \(Z_1, (\sqrt{\beta-\alpha ^\T \Sigma ^{-1}\alpha} \tilde{Z}_2 +  Z_1 \Sigma ^{-1} \alpha)  \) is an iid draw from \( \Nscr(0,\overline{\Sigma}) \). 
\end{proof}

\begin{lemma}[Extended Stein's Lemma]\label{lem:apd.stein}
  Suppose we have a matrix of jointly normal elements, \( Z\in\R^{N\times(m_1+m_2)} \), with two blocks \( Z = [Z_1|Z_2] \), \( Z_1\in\R^{N\times m_1} \) and \( Z_2\in\R^{N\times m_2} \). There is some covariance matrix \( \Sigma \in \R^{(m_1+m_2)\times (m_1+m_2)} \), s.t. \( Z^{\T} \) has iid column vectors drawn from the \( (m_1+m_2) \)-dimensional normal distribution, \( \Nscr(0,\Sigma) \). In particular, on the same row, the row vectors of \( Z_1 \)  and \( Z_2 \) have a covariance matrix, \( \Sigma _{1,2} \in\R^{m_1\times m_2} \). 
  
  Suppose we have a Lipschitz function \( f:\R^{N\times m_2}\to \R^{N\times q} \), \( f^{\T} = [f_1|f_2|\cdots |f_N] \), where $\forall\,j\in[N]$, \( f_j :\R^{N\times m_2}\to \R^{q} \) is also Lipschitz, with a Jacobian matrix taken only w.r.t. its \( j \)-th row denoted as \( \frac{\partial f_j}{\partial  x_{j,\cdot }}(\cdot ) \in\R^{q \times m_2} \). 

  Then we have \begin{equation*}
    \E\left[ f(Z_2)^{\T} Z_1 \right] = \sum_{j\in[N]} \E \left[  \frac{\partial f_j}{\partial  x_{j,\cdot }}( Z_2 ) \right] \Sigma _{1,2}^{\T} \in \R^{q\times m_1}.
  \end{equation*}
\end{lemma}
The proof is straightforward and thus omitted. In the proof, it is useful to condition on all other rows when studying one specific row of both \( Z_1 \) and \( Z_2 \). 

\begin{lemma}[Singular Values of Gaussian Matrices]\label{lem:apd.sing.val.limit.concent}
  Let \( \sigma _{\min }(\cdot ) \) and \( \sigma _{\max }(\cdot ) \) denote the smallest and largest non-zero singular values of a matrix. 
  Suppose \( X\in\R^{N\times p} \), \( X_{i,j}\iidsim \Nscr(0,\frac{1}{N}) \). Consider the asymptotic regime where \( N,p\to \infty \) and \( \kappa=\lim_{N,p\to \infty}\frac{p}{N}\in(0,\infty) \). Then we have the following two facts:
  \begin{enumerate}
    \item[(a)] \( \lim_{p} \sigma _{\max }(X) = 1 + \sqrt{\kappa },a.s.,\quad \lim_{p} \sigma _{\min }(X) = \abs{1 - \sqrt{\kappa }}, a.s. \)
    \item[(b)] There exists \( p_0 \) s.t. \( \forall\,p\geq p_0 \), \( \PR\left[ \abs{1-\sqrt{\kappa }} - \epsilon \leq \sigma _{\min }(X) \leq \sigma_{\max }(X) \leq 1+\sqrt{\kappa } + \epsilon  \right] \geq 1 - 2 \exp\left(- \frac{\epsilon^{2}p}{8 \kappa }\right). \) 
  \end{enumerate}
\end{lemma}
Lemma~\ref{lem:apd.sing.val.limit.concent}(a) is cited from Theorem~2 of \citet{bai2008limit}. Lemma~\ref{lem:apd.sing.val.limit.concent}(b) is cited from Corollary 5.35 of \citet{vershynin2010introduction}, slightly adapted by re-writing the bound in terms of \( \kappa =\lim_{p} p/N \).

\begin{lemma}[Kashin's Theorem, 1977]\label{thm:kashin}
For any positive number $v$ there exist a universal constant $c_v$ such that for any $n \geq 1$, with probability at least $1 - 2^{-n}$, for a uniformly random subspace $V_{n,v}$ of dimension $\lfloor n(1 - v)\rfloor$, $\forall x \in V_{n,v} : c_v \|x\|_2 \leq \frac{1}{\sqrt{n}} \|x\|_1$.
\end{lemma}

\begin{lemma}[Non-Asymptotic Singular Value Bounds]\label{lem:apd.sing.val.bounds}
  Let \( \sigma _{\min }(\cdot ) \) and \( \sigma _{\max }(\cdot ) \) denote the smallest and largest non-zero singular values of a matrix. 
  Suppose \( X\in\R^{N\times p} \), \( \{X_{i,j}\}\) are iid Gaussian with the same variance for all \( i\in[N],j\in[p] \). Let \( \{ \Sigma _{i}\in\R^{p\times p}:i\in[N] \}  \) be a collection of covariance matrices satisfying \( \min _{i\in[N]} \lambda _{\min }(\Sigma _{i}) > 0 \). Let \( X_{\Sigma}\in\R^{N\times p} \) be defined such that \( (x_{\Sigma})^{(i)} = x^{(i)} \Sigma_{i}^{1/2} \), where \( (x_{\Sigma})^{(i)}, x^{(i)} \) are respectively the \( i \)-th rows of \( X_{\Sigma} \) and \( X \). Then: 
  \begin{enumerate}
    \item[(a)] There exists  \( \widetilde{X} \in\R^{N\times p} \), s.t. \( \widetilde{X} \disteq X \) and \( \sigma _{\min }(X_{\Sigma}) \geq ({\min _{i\in[N]} \lambda _{\min }(\Sigma _{i})})^{1/2} \sigma_{\min} ({\widetilde{X} }) \);
    \item[(b)]  There exists  \( \widetilde{X} \in\R^{N\times p} \), s.t. \( \widetilde{X} \disteq X \) and \( \sigma _{\max }(X_{\Sigma}) \leq ({\max _{i\in[N]} \lambda _{\max }(\Sigma _{i})})^{1/2} \sigma_{\max} ({\widetilde{X} }) \). 
  \end{enumerate}
\end{lemma}
\begin{proof}[\textbf{Proof of Lemma~\ref{lem:apd.sing.val.bounds}}]
  Without loss of generality we let \( \Var(X_{i,j})=1 \). Let \( z_1,\cdots ,z_N \in\R \) be iid standard normal variables to be used later. For any \( \Sigma _i,i\in[N] \), let its spectral decomposition be \( \Sigma _i = S^{T}_i \Lambda _{i} S_i \), where \( S \) is orthogonal and \( \Lambda _{i} \) is diagonal. 
  
  We prove \emph{(a)} of Lemma~\ref{lem:apd.sing.val.bounds}, and \emph{(b)} can be proved analogously. Let \( c := \min _{i\in[N]} \lambda _{\min }(\Sigma _{i}) >0 \).
  For any vector \( v\in\R^{p} \), 
  \begin{equation*}
\begin{aligned}
  \norm{X^{\Sigma}v}_2^{2} =&~ \sum_{i\in[N]} (x^{(i)}\Sigma _{i}^{1/2}v)^{2} \disteq \sum_{i\in[N]} (X_{i,\cdot }\Lambda _{i}^{1/2}S_i v)^{2} = \sum_{i\in[N]} \left[\sum_{j\in[p]} X_{i,j} \Lambda^{1/2}_{i,j} (S_i v)_j \right]^{2} 
  \disteq \sum_{i\in[N]} z^{2}_{i} \cdot \left[ \sum_{j\in[p]} \Lambda _{i,j} (S_i v)_j^{2} \right] \\
  \geq &~ c \cdot \sum_{i\in[N]} z_{i}^{2} \cdot \left[ \sum_{j\in[p]} (S_i v)_j^{2} \right] = c \cdot \sum_{i\in[N]} z^{2}_{i} \cdot \norm{v}_2^{2} \disteq c \cdot \sum_{i\in[N]} \left[\sum_{j\in[p]} \widetilde{X}_{i,j} v_j \right]^{2} = c\cdot \sum_{i\in[N]} \left(\tilde{x}^{(i)} v \right)^{2} = c\,\norm{\widetilde{X} v}_2^{2}.
\end{aligned}
\end{equation*}
As a direct corollary, we know $\|\widetilde{X} v\|_2=0$ whenever $\|X^{\Sigma} v\|_2 = 0$, which means the null space of $X^{\Sigma}$ and $\widetilde{X}$ coincide.
Thus the bound holds for smallest non-zero singular value, and we have proved Lemma~\ref{lem:apd.sing.val.bounds} (a). 
\end{proof}

\end{document}